\documentclass[12pt]{article}
\usepackage{fullpage,amsmath,amsthm,amsfonts,amssymb,graphicx,setspace,amscd,float}
\DeclareMathOperator\tni{int}

\newtheorem{thm}{Theorem}[section]
\newtheorem{lem}[thm]{Lemma}
\newtheorem{qst}[thm]{Question}
\newtheorem{prop}[thm]{Proposition}
\newtheorem{cor}[thm]{Corollary}
\theoremstyle{definition}
\newtheorem{df}[thm]{Definition}
\newtheorem{rk}[thm]{Remark}
\newtheorem{ex}[thm]{Example}

\begin{document}

\title{Constructing and Classifying Fully Irreducible Outer Automorphisms of Free Groups}
\author{Catherine Pfaff}
\date{\today}
\maketitle

\abstract{The main theorem of this document emulates, in the context of $Out(F_r)$ theory, a mapping class group theorem (by H. Masur and J. Smillie) that determines precisely which index lists arise from pseudo-Anosov mapping classes.  Since the ideal Whitehead graph gives a finer invariant in the analogous setting of a fully irreducible $\phi \in Out(F_r)$, we instead focus on determining which of the 21 connected, loop-free, 5-vertex graphs are ideal Whitehead graphs of ageometric, fully irreducible $\phi \in Out(F_3)$.  Our main theorem accomplishes this by showing that there are precisely 18 graphs arising as such.  We also give a method for identifying certain complications called periodic Nielsen paths, prove the existence of conveniently decomposed representatives of ageometric, fully irreducible $\phi\in Out(F_r)$ having connected, loop-free, ($2r-1$)-vertex ideal Whitehead graphs, and prove a criterion for identifying representatives of ageometric, fully irreducible $\phi\in Out(F_r)$.  The methods we use for constructing fully irreducible outer automorphisms of free groups, as well as our identification and decomposition techniques, can be used to extend our main theorem, as they are valid in any rank.  Our methods of proof rely primarily on Bestvina-Feighn-Handel train track theory and the theory of attracting laminations.}

\newpage

\tableofcontents

\newpage

\section{Introduction}
\noindent The main theorem of this document (Theorem \ref{T:MainTheorem}) is motivated by a theorem in mapping class group theory.  The \emph{mapping class group $MCG(S)$} of a compact surface $S$ is the group of homotopy classes of homeomorphisms $h: S\to S$.  The most common mapping classes are called \emph{pseudo-Anosov}.  (One characterization of a pseudo-Anosov mapping class is that some representative of a pseudo-Anosov mapping class expands and contracts a pair of transverse singular measured foliations on the surface).  Because of their fundamental importance in topology and geometry, both mapping class groups and pseudo-Anosov mapping classes have been objects of extensive research.  The list of singularity indices associated to a pseudo-Anosov mapping class is an important invariant of the class.  (Each foliation singularity for the pair of transverse singular measured foliations has an associated index).  H. Masur and J. Smillie (see [MS93]) proved precisely which lists of singularity indices arise from pseudo-Anosov mapping classes.  This document is the first step to proving an analogous theorem for outer automorphism groups of free groups.

For a free group of rank $r$, \emph{$F_r$}, the \emph{outer automorphism group}, $Out(F_r)$, consists of equivalence classes of automorphisms $\Phi: F_r\to F_r$, where two automorphisms are equivalent when they differ by an inner automorphism, i.e. a map $\Phi_b$ defined by $\Phi_b(a)=b^{-1}ab$ for all $a \in F_r$.  Outer automorphisms can be described geometrically as follows.  Let $R_r$ be the $r$-petaled rose (graph having $r$ edges and a single vertex $v$).  Given a graph $\Gamma$  with no valence-one vertices, we can assign to $\Gamma$ a $\emph{marking}$ (identification of the fundamental group with the free group $F_r$) via a homotopy equivalence $R_r \to \Gamma$.  We call such a graph $\Gamma$, together with its marking $R_r \to \Gamma$, a \emph{marked graph}.  Each element $\phi$ of $Out(F_r)$ can be represented geometrically by a homotopy equivalence $g: \Gamma \to \Gamma$ of a marked graph, where $\phi=g_*$ is the induced map of fundamental groups.  For the proof of our analogue to the mapping class group theorem, we will focus on constructing representatives of ageometric, fully irreducible outer automorphisms with particular ideal Whitehead graphs.  An ideal Whitehead graph is a strictly finer invariant than a singularity index list and encodes information about the attracting lamination for a fully irreducible outer automorphism.

There is potential for an $Out(F_r)$ analogue to the Masur-Smillie theorem because of deep connections between outer automorphism groups of free groups and mapping class groups of surfaces.  When $r=2$, we have $Out(F_2) \cong Out(\Pi_1(\Sigma_{1,1})) \cong MCG(\Sigma_{1,1}) $, where $\Sigma_{1,1}$ denotes a genus-1 torus with a single punture. Furthermore, elements $\phi \in Out(F_2)$ are induced by homeomorphisms of $\Sigma_{1,1}$ and ``fully irreducible outer automorphisms'' (Subsection \ref{S:Reducibility}) are induced by pseudo-Anosov homeomorphisms.  While we do not have such exact correspondences for $r>2$, there are still strong similarities between all of the outer automorphism groups $Out(F_r)$ and mapping class groups $MCG(S)$, as well as between the fully irreducible $\phi \in Out(F_r)$ and pseudo-Anosov $\psi \in MCG(S)$.  In fact, some $\phi \in Out(F_r)$ with $r>2$ are even still induced by homeomorphisms of a compact surface with boundary (such $\phi$ are called \emph{geometric}).

There is a large group of mathematicians exploring the parallel properties between the $Out(F_r)$ groups and the $MCG(S)$ groups. They have made significant progress to this affect.  We use some of their definitions and machinery (including the definitions of singularities, indices, and ideal Whitehead graphs for outer automorphisms of free groups, as defined in [GJLL98] and [HM11]), in order to understand an appropriate $Out(F_r)$-analogue to the Masur-Smillie theorem for mapping class groups (as described in the next section).

\subsection{The Question}

Let $i(\phi)$ denote the sum of the singularity indices for a fully irreducible $\phi \in Out(F_r)$ as defined in [GJLL98] and in Subsection \ref{S:IWGs} below.  An index can in some sense be thought of as recording the number of germs of initial edge segments (\emph{directions}) emanating from a vertex that are fixed by a given geometric representative of $\phi$.  [GJLL98] gives an inequality $i(\phi) \geq 1-r$ bounding the index sum for a fully irreducible outer automorphism $\phi \in Out(F_r)$, in contrast with the equality $i(\psi)=\chi(S)$, for a pseudo-Anosov $\psi$ on a surface $S$, dictated by the Poincare-Hopf Theorem.  With this in mind, one can ask whether every index list who's sum satisfies this inequality is achieved.  M. Handel and L. Mosher pose the question in [HM11]:

\begin{qst}{\label{Q:IndexIndequality}}  What possible index types, which satisfy the index inequality $i(\phi) \geq 1-r$, are achieved by a nongeometric, fully irreducible element of $Out(F_r)$?  \end{qst}

What we present here focuses on several constructions that may eventually allow us to attack this question directly and that, in the meantime, give us a stronger result for the rank-3 case when restricting to ``ageometric'' fully irreducible outer automorphisms and connected, $(2r-1)$-vertex ideal Whitehead graphs with no single-vertex edges (see Theorem \ref{T:MainTheorem}).  For an ageometric, fully irreducible $\phi \in Out(F_r)$ and TT representative $g: \Gamma \to \Gamma$ on a rose having $2r-1$ fixed directions at the unique vertex and no \emph{periodic Nielsen paths (PNPs)}, i.e. paths $\rho$ in $\Gamma$ such that $g^k(\rho) \cong \rho$ for some $k$, the ideal Whitehead graph $IW(\phi)$ is the graph with one vertex for each fixed direction of $\Gamma$ and an edge between two such vertices when there exists some $k >0$ and edge $e$ of $\Gamma$ such that $g^k(e)$ crosses over the turn formed by the directions corresponding to the vertices.

Having $2r-1$ vertices is maximal when we refine our search to ageometric, fully irreducible outer automorphisms and connected, loop-free graphs.  We focus on ageometric, fully irreducible outer automorphisms, as they are far more common and better understood than parageometrics, the only other kind of nongeometric, fully irreducible outer automorphism (geometric outer automorphisms are induced by surface homeomorphisms, thus are already understood).  We focus on connected graphs, as this allows us to only look at homotopy equivalences of roses (see Proposition \ref{P:IdealDecomposition}).

In the mapping class group case, one only sees circular ideal Whitehead graphs, making singularity index lists the best possible invariant.  However, this is not true for fully irreducible outer automorphisms.  Thus, a better analogue to the Masur-Smillie theorem would record possible ``ideal Whitehead graphs,'' instead of just singularity indices.  For an ageometric, fully irreducible $\phi \in Out(F_r)$, the index of a component in $IW(\phi)$ is simply $1-\frac{k}{2}$, where $k$ is the number of vertices of the component.  Ideal Whitehead graphs for outer automorphisms of free groups are defined in [HM11] and discussed in Subsection \ref{S:IWGs} below.  In the spirit of focusing on ideal Whitehead graphs, the question we give a partial answer to in this document (see Theorem \ref{T:MainTheorem}) is that also posed by L. Mosher and M. Handel in [HM11]:

\begin{qst}  Which ideal Whitehead graphs arise from ageometric, fully irreducible outer automorphisms of rank-3 free groups?  \end{qst}

We call graphs with no single-vertex edges (i.e. no edges are loops) \emph{loop-free} and a connected, ($2r-1$)-vertex, loop-free graph a \emph{Type (*) potential ideal Whitehead graph} or \emph{Type (*) pIW graph (pIWG)}.  The partial answer (Theorem \ref{T:MainTheorem}) we give completely answers the following subquestion posed in person by L. Mosher and M. Feighn:

\begin{qst}  Which of the twenty-one five-vertex Type (*) pIW graphs are ideal Whitehead graphs for ageometric, fully irreducible $\phi \in Out(F_3)$?  \end{qst}

\noindent What we state here is Theorem \ref{T:MainTheorem} and is the complete answer to Questions \ref{Q:IndexIndequality}.

\smallskip

\begin{thm} Precisely eighteen of the twenty-one five-vertex Type (*) pIW graphs are ideal Whitehead graphs for ageometric, fully irreducible $\phi\in Out(F_3)$. \end{thm}

As mentioned above, we chose to look at $5$-vertex graphs because, with the restriction that $\phi\in Out(F_3)$ must be ageometric and fully irreducible and the restriction that $IW(\phi)$ is loop-free and connected, five vertices is maximal.  We focused on connected graphs as this allowed us to focus on representatives on the rose (see Proposition \ref{P:IdealDecomposition}).

As they will be used throughout, we list now the $21$ $5$-vertex Type (*) pIWGs.

\begin{figure}[H]
\centering
\includegraphics[width=4.5in]{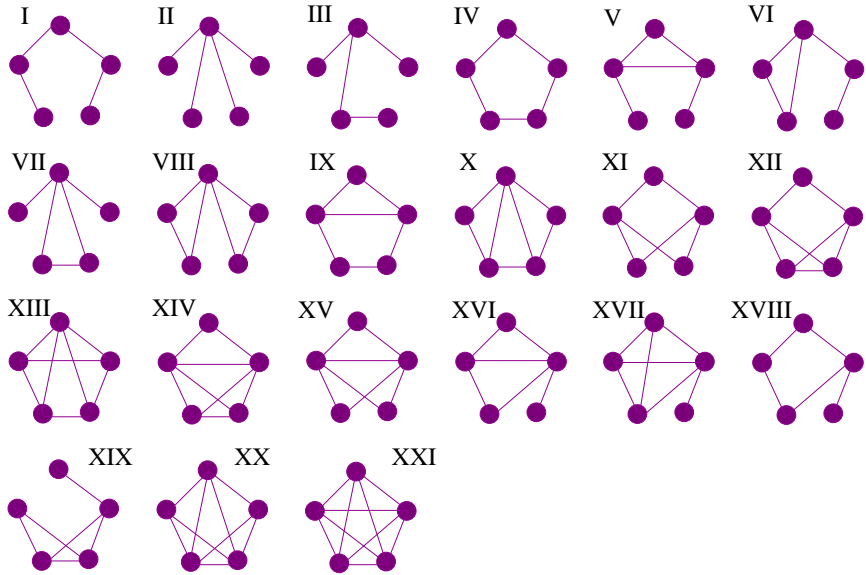}
\caption{{\small{\emph{The 21 $5$-vertex Type (*) pIW graphs (up to graph isomorphism) (See [CP84])}}}}
\label{fig:21Graphs}
\end{figure}

\begin{rk}
The $5$-vertex Type (*) pIWGs that are not ideal Whitehead graphs for ageometric, fully irreducible $\phi\in Out(F_3)$ are:

\begin{figure}[H]
\centering
\includegraphics[width=3.3in]{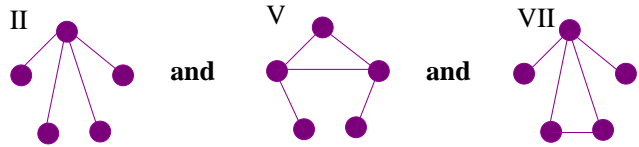}
\caption{{\small{\emph{Unachievable Graphs}}}}
\label{fig:UnachievableGraphs}
\end{figure}
\end{rk}

\subsection{Outline of Document}

\noindent The first step to proving the main theorem, Theorem \ref{T:MainTheorem}, is Proposition \ref{P:IdealDecomposition}:

\begin{prop} Let an ageometric, fully irreducible $\phi \in Out(F_r)$ be such that $IW(\phi)$ is a Type (*) pIWG.  Then there exists a PNP-free, rotationless representative of a power $\psi = \phi^R$ of $\phi$ on the rose.  Furthermore, the representative can be decomposed as a sequence of proper full folds of roses. \newline \emph{(The point of $\psi$ being rotationless is its representative fixing the periodic directions.)}
\end{prop}

\smallskip

The significance, for our purposes, of the proposition is its refining our search for representatives with desired ideal Whitehead graphs to those ``ideally decomposed,'' as in the proposition.

The next step to proving Theorem \ref{T:MainTheorem} is to define, as we do in Section \ref{Ch:LTT}, ``lamination train track structures '' (LTT structures). We ``build'' or ``construct'' portions of desired ideal Whitehead graphs by determining ``construction compositions'' from smooth paths in LTT structures.  ``Construction compositions'' are in ways analogues to Dehn twists (mapping class group elements used in pseudo-Anosov construction methods, including those of Penner in [P88]).  We appropriately compose construction compositions to construct (for the Theorem \ref{T:MainTheorem} proof) the representatives of outer automorphisms with particular ideal Whitehead graphs.  On the other hand, we use that LTT structures for ageometric, fully irreducible outer automorphisms are ``birecurrent'' to show in Proposition \ref{P:Birecurrent} that the ideal Whitehead graph of an ageometric, fully irreducible outer automorphism cannot be of a certain type.  While stated in the restricted form it is used for (and we give definitions for) in this document, what the Proposition \ref{P:Birecurrent} proof really says is:

\begin{prop} The LTT structure for a train track representative of a fully irreducible outer automorphism is birecurrent.
\end{prop}

However, Proposition \ref{P:Birecurrent} only explains one of the three graphs Theorem \ref{T:MainTheorem} deems unachievable.  We explain in Section \ref{Ch:Procedure} methods used to show the remaining two graphs are also unachievable.  We prove their unachievability in Section \ref{Ch:UnachievableIWGs}.

To determine which construction compositions to use and how to appropriately compose them, we define ``AM Diagrams'' in Section \ref{Ch:AMDiagrams}.   If there is an ageometric, fully irreducible $\phi \in Out(F_r)$ with a particular Type (*) pIWG, $\mathcal{G}$, as its ideal Whitehead graph, then there is a loop in the AM Diagram for $\mathcal{G}$ that corresponds to an ``ideal decomposition'' (defined in Section \ref{Ch:IdealDecompositions}) of a representative $g$ for a power of $\phi$.  This fact is proved in Proposition \ref{P:RepresentativeLoops} and helps us rule out unobtainable ideal Whitehead graphs.  Additionally, Sections \ref{Ch:AMDiagrams} and \ref{Ch:Procedure} tell us how to construct representatives yielding the obtainable Type (*) pIWGs.

Finally, in order to prove that our representatives are representatives of ageometric fully irreducible outer automorphisms, we proved in Section \ref{Ch:FIC} the ``Full Irreducibility Criterion''  or ``FIC'' (Lemma \ref{L:FIC}).  We need three conditions to apply the criterion.  First, the FIC requires that a representative is PNP-free.  Proposition \ref{P:NPIdentification} of Section \ref{Ch:NPIdentification} offers a method for identifying PNPs for an ideally decomposed train track representative $g$.  Second, the criterion includes a condition that the local Whitehead graph at a vertex be connected (a condition satisfied in our case by the ideal Whitehead graph being connected).  A local Whitehead graph records how images of edges enter and exit a particular vertex. For a representative $g$ at a vertex $x$ the local Whitehead graph will be denoted by $LW(g; x)$.  Finally, the criterion includes a condition on the transition matrix for $g: \Gamma \to \Gamma$ satisfied when there exists some $k >0$ such that $g^k$ maps each edge of $\Gamma$ over each other edge of $\Gamma$.

\begin{lem} (The Full Irreducibility Criterion) Let $g: \Gamma \to \Gamma$ be an irreducible train track representative of $\phi \in Out(F_r)$.  Suppose that $g$ has no PNPs, that the transition matrix for $g$ is Perron-Frobenius, and that all $LW(g; x)$ for $g$ are connected.  Then $\phi$ is fully irreducible. \end{lem}

For our proof of the criterion we appeal to the train track machinery of M. Bestvina, M. Feighn, and M. Handel.  The proof uses several different revised versions (defined in [BFH00] and [FH09]) of ``relative train track representatives.''  Definitions of the relative train track representatives relevant to our situation are also given in Section \ref{Ch:FIC}.  Outside of Section \ref{Ch:FIC} we restrict our discussions to train track representatives.

As a final note, we comment that, while our methods have been designed for constructing ageometric, fully irreducible outer automorphisms with particular ideal Whitehead graphs, and thus are particularly well-suited for this purpose, there are other methods for constructing fully irreducible outer automorphisms, such as the recent one described in [CP10].

\subsection{Summary of Results}

\indent The main theorem of this document (Theorem \ref{T:MainTheorem}) lists precisely which Type (*) pIWGs arise as ideal Whitehead graphs for ageometric, fully irreducible $\phi\in Out(F_3)$.  Of independent interest and use are two propositions and a lemma used in the proof of the main theorem.  The propositions and lemma show the existence of a useful decomposition of a certain class of fully irreducible $\phi \in Out(F_r)$ (Proposition \ref{P:IdealDecomposition}), a method for identifying pNps (Proposition \ref{P:NPIdentification}), and a criterion for identifying representatives of ageometric fully irreducible $\phi\in Out(F_r)$ (Lemma \ref{L:FIC}).  Finally, our Theorem \ref{T:MainTheorem} proof outlines fully irreducible representative construction methods that have use beyond proving Theorem \ref{T:MainTheorem}, as they can, for example, they could be used in any rank.

\subsection{Acknowledgements}

\indent The author would like to thank her thesis advisor Lee Mosher for his incredible patience, suggestions, ideas, corrections, and answers to her endless questions.  The author would like to thank Mark Feighn for his numerous meetings with her and answers to her questions, as well as his recommendations for uses of the methods she developed.  The author would like to thank Michael Handel for meeting with her and his recommendations, particularly on how to finish our Full Irreducibility Criterion proof.

 The author would like to thank Mladen Bestvina, Alexandra Pettet, Andy Putman, Martin Lustig, Benson Farb, and Ilya Kapovich for their conversations and support.

 The author would like to thank Yael Algom-Kfir and Johanna Mangahas for the knowledge and understanding they imparted on her during their numerous conversations, as well as the inspiration their strength and courage have always provided.

 Finally, the author thanks Bard College at Simon's Rock for their financial support.

\section{Preliminary Definitions}{\label{Ch:PrelimDfns}}

In this section we give some definitions used throughout the document.  We continue with the notation established in the introduction.

\subsection{Train Track Representatives}

Thurston defined a homotopy equivalence $g:\Gamma \to \Gamma$ of marked graphs to be a \emph{train track map} if, for all $k>0$, the restriction of $g^k$ to the interior of each edge of $\Gamma$ is locally injective (there is no ``backtracking'' in edge images).  If $g$ induces a $\phi \in Out(F_r)$ (as a map of fundamental groups) and $g(\mathcal{V}) \subset \mathcal{V}$ (where $\mathcal{V}$ is the vertex set of  $\Gamma$) then $g$ is called a \emph{topological (or train track) representative} for $\phi$ [BH92].  Train track representatives are in many ways the most natural representatives to work with and [BH92] gives an algorithm for finding a train track representative of any irreducible $\phi \in Out(F_r)$.  In this document we focus on train track representatives and several versions of their more general ``relative train track'' representatives, as defined in Section \ref{Ch:FIC}.  Unless otherwise stated, one should assume throughout this document that a representative of an outer automorphism is a train track representative.

\subsection{Reducibility}{\label{S:Reducibility}}

``Fully irreducible'' outer automorphisms of free groups are induced by pseudo-Anosov homeomorphisms of surfaces in rank $2$ and still resemble pseudo-Anosovs in higher rank.  They are our main focus and can be defined either algebraically or geometrically.  We algebraically define fully irreducible outer automorphisms, but geometrically define representative reducibility and irreducibility.

Again let $F_r$ denote the free group of rank $r$.  An outer automorphism $\phi \in Out(F_r)$ is \emph{reducible} if there are proper free factors $F^1, \dots , F^k$ of $F_r$ such that $\phi$ permutes the conjugacy classes of the $F^i$'s and such that $F^1\ast \dots \ast F^k$ is a free factor of $F_r$ (i.e. there is a free group $F_l$ such that ($F^1\ast \dots \ast F^k$)$\ast F_l=F_r$).  If $\phi$ is not reducible, then we say $\phi$ is \emph{irreducible}.  If every power of $\phi$ is irreducible, then we say $\phi$ is \emph{fully irreducible}.

A train track representative $g: \Gamma \to \Gamma$ of a $\phi \in Out(F_r)$ is \emph{reducible} if it has a nontrivial invariant subgraph $\Gamma_0$ (meaning $g(\Gamma_0) \subset \Gamma_0$) with at least one noncontractible component. The representative $g$ is otherwise called \emph{irreducible}. [BH92, BFH97]

It is important to note that a reducible outer automorphism may have irreducible representatives.  It is only necessary that it have at least one reducible representative in order for it to be reducible.  Thus, a fully irreducible outer automorphism is an outer automorphism such that no power has a reducible representative.

\subsection{Turns, Paths, Circuits, and Tightening}{\label{S:Paths}}

We remind the reader here of a few definitions (unless otherwise indicated) from [BH92].  These definitions are important in establishing notions of ``legality,'' prevalent in train track theory, and are needed to define ideal Whitehead graphs, the outer automorphism invariant finer than a singularity index list.  First we establish notation.

Let $g: \Gamma \to \Gamma$ be a train track representative of $\phi \in Out(F_r)$ and $\mathcal{E}^+(\Gamma)= \{E_1, \dots, E_{n}\}$ the set of edges in $\Gamma$ with some prescribed orientation.  For any edge $E \in \mathcal{E}^+(\Gamma)$, let $\overline{E}$ denote $E$ oriented oppositely as $E$ and then let $\mathcal{E}(\Gamma)=\{E_1, \overline{E_1}, \dots, E_n, \overline{E_n} \}$.  Finally, let $\mathcal{V}(\Gamma)$ denote the set of vertices of $\Gamma$ (or just $\mathcal{V}$, when $\Gamma$ is clear).  We continue with this notation throughout the document.

Next we state the definition versions used here for paths and circuits. These notions are important for analyzing images of edges under train tracks and in discussions of RTTs and laminations.

Let $\Gamma$ be a marked graph with universal cover $\tilde{\Gamma}$ and projection map $p: \tilde{\Gamma} \to \Gamma$.  A \emph{path} in $\tilde{\Gamma}$ is either a proper embedding $\tilde{\gamma}: I \to \tilde{\Gamma}$, where $I$ is a (possibly infinite) interval, or a map of a point $\tilde{\gamma}: x \to \tilde{\Gamma}$.  \emph{Paths} in $\Gamma$ are projections $p \circ \tilde{\gamma}$, where $\tilde{\gamma}$ is a path in $\tilde{\Gamma}$.  Paths differing by an orientation-preserving change of parametrization are considered to be the same path.  [BFH00]

Let $\tilde{\gamma}$ be a path in $\tilde{\Gamma}$.  Then $\tilde{\gamma}$ can be written as a concatenation of subpaths, each of which is an oriented edge of $\tilde{\Gamma}$ (with the exception that the first and last subpaths of $\tilde{\gamma}$ may actually only be partial edges).  We call this sequence of oriented edges (and partial edges) the \emph{edge path associated to $\tilde{\gamma}$}.  Its projection gives a decomposition of $\gamma$ as a concatenation of oriented edges in $\Gamma$.  We will call this sequence of oriented edges in $\Gamma$ the \emph{edge path associated to $\gamma$}.  [BFH00]

A \emph{circuit} in $\Gamma$ is an immersion $\alpha : S^1 \to \Gamma$ of the circle into $\Gamma$.  Edge paths for circuits in $\Gamma$ are defined the same as paths in $\Gamma$ except that the edges for only one period of the edge path are listed. [BFH00]  This will mean that there can be multiple edge paths representing the same circuit.

Directions will be important for defining ideal Whitehead graphs and will be prevalent throughout the proofs of this document.

The \emph{directions} at a point $x\in \Gamma$ are the germs of initial segments of edges emanating from $x$.  Let $\mathcal{D}(x)$ denote the set of directions at $x$ and $\mathcal{D}(\Gamma)=\underset{v \in \mathcal{V}(\Gamma)}{\cup} \mathcal{D}(v)$.  For an edge $e \in \mathcal{E}(\Gamma)$, let $D_0(e)$ denote the initial direction of $e$ (the germ of initial segments of e).  For a path $\gamma=e_1 \dots e_k$, define $D_0 \gamma = D_0(e_1)$.  We denote the map of directions induced by $g$ as \emph{$Dg$}, i.e. $Dg(d)=D_0(g(e))$ for $d=D_0(e)$.  (Note that $D(f \circ g)=Df \circ Dg$).  $d \in \mathcal{D}(\Gamma)$ is \emph{periodic} if $Dg^k(d)=d$ for some $k>0$ and \emph{fixed} if $k=1$.  We denote the set of periodic directions at a $x \in \Gamma$ by $Per(x)$ and the set of fixed points by $Fix(x)$.

The following notions of turns, legality, and tightening will be important for stating the properties of the different RTT variants and laminations.  They will also be prevalent throughout the proofs of this document.

A \emph{turn} in $\Gamma$ is defined as an unordered pair of directions $\{d_1, d_2\}$ at a vertex $v \in \Gamma$.  Let $\mathcal{T}(v)$ denote the set of turns at $v$.  For a vertex $v \in \Gamma$, $Dg$ induces a map of turns $D^tg$ on $\mathcal{T}(v)$, defined by $D^tg(\{d_1, d_2\}) = \{Dg(d_1), Dg(d_2)\}$ for each $\{d_1, d_2\} \in \mathcal{T}(v)$. A turn $\{d_i, d_j \}$ is  \emph{degenerate} if $d_i = d_j$ and \emph{nondegenerate} otherwise.  The turn is \emph{illegal} with respect to $g$: $\Gamma \to \Gamma$ if some $D^tg^k(\{d_1, d_2 \})$ is degenerate and is otherwise \emph{legal}.

It is an important property of any train track representative $g: \Gamma \to \Gamma$ that one never has $g^k(e)= \dots \overline{e_i} e_j \dots$, where $D_0(e_i)=d_i$, $D_0(e_j)=d_j$, $\{d_i, d_j\}$ is an illegal turn for $g$, and $e, e_i, e_j \in \mathcal{E}(\Gamma)$.  (In other words, for a train track representative $g$, no iterate of $g$ maps an edge over an illegal turn).

The set of \emph{gates} with respect to $g$ at a vertex $v \in \Gamma$ is the set of equivalence classes in $\mathcal{D}(v)$ where $d \sim d'$ if and only if $(Dg)^k(d)=(Dg)^k(d')$ for some $k \geq 1$.  In other words, pairs of directions in the same gate form illegal turns and pairs of directions in different gates form legal turns.

For an edge path $e_1e_2 \dots e_{k-1}e_k$ associated to a path $\gamma$ in $\Gamma$, we say that $\gamma$ \emph{contains (or crosses over)} the turn $\{\overline{e_i}, e_{i+1}\}$ for each $1 \leq i < k$.  A path $\gamma \in \Gamma$ is called \emph{legal} if it does not contain any illegal turns and \emph{illegal} if it contains at least one illegal turn.

Every map of the unit interval $\tilde{\alpha}: I \to \tilde{\Gamma}$ is homotopic rel endpoints to a unique path in $\tilde{\Gamma}$, which we denote by $[\tilde{\alpha}]$.  We then say that $[\tilde{\alpha}]$ is obtained from $\tilde{\alpha}$ by \emph{tightening}.  ([$\tilde{\alpha}$] is obtained from $\tilde{\alpha}$ by removing all ``backtracking'').  If $\alpha$ is the projection to $\Gamma$ of $\tilde{\alpha}$, then the projection $[\alpha]$ of $[\tilde{\alpha}]$ is said to be obtained from $\alpha$ by \emph{tightening}.

A homotopy equivalence $g: \Gamma \to \Gamma$ is \emph{tight} if, for each edge $e \in \mathcal{E}(\Gamma)$, either $g(e) \in \mathcal{V}$ or $g$ is locally injective on $int(e)$.  Any homotopy equivalence can be tightened to a unique tight homotopy equivalence by a homotopy rel $\mathcal{V}$.  For a train track representative $g: \Gamma \to \Gamma$, we define $g_{\#}$ by $g_{\#}(\alpha)=[g(\alpha)]$ for each path $\alpha$ in $\Gamma$.

\subsection{Lines}

We give here several definitions from [BFH00].  These definitions will be important for defining the laminations analogous to the attracting laminations for pseudo-Anosovs that we use in the proofs of Proposition \ref{P:Birecurrent} and Lemma \ref{L:FIC}.

We start by establishing the notion of a line in a marked graph and in its universal cover.  Again let $\Gamma$ be a marked graph with universal cover $\tilde{\Gamma}$ and projection map $p: \tilde{\Gamma} \to \Gamma$.  A \emph{line} in $\tilde{\Gamma}$ is the image of a proper embedding of the real line $\tilde{\lambda}: \textbf{R} \to \tilde{\Gamma}$.  We denote by $\tilde{\mathcal{B}}(\Gamma)$ the space of lines in $\tilde{\Gamma}$ with the compact-open topology (one can define a basis for $\tilde{\mathcal{B}}(\Gamma)$ where an open set consists of all lines sharing a given line segment).  A \emph{line} in $\Gamma$ is the image of a projection $p \circ \tilde{\lambda}$ of a line $\tilde{\lambda}$ in $\tilde{\Gamma}$.  We denote by $\mathcal{B}(\Gamma)$ the space of lines in $\Gamma$ with the quotient topology induced by the natural projection map from $\tilde{\mathcal{B}}(\tilde{\Gamma})$ to $\mathcal{B}(\Gamma)$.

Now we give a characterization of lines as pairs of points in the space of ends of $\tilde{\Gamma}$ (viewed as $\partial F_r$).  We then relate this characterization back to the definitions just given. The characterization of lines as pairs of points in the space of ends of $\tilde{\Gamma}$ is used to discuss laminations.  Let $\Delta$ be the diagonal in $\partial F_r$ x $\partial F_r$.  $\tilde{\mathcal{B}}$ is obtained from ($\partial F_r$ x $\partial F_r$) - $\Delta$ by quotienting out by the action that interchanges the factors of $\partial F_r$ x $\partial F_r$.  We denote by $\mathcal{B}$ the quotient of $\tilde{\mathcal{B}}$ under the diagonal action of $F_r$ on $\partial F_r$ x $\partial F_r$. We can identify the Cantor Set $\partial F_r$ with the space of ends of $\tilde{\Gamma}$.  In particular, if ($b_1, b_2$) $\in \partial F_r$ x $\partial F_r$ is an unordered pair of distinct elements of $\partial F_r$, then there exists a unique line $\tilde{\gamma} \in \tilde{\Gamma}$ with endpoints corresponding to $b_1$ and $b_2$.  This defines a homeomorphism between $\tilde{\mathcal{B}}$ and $\tilde{\mathcal{B}}(\tilde{\Gamma})$ that projects to a homeomorphism between $\mathcal{B}$ and $\mathcal{B}(\Gamma)$ (see [BFH00]).  For a path $\beta \in \mathcal{B}$, we say that \emph{$\gamma \in \mathcal{B}(\Gamma)$ realizes $\beta$ in $\Gamma$} if $\gamma$ corresponds to $\beta$ under the projection of the homeomorphism between $\tilde{\mathcal{B}}$ and $\tilde{\mathcal{B}}(\Gamma)$.

As we use it in Subsection \ref{SS:Laminations}, we give one last definition here.  For a marked graph $\Gamma$, we say that a line $\tilde{\gamma}$ in $\tilde{\Gamma}$ is \emph{birecurrent} if every finite subpath of $\tilde{\gamma}$ occurs infinitely often as an unoriented subpath in each end of $\tilde{\gamma}$.  A line $\gamma$ in $\Gamma$ representing a birecurrent line $\tilde{\gamma} \in \tilde{\Gamma}$ (with either choice of orientation) is called \emph{birecurrent}. [BFH00]

\subsection{Laminations}{\label{SS:Laminations}}

The following two definitions are required to state the attracting lamination definition for a $\phi \in Out(F_r)$.  Attracting laminations are used in Section \ref{Ch:AMProperties} to prove a necessary condition for a Type (*) pIWG to be the ideal Whitehead graph of a fully irreducible $\phi \in Out(F_r)$ and are used in the the Full Irreducibility Criterion proof (Section \ref{Ch:FIC}).  To avoid reading about laminations, one can simply skip the proofs requiring them in Sections \ref{Ch:AMProperties} and \ref{Ch:FIC}.  All definitions in this subsection are from [BFH00].

\begin{df} An \emph{attracting neighborhood} of $\beta \in \mathcal{B}$ for the action of $\phi$ is a subset $U \subset \mathcal{B}$ such that $\phi_{\#}(U) \subset U$ and $\{\phi^k_{\#}(U) : k \geq 0 \}$ is a neighborhood basis for $\beta$ in $\mathcal{B}$. \end{df}

\begin{df}  For a free factor $F^i$ of $F_r$, [[$F^i$]] will denote the conjugacy class of $F^i$.  Consider the set of circuits in $\mathcal{B}$ determined by the conjugacy classes in $F_r$ of $F^i$.  Lines $\beta \in \mathcal{B}$ in the closure of this set of circuits are said to be \emph{carried by} [[$F^i$]]. \end{df}

\begin{rk}  The notion of ``closure'' in this context can be understood by recognizing that the appropriate notion of convergence in $\mathcal{B}$ is ``weak convergence.''  Suppose that $g: \Gamma \to \Gamma$ represents $\phi$.  If $\gamma' \in \Gamma$ realizes $\beta' \in \mathcal{B}$ and $\gamma \in \Gamma$ realizes $\beta \in \mathcal{B}$, then $\beta'$ is \emph{weakly attracted} to $\beta$ if, for each subpath $\alpha \in \gamma$, $\alpha \subset g^k_{\#}(\gamma')$ for all sufficiently large $k$. \end{rk}

\noindent We are now ready to give the definition of an attracting lamination.

\begin{df}  An \emph{attracting lamination} $\Lambda$ for $\phi \in Out(F_r)$ is a closed subset of $\mathcal{B}$ that is the closure of a single point $\lambda$ which:

(1) is birecurrent,

(2) has an attracting neighborhood for the action of some $\phi^k$, and

(3) is not carried by a $\phi$-periodic rank-1 free factor. \newline
\noindent In such a circumstance we say that $\gamma$ is \emph{generic} for $\Lambda$ or \emph{$\Lambda$-generic}.  $\mathcal{L}(\phi)$ will denote the set of attracting laminations for $\phi$.
\end{df}

\begin{rk} It is proved in [BFH00] that a fully irreducible outer automorphism $\phi \in Out(F_r)$ has a unique attracting lamination associated to it (in fact, any irreducible train track representative having a Perron-Frobenius transition matrix, as defined below, has a unique attracting lamination associated with it).

The notation in the literature for this unique attracting lamination varies in ways possibly confusing to the reader unaware of this fact.  For example, in [BFH97] and [BFH00] it is denoted by $\Lambda^+_{\phi}$, or just $\Lambda^+$, while the authors of [HM11] used the notation $\Lambda_-$, more consistent with the terminology of dynamical systems ($\Lambda_-$ turns out to be the dual lamination to the tree $T_-$).  To avoid confusion, we simply denote the unique attracting lamination associated to the fully irreducible outer automorphism $\phi$ by $\Lambda(\phi)$ (or just $\Lambda$ when we believe $\phi$ to be clear).

In addition to the notational variance, there is also variance in the name assigned to $\Lambda(\phi)$.  An attracting lamination is called a \emph{stable lamination} in [BFH97].  It is also referred to in the literature, at times, as an \emph{expanding lamination}.
\end {rk}

\subsection{Periodic Nielsen Paths and Geometric, Parageometric, and Ageometric Fully Irreducible Outer Automorphisms}

Recall that ``periodic Nielsen paths'' are important for determining fully irreducibility (see the Full Irreducibility Criterion) and are used to identify ageometric outer automorphisms, the type of outer automorphisms we focus on.

\begin{df} A nontrivial path $\rho$ between fixed points $x,y \in \Gamma$ is called a \emph{periodic Nielsen Path (PNP)} if, for some $k$, $g^k(\rho) \simeq \rho$ rel endpoints.  If $k=1$, then $\rho$ is called a \emph{Nielsen Path (NP)}.  $\rho$ is called an \emph{indivisible Nielsen Path (iNP)} if it cannot be written as a nontrivial concatenation $\rho=\rho_1 \cdot \rho_2$, where $\rho_1$ and $\rho_2$ are NPs.  A particularly nice property of an iNP for an irreducible train track representative [Lemma 3.4, BH97] is that there exist unique, nontrivial, legal paths $\alpha$, $\beta$, and $\tau$ in $\Gamma$ such that $\rho=\bar{\alpha}\beta$, $g(\alpha)= \tau\alpha$, and $g(\beta)=\tau\beta$.

In [BF94], immersed paths $\alpha_1, \dots \alpha_k \in \Gamma$ are said to form an \emph{orbit of periodic Nielsen paths} if $g^k(\alpha_i) \simeq \alpha_{\text{i+1 mod k}}$ rel endpoints, for all $1 \leq i \leq k$.  The orbit is called \emph{indivisible} if $\alpha_1$ is not a concatenation of subpaths belonging to orbits of PNPs.  We call each $\alpha_i$ in an indivisible orbit an \emph{indivisible periodic Nielsen path} (\emph{iPNP}).
\end{df}

\noindent In order to define geometric, ageometric, and parageometric fully irreducible outer automorphisms, we first remind the reader of the following definitions.

\begin{df} Let $CV_r$ denote \emph{Outerspace}, defined in [CV86] to be the set of projective equivalence classes of marked graphs (where the equivalence is up to marking-preserving isometry).  We remind the reader that Outerspace can also be defined in terms of free, simplicial $F_r$-trees up to isometric conjugacy and that elements of the compactification are represented by equivalence classes of actions of $F_r$ on $\textbf{R}$-trees (sometimes called \emph{$F_r$-trees}) that are: \newline
\indent (1) minimal: there exists no proper, nonempty, $F_r$-invariant subtree and \newline
\indent (2)very small: \newline
\indent \indent (a) the stabilizer of every nondegenerate arc is either trivial or a cyclic subgroup generated by a primitive element of $F_r$ and \newline
\indent \indent (b) the stabilizer of every triod is trivial. \newline
\noindent In the tree definition, elements in an equivalence class differ by $F_r$-equivariant bijections that multiply their metrics by a constant.

Let $\phi \in Out(F_r)$ be fully irreducible.  \emph{$T_+$} is defined as the unique point in $\partial CV_r$ which is the attracting point for every forward orbit of $\phi$ in $CV_r$.  The point's uniqueness is proved in [LL03].
\end{df}

It is proved in [BF94, Theorem 3.2] that, for a fully irreducible $\phi \in Out(F_r)$, $T_+$ is a geometric $\textbf{R}$-tree if and only if every TT representative of every positive power of $\phi$ has at least one iPNP.  Recall that a fully irreducible $\phi \in Out(F_r)$ is called \emph{geometric} if it is induced by a homeomorphism of a compact surface with boundary.  A defining characteristic of geometric fully irreducible outer automorphisms is that they have a power with a representative having only a single closed iPNP (and no other iPNPs) [BH92].  In fact, such a $\phi$ can be realized as a pseudo-Anosov homeomorphism of the surface obtained from $\Gamma$ by gluing a boundary component of an annulus along this loop [BFH, Proposition 4.5].  In the remaining circumstances where $T_+$ is a geometric $\textbf{R}$-tree, but $\phi$ is not geometric, every representative of a positive power of $\phi$ has at least one indivisible periodic Nielsen path that is not closed.  This type of outer automorphism was defined by M. Lustig and is called \emph{parageometric} [GJLL98].

In the case where $T_+$ is nongeometric, $\phi$ is called \emph{ageometric}.  In other words, a fully irreducible $\phi \in Out(F_r)$ is ageometric if and only if there exists a representative of a power of $\phi$ having no PNPs (closed or otherwise).

Since our question was answered in the geometric case by the work of H. Masur and J. Smillie, we do not focus on geometric outer automorphisms in this document.  We also ignore the parageometric case and instead focus on ageometric fully irreducible outer automorphisms.

\subsection{Rotationlessness}{\label{SS:Rotationless}}

M. Feighn and M. Handel defined rotationless outer automorphisms and rotationless train track representatives in [FH09].  The following (from [HM11]) is the description of a rotationless train track map that we will use.  A vertex is called \emph{principal} if it is either an endpoint of an iPNP or has at least three periodic directions.

\begin{df} A TT map $g: \Gamma \to \Gamma$ is called \emph{(forward) rotationless} if it satisfies: \newline
\indent (1) every principal vertex is fixed and \newline
\indent (2) every periodic direction at a principal vertex is fixed.
\end{df}

\noindent The property of being rotationless is an outer automorphism invariant and so it suffices to have a definition of a rotationless representative, as above. That is, $\phi$ is rotationless if and only if some (every) RTT representative is rotationless [FH09, Proposition 3.29].

\begin{rk}  An important fact proved in [FH09, Lemma 4.43] is that there exists a $K_r > 0$, depending only on $r$, such that $\phi^{K_r}$ is forward rotationless for all $\phi \in Out(F_r)$.  (Thus all representatives of a given $\phi \in Out(F_r)$ have a rotationless power).
\end{rk}

\subsection{Local Whitehead Graphs, Local Stable Whitehead Graphs, Ideal Whitehead Graphs, and Singularity Indices}{\label{S:IWGs}}

In order to define singularity indices (the weaker outer automorphism invariant), we first give a special case definition for ideal Whitehead graphs (the finer outer automorphism invariant).  It is important to notice that these ideal Whitehead graphs, local Whitehead graphs, and local stable Whitehead graphs given here are as defined in [HM11] differ from the Whitehead graphs mentioned elsewhere in the literature.  As this has been a reoccurring point of confusion, we clarify a difference here.  In general, Whitehead graphs come from looking at the turns taken by immersions of 1-manifolds into graphs.  In our case the 1-manifold is a set of lines, the attracting lamination.  In much of the literature the 1-manifolds are circuits representing conjugacy classes of free group elements.  For example, for the Whitehead graphs referred to in [CV86], the images of edges are viewed as cyclic words.  This is not the case for local Whitehead graphs, local stable Whitehead graphs, or ideal Whitehead graphs, as we define them.

The following set of definitions is taken from [HM11], though it is not their original source.  We start by defining the Whitehead graph variants in a way more user-friendly for the purposes of our document and then give the definitions, involving singular leaves and points in $\partial F_r$, found in [HM11].  The definitions we begin with involve turns taken by a given representative of $\phi \in Out(F_r)$.

\begin{df} Let $v \in \Gamma$ be a vertex of the connected marked graph $\Gamma$ and let $g: \Gamma \to \Gamma$ be a train track representative of $\phi \in Out(F_r)$. Then the \emph{local Whitehead graph} for $g$ at $v$ (denoted \emph{$LW(g; v)$}) has:

(1) a vertex for each direction $d \in \mathcal{D}(v)$ and

(2) an edge connecting vertices corresponding to directions $d_1, d_2 \in \mathcal{D}(v)$ if the turn $\{d_1, d_2 \} \in \mathcal{T}(v)$ is taken by some $g^k(e)$, where $e \in \mathcal{E}(\Gamma)$.

The \emph{local Stable Whitehead graph} for $g$ at $v$, $SW(g; v)$, is the subgraph of $LW(g; v)$ obtained by restricting to precisely the vertices with labels $d \in Per(v)$, i.e. vertices corresponding to periodic directions at $v$. If $\Gamma$ is a rose with vertex $v$, then we denote the single local stable Whitehead graph $SW(g; v)$ by $SW(g)$ and the single local Whitehead graph $LW(g; v)$ by $LW(g)$.

If $g$ has no PNPs (which is the only case we consider in this document), then the ideal \emph{Whitehead graph of $\phi$}, \emph{$IW(\phi)$}, is isomorphic to $\underset{\text{singularities v} \in \Gamma}{\bigsqcup} SW(g;v)$, where a \emph{singularity} for $g$ in $\Gamma$ is a vertex with at least three periodic directions.  In particular, when $\Gamma$ has only one vertex $v$ (and no PNPs), $IW(\phi)=SW(g;v)$.

Let $g: \Gamma \to \Gamma$ be a PNP-free representative of an ageometric, fully irreducible $\phi \in Out(F_r)$.  We take from [HM11] the definition of the index for a singularity $v$ to be $i(g;v)= 1-\frac{k}{2}$, where $k$ is the number of vertices of $SW(g;v)$. The index of $\phi$ is then the sum $i(\phi)=\underset{\text{singularities v} \in \Gamma}{\sum i(g;v)}$.  When $\Gamma$ has only a single vertex $v$, $i(\phi)= i(g;v)$. The \emph{index type} of $\phi$ is the list of indices of the components of $IW(\phi)$, written in increasing order.  Since the index type can be determined by counting the vertices in the components of the ideal Whitehead graph, one can ascertain that the ideal Whitehead graph is indeed a finer invariant than the index type for a fully irreducible outer automorphism.

While we took the definition from [HM11], the index sum of a fully irreducible $\phi \in Out(F_r)$ was studied much before [HM11], in papers including [GJLL98].  The papers written by D. Gaboriau, A. Jaeger, G. Levitt, and M. Lustig take a perspective of studying outer automorphisms via \textbf{R}-trees.  We focus instead on TT representatives.
\end{df}

\begin{ex}{\label{Ex:G(g)}}
Let $g: \Gamma \to \Gamma$, where $\Gamma$ is a rose and $g$, defined by
$$g =
\begin{cases}
a \mapsto abacbaba\bar{c}abacbaba \\
b \mapsto ba\bar{c} \\
c \mapsto c\bar{a}\bar{b}\bar{a}\bar{b}\bar{a}\bar{b}\bar{c}\bar{a}\bar{b}\bar{a}c
\end{cases},
$$
is a train track representative of an ageometric, fully irreducible $\phi \in Out(F_r)$. We will see in Section \ref{Ch:NPIdentification} that $g$ has no PNPs and in Section \ref{Ch:Procedure} that $\phi$ is fully irreducible.

In Section \ref{Ch:LTT} we officially define the lamination train track structure (LTT structure) $G(g)$ for a Type (*) representative $g$.  We will end this example by giving the LTT structure for $g$ as a warm up for the definitions of Section \ref{Ch:LTT}.  Since $G(g)$ will encapsulate the information of $SW(g)$ and $LW(g)$ into a single graph (along with information about the marked graph $\Gamma$), we first determine $SW(g)$ and $LW(g)$.

The periodic (actually fixed) directions for $g$ are $\{a, \bar a, b, c, \bar c \}$. $\bar b$ is not periodic since $Dg(\bar b)=c$, which is a fixed direction, meaning that $Dg^k(\bar b)=c$ for all $k \geq 1$ and thus $Dg^k(\bar{b})$ does NOT equal $\bar{b}$ for any $k \geq 1$.  The vertices for $LW(g)$ are $\{a, \bar a, b, \bar b, c, \bar c \}$ and the vertices of $SW(g)$ are $\{a, \bar a, b, c, \bar c \}$.

The turns taken by the $g^k(E)$ where $E \in \mathcal{E}(\Gamma)$ are $\{a,\bar{b}\}$, $\{\bar{a},\bar{c}\}$, $\{b,\bar{a}\}$, $\{b,\bar{c}\}$, $\{c,\bar{a}\}$, and $\{a, c\}$.  Since $\{a,\bar{b}\}$ contains the nonperiodic direction $\bar{b}$, this turn is not represented by an edge in $SW(g)$, though is represented by an edge in $LW(g)$.  All of the other turns listed are represented by edges in both $SW(g)$ and $LW(g)$.

There will be a vertex in $G(g)$ and $LW(g)$ for each of the directions $a, \bar a, b, \bar b, c, \bar c$.  The vertex in $G(g)$ corresponding to $\bar b$ is red and all others are purple.  There are purple edges in $G(g)$ for each edge in $SW(g)$.  And $G(g)$ has a single red edge for the turn $\{a,\bar{b}\}$ (the turn represented by an edge in $LW(g)$, but not in $SW(g)$.  $G(g)$ is obtained from $LW(g)$ by adding black edges connecting the pairs of vertices $\{a,\bar{a}\}$, $\{b,\bar{b}\}$, and $\{c,\bar{c}\}$ (these black edges correspond precisely to the edges $a, b,$ and $c$ of $\Gamma$).

$SW(g)$, $LW(g)$, and $G(g)$ respectively look like (in these figures, A will be used to denote $\bar{a}$, B will be used to denote $\bar{b}$, and $C$ will be used to denote $\bar{c}$):

\begin{figure}[H]
\centering
\includegraphics{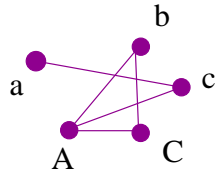}
\caption{{\small{\emph{Stable Whitehead Graph $SW(g)$ for $g$}}}}
\label{fig:SWExample}
\end{figure}

\begin{figure}[H]
\centering
\includegraphics{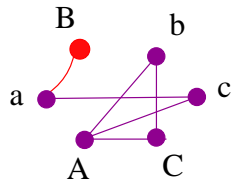}
\caption{{\small{\emph{Local Whitehead Graph $LW(g)$ for $g$}}}}
\label{fig:LWExample}
\end{figure}

\begin{figure}[H]
\centering
\noindent \includegraphics{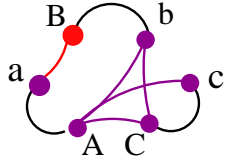}
\caption{{\small{\emph{LTT Structure $G(g)$ for $g$}}}}
\label{fig:LTTExample}
\end{figure}
\end{ex}

\noindent We now relate the definitions of ideal Whitehead graphs, etc, given above to those only relying on the attracting lamination for a fully irreducible outer automorphism.  The purpose will be to show that an ideal Whitehead graph is indeed an outer automorphism invariant.  Each of the following definitions can be found in [HM11].

\begin{df}{\label{D:IWG1}}
A fixed point $x$ is \emph{repelling} for the action of $f$ if it is an attracting fixed point for the action of $f^{-1}$, i.e. if there exists a neighborhood $U$ of $x$ such that, for each neighborhood $V \subset U$ of $x$, there exists an $N>0$ such that $f^{-k}(y) \in V$ for all $y \in U$ and $k \geq N$.  \newline
\indent Let $g: \Gamma \to \Gamma$ be a train track representative of a nongeometric, fully irreducible, rotationless $\phi \in Out(F_r)$ and let $\tilde{g}:\tilde{\Gamma} \to \tilde{\Gamma}$ be a \emph{principal lift} of $g$, i.e. a lift to the universal cover such that the boundary extension has at least three nonrepelling fixed points. We denote the boundary extension of $g$ by $\hat{g}$.  $\tilde{\Lambda}(\phi)$ will denote the lift of the attracting lamination to the universal cover $\tilde{\Gamma}$ of $\Gamma$.  The ideal Whitehead graph, $W(\tilde{g})$, for $\tilde{g}$ is defined to be the graph where: \newline
\indent (1) The vertices correspond to nonrepelling fixed points of $\hat{g}$. \newline
\indent (2) The edges connect vertices corresponding to $P_1$ and $P_2$ precisely when $P_1$ and $P_2$ are the ideal (boundary) endpoints of some leaf in $\tilde{\Lambda}(\phi)$. \newline
\indent We then define the \emph{ideal Whitehead graph} for $g$ by $W(g)=\sqcup W(\tilde{g})$, leaving out components having two or fewer vertices.
\end{df}

\indent Since the attracting lamination is an outer automorphism invariant (and, in particular, the properties of leaves having nonrepelling fixed point endpoints and sharing an endpoint are invariant), the definition we just gave does not rely on the choice of representative $g$ for a given $\phi \in Out(F_r)$.  Thus, once we establish equivalence between this definition and that given at the beginning of this subsection, it should be clear that the ideal Whitehead graph is an outer automorphism invariant.  \newline
\indent Corollary \ref{C:WGfromSW} below is Corollary 3.2 of [HM11].  It relates the definition of an ideal Whitehead graph that we gave above Example \ref{Ex:G(g)} to that given in Definition \ref{D:IWG1}.  \newline
\indent For Corollary \ref{C:WGfromSW} to actually make sense, one needs the following definitions and identification from [HM11].  A \emph{cut vertex} of a graph is a vertex separating a component of the graph into two components.  $SW(\tilde{v}; \tilde{\Gamma})$ denotes the lift of $SW(v;\Gamma)$ to the universal cover $\tilde{\Gamma}$ of $\Gamma$ (having countably many disjoint copies of $SW(v;\Gamma)$, one for each lift of $v$).  \newline
\indent Let $g: \Gamma \to \Gamma$ be an irreducible train track representative of an iterate of $\phi \in Out(F_r)$ such that: \newline
\indent (1) each periodic vertex $v \in \Gamma$ is fixed and \newline
\indent (2) each periodic direction at such a $v$ is fixed. \newline
\noindent Choose one of these fixed vertices $v$.  Suppose $\tilde{v} \in \tilde{\Gamma}$ is a lift of $v$ to the universal cover, $\tilde{g}: \tilde{\Gamma} \to \tilde{\Gamma}$ is a lift of $g$ fixing $\tilde{v}$, and $d$ is a direction at $\tilde{v}$ fixed by $D\tilde{g}$.  Furthermore, let $\tilde{E}$ be the edge at $\tilde{v}$ whose initial direction is $d$.  The \emph{ray determined by $d$} (or \emph{by $\tilde{E}$}) is defined as $\tilde{R}= \bigcup\limits_{j=0}^{j=\infty} \tilde{g}^j(\tilde{E})$.  This is a ray in $\tilde{\Gamma}$ converging to a nonrepelling fixed point for $\hat{g}$.  Such a ray is called \emph{singular} if the vertex $\tilde{v}$ it originates at is principal (i.e. $v$ is principal). With these definitions: \newline
\indent (1) the vertices of $SW(\tilde{v}; \tilde{\Gamma})$ correspond to singular rays $\tilde{R}$ based at $\tilde{v}$ and \newline
\indent (2) directions $d_1$ and $d_2$ represent endpoints of an edge in $SW(\tilde{v}; \tilde{\Gamma})$ if and only if $\tilde{l}= \tilde{R_1} \cup \tilde{R_2}$ is a singular leaf of $\tilde{\Lambda}$ realized in $\tilde{\Gamma}$, where $\tilde{R_1}$ and $\tilde{R_2}$ are the rays determined by $d_1$ and $d_2$ respectively. \newline
\indent Noticing that the ideal (boundary) endpoints of singular rays are precisely the nonrepelling fixed points at infinity for the action of $\tilde{g}$, combining this with what has already been said, as well as Corollary \ref{C:WGfromSW} and what follows, we have the correspondence proving ideal Whitehead graph invariance.

\begin{cor}{\label{C:WGfromSW}} [HM11]
Let $\tilde{g}$ be a principal lift of $g$. Then: \newline
\indent (1) $W(\tilde{g})$ is connected. \newline
\indent (2) $W(\tilde{g}) = \underset{\tilde{v} \in Fix(\tilde{g}) \in \Gamma}{\cup} SW(\tilde{v})$. \newline
\indent (3) For $i \neq j$, $SW(\tilde{v_i})$ and $SW(\tilde{v_j})$ intersect in at most one vertex.  If they do intersect at a vertex $P$, then $P$ is a cut point of $W(\tilde{g})$, in fact $P$ separates $SW(\tilde{v_i})$ and $SW(\tilde{v_j})$ in $W(\tilde{g})$.
\end{cor}

By [Lemma 3.1, HM11], in our case (where there are no PNPs), there is in fact only one $\tilde{v} \in Fix(\tilde{g})$ and so the above corollary gives that $W(\tilde{g}) = SW(\tilde{v})$.

This concludes our justification of how an ideal Whitehead graph is an outer automorphism invariant. Consult [HM11] for clarification of the relationship between ideal Whitehead graphs and \textbf{R}-trees or for other ideal Whitehead graph characterizations.

\subsection{Folds, Decompositions, and Generators}{\label{S:Folds}}

John Stallings introduced ``folds'' in [St83].  Bestvina and Handel use in [BH92] several versions of folds in their construction of TT representatives of irreducible $\Phi \in Out(F_r)$.  We use folds in Section \ref{Ch:IdealDecompositions} for defining and proving ideal decomposition existence.

Let $g: \Gamma \to \Gamma$ be a homotopy equivalence of marked graphs.  Suppose that $g(e_1)=g(e_2)$ as edge paths, where the edges $e_1, e_2 \in \mathcal{E}(\Gamma)$ emanate from a common vertex $v \in \mathcal{V} (\Gamma)$.  One can obtain a graph $\Gamma_1$ by identifying $e_1$ and $e_2$ in such a way that $g:\Gamma \to \Gamma$ projects to $g_1: \Gamma_1 \to \Gamma_1$ under the quotient map induced by the identification of $e_1$ and $e_2$.  $g_1$ is also a homotopy equivalence and one says that $g_1$ and $\Gamma_1$ are obtained from $g: \Gamma \to \Gamma$ by an \emph{elementary fold} of $e_1$ and $e_2$. [St83, BH92]

One can generalize this definition by only requiring that $e_1' \subset e_1$ and $e_2' \subset e_2$ be maximal, initial, nontrivial subsegments of edges emanating from a common vertex $v \in \mathcal{V} (\Gamma)$ such that $g(e_1')=g(e_2')$ as edge paths and such that the terminal endpoints of $e_1$ and $e_2$ are in $g^{-1}(\mathcal{V}(\Gamma))$.  Possibly redefining $\Gamma$ to have vertices at the endpoints of $e_1'$ and $e_2'$, one can fold $e_1'$ and $e_2'$ as $e_1$ and $e_2$ were folded above.  If both $e_1'$ and $e_2'$ are proper subedges then we say that $g_1:\Gamma_1 \to \Gamma_1$ is obtained by a \emph{partial fold} of $e_1$ and $e_2$.  If only one of $e_1'$ and $e_2'$ is a proper subedge (and the other is a full edge), then we call the fold a \emph{proper full fold} of $e_1$ and $e_2$.  In the remaining case where $e_1'$ and $e_2'$ are both full edges, we call the fold an \emph{improper full fold}. [St83, BH92]

Now let $S=<x_1, \dots, x_r>$ be a free basis for the free group $F_r$.  From [N86] we know that any $\Phi \in Aut(F_r)$ can be written as a composition of ``Nielsen generators'' having one of the following two forms (Nielsen gave a longer list, but these suffice):
\begin{itemize}
\item [(1)] $\Phi(x)=xy$ for some $x,y \in S \cup S^{-1}$ (and $\Phi(z)=z$ for all $z \in S \cup S^{-1}$ with $z \neq x^{\pm 1}$)
\item [(2)] a permutation $\sigma$ of $S \cup S^{-1}$ preserving inverses (if $\sigma(x)=y$, then $\sigma(x^{-1})=y^{-1}$).
\end{itemize}

\begin{df}
In general, we will call an automorphism such as $\Phi$ in (1) the \emph{Nielsen generator} (or just \emph{generator}) $x \mapsto xy$.
\end{df}

Consider two metric roses $R_r$ and $R'_r$ with respective edge-labelings $\{a_1, a_2, a_3, \dots \}$ and \newline
\noindent $\{A_1, A_2, A_3, \dots \}$ and markings where the homotopy class of each $a_i$ in $\pi_1(R_r)$ and each $A_i$ in $\pi_1(R_r')$ are identified with the free basis element $x_i$ under the respective markings.  Consider a homotopy equivalence $g: R_r \to R'_r$ that linearly maps $a_i$ over $A_i \cup A_j$ and, for each $k \neq i$, linearly maps $a_k$ over $A_k$. Let $a_i=a_i' \cup a_i''$ where $a_i'$ is mapped by $g$ over $A_i$ and $a_i''$ is mapped by $g$ over $A_j$.  Now consider a quotient map (a ``proper full fold'') $q: R_r \to R^q_r$ identifying $a_i''$ with $a_j$.  There exists a homeomorphism $h: R^q_r \to R'_r$ such that $g=h \circ q$, i.e. $g$ and $h \circ q$ give the same induced map of fundamental groups.  In fact, the homeomorphism linearly maps $a_i'$ over $A_i$ and linearly maps each other $a_k$ over $A_k$.  Sometimes we instead call $g$ the proper full fold.

\begin{figure}[H]
\centering
\noindent \includegraphics[width=5in]{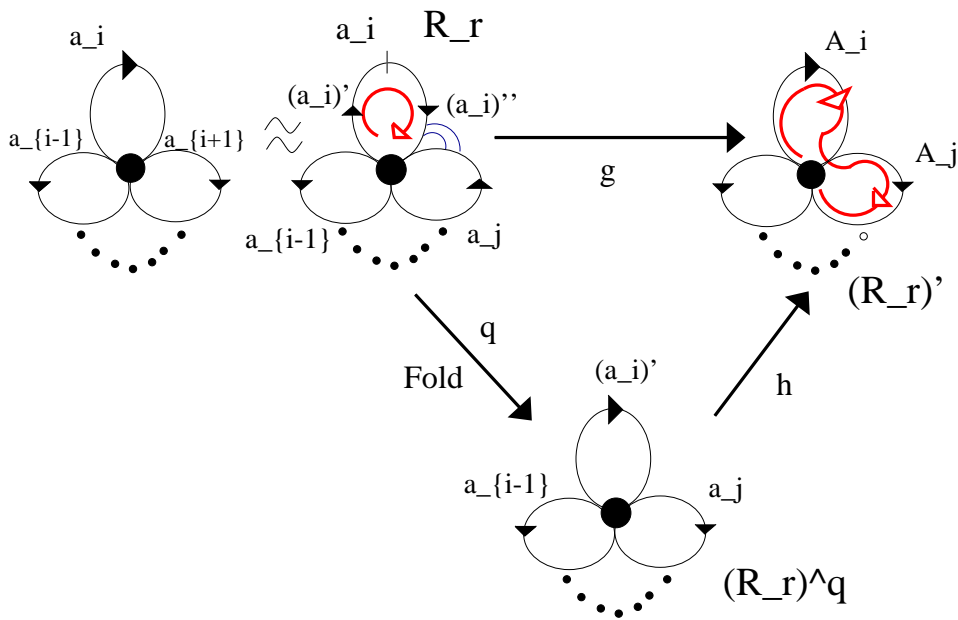} \caption{{\small{\emph{Proper full fold}}}}
\label{fig:ProperFullFold}
\end{figure}

Under the identification, for each $i$, of the homotopy classes of the $a_i$ and $A_i$ with the free basis element $x_i$, the \emph{induced} automorphism $\Phi \in Aut(F_r)$ is the automorphism where $\Phi(x_i)=x_ix_j$ and $\Phi(x_k)=x_k$ for all $k \neq i$. We say that $g$ in the pervious paragraph \emph{corresponds} to the Nielsen generator $x_i \to x_ix_j$. We have similar situations for cases where $g: R_r \to R'_r$ maps $a_i$ linearly over $A_i \cup \overline{A_j}$ and $\Phi(x_i)=x_i(x_j)^{-1}$, or $A_i$ is replaced by its inverse, or both $A_i$ and $A_j$ are replaced by their inverses. \newline
\indent Stallings showed in [St83] that one can decompose a tight (in the sense defined in Subsection \ref{S:Paths} above) homotopy equivalence of graphs as a composition of elementary folds together with a final homeomorphism.  In the circumstance where the elementary folds are proper full folds of roses, the elements of this decomposition have induced Nielsen generators, as described above.\newline
\indent Ideally, the Nielsen generators in a decomposition of $\Phi \in Aut(F_r)$ would all be of form (1) above and there would be a representative $g:\Gamma \to \Gamma$ of $\phi$ where
\begin{itemize}
\item $(\Gamma, \pi: R_r \to \Gamma)$ is a marked rose,
\item $\Phi= \pi^{-1} \circ g \circ \pi$ ($g$ \emph{corresponds to} $\Phi$),
\item the Stallings fold decomposition $\Gamma = \Gamma_0 \xrightarrow{g_1} \Gamma_1 \xrightarrow{g_2} \cdots \xrightarrow{g_{n-1}} \Gamma_{n-1} \xrightarrow{g_n} \Gamma_n = \Gamma$
\emph{corresponds} to the Nielsen generator decomposition $\Phi = \phi_n \circ \phi_{n-1} \circ \cdots \circ \phi_1 \circ \phi_0$ in the sense that $\phi_i= \pi_i^{-1} \circ g_i \circ \pi_i$ (where $\pi_i$ is the marking on $\Gamma_i$) for each $i$, and
\item each $g_i$ is a proper full fold of a rose.
\end{itemize}

\noindent In Section \ref{Ch:IdealDecompositions} we prove that such is the case in the scenario we want it for.

\section{Ideal Decompositions}{\label{Ch:IdealDecompositions}}

As mentioned in Subsection \ref{S:Folds}, TT representatives are composed of elementary folds, followed by a homeomorphism.  The ideal situation would be to have a TT representative $g$ for each $\phi \in Out(F_r)$ composed only of proper full folds of roses and a homeomorphism inducing a trivial permutation.  We would further like for the automorphism $\Phi$ corresponding to $g$ to have several properties.

In this section we show that for an ageometric, fully irreducible $\phi \in Out(F_r)$ such that $IW(\phi)$ is Type (*) pIW we have a representative of a rotationless power decomposed as desired.

For Proposition \ref{P:IdealDecomposition}, we need the following from [HM11]: Let an ageometric, fully irreducible $\phi \in Out(F_r)$ be such that $IW(\phi)$ is a Type (*) pIWG $\mathcal{G}$, then $\phi$ is \emph{rotationless} if and only if the vertices of $IW(\phi)$ are fixed by the action of $\phi$. \newline
\indent We will also need the following lemmas:

\begin{lem}{\label{L:PNPFreePreserved}} Let $g: \Gamma \to \Gamma$ be a PNP-free TT representative of a fully irreducible $\phi \in Out(F_r)$ and let $\Gamma = \Gamma_0 \xrightarrow{g_1} \Gamma_1 \xrightarrow{g_2} \cdots \xrightarrow{g_{n-1}} \Gamma_{n-1} \xrightarrow{g_n} \Gamma_n = \Gamma$ be a decomposition of $g$ into homotopy equivalence of marked metric graphs.  Let $f_k: \Gamma_k \to \Gamma_k$ denote the composition $\Gamma_k \xrightarrow{g_{k+1}} \Gamma_{k+1} \xrightarrow{g_{k+2}} \cdots \xrightarrow{g_{k-1}} \Gamma_{k-1} \xrightarrow{g_k} \Gamma_k$.
Then $f_k$ is also a PNP-free TT representative of $\phi$ (and, in particular, $IW(f_k) \cong IW(g)$).
\end{lem}

\noindent Proof: Suppose that $h=f_k$ had a PNP $\rho$ and let $h^p$ be such that the path is fixed (up to homotopy rel endpoints), i.e. $h^p(\rho) \simeq \rho$ rel endpoints.  Let $\rho_1=g_n \circ g_{k+1}(\rho)$.  First notice that $\rho_1$ cannot be trivial or $h_p(\rho)=(g_k \circ g_1 \circ g^{p-1})(g_n \circ g_{k+1}(\rho))=(g_k \circ g_1 \circ g^{p-1})(\rho_1)$ would be trivial, contradicting that $\rho$ is a PNP.

$g^p(\rho_1)=g^p((g_k \circ g_1)(\rho))=(g_n \circ g_{k+1}) \circ h^p(\rho)$.  Now, $h^p(\rho) \simeq \rho$ rel endpoints and so $(g_n \circ g_{k+1}) \circ h^p(\rho) \simeq (g_n \circ g_{k+1})(\rho)$ rel endpoints (just compose the homotopy with $g_n \circ g_{k+1}$)  But that means that $g^p(\rho_1)=g^p((g_k \circ g_1)(\rho))=(g_n \circ g_{k+1}) \circ h^p(\rho)$ is homotopic to $(g_n \circ g_{k+1})(\rho)=\rho_1$ relative endpoints.  Making $\rho_1$ a PNP for $g$ contradicting that $g$ is PNP-free.  Thus, $h=f_k$ must be PNP-free, as desired.

Let $\pi: R_r \to \Gamma$ be the marking on $\Gamma_1$.  Since $g_1$ is a homotopy equivalence, $g_1 \circ \pi$ gives a marking on $\Gamma$ and $g$ and $h$ just differ by a change of marking.  Thus, $g$ and $h$ are representatives of the same outer automorphism $\phi$. \newline
\noindent QED.

\begin{lem}{\label{L:GateCollapsing}} Let $h: \Gamma \to \Gamma$ be a PNP-free train track representative of a fully irreducible $\phi \in Out(F_r)$ such that $h$ has $2r-1$ fixed directions.  Let \newline
\noindent $\Gamma = \Gamma_0 \xrightarrow{h_1} \Gamma_1 \xrightarrow{h_2} \cdots \xrightarrow{h_{n-1}} \Gamma_{n-1} \xrightarrow{h_n} \Gamma_n = \Gamma$ be the Stallings fold decomposition for $h$.  Let $h^i$ be such that $h=h^i \circ h_i \circ \dots \circ h_1$. Let $d_{(1,1)}, \dots, d_{(1,2r-1)}$ be the fixed directions for $Df$ and let $d_{j,k}=D(h_j \circ \dots \circ h_1)(d_{1,k})$ for each $1 \leq j \leq n$ and $1 \leq k \leq 2r-1$. Then $D(h^i)$ cannot identify any of the directions $d_{(i,1)}, \dots, d_{(i,2r-1)}$.
\end{lem}

\noindent Proof: Let $d_{(1,1)}, \dots, d_{(1,2r-1)}$ be the fixed directions for $Df$ and let $d_{j,k}=D(h_j \circ \dots \circ h_1)(d_{1,k})$ for each $1 \leq j \leq n$ and $1 \leq k \leq 2r-1$. Suppose that $D(h^i)$ identified any of the directions $d_{(i,1)}, \dots, d_{(i,2r-1)}$, then we would have that $Df$ had fewer than $2r-1$ directions in its image, contradicting that it has $2r-1$ fixed directions. \newline
\noindent QED.

\begin{prop}{\label{P:IdealDecomposition}} Let an ageometric, fully irreducible $\phi \in Out(F_r)$ be such that $IW(\phi)$ is a Type (*) pIW graph.  Then there exists a PNP-free, rotationless representative of a power $\psi=\phi^R$ of $\phi$ on the rose.  Further, the representative can be decomposed as a sequence of proper full folds of roses.
\end{prop}

\noindent Proof: Suppose $\phi$ is as described in the proposition.  Since $\phi$ is ageometric, there exists a PNP-free TT representative $g$ of a power of $\phi$.  Let $h=g^k: \Gamma \to \Gamma$ be such that $h$ fixes all of $g$'s periodic directions ($h$ is rotationless).  Then $h$ is also a PNP-free TT representative of some $\phi^l$.  Since $h$ has no PNPs (meaning $IW(\phi^R) \cong \underset{\text{singularities v} \in \Gamma}{\bigsqcup} SW(h;v)$), since $h$ fixes all of its periodic directions (in particular $SW(h;v) \cong LW(h;v)$ for all $V$), and since the ideal Whitehead graph of $\phi$ (hence $\phi^R$) is a Type (*) pIW graph, $\Gamma$ must have a vertex with $2r-1$ fixed directions.  Thus, $\Gamma$ must be one of the only three graphs of rank $r$ with a valence $2r-1$ or higher vertex:

\begin{figure}[H]
\centering
\noindent \includegraphics[width=3.8in]{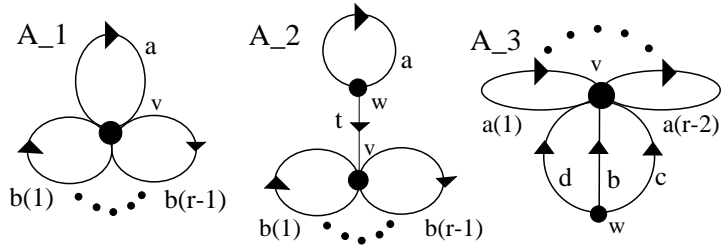}
\caption{{\small{\emph{Graphs of rank $r$ with valence $2r-1$}}}}
\label{fig:HigherRankGraphChoices}
\end{figure}

If $\Gamma=A_1$, then $h$ will be the desired representative.  We will show that, in this case we have the desired decomposition.  However, first we will show that $\Gamma$ cannot be $A_2$ or $A_3$ by ruling out all possibilities for folds in the Stallings fold decomposition for $h$ in the cases of $\Gamma=A_2$ and $\Gamma=A_3$.

If we had $\Gamma=A_2$, then the vertex with $2r-1$ fixed directions could only be $v$.  $h$ must have an illegal turn unless it were a homeomorphism, which it could not be and still be irreducible.  Notice that $w$ could not be mapped to $v$ in a way not forcing an illegal turn at $w$, as this would mean that either we would have an illegal turn at $v$ (if $t$ were wrapped around some $b_i$) or we would have backtracking on $t$ (contradicting that $g$ is a TT map, so must be locally injective on edge interiors).  Because all $2r-1$ of the directions at $v$ are fixed by $h$, if $h$ had an illegal turn, it would have to occur at $w$ (as no two of the $2r-1$ fixed directions can share a gate).

The turns at $w$ are $\{a, \bar{a}\}$, $\{a, t\}$, and $\{\bar{a}, t\}$.  However, we only need to rule out illegal turns at $\{a, \bar{a}\}$ and $\{a, t\}$, as the situations with $\{\bar{a}, t\}$ and $\{a, t\}$ are identical.

First, suppose that the illegal turn at $w$ were $\{a, \bar{a}\}$, so that the first fold in the Stallings decomposition would have to be of $\{a, \bar{a}\}$.  Fold $\{a, \bar{a}\}$ maximally to obtain $(A_2)_1$.  The fold cannot completely collapse $a$, as this would change the homotopy type of $A_2$.  Also, we can assume the fold is maximal or the next fold in the sequence would just anyway be a continuation of the same fold.

\begin{figure}[H]
\centering
\noindent \includegraphics[width=4.5in]{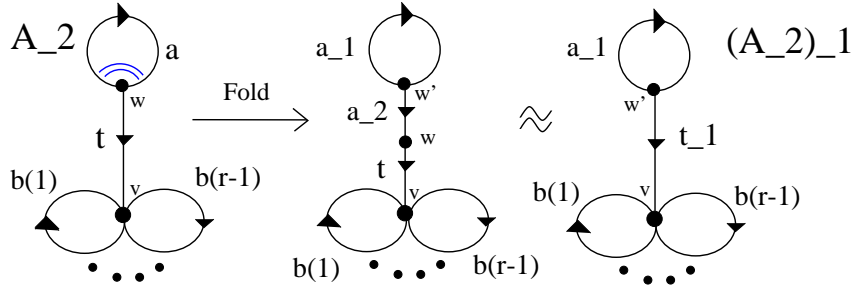}
\caption{{\small{\emph{$a_1$ is the portion of $a$ not folded, $a_2$ is the edge created by the fold, $w'$ is the vertex created by the fold, and $t_1$ is $a_2 \cup t$ without the (now unnecessary) vertex $w$}}}}
\label{fig:RosewithStemPetalFold}
\end{figure}

Let $h_1: (A_2)_1 \to (A_2)_1$ be the induced map as in [BH92] and explained in Subsection \ref{S:Folds} above.  Since the fold of $\{a, \bar{a}\}$ was maximal, $\{a_1, \overline{a_1}\}$ must be legal.  Since $h$ was a TT map, and thus was locally injective on edges of $\Gamma$ (on the edge $a$ in particular), $\{t_1, a_1\}$ and $\{t_1, \overline{a_1}\}$ must also be legal.  But then $h_1$ would fix all directions at both vertices of $\Gamma_1$ (since it still must fix all directions at $v$).  This would make $h_1$ a homeomorphism, again contradicting irreducibility.  Thus, $\{a, \bar{a}\}$ could not have been the illegal turn at $w$.  This leaves us to rule out $\{a, t\}$.

Suppose that the illegal turn at $w$ were $\{a, t\}$, so that the first fold in the Stallings decomposition would have to be of $\{a, t\}$.  Fold $\{a, t\}$ maximally (we can again assume the fold is maximal).  Let $h_1': (A_2)'_1 \to (A_2)'_1$ be the induced map of [BH92] and Subsection \ref{S:Folds}.  Either \newline
\indent \indent \indent I. all of $t$ is folded with all of $a$ or a power of $a$; \newline
\indent \indent \indent II. all of $t$ is folded with part of $a$ or a power of $a$; \newline
\indent \indent \indent III. part of $t$ is folded with all of $a$ or a power of $a$; or \newline
\indent \indent \indent IV. part of $t$ is folded with part of $a$ or a power of $a$.

If (I) or (II) held, $(A_2)_1'$ would be a rose and $h_1'$ would give a representative on the rose, returning us to the case of $A_1$.  So we just need to analyze (III) and (IV).

Consider first (III), i.e. suppose part of $t$ is folded with part of $a$ or a power of $a$:

\begin{figure}[H]
\centering
\noindent \includegraphics[width=3.7in]{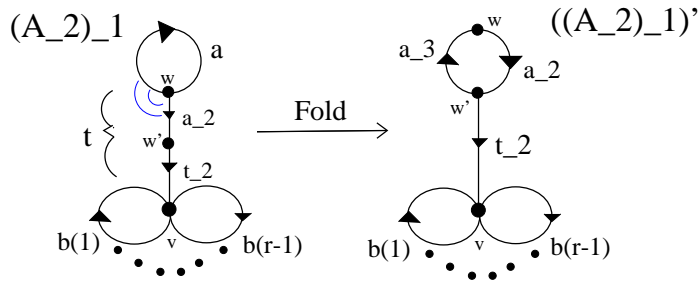}
\caption{{\small{\emph{$a_2$ is the portion of $a$ folded with the portion of $t$, $a_3$ is the portion of $a$ not folded with $t$, and $t_2$ is the portion of $t$ not folded with $a$}}}}
\label{fig:RoseWithStemChoices2}
\end{figure}

Let $h^1$ be such that $h=h^1 \circ g_1$ with $g_1$ being the single fold performed thus far.  $h^1$ could not identify any of the directions at $w'$: $h^1$ could not identify $a_2$ and $t_2$ or $h$ would have had back-tracking on $t$; $h^1$ could not identify $a_2$ and $\overline{a_3}$ or $h$ would have had back-tracking on $a$; and $h^1$ could not identify $t_2$ and $\overline{a_3}$ because the fold was maximal.  But then all the directions of $(A_2)_{1}'$ would be fixed by $h^1$, making $h^1$ a homeomorphism and the decomposition complete.  However, this would make $h$ consist of the single fold $g_1$ and a homeomorphism, contradicting that $h$'s irreducibility.   So, in analyzing what might happen if $\{a, t\}$ were the illegal turn for $h$ at $w$, we are left to analyze (IV):

\begin{figure}[H]
\centering
\noindent \includegraphics[width=3.5in]{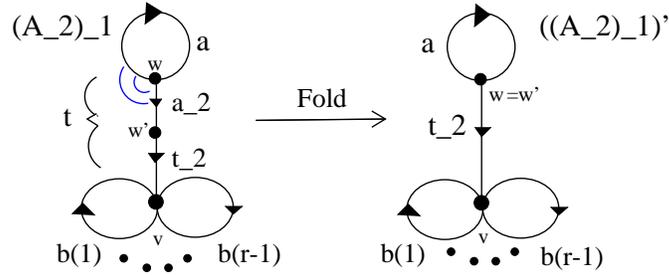}
\caption{{\small{\emph{$t_2$ is the portion of $t$ that was not folded with $a$}}}}
\label{fig:RoseWithStemNonfold}
\end{figure}

Suppose part of $t$ were folded with all of $a$ (or a power of $a$).  Again let $h^1$ be such that $h=h^1 \circ g_1$, where $g_1$ is the single fold performed thus far.  If $h^1$ did not identify any of the directions at $w'$, it would be a homeomorphism, causing the same contradiction with $h$'s irreducibility as above.  But $h^1$ could not identify any of the directions at $w'$: $h^1$ could not identify $a$ and $\bar{a}$ for the same reasons as above; $h^1$ could not identify $a$ and $t_2$, as the fold of $a$ and $t$ was maximal; and $h^1$ cannot identify $t_2$ and $\bar{a}$ of $h$ would have had backtracking on $t$.  We have thus shown that all cases where $\Gamma = A_2$ are either impossible or yield the desired representative on the rose.

We are left to analyze when $\Gamma=A_3$.  In this case, $v$ must be the vertex with $2r-1$ fixed directions.  As with $A_2$, because $h$ must fix all $2r-1$ directions at $v$, if $h$ had an illegal turn (which it still has to) the turn would be at $w$.   Without loss of generality assume the illegal turn is $\{b, d\}$.  Maximally fold $\{b, d\}$.  If all of $b$ and $d$ were folded, this would change the homotopy type.  Thus also assume (again without losing generality) that either 1. all of $b$ is folded with part of $d$ or 2. only proper initial segments of $b$ and $d$ are folded with each other.  If (1) holds, we get a PNP-free TT representative on the rose.  So suppose (2) holds.  Let $h_1: (A_3)_1 \to (A_3)_1$ be the induced map of [BH92].  The fold and $(A_3)_1$ look like:

\begin{figure}[H]
\centering
\noindent \includegraphics[width=5in]{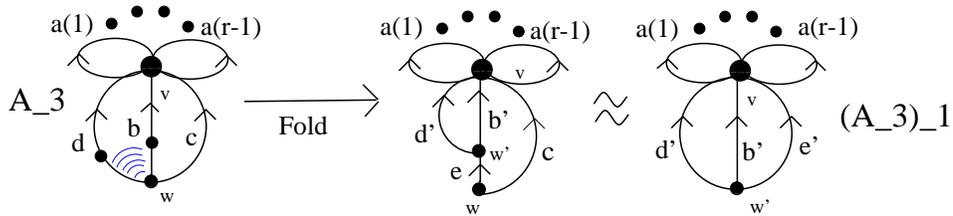}
\caption{{\small{\emph{$e$ is the edge created by the fold and $e'$ is $\bar{e} \cup c$ without the (now unnecessary) vertex $w$}}}}
\label{fig:TheFold}
\end{figure}

The new vertex $w'$ has 3 distinct gates ($\{b', d'\}$ is legal since the fold was maximal and $\{b', \bar{e}\}$ and $\{d', \bar{e}\}$ must also be legal or $h$ would have had back-tracking on $b$ or $d$, respectively).  This leaves the situation where the entire decomposition is a single fold with a homeomorphism, again leading to the contradiction of $h$ being reducible.

Having ruled out all cases, we have completed the analysis of $A_3$ and thus proved that we have a PNP-free representative of a power $\psi=\phi^R$ of $\phi$ on the rose.

Let $h: \Gamma \to \Gamma$ be the PNP-free train track representative of $\phi^R$ on the rose.  Let \newline
\noindent $\Gamma = \Gamma_0 \xrightarrow{g_1} \Gamma_1 \xrightarrow{g_2} \cdots \xrightarrow{g_{n-1}} \Gamma_{n-1} \xrightarrow{g_n} \Gamma_n = \Gamma$ be the Stallings decomposition for $h$.  Each $g_i$ is either an elementary fold or locally injective (in which case it would be a homeomorphism).  We can assume that $g_n$ is the only homeomorphism.  Let $h^i=g_n \circ \dots \circ g_{i+1}$.  Since $h$ has precisely $2r-1$ gates, $h$ has precisely one illegal turn.  We first determine what $g_1$ could be.  $g_1$ cannot be a homeomorphism or we would have $h=g_1$, making $h$ reducible.  Thus, $\Gamma$'s vertex contains an illegal turn for $h$.  We maximally fold the illegal turn.  Suppose first that the fold is a proper full fold.  (If the fold is not a proper full fold, then see the analysis below about what would happen with an improper or partial fold.)

\begin{figure}[H]
\centering
\noindent \includegraphics[width=4.5in]{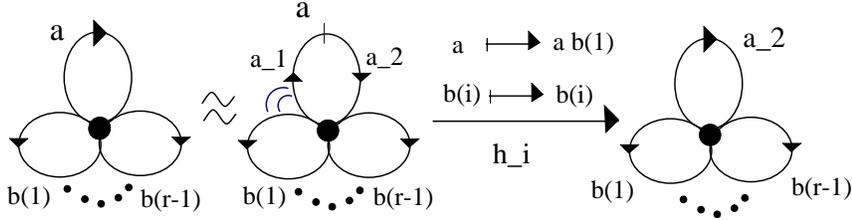} \caption{{\small{\emph{Proper Full Fold}}}}
\label{fig:HigherRankRoseProperFullFold}
\end{figure}

By Lemma \ref{L:GateCollapsing}, $h^1$ cannot have more than one turn $\{d_1, d_2\}$ such that $Dh^1(\{d_1, d_2\})$ is degenerate (we will call such a turn an \emph{order one illegal turn} for $h^1$).  If it has no order one illegal turn, then $h^1$ must be a homeomorphism and we have determined the entire decomposition.  So suppose that $h^1$ has an order one illegal turn and maximally fold this illegal turn.  Continue as such until either $h$ is obtained, in which case the desired decomposition has been found, or until the next fold is not a proper full fold.

The next fold cannot be an improper full fold because this would change the homotopy type of the rose.  So suppose, without loss of generality, that this first fold other than a proper full fold is a partial fold of $b'$ and $c'$.

\begin{figure}[H]
\centering
\noindent \includegraphics[width=3.6in]{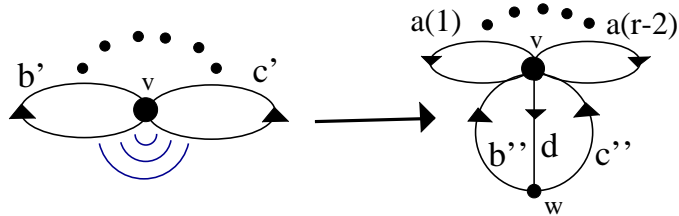}
\caption{{\small{\emph{$d$ is the edge created by folding initial segments of $b'$ and $c'$, $b''$ is the terminal segment of $b'$ not folded, and $c''$ is the terminal segment of $c'$ not folded}}}}
\label{fig:ImproperFold}
\end{figure}

As in the analysis of the case of $\Gamma=A_3$ above, the next fold has to be at $w$ or the next generator would be a homeomorphism, which does not make sense since $A_3$ is not a rose and the image of $h$ is a rose.  Since the previous fold was maximal, the illegal turn cannot be $\{b'', c''\}$.  The illegal turn also cannot be $\{b'', \bar{d}\}$ or $\{c'', \bar{d}\}$, since this would imply $h$ had backtracking on edges, contradicting that $h$ is a train track.  Thus, $h_i$ was not possible in the first place, meaning that all of the folds in the Stallings decomposition must be proper full folds between roses.

Since all folds in the Stallings decomposition are proper full folds of roses, it is possible to index the edge sets $\mathcal{E}_k = \mathcal{E}(\Gamma_k)$ as \newline
\noindent $\{E_{(k,1)},\overline{E_{(k,1)}}, E_{(k,2)}, \overline{E_{(k,2)}}, \dots , E_{(k,r)}, \overline{E_{(k,r)}} \} = \{e_{(k,1)}, e_{(k,2)}, \dots, e_{(k,2r-1)}, e_{(k,2r)}\}$ so that \newline
\indent(a) $g_k: e_{k-1,j_k} \mapsto e_{k,i_k} e_{k,j_k}$ where $e_{k-1,j_k} \in \mathcal{E}_{k-1}$, $e_{k,i_k}, e_{k,j_k} \in \mathcal{E}_k$, and \newline
\indent (b) $g_k(e_{k-1,i})=e_{k,i}$ for all $e_{k-1,i} \neq e_{k-1,j_k}^{\pm 1}$. \newline
\noindent Suppose we similarly index the directions $D(e_{k,i}) = d_{k,i}$.

Let $g_n=h'$ be the homeomorphism in the Stalling's decomposition and suppose that $Dh'$ gave a nontrivial permutation of the second indices of the directions.  Some power $p$ of the permutation would be trivial.  Replace $h$ by $h^p$. We rewrite the decomposition of $h^p$ as follows.  Let $\sigma$ denote the permutation of second indices defined by $Dh'$.  Then, for $n \leq k \leq 2n-p$ define $g_k$ by $g_k: e_{k-1,\sigma^{-s+1}(j_t)} \mapsto e_{k,\sigma^{-s+1}(i_t)} e_{k,\sigma^{-s+1}(j_t)}$ where $k=sp+t$ and $0 \leq t \leq p$.  Adjust the corresponding proper full folds accordingly.  This decomposition still gives $h$, but now the permutation caused by the final homeomorphism is trivial.

This concludes the proof that the PNP-free TT representative of a power of $\phi$ on the rose can be decomposed as a sequence of proper full folds between roses and thus the proof of the proposition. \newline
\noindent QED

\vskip8pt

Representatives such as that given in Proposition \ref{P:IdealDecomposition} will be called \emph{ideally decomposable} and the decomposition an \emph{ideal decomposition}.  We establish here notation used for discussing ideally decomposed representatives and restate the results in that notation.  That (4) and (3a) hold (i.e that $g_n$ is not a homeomorphism permuting the direction second indices) can be ascertained from the Proposition \ref{P:IdealDecomposition} proof.

\vskip5pt

\noindent \textbf{Ideal Decomposition Standard Statement and Notation:} \newline
\emph{Let $\phi \in Out(F_r)$ be an ageometric, fully irreducible outer automorphism such that $IW(\phi)$ has $2r-1$ vertices.  Then Proposition \ref{P:IdealDecomposition} gives a PNP-free representative $g:\Gamma \to \Gamma$ of a rotationless power $\psi=\phi^R$ of $\phi$ with a decomposition \newline
$\Gamma = \Gamma_0 \xrightarrow{g_1} \Gamma_1 \xrightarrow{g_2} \cdots \xrightarrow{g_{n-1}}
\Gamma_{n-1} \xrightarrow{g_n} \Gamma_n = \Gamma$, where: \newline
\noindent (1) the index set $\{1, \dots, n \}$ is viewed as the set $\mathbf {Z}$/$n \mathbf {Z}$ with its natural cyclic ordering; \newline
\noindent (2) each $\Gamma_k$ is a rose; \newline
\noindent (3) we can index the sets $\mathcal{E}_k^+ =\mathcal{E}^+(\Gamma_k)$ as \newline
\indent $\{E_{(k,1)}, E_{(k,2)}, \dots , E_{(k,r)} \}$ and $\mathcal{E}_k =\mathcal{E}(\Gamma_k)$ as \newline
\indent $\{E_{(k,1)}, \overline{E_{(k,1)}}, E_{(k,2)}, \overline{E_{(k,2)}}, \dots , E_{(k,r)}, \overline{E_{(k,r)}} \} = \{e_{(k,1)}, e_{(k,2)}, \dots, e_{(k,2r-1)}, e_{(k,2r)}\}$ \newline
\noindent ($e_{k, 2j-1}=E_{k,j}$ and $e_{k, 2j}=\overline{E_{k,j}}$ for each $E_{k,j} \in \mathcal{E}^+_k$) so that} \newline
\indent \emph{(a) $g_k: e_{k-1,j_k} \mapsto e_{k,i_k} e_{k,j_k}$ where $e_{k-1,j_k} \in \mathcal{E}_{k-1}$, $e_{k,i_k}, e_{k,j_k} \in \mathcal{E}_k$, and $e_{k,i_k} \neq (e_{k,j_k})^{\pm 1}$ and} \newline
\indent \emph{(b) $g_k(e_{k-1,i})=e_{k,i}$ for all $e_{k-1,i} \neq e_{k-1,j_k}^{\pm 1}$ \newline
\noindent (we denote $e_{k-1,j_k}$ by $e^{pu}_{k-1}$, $e_{k,j_k}$ by $e^u_k$, $e_{k,i_k}$ by $e^a_k$, and $e_{k-1,i_{k-1}}$ by $e^{pa}_{k-1}$); and} \newline
\emph{\noindent (4) for all $e \in \mathcal{E}(\Gamma)$ such that $e \neq e^u_n$, we have $Dg(d)=d$, where $d=D_0(e)$.  ($g$ fixes every direction except for $D_0(e^u_n)=d^u_n$).} \newline

\noindent Additionally,
\begin{itemize}
\item $\mathcal{D}_k$ will denote the set of directions corresponding to $\mathcal{E}_k$.
\item $f_k= g_k \circ \cdots \circ g_1 \circ g_n \circ  \cdots \circ g_{k+1}: \Gamma_k \to \Gamma_k$.
\item $g_{k,i}=g_k \circ \cdots \circ g_i: \Gamma_{i-1} \to \Gamma_k$ if $k>i$ and $g_{k,i}=g_k \circ \cdots \circ g_1 \circ g_n \circ  \cdots \circ g_i$ if $k<i$.
\item $d^u_k$ will denote $D_0(e^u_k)$, which will sometimes be called the \emph{unachieved direction} for $g_k$ because it is not in the image of $Dg_k$ (the ``$u$'' in $d^u_k$ is for ``unachieved'').
\item $d^a_k$ will denote $D_0(e^a_k)$ and sometimes be called the \emph{twice-achieved direction} for $g_k$, as it is the image of both $d^{pu}_{k-1}$ ($=D_0(e_{k-1,j_k})$) and $d^{pa}_{k-1}$ ($=D_0(e_{k-1,i_k})$) under $Dg_k$ (the ``$a$'' in $d^a_k$ is for ``(twice) achieved'' and the ``$p$'' in $d^{pu}_{k-1}$ and $d^{pa}_{k-1}$ is for ``pre'').
\item $G_k$ will denote $G(f_k)$
\item $G_{k,l}$ will denote the subgraph of $G_l$ containing
{\begin{itemize}
\item[(1)] \emph{all black edges and vertices (given the same colors and labels as in $G_l$) and}
\item[(2)] \emph{all colored edges representing turns in $g_{k,l}(e)$ taken by $e \in \mathcal{E}_{k-1}$.}
\end{itemize}}
\item If we additionally require that $\phi \in Out(F_r)$ is ageometric and fully irreducible and that $IW(\phi)$ is a Type (*) pIW graph, then we will say $g$ is of \emph{Type (*)}. (By saying $g$ is of Type (*), it will be implicit that, not only is $\phi$ ageometric and fully irreducible, but $\phi$ is ideally decomposed, or at least ideally decomposable).
\end{itemize}

\begin{rk} We refer to $E_{k,i}$ as $E_i$ for all $k$ in circumstances where we believe it will not cause confusion.  In these circumstances we may also refer to $\Gamma_k$ as $\Gamma$.

While we may abuse notation by writing $E_i$ instead of $E_{(j,i)}$, unless otherwise specified, $e_i$ will always denote an element of $\mathcal{E}(\Gamma)$ (or $\mathcal{E}(\Gamma_k)$ when specified), where the index of $e_i$ will not necessarily match the index of the corresponding element of $\mathcal{E}(\Gamma)$ (or $\mathcal{E}(\Gamma_k)$).  $d_i$ will still denote $D_0(e_i)$.

Additionally, for reasons of typographical clarity we sometimes put parantheses around the subscripts.
\end{rk}

\section{Lamination Train Track (LTT) Structures}{\label{Ch:LTT}}

This section contains our definitions for several different abstract and specific notions of ``lamination train track (LTT) structures.''  M. Bestvina, M. Feighn, and M. Handel discussed in their papers slightly different notions of train track structures than the notions we describe here.  However, those we describe in this section contain as smooth paths realizations of leaves of the attracting lamination for the outer automorphism. This fact makes them useful for ruling out the achievability of particular ideal Whitehead graphs and for constructing the particular representatives we seek.  The need for all of the properties included in the LTT definitions should become clear in Section \ref{Ch:AMProperties} when we prove the necessity of the ``Admissible Map Properties.''

\subsection{Abstract Lamination Train Track Structures}

\begin{df} A \emph{smooth train track graph} is a finite graph $G$ satisfying: \begin{description}
\item[STTG1:] $G$ has no valence-1 vertices;
\item[STTG2:] each edge of $G$ has 2 distinct vertices (single edges are never loops); and
\item[STTG3:] the set of edges of $G$ can be partitioned into two subsets, $\mathcal{E}_b$ (the ``black'' edges) and $\mathcal{E}_c$ (the ``colored'' edges), such that each vertex is incident to at least one $E_b \in \mathcal{E}_b$ and at least one $E_c \in \mathcal{E}_c$.
\end{description}

\noindent Two train track graphs will be considered \emph{equivalent} if they are isomorphic as graphs.
\end{df}

\begin{df}
\noindent A path in a smooth train track graph is \emph{smooth} if no two consecutive edges of the path are in $\mathcal{E}_b$ and no two consecutive edges of the path are in $\mathcal{E}_c$.
\end{df}

\vskip15pt

\noindent We now give our first abstract notion of a lamination train track (LTT) structure.

\begin{df} A \emph{Lamination Train Track (LTT) Structure $G$} is a smooth colored train track graph (black edges will be included but not considered colored) satisfying each of the following:
\begin{description}
\item[LTT1:] Vertices are either purple or red.
\item[LTT2:] There are an even number of vertices and they are labeled via a one-to-one correspondence with a set $\{d_1, \dots, d_k, \overline{d_1}, \dots, \overline{d_k}\}$.  (In the case of an LTT structure for a Type (*) representative $g: \Gamma \to \Gamma$, the $d_i$ and $\overline{d_i}$ will be the initial and terminal directions of the edges $e_i$ of $\Gamma$).
\item[LTT3:] No pair of vertices is connected by two distinct colored edges.
\item[LTT4:]Edges of $G$ are of the following 3 types:
\indent \indent {\begin{description}
\item[(Black Edges):] A single black edge connects each pair of vertices of the form $\{d_i, \overline{d_i} \}$.  There are no other black edges. In particular, there is precisely one black edge containing each vertex.
\item[(Red Edges):] A colored edge is red if and only if at least one of its endpoint vertices is red.
\item[(Purple Edges):] A colored edge is purple if and only if both endpoint vertices are purple.
    \end{description}}
\item[LTT5:] The partition of the set of edges of $G$, where $\mathcal{E}_b$ is the set of black edges of $G$ and $\mathcal{E}_c$ is the set of colored edges of $G$, satisfies (STTG3).
    \end{description}
\noindent $G_A$ is an \emph{augmented} LTT structure with \emph{legal structure} $G$ if it is obtained from $G$ by adding some green edges to the colored edges of $G$ and the green edges satisfy: \begin{description}
\item[LTT6:] At least one vertex of each green edge is red.
\item[LTT7:] A green edge and a nongreen edge never connect the same vertex pair.
\end{description}
\end{df}

\begin{df}
\noindent  Two LTT structures differing by an ornamentation-preserving (label and color preserving), vertex-preserving homeomorphism will be considered \emph{equivalent}.
\end{df}

\vskip10pt

\noindent \textbf{Standard LTT Structure Notation and Terminology:} In the context of an LTT Structure $G$:
\begin{itemize}
\item An edge connecting a vertex pair $\{d_i, d_j \}$ will be denoted [$d_i, d_j$].
\item The interior of [$d_i, d_j$] will be denoted ($d_i, d_j$).  \newline
\indent (While the notation [$d_i, d_j$] may be ambiguous when there is more than one edge connecting the vertex pair $\{d_i, d_j \}$, we will be clear in such cases as to which edge we refer to.)
\item $[e_i]$ will denote [$d_i, \overline{d_i}$]
\item Red vertices will be called \emph{nonperiodic (direction) vertices}.
\item Red edges will be called \emph{nonperiodic (turn) edges}.
\item Purple vertices will be called \emph{periodic (direction) vertices}.
\item Purple edges will be called \emph{periodic (turn) edges}.
\item The purple subgraph of an LTT structure $G$ will be called the \emph{potential ideal Whitehead graph associated to $G$} and will be denoted \emph{$PI(G)$}.  For a finite graph $\mathcal{G} \cong PI(G)$, we will say that $G$ \emph{is an LTT Structure for $\mathcal{G}$}.
\item \emph{$C(G)$} will denote the colored subgraph of the LTT structure $G$ and will be called the \emph{colored subgraph associated to} (or \emph{of}) $G$.
\item We say that the LTT structure $G$ is \emph{admissible} if $G$ is additionally \emph{birecurrent} as a graph, i.e. if there exists a smooth line in $G$ such that each edge of $G$ occurs infinitely often in each end of the line.
\item For an augmented LTT structure $G_A$ with legal structure $G$, we denote the set of green edges of $G_A$ by $\mathcal{E}_g(G_A)$ and call elements of $\mathcal{E}_g(G_A)$ \emph{green illegal turn edges} (or say that they \emph{correspond to illegal turns}).
\end{itemize}

\subsubsection{Type (*) LTT Structures for Type (*) pIWGs}

The following specialized abstract LTT structure is tailored for the case of a fully irreducible, ageometric $\phi \in Out(F_r)$ such that $IW(\phi) \cong \mathcal{G}$ is a Type (*) pIWG.  For this definition, a (potential) ideal Whitehead graph must be designated, but the structure does not use or record any other information about $\phi \in Out(F_r)$.

\begin{df} A \emph{Type (*) Lamination Train Track Structure} is an LTT structure a Type (*) pIW graph $G$ for $\mathcal{G}$ such that:
\begin{description}
\item[LTT(*)1:] $G$ has only a single red vertex (all other vertices are purple).
\item[LTT(*)2:] $G$ has a unique red edge.
\item[LTT(*)3:] $PI(G) \cong \mathcal{G}$.
\end{description}
\end{df}

\vskip10pt

\noindent \textbf{Standard Type (*) Notation and Terminology:} In the context of a Type (*) LTT structure $G$ for $\mathcal{G}$:
\begin{itemize}
\item The label on the unique red vertex will sometimes be written $d^u$ and called the \emph{unachieved direction}.
\item The unique red edge is denoted $e^R$, or [$t^R$], and the label on its purple vertex is denoted $\overline{d^a}$.
\item $\overline{d^a}$ is contained in a unique black edge, which we call the \emph{twice-achieved edge}.
\item The other twice-achieved edge vertex will be labeled by $d^a$ and called the \emph{twice-achieved direction}.
\item For a Type (*) LTT structure to be \emph{admissible}, we will require that it is admissible as an LTT structure and, in particular:
\end{itemize}
\begin{description}
\item[LTT(*)4:] $PI(G) \cup [t^R]$ has no valence-1 vertices contained in purple or red edges of the form [$d, \bar{d}$].
\end{description}

\subsubsection{Based Lamination Train Track Structures}

Instead of relying on information about $IW(\phi)$ for $\phi \in Out(F_r)$, the following LTT structure focuses on the graph $\Gamma$ for a representative $g: \Gamma \to \Gamma$ of $\phi$.

\begin{df} {\label{D:BasedLTTStructure}} Let $\Gamma$ be a connected marked graph with no valence-one vertices.  An \emph{LTT Structure} with \emph{Base Graph} $\Gamma$ is an LTT structure $G$ such that:
\begin{description}
\item[LTT(Based)1:] For each vertex $v \in \Gamma$ and direction $d \in \mathcal{D}(v)$, there exists a vertex in $G$ labeled by $d$.  In particular, the vertices of $G$ are in one-to-one correspondence with $\mathcal{D}(\Gamma)$
\item[LTT(Based)2:] Edges of $G$ are of the following 3 types:
\indent \indent {\begin{description}
\item[(Purple Edges)] connect the purple vertices in $G$ corresponding to certain distinct pairs $\{d_1, d_2 \}$ of directions at a common vertex $v$ of $\Gamma$ (we will call such pairs of directions \emph{periodic turns});
\item [(Red Edges)]  connect the vertices in $G$ corresponding to certain distinct pairs $\{d_1, d_2 \}$ of directions at a common vertex $v$ of $\Gamma$ such that at least one direction in the pair is represented by a red vertex in $G$.  (We will call such pairs of directions \emph{nonperiodic turns});
\item [(Black Edges)]  connect precisely vertex pairs $\{D_0(e_i), D_0(\overline{e_i}) \}$ where $e_i \in \mathcal{E}(\Gamma)$.
\end{description}}
\end{description}
\end{df}

\vskip8pt

\noindent \textbf{Based LTT Structure Terminology:} In the context of Definition \ref{D:BasedLTTStructure}, for a given vertex $v \in \Gamma$:
\begin{itemize}
\item We call the union of the purple edges $[d_1, d_2]$, where $d_1, d_2 \in \mathcal{D}(v)$, the \emph{stable Whitehead graph $SW(v, \Gamma, G)$ at $v$}.
\item We call the union of the purple and red edges [$d_1, d_2$] corresponding to turns $\{d_1, d_2 \}$, where $d_1, d_2 \in \mathcal{D}(v)$, the \emph{local Whitehead graph $LW(v, \Gamma, G)$ at $v$}.
\end{itemize}

\vskip10pt

\begin{df} Let $\Gamma$ be an r-petaled rose with vertex $v$.  A \emph{Type (*) LTT Structure $G$} with \emph{base graph} $\Gamma$ is an LTT structure with base graph $\Gamma$ additionally satisfying:
\begin{description}
\item[LTT(*)(Based)1:] $SW(v, \Gamma, G)$ is a Type (*) pIW graph (and is actually PI($G$)).
\end{description}
An \emph{Augmented Type (*) Lamination Train Track Structure for $G$ with Base Graph $\Gamma$} is an augmented LTT structure $G_A$ with legal structure $G$ additionally satisfying:
\begin{description}
\item[LTT(*)(Based)2:] $\overline{G_A-\mathcal{E}_G}$ is a Type (*) LTT structure with base graph $\Gamma$.
\item[LTT(*)(Based)3:] $\mathcal{E}_G$ contains only a single edge, which we denote $T(G)$, or just $T$, and call the $\emph{green illegal turn edge}$ of $G$ or \emph{edge corresponding to the illegal turn};
\end{description}
\end{df}

\vskip5pt

\noindent We describe here what it means for two based LTT structures to be equivalent.

\begin{df}{\label{D:EquivalentBasedLTTStructures}}
Suppose $G$ and $G'$ are LTT structures with respective base graphs $\Gamma$ and $\Gamma'$.  A homeomorphism $H: \Gamma_i \to \Gamma_i'$ \emph{extends} to an ornamentation-preserving homeomorphism $H^T: G_i \to G_i'$ if, for each edge $e \in \mathcal{E}(\Gamma)$ there exists a homeomorphism $i_e: \tni(e) \to \tni([e])$ and for each edge $e'=H(e) \in \mathcal{E}(\Gamma)$ there exists a homeomorphism $i_e': H(\tni(e)) \to H(\tni([e]))$ such that the following commutes:

{\large{\[ \begin{CD}
\tni([e]) @>{H^T_{\tni([e])}}>>  G_i'\\
@A{i_e}AA           @AA{i_e'}A \\
\tni(e) @>{H_{\tni(e)}}>>  H(\tni(e))
\end{CD} \]}}

\noindent One would also say in such a circumstance that $H^T$ \emph{restricts} to $H$.

LTT structures $G$ and $G'$ with respective bases $\Gamma$ and $\Gamma'$ are equivalent if there exists an ornamentation-preserving homeomorphism $H^T: G_i \to G_i'$ of based LTT structures restricting to a label-preserving homeomorphism $H: \Gamma \to \Gamma'$ of marked graphs.
\end{df}

\subsubsection{Maps of Based Lamination Train Track Structures}

\indent Let $G$ and $G'$ be LTT structures with respective base graphs $\Gamma$ and $\Gamma'$.  Let $g:\Gamma \to \Gamma'$ be a homotopy equivalence.  Recall that $Dg$ induces a map of turns $D^tg: \{a,b\} \mapsto \{Dg(a), Dg(b)\}$. $Dg$ additionally induces a map on the corresponding edges of $C(G)$ and $C(G')$ (if the appropriate edges exist in $C(G')$):

\begin{df}  Let $G$ and $G'$ be LTT structures with respective base graphs $\Gamma$ and $\Gamma'$.  We say $D^C(g): C(G) \to C(G')$ is \emph{a map of colored subgraphs induced by $g$} if:
\begin{itemize}
\item [1.] $D^C(g)$ sends the vertex labeled $d$ in $G$ to that labeled by $Dg(d)$ in $G'$;
\item [2.] $D^C(g)$ sends the edge [$d_i, d_j$] in $C(G)$ to the edge [$Dg(d_i), Dg(d_j)$] in $C(G')$;
\item [3.] $D^C(g)$ maps each $LW(\Gamma, v)$ into $LW(\Gamma', g(v))$;
\item [4.] $D^C(g)$ maps each $SW(\Gamma, v)$ isomorphically onto $SW(\Gamma', g(v))$
\end{itemize}
\end{df}

\noindent We now describe what it means for a continuous map of marked graphs to extend to a map of based LTT structures.

\begin{df}{\label{D:GeneratorExtendsToLTTStructures}}
Suppose that $G$ and $G'$ are LTT structures with respective base graphs $\Gamma$ and $\Gamma'$.  Let $g: \Gamma \to \Gamma'$ be a homotopy equivalence such that $g(\mathcal{V}(\Gamma)) \subset \mathcal{V}(\Gamma')$. A map $D^T(g): G \to G'$ \emph{induced by $g$} is an extension of $D^C(g): C(G) \to C(G')$ taking the interior of the black edge of $G$ corresponding to the edge $E \in \mathcal{E}(\Gamma)$ to the interior of the smooth path in $G'$ corresponding to $g(E)$.  In this case we say that $g^T=D^T(g): G \to G'$ is the \emph{extension} of $g: \Gamma \to \Gamma'$ to the continuous map of based LTT structures and that $g: \Gamma \to \Gamma'$ \emph{extends} to $D^T(g)$.
\end{df}

\begin{df}
A vertex-preserving homeomorphism $H: \Gamma_i \to \Gamma_i'$ \emph{extends} to an ornamentation-preserving homeomorphism $H^T: G_i \to G_i'$ if there exist homeomorphisms $i_e: \tni(e) \to \tni([e])$ and homeomorphisms $i_e': H(\tni(e)) \to H(\tni([e]))$ such that the following diagram commutes:

{\large{\[ \begin{CD}
\tni([e]) @>{H^T_{\tni([e])}}>>  G_i'\\
@A{i_e}AA           @AA{i_e'}A \\
\tni(e) @>{H_{\tni(e)}}>>  H(\tni(e))
\end{CD} \]}}

\smallskip

\noindent One would also say in this circumstance that $H^T$ \emph{restricts} to $H$.
\end{df}

\noindent We now describe what it means for two based LTT structures to be equivalent.

\begin{df}
 Suppose that $G$ and $G'$ are LTT structures with respective base graphs $\Gamma$ and $\Gamma'$.  If there exists a vertex-preserving homeomorphism $H: \Gamma_i \to \Gamma_i'$ extending to an ornamentation-preserving homeomorphism $H^T: G_i \to G_i'$ then we say that $H: \Gamma_i \to \Gamma_i'$ \emph{induces an equivalence of $G_i$ and $G_i'$} and that $G_i$ and $G_i'$ are \emph{equivalent} based LTT structures.
\end{df}

\subsubsection{Generating Triples}

\noindent Since we deal with representatives decomposed into Nielsen generators, we will use an abstract notion of a ``generating triple.''

\begin{df}  A \emph{generating triple} is a triple $(g_i, G_{i-1}, G_i)$ such that
\begin{itemize}
\item [(1)]  $g_i: \Gamma_{i-1} \to \Gamma_i$ is a Nielsen generator of marked graphs defined by $g_i: e_{i-1,j_i} \mapsto e_{i,k_i} e_{i,j_i}$ for some $e_{i-1,j_i} \in \mathcal{E}_{i-1}$, $e_{i,k_i}, e_{i,j_i} \in \mathcal{E}_i$, and $e_{i,k_i} \neq (e_{i,j_i})^{\pm 1}$;
\item [(2)] $G_{i-1}$ is an LTT structure with base graph $\Gamma_{i-1}$;
\item [(3)] $G_i$ is an LTT structure with base graph $\Gamma_i$;
\item [(4)] $D^T(g_i): G_{i-1} \to G_i$ is the induced map of based LTT structures;
\item [(5)] $G_i$ contains the red edge $[\overline{d^a_i}, d^u_i]$ where $\overline{d^a_i}=D_0(\overline{e_{i,k_i}})$ and $\overline{d^u_i}=D_0(e_{i,j_i})$; and
\item [(6)] either $D_0(e_{i-1,j_i})$ or $D_0(e_{i-1,k_i})$ is a red vertex in $G_{i-1}$.
\end{itemize}
\end{df}

\noindent The following establishes equivalences for generating triples of based LTT structures.

\begin{df}{\label{D:GeneratorExtendsToLTTStructures}}
Suppose that $(g_i, G_{i-1}, G_i)$ and $(g_i', G_{i-1}', G_i)'$ are generating triples.  Let $g_i^T: G_{i-1} \to G_i$ be the map of LTT structures induced by  $g_i: \Gamma_{i-1} \to \Gamma_i$ and let $g_i^T: G_{i-1}' \to G_i'$ be the map of LTT structures induced by $g_i: \Gamma_{i-1}' \to \Gamma_i'$.

We say that $(g_i, G_{i-1}, G_i)$ and $(g_i', G_{i-1}', G_i')$ are \emph{equivalent generating triples} if there exist homeomorphisms $H_{i-1}: \Gamma_{i-1} \to \Gamma_{i-1}'$ and $H_i: \Gamma_i \to \Gamma_i'$ such that
\begin{itemize}
\item $H_i: \Gamma_i \to \Gamma_i'$ induces an equivalence of $G_i$ and $G_i'$ as based LTT structures,
\item $H_{i-1}: \Gamma_{i-1} \to \Gamma_{i-1}'$ induces an equivalence of $G_{i-1}$ and $G_{i-1}'$ as based LTT structures,
\item and the following diagram commutes:

{\Large{\[ \begin{CD}
\Gamma_i @>{H_i}>>  \Gamma_i'\\
@A{g_i}AA           @AA{g_i'}A \\
\Gamma_{i-1} @>{H_{i-1}}>>  \Gamma_{i-1}'
\end{CD} \]}}
\end{itemize}
\end{df}

\vskip10pt

\subsection{LTT Structures of Type (*) Representatives}{\label{SS:RealLTTs}}

We now give a few definitions that will enable us to apply the abstract definitions given earlier to the setting of Type (*) representatives, as defined in Section \ref{Ch:IdealDecompositions}.

\begin{df} Let $g: \Gamma \to \Gamma$ be ideally decomposable (with the condition on directions being fixed dropped), let $\mathcal{E}^+(\Gamma)= \{E_1, E_2, \dots , E_r \}$ be the set of edges of $\Gamma$ with a particular orientation, and let $\mathcal{E}(\Gamma)=\{E_1, \overline{E_1}, E_2, \overline{E_2}, \dots , E_r, \overline{E_r} \}$. For an oriented edge $E_i$, we let $d_i=D_0(E_i)$ and $\overline{d_i}= D_0(\overline{E_i})$ (the terminal direction of $E_i$).

The \emph{Colored local Whitehead graph at the vertex $v \in \Gamma$}, \emph{$CW(g;v)$}, is the uncolored graph $LW(g;v)$ but with the subgraph $SW(g;v)$ colored purple and $LW(g;v)-SW(g;v)$ colored red (including the nonperiodic vertices).

Let $\Gamma_N$ be the graph obtained from $\Gamma$ by removing a contractible neighborhood, $N(v)$, of the vertex $v$ of $\Gamma$ and adding vertices $d_i$ and $\overline{d_i}$ at the corresponding boundary points of each partial edge $E_i-(N(v) \cap E_i)$, for each $E_i \in \mathcal{E}^+$.  A \emph{Lamination Train Track Structure} \emph{$G(g)$} for $g$ is formed from $\Gamma_N \bigsqcup CW(g;v)$ by identifying the vertex labeled $d_i$ in $\Gamma_N$ with the vertex labeled $d_i$ in $CW(g;v)$.  The vertices for nonperiodic directions are red, the edges of $\Gamma_N$ remain black, and all periodic vertices remain purple.

An LTT structure $G(g)$ is given a \emph{smooth structure} via a partition of the edges at each vertex into the two sets: $\mathcal{E}_b$ (containing all black edges of $G(g)$) and $\mathcal{E}_c$ (containing all colored edges of $G(g)$).  A \emph{smooth path} will be a path alternating between colored and black edges.

The \emph{Augmented Lamination Train Track Structure} for $g$, \emph{$G_{A}(g)$}, is formed from $G(g)$ by adding a green edge for each illegal turn of $g$.   \end{df}

\begin{rk} We record here the following remarks about LTT structures:
\begin{itemize}
\item [(1)] $G(g)$ could also be built from $\underset{\text{vertices v} \in \Gamma} {\bigsqcup} CW(g;v)$ by adding a black edge connecting each vertex pair $\{D_0(e_i), \overline{D_0(e_i)} \}$.
\item [(2)]  Each edge image path $g(e_i)=e_{j_1}...e_{j_k}$ determines a smooth path in $G(g)$ that transverses the black edge $[d_{j_1}, \overline{d_{j_1}}]$, then the colored edge $[\overline{d_{j_1}}, d_{j_2}]$, then the black edge $[d_{j_2}, \overline{d_{j_2}}]$, and so on, until it ends with the black edge $[d_{j_k}, \overline{d_{j_k}}]$.  This observation is related to one of the most important properties of LTT structures for fully irreducible representatives, i.e. they contain leaves of the attracting lamination as locally smoothly embedded lines.
\item [(3)]  The train track structures we define are not quite the same as those in [BH97].
\end{itemize}
\end{rk}

\begin{df} Let $g:\Gamma \to \Gamma$ be an ideally decomposed Type (*) representative of $\phi \in Out(F_r)$ with the standard ideal decomposition notation.  Then $G_k$ will denote the LTT structure $G(f_k)$ and $G_{k,l}$ will denote the subgraph of $G_l$ containing \newline
\indent (1) all black edges and vertices (given the same colors and labels as in $G_l$) and \newline
\indent (2) all colored edges representing turns in $g_{k,l}(e)$ for some $e \in \mathcal{E}_{k-1}$. \newline
\noindent For any $k,l$, we have a direction map $Dg_{k,l}$ and an induced map of turns $Dg_{k,l}^t$.  The \emph{induced map of LTT Structures} $Dg_{k,l}^T: G_{l-1} \mapsto G_k$ is such that \newline
\indent (1) the vertex corresponding to a direction $d$ is mapped to the vertex corresponding to the direction $Dg_{k,l}(d)$, \newline
\indent (2) the colored edge [$d_1, d_2$] is mapped linearly as an extension of the vertex map to the edge [$Dg_{k,l}^t(\{d_1, d_2 \})$] $=$ [$Dg_{k,l}(d_1), Dg_{k,l}(d_2)$], and \newline
\indent (3) the interior of the black edge of $G_{l-1}$ corresponding to the edge $E \in \mathcal{E}(\Gamma_{l-1})$ to the interior of the smooth path in $G_k$ corresponding to $g(E)$.
\end{df}

\begin{rk} It still makes sense to define $G_k$ when $\phi$ is only irreducible (not fully irreducible) and possibly even is not ageometric.  The difference will be that, while the purple subgraph will be $SW(g)$, it will not necessarily be $IW(g)$.  If $\Gamma$ had more than one vertex, one would define $G(g)$ by creating a colored graph $CW(g;v)$ for each vertex, removing an open neighborhood of each vertex when forming $\Gamma_N$, and then continuing with the identifications as above in $\Gamma_N \bigsqcup (\cup CW(g;v))$. Dropping the condition on $g$ having $2r-1$ fixed directions more drastically changes what definitions actually make sense or what they look like if they do make sense.
\end{rk}

\section{Admissible Map Properties}{\label{Ch:AMProperties}}

The aim of this section is establishing additional properties held by any Type (*) representative.  In particular, we determine several necessary characteristics of LTT structures $G_k$ arising in an ideal decomposition of a Type (*) representative and give the background to identify (as described in Section \ref{Ch:Peels}) the only two possible types of (fold/peel) relationships between any LTT structures $G_{k-1}$ and $G_k$ in an ideal decomposition.  The properties proved necessary in this section will be called ``Admissible Map Properties.''  They are summarized in the final subsection, Subsection \ref{S:AMPropertiesSummarized}.

In subsequent sections, we will define and outline a method for associating, a diagram (the ``AM Diagram'') to a Type (*) pIW graph $\mathcal{G}$.  This diagram will contain a loop for each map having the Admissible Map Properties we establish in this section.  Thus, in particular, the diagram contains a loop for any Type (*) representative $g$ with $IW(g)=\mathcal{G}$.  If no loop in the diagram gives an irreducible, PNP-free representative $g$ with $IW(g)=\mathcal{G}$, then we know that $\mathcal{G}$ does not occur as $IW(\phi)$ for any ageometric, fully irreducible $\phi \in Out(F_r)$.  We will use this fact to rule out the possibility of achieving certain graphs in Sections \ref{Ch:UnachievableIWGs} and \ref{Ch:Achievable}.  We also give in subsequent sections methods for constructing the Type (*) representatives, if they do exist.

The conditions of an ideal decomposition will be relaxed slightly for many of the subsections of this section in order to highlight the necessity of certain properties (we will make it clear when representatives will also be required to be ideally decomposed of Type (*)).  \emph{For each of these subsections, $g:\Gamma \to \Gamma$ will be an irreducible TT representative of $\phi\in Out(F_r)$ \emph{semi-ideally decomposed} as: $\Gamma = \Gamma_0 \xrightarrow{g_1} \Gamma_1 \xrightarrow{g_2} \cdots \xrightarrow{g_{n-1}}\Gamma_{n-1} \xrightarrow{g_n} \Gamma_n = \Gamma$, where the decomposition differs from an ideal decomposition (as described in the end of Section \ref{Ch:IdealDecompositions}) in that the periodic directions may not be fixed.  The notation for the decomposition will also be as in the end of Section \ref{Ch:IdealDecompositions} and saying that $g$ is of Type (*) will still mean that $IW(\phi)$ is a Type (*) pIWG.  We will use the standard ideal decomposition notation, even when our decomposition in only semi-ideal.}

\vskip7pt

\noindent We begin this section by proving a preliminary lemma that will be used later in this section to prove the necessity of
the ``Admissible Map Properties.''

\subsection{Cancellation and a Preliminary Lemma}

\noindent Before stating the lemma, we clarify for the reader what is meant by ``cancellation''.

\begin{df}  We say that an edge path $\gamma=e_1 \dots e_k$ in $\Gamma$ has \emph{cancellation} if $\overline{e_i} =e_{i+1}$ for some $1 \leq i \leq k-1$.  We say that $g$ has \emph{no cancellation on edges} if for no $l>0$ and edge $e \in \mathcal{E}(\Gamma)$ does the edge path $g^l(e)$ have cancellation.
\end{df}

\noindent We are now ready to state and prove the lemma.

\begin{lem}{\label{L:PreLemma}} Suppose that $g: \Gamma \to \Gamma$ is a semi-ideally decomposed TT representative of $\phi \in Out(F_r)$ with the standard ideal decomposition notation.  For this lemma we index the generators in the decomposition of all powers $g^p$ of $g$ so that $g^p=g_{pn} \circ g_{pn-1} \circ \dots \circ g_{(p-1)n} \circ \dots \circ g_{(p-2)n} \circ \dots \circ g_{n+1} \circ g_n \circ \dots \circ g_1$ ($g_{mn+i}=g_i$, but we want to use the indices to keep track of a generator's place in the decomposition of $g^p$).  With this notation, $g_{k,l}$ will mean $g_k \circ \dots \circ g_l$.  Then: \newline
\indent (1)  for each $e \in \mathcal{E}(\Gamma_{l-1}$), no $g_{k,l}(e)$ has cancellation; \newline
\indent (2) for each $0 \leq l \leq k$ and each edge $E_{l-1,i} \in \mathcal{E}^+(\Gamma_{l-1})$, the edge $E_{k,i}$ is contained in the edge path $g_{k,l} (E_{l-1,i})$; and \newline
\indent (3) if $e^u_k=e_{k,j}$, then the turn $\{\overline{d^a_k}, d^u_k \}$ is in the edge path $g_{k,l}(e_{l-1,j})$, for all $0 \leq l \leq k$.
\end{lem}

\noindent \emph{Proof}: Let $s$ be minimal so that some $g_{s,t} (e_{t-1,j})$ has cancellation. Before continuing with our proof of (1), we first proceed by induction on $k-l$ to show that (2) holds for $k<s$.  For the base case observe that $g_{l+1}(e_{l,j})=e_{l+1,j}$ for all $e_{l+1,j} \neq (e^{pu}_l)^{\pm 1}$. Thus, if $e_{l,j} \neq e^{pu}_l$ and $e_{l,j} \neq \overline{e^{pu}_l}$ then $g_{l+1}(e_{l,j})$ is precisely the path $e_{l+1,j}$ and so we are only left for the base case to consider when $e_{l,j} = (e^{pu}_l)^{\pm 1}$.  If $e_{l,j} = e^{pu}_l$, then $g_{l+1}(e_{l,j})=e^a_{l+1} e_{l+1,j}$ and so the edge path $g_{l+1}(e_{l,j})$ contains $e_{l+1,j}$, as desired.  If $e_{l,j} =\overline{e^{pu}_l}$, then $g_{l+1}(e_{l,j})=e_{l+1,j} \overline{e^a_{l+1}}$ and so the edge path $g_{l+1}(e_{l,j})$ also contains $e_{l+1,j}$ in this case.  Having considered all possibilities, the base case is proved.

For the inductive step, we assume that $g_{k-1,l+1} (e_{l,j})$ contains $e_{k-1,j}$ and show that $e_{k,j}$ is in the edge path $g_{k,l+1}(e_{l,j})$.  Let $g_{k-1,l+1}(e_{l,j})= e_{i_1}\dots e_{i_{q-1}} e_{k-1,j} e_{i_{q+1}}\dots e_{i_r}$ for some edges $e_i \in \mathcal{E}_{k-1}$.  As in the base case, for all $e_{k-1,j} \neq (e^u_k)^{\pm 1}$, $g_k(e_{k-1,j})$ is precisely the edge path $e_{k,j}$.
Thus (since $g_k$ is an automorphism and since there is no cancellation in $g_ {j_1,j_2}(e_{j_1,j_2})$ for $1 \leq j_1 \leq j_2 \leq k$), $g_{k,l+1} (e_{l,j})=\gamma_1...\gamma_{q-1}  (e_{k,j}) \gamma_{q+1}...\gamma_m$ where each $\gamma_{i_j}= g_l({e_{i_j}})$ and where no $\{\overline{\gamma_i}, \gamma_{i+1} \}$, $\{\overline{e_{k,j}}, \gamma_{q+1} \}$, or $\{\overline{\gamma_{q-1}}, e_{k,j} \}$ is an illegal turn.
So each $e_{k,j}$ is in $g_{k,l+1} (e_{l,j})$, as desired.  We are only left to consider for the inductive step the cases where $e_{k-1,j}= e^{pu}_k$ and where $e_{k-1,j} =\overline {e^{pu}_k}$.

If $e_{k-1,j} = e^{pu}_k$, then $g_k(e_{k-1,j})=e^a_k e_{k,j}$, and so $g_{k,l+1}(e_{l,j}) = \gamma_1...\gamma_{q-1} e^a_k e_{k,j} \gamma_{q+1} \dots \gamma_m$ (where no $\{\overline{\gamma_i}, \gamma_{i+1} \}$, $\{\overline{e_{k,j}}, \gamma_{q+1} \}$, or $\{\overline{\gamma_{q-1}}, e^a_k \}$ is an illegal turn), which contains $e_{k,j}$, as desired.  If instead $e_{k-1,j} =\overline {e^{pu}_k}$, then $g_k(e_{k-1,j})=e_{k,j} \overline{e^a_k}$ and so $g_{k,l+1}(e_{l,j})=\gamma_1 \dots \gamma_{q-1} e_{k,j} \overline{e^a_k} \gamma_{q+1} \dots\gamma_m$, which also contains $e_{k,j}$.  Having considered all possibilities, the inductive step is now also proven and the proof is complete for (2) in the case of $k<s$.

We now finish our proof of (1).  We are still assuming that $s$ is minimal so that $g_{s,t}(e_{t-1,j})$ has cancellation for some $e_{t-1,j} \in \mathcal{E}_j$.  Let $t$ be such that $g_{s,t}(e_{t-1,j})$ has cancellation.  Let $\alpha_j$, for $1 \leq j \leq m$, be edges in $\Gamma_{s-1}$ so that $g_{s-1,t}(e_{t-1,j})= \alpha_1 \cdots \alpha_m$.  Since $s$ was minimal, either $g_s(\alpha_i)$ has cancellation for some $1 \leq i \leq m$ or $Dg_s(\overline{\alpha_i})=Dg_s(\alpha_{i+1})$ for some $1 \leq i < m$. Since each $g_s$ is a generator, no $g_s(\alpha_i)$ has cancellation.  Thus, there exists an $i$ such that $Dg_s(\overline{\alpha_i})= Dg_s(\alpha_{i+1})$.  Since we have already proved (1) for all $k<s$, we know that the edge path $g_{t-1,1}(e_{0,j})$ contains $e_{t-1,j}$.  Then $g_{s,1}(e_{0,j})= g_{s,t}(g_{t-1,1}(e_{0,j}))$ contains cancellation, which implies that $g^p(e_{0,j})= g_{pn,s+1}(g_{s,1}(e_{0,j}))= g_{s,t}(\dots e_{t-1,j} \dots)$ for some $p$ (with $pn>s+1$) contains cancellation, which contradicts that $g$ is a train track map.

We now prove (3).  Let $e^u_k=e_{k,l}$ By (2) we know that the edge path $g_{k-1,l}(e_{l-1,j})$ contains $e_{k-1,j}$.  Let $e_1, \dots e_m \in \mathcal{E}_{k-1}$ be such that $g_{k-1,l}(e_{l-1,j})=e_1 \dots e_{q-1} e_{k-1,j} e_{q+1} \dots e_m$.  Then $g_{k,l}(e_{l-1,j})=\gamma_1 \dots \gamma_{q-1} e^a_k e^u_k \gamma_{q+1} \dots \gamma_r$ where $\gamma_j=g_k(e_j)$ for all $j$.  Thus $g_{k,l}(e^{pu}_{k-1})$ contains $\{ \overline{d^a_k}, d^u_k \}$, as desired. \newline
\noindent QED.

\subsection{LTT Structures, Birecurrency, and AM Property I}

\vskip5pt

\indent LTT structures were defined in Section \ref{Ch:LTT} and the ``birecurrency''(defined below) of each LTT structure $G_k$ in a semi-ideal decomposition is the first property we will prove necessary for a Type (*) representative, i.e. \emph{AM Property I}.

\begin{df} We will say that a smooth train track graph $G$ is \emph{birecurrent} if there exists a locally smoothly embedded line in $G$ that crosses each edge of $G$ infinitely many times as $\bold{R}\to \infty$ and as $\bold{R}\to -\infty$. \end{df}

\begin{prop}{\label{P:Birecurrent}} Let $g: \Gamma \to \Gamma$ be a Type (*) representative of $\phi \in Out(F_r)$.  Then $G(g)$ is birecurrent. \end{prop}

Our proof of this proposition will require the following lemmas recording the relationship between the local Whitehead graph for $g$, $LW(g)$, and the realization of the leaves of the attracting lamination, $\Lambda_{\phi}$, for $\phi$.  The proofs will use facts about laminations that can be found in [BFH97] and [HM11], but will not be recorded here.

\begin{lem}{\label{L:LeafTurns}} Let $g: \Gamma \to \Gamma$ be a Type (*) representative of $\phi \in Out(F_r)$.  The only turns possible in the realization in $\Gamma$ of a leaf of the attracting lamination $\Lambda_{\phi}$ for $\phi$ are those corresponding to edges in $LW(g)$.  Conversely, each turn represented by an edge of $LW(g)$ is a turn of some (hence all) leaves of $\Lambda_{\phi}$ (as realized in $\Gamma$).
\end{lem}

\noindent \emph{Proof}: To prove the forward direction, we first notice, as follows, that each edge $E_i \in \mathcal{E}(\Gamma)$ has a fixed point in its interior.  Since $g$ is irreducible, some $g^k(E_i)$ contains a path with at least three edges (some $g^k(E_i)$ contains at least two edges of $\Gamma$, including $E_i$ and then $g^{2k}(E_i)$ contains at least three edges).  Let $g^k(E_i)=e_1e_2...e_m$, with each $e_i \in \mathcal{E}(\Gamma)$.  Again, since $g$ is irreducible, for some $l$, the edge path $g^l(e_2)$ contains either $E_i$ or $\overline{E_i}$. Thus, $g^{k+l}(E_i)$ contains either $E_i$ or $\overline{E_i}$ in its interior, implying that $E_i$ has a fixed point in its interior.  This then tells us that, for each edge $E_i \in \mathcal{E}(\Gamma)$, there is a periodic leaf of $\Lambda_{\phi}$ obtained by iterating a neighborhood of a fixed point of $E_i$.

Consider any turn $\{d_1, d_2\}$ taken by the realization in $\Gamma$ of a leaf $L$ of $\Lambda_{\phi}$. Since periodic leaves are dense in the lamination, either $\overline{e_1} e_2$ or $\overline{e_2} e_1$ (where $D_0(e_1)=d_1$ or $D_0(e_2)=d_2$) is a subpath of any periodic leaf of the lamination. In particular, either $\overline{e_1} e_2$ or $\overline{e_2} e_1$ is a subpath of the leaf obtained by iterating a neighborhood of a fixed point of $e$ for any $e \in \mathcal{E}(\Gamma)$, so $\overline{e_1} e_2$ is contained in some $g^k(e)$, for each $e \in \mathcal{E}(\Gamma)$.  Thus, $\{d_1, d_2\}$ is represented by an edge in $LW(g)$, as desired. This concludes the forward direction.

We now prove the converse.  The presence of the turn $\{d_1, d_2\}$ as an edge of $LW(g)$ indicates that, for some $i$ and $k$, $\overline{e_1} e_2$ is a subpath of $g^k(E_i)$.  We showed above that each $E_i \in \mathcal{E}(\Gamma)$ has a fixed point in its interior and hence that there is a periodic leaf of $\Lambda_{\phi}$ obtained by iterating a neighborhood of the fixed point of $E_i$.  $g^k(E_i)$ is a subpath of this periodic leaf and (since periodic leaves are dense) of every leaf of $\Lambda_{\phi}$.   Since the leaves contain $g^k(E_i)$ as a subpath, they contain $\overline{e_1} e_2$ as a subpath, and thus the turn $\{d_1, d_2\}$.  This concludes the proof of the converse, and hence lemma. \newline
\noindent QED.

\vskip8pt

\noindent We will need one more definition for the proof of the second lemma.

\begin{df}  Let $g: \Gamma \to \Gamma$ be a train track representative of a fully irreducible, ageometric $\phi \in Out(F_r)$.  Let $\gamma$ be a smooth (possibly infinite) path in $G(g)$. The \emph{path (or line) in $\Gamma$ corresponding to $\gamma$} is $\dots e_{-j}e_{-j+1} \dots e_{-1}e_0e_1 \dots e_j \dots$, where \newline
\indent$\gamma= \dots [d_{-j}, \overline{d_{-j}}][\overline{d_{-j}}, d_{-j+1}] \dots [d_{-1}, \overline{d_{-1}}][\overline{d_{-1}}, d_0][d_0, \overline{d_0}][\overline{d_0}, d_1][d_1, \overline{d_1}] \dots [d_j, \overline{d_j}] \dots$, \newline
\noindent where each $d_i=D_0(e_i)$, each $[d_i, \overline{d_i}]=[e_i]$ is the black edge of $G$ corresponding to the edge $e_i \in \mathcal{E}(\Gamma)$, and each $[d_i, \overline{d_{i+1}}]$ is a colored edge. \end{df}

\begin{lem}{\label{L:SmoothPathsforLeaves}} Let $g: \Gamma \to \Gamma$ be a TT representative of a fully irreducible, ageometric $\phi \in Out(F_r)$.  Then $G(g)$ contains smooth paths corresponding to the realizations in $\Gamma$ of the leaves of $\Lambda_{\phi}$. \end{lem}

\noindent \emph{Proof of Lemma}: Consider the realization $\lambda$ of a leaf of $\Lambda_{\phi}$ and any single subpath $\sigma=e_1 e_2 e_3$ in $\lambda$.  If it exists, the representation in $G(g)$ of $\sigma$ would be by the path \newline
\noindent $[d_1, \overline{d_1}] [\overline{d_1}, d_2] [d_2, \overline{d_2}] [\overline{d_2}, d_3] [d_3,\overline{d_3}]$, as above.  Lemma \ref{L:LeafTurns} above tells us that $[\overline{d}_1, d_2]$ and $[\overline{d_2}, d_3]$ are edges of $LW(g)$ and hence are colored edges in $G(g)$.  The path representing $\sigma$ in $G(g)$ thus exists and alternates between colored and black edges.  By looking at overlapping subpaths, we can see that the path in $G(g)$ corresponding to $\lambda$ has no consecutive colored or black edges and so is smooth.  We have proved the lemma. \newline
\noindent QED.

\vskip10pt

\noindent We are now ready for the proof of the proposition.

\bigskip

\noindent \textbf{Proof of Proposition \ref{P:Birecurrent}:} We need that $G(g)$ contains a locally smoothly embedded line crossing over each edge of $G(g)$ infinitely many times as $\bold{R}\to \infty$ and as $\bold{R}\to -\infty$.  We will show that the path $\gamma$ corresponding to the realization $\lambda$ of a leaf of $\Lambda_{\phi}$ is such a line.  We first consider any colored edge [$d_i, d_j$] in $G(g)$.  By Lemma \ref{L:LeafTurns}, $\lambda$ must contain either $\overline{e_i} e_j$ or $\overline{e_2} e_1$ as a subpath.  Birecurrency of the lamination leaves of a fully irreducible $\phi \in Out(F_r)$ implies that $\gamma$ must cross the subpath $\overline{e_i} e_j$ or $\overline{e_j} e_i$ infinitely many times as $\bold{R}\to \infty$ and as $\bold{R}\to -\infty$.  We showed in Lemma \ref{L:SmoothPathsforLeaves} above that this means that $\lambda$ contains either $\overline{e_i} e_j$ or $\overline{e_j} e_i$ infinitely many times as $\bold{R}\to \infty$ and as $\bold{R}\to -\infty$.  This concludes the proof for a colored edge.

Now consider a black edge $[d_l, \overline{d_l}]=[e_l]$.  Each vertex is shared with a colored edge.  Let [$d_l, \overline{d_m}$] be such an edge.  As shown above, $\overline{e_l} e_m$ or $\overline{e_m} e_l$ occur in realizations $\lambda$ infinitely many times as $\bold{R}\to \infty$ and as $\bold{R}\to -\infty$. In particular, it crosses over $e_l$ infinitely many times as $\bold{R}\to \infty$ and as $\bold{R}\to -\infty$.  And so $\gamma$ crosses over $[d_l, \overline{d_l}]=[e_l]$ infinitely many times as $\bold{R}\to \infty$ and as $\bold{R}\to -\infty$.  This concludes the proof. \newline
\noindent QED.

\vskip8pt

In combination with Proposition \ref{P:Birecurrent}, the second of the following two lemmas proves the necessity of AM Property I.  The first (Lemma \ref{L:fkGates}) is used in the proof of the second (Lemma \ref{L:fk}) .

\begin{lem}{\label{L:fkGates}} Let $g$ be a semi-ideally decomposed train track representative.  Each $f_k$ has the same number of gates (and thus periodic directions).
\end{lem}

\noindent \emph{Proof}: Suppose, for the sake of contradiction, that $f_k$ had more gates than $f_l$.  Let $p_k$ be such that $D(f^{p_k}_k)$ maps each gate of $f_k$ to a single direction and let $p_l$ be such that $D(f^{p_l}_l)$ maps each gate of $f_l$ to a single direction.  Let $\{\mathcal{G}_1, \dots , \mathcal{G}_s \}$ be the set of gates for $f_k$, let $\alpha_i$ be the periodic direction of $\mathcal{G}_i$ for each $1 \leq i \leq s$, let $\{\mathcal{G}'_1, \dots, \mathcal{G}'_{s'} \}$ be the set of gates for $f_l$, and let $\alpha'_i$ be the periodic direction of $\mathcal{G}'_i$ for each $1 \leq i \leq s'$.  Consider $f_k^{p_k+p_l+1}=f_{k,l+1} \circ f_l^{p_l} \circ f_{l, k+1} \circ f_k^{p_k}$.  Let $\{d_1, \dots , d_t \}=D(f_{l, k+1} \circ f_k^{p_k})(\mathcal{D}_k)$.  Then $\{d_1, \dots , d_t \}$ is mapped by $D(f_l^{p_l})$ into $\{\alpha'_1 \dots \alpha'_{s'} \}$ and, consequently, $D(f_l^{p_l} \circ  f_{l, k+1} \circ f_k^{p_k}) (\mathcal{D}_k) \subset \{\alpha'_1 \dots \alpha'_{s'} \}$.  This implies that $D(f_{k,l+1})(D(f_l^{p_l} \circ f_{l, k+1} \circ f_k^{p_k})(\mathcal{D}_k)) =D(f_k^{p_k+p_l+1})(\mathcal{D}_k) \subset D(f_{k,l+1})(\{\alpha'_1 \dots \alpha'_{s'} \})$, which has at most $s'$ elements.  But this contradicts $f_k$ having more gates that $f_l$.  Thus, all $f_k$ have the same number of gates.  \newline
\noindent QED.

\vskip10pt

\begin{rk} If $g$ is an ideally decomposed Type (*) representative, then the above lemma shows that each $G_k$ has the same number of purple periodic vertices.
\end{rk}

\begin{lem}{\label{L:fk}} If $g:\Gamma \to \Gamma$ is an ideally decomposed Type (*) representative of $\phi \in Out(F_r)$ with the standard ideal decomposition notation, then each $f_k$ is also an ideally decomposed Type (*) representative of $\phi^p$ for some $p$.
\end{lem}

\noindent \emph{Proof of Lemma}: If $\Gamma = \Gamma_0 \xrightarrow{g_1} \Gamma_1 \xrightarrow{g_2} \cdots \xrightarrow{g_{n-1}}\Gamma_{n-1} \xrightarrow{g_n} \Gamma_n = \Gamma$ is an ideal decomposition of $g$, then $f_k$ can be decomposed as $\Gamma_k \xrightarrow{g_{k+1}} \Gamma_{k+1} \xrightarrow{g_{k+2}} \cdots \xrightarrow{g_{k-1}} \Gamma_{k-1} \xrightarrow{g_k} \Gamma_k.$

What we need to show is that this decomposition of $f_k$ is an ideal decomposition and that $f_k$ is a representative of $\phi^p$ for some $p$ (since we already know that $\phi$ is ageometric and fully irreducible, as well as that $IW(\phi)$ is a Type (*) pIW graph, we know that all of these things must also be true of any $\phi^p$).  Properties (1)-(3) of an ideal decomposition hold for the decomposition of $f_k$ because they hold for the decomposition of $f$ and the decompositions have the same $\Gamma_i$ and $g_i$ (just renumbered).  By the previous lemma we know that $f_k$ also has $2r-1$ gates.  Thus, some $Df^p_k$ fixes $2r-1$ directions.  Since $d^u_k$ is not in the image of $Dg_k$, it cannot be in the image of $Df^p_k$.  We thus also know that (4) holds for the decomposition of $f^k$ and we are only left to prove that $f_k$ is a representative of $\phi^p$.  Now, $g^p$ is a representative of $\phi^p$ and $g^p=(g_{1,k})^{-1}f^pg_{1,k}$.  Let $\pi:R_r \to \Gamma$ be the marking on $\Gamma$.  Since $g_{1,k}$ is a homotopy equivalence, $g_{1,k} \circ \pi$ gives a marking on $\Gamma_k$ and $g^p$ and $f^p$ just differ by a change of marking.  Thus, $g^p$ and $f^p$ are representatives of the same outer automorphism, i.e. $\phi^p$.  This concludes the proof.  \newline
\noindent QED.

\begin{rk} The same proof works to prove the statement where ``ideally'' is replaced by ``semi-ideally.''
\end{rk}

\vskip5pt

\noindent We have thus shown that every ideally decomposed Type (*) representative satisfies
\begin{description}
\item[AM Property I:]Each LTT structure $G_k$ is birecurrent.
\end{description}

\vskip10pt

\subsection{Periodic Directions and AM Property II}

\noindent The main goal of this subsection is Proposition \ref{P:PeriodicDirections}, giving \emph{AM Property II} for a Type (*) representative of $\phi \in Out(F_r)$ such that $IW(\phi)=\mathcal{G}$.

\emph{We add to the notation already established that $t^R_k=\{\overline{d^a_k},d^u_k\}$ and $T_{k+1}=\{d^{pa}_k, d^{pu}_k \}$.}

\smallskip

\noindent The following lemma is used in the proof of Proposition \ref{P:PeriodicDirections}.

\smallskip

\begin{lem} $T_k$ is an illegal turn for $f_{k-1}$.
\end{lem}

\noindent \emph{Proof}: Recall that $T_k=\{d^{pa}_{k-1}, d^{pu}_{k-1} \}$. \newline
\indent Since $D^tg_k(\{d^{pa}_{k-1}, d^{pu}_{k-1} \})= \{Dg_k(d^{pa}_{k-1}), Dg_k(d^{pu}_{k-1}) \} = \{d^a_k, d^a_k \}$, \newline
\noindent $D^tf_{k-1}(\{d^{pa}_{k-1}, d^{pu}_{k-1}\})=D^t(g_{k-1,k+1} \circ g_k) (\{d^{pa}_{k-1}, d^{pu}_{k-1}\})=D^t(g_{k-1,k+1}) (D^tg_k (\{d^{pa}_{k-1}, d^{pu}_{k-1}\}))$  \newline
\noindent $= D^tg_{k-1,k+1} (\{d^a_k, d^a_k\}) = \{D^tg_{k-1,k+1}(d^a_k), D^tg_{k-1,k+1}(d^a_k)\}$,
which is degenerate.  So $T_k$ is an illegal turn for $f_{k-1}$, as desired. \newline
\noindent QED.

\begin{df} As a result of the previous lemma, for the generator $g_k: e^{pu}_{k-1} \mapsto e^a_k e^u_k$, we will sometimes call $T_k = \{d^{pa}_{k-1}, d^{pu}_{k-1} \}$ the \emph{green illegal turn in $G_{k-1}$} (even though $T_k$ does not technically live in $G_{k-1}$, but in the augmented LTT structure for $f_{k-1}$).
\end{df}

\noindent We are now ready to prove the proposition.

\begin{prop}{\label{P:PeriodicDirections}} Let $g: \Gamma \to \Gamma$ be semi-ideally decomposed (though not necessarily irreducible).  $g$ has $2r-1$ periodic directions if and only if, for each $k$, the illegal turn $T_{k+1}=\{d^{pa}_k, d^{pu}_k \}$ contains $d^{u}_{k}$, ie, either $d^{pu}_{k}=d^{u}_{k}$ or $d^{pa}_{k}=d^{u}_{k}$.  In fact, if each $T_{k+1}$ contains $d^{u}_{k}$, the image of $Dg$ contains all directions in $\mathcal{D}(\Gamma)$ except $d^u_n$.
\end{prop}

\noindent \emph{Proof}: We start by proving the forward direction.  Suppose that our map $g$ has $2r-1$ periodic directions and, for the sake of contradiction, that the illegal turn $T_{k+1}= \{d^{pa}_k, d^{pu}_k \}$ does not contain $d^u_k=d_{k,i}$.  Let $d^u_{k+1}=d_{k+1,s}$ and $d^a_{k+1}=d_{k+1,t}$.  Then $Dg_k(d_{k-1,s})=d_{k,s}$ and $Dg_k(d_{k-1,t})=d_{k,t}$, which means that \newline
\noindent $D^t(g_{k+1} \circ g_k)(\{d_{(k-1,s)}, d_{(k-1,t)} \})=\{D(g_{k+1} \circ g_k)(d_{(k-1,s)}), D(g_{k+1} \circ g_k)(d_{(k-1,t)})\}=$ \newline
\noindent $\{ Dg_{k+1}(d_{k,s}= d^{pu}_k), Dg_{k+1}(d_{k,t}= d^{pa}_k)\}= \{ d^a_{k+1}, d^a_{k+1} \}$
and so $d_{k-1,s}$ and $d_{k-1,t}$ share a gate.  But $d_{k-1,i}$ is already in a gate with more than one element and we already established that $d_{k-1,i} \neq d_{k-1,s}$ and $d_{k-1,i} \neq d_{k-1,t}$.  So $f_{k-1}$ has a maximum of $2r-2$ gates.  Since each $f_k$ has the same number of gates, this would imply that $g$ has a maximum of $2r-2$ gates, giving a contradiction.  The forward direction is thus proved.

Now suppose that, for each $1 \leq k \leq n$, the illegal turn $T_{k+1}$ for the generator $g_{k+1}$ always contained the unachieved direction $d^u_k$ for the generator $g_k$.  We will proceed by induction to prove that $g$ would then have $2r-1$ distinct gates.  In fact, we will show that the image of $Dg$ is missing precisely $d^u_n$, where $g=g_n \circ \cdots \circ g_1$.

For the base case we need that $g_1$ has $2r-1$ distinct gates.  By our assumptions, $g_1: e^{pu}_0 \mapsto e^a_1 e^u_1$.  The direction map for $g_1$, $Dg_1$, is defined by $Dg_1(d^{pu}_0) = d^a_1$ and $Dg_1(d_{0,t})=d_{1,t}$ for all $t$ with $d_{0,t} \neq d^{pu}_0$.  Thus, the image of $Dg_1$ includes $2r-1$ distinct directions and is missing precisely $d^u_1$.  Also, the only direction with two preimages is $d^a_1$.  This concludes the proof of the entire base case.

For the inductive step assume that $g_{k-1,1}$ has $2r-1$ distinct gates (there are $2r-1$ distinct directions (and second indices) in the image of $Dg_{k-1,1}$) and that $d^u_{k-1}$ is the only direction not in the image of $Dg_{k-1,1}$.  We also assume that $g_k$ is defined by $g_k: e^{pu}_{k-1} \mapsto e^a_k e^u_k$ where either (1) $d^u_{k-1} = e^{pu}_{k-1} $ or (2) $d^u_{k-1} = e^{pa}_{k-1}$.

Consider Case (1) where $d_{(k-1,i_1)}, d_{(k-1,i_2)}, \dots, d_{(k-1,i_{2r-1})}$ are the $2r-1$ directions in the image of $Dg_{k-1,1}$ (none of which is $d^{pu}_{k-1}=d^u_{k-1}=d_{k-1,j}$). $Dg_k: d^u_{k-1}=d^{pu}_{k-1} \mapsto d^a_k$ and preserves the second indices of all of other directions.  Since none of the \newline
$d_{(k-1,i_1)}, d_{(k-1,i_2)}, \dots, d_{(k-1,i_{2r-1})}$ are $d^{pu}_{k-1} = d^u_{k-1} = d_{k-1,j}$, $Dg_k$ acts as the identity on the second indices of $d_{(k-1,i_1)}, d_{(k-1,i_2)}, ..., d_{(k-1,i_{2r-1})}$, leaving \newline
\noindent $d_{(k,i_1)}, d_{(k,i_2)}, ..., d_{(k,i_{2r-1})}$
as $2r-1$ distinct directions in the image of $Dg_k$ (still none of which is equal to $d^u_k=d_{k,j}$) and the second indices in the image of $Dg_{k,1}$ and $Dg_{k-1,1}$ are the same.  Since $d_{k-1,i_1}, d_{k-1,i_2}, ..., d_{k-1,i_{2r-1}}$ were the only directions in the image of $Dg_{k-1,1}$, their images are the only directions in the image of $Dg_{k,1}$, meaning that $d^u_k=d_{k,j}$ is also not in the image of $Dg_{k,1}$.  Thus, $g_{k,1}$ has precisely $2r-1$ distinct gates and $d^u_k=d_{k,j}$ is not in the image of $Dg_{k,1}$, which were our two desired conclusions.

Now, consider Case (2), i.e $d^u_{k-1} = d^{pa}_{k-1}(=d_{k-1,j})$, where $d^{pu}_{k-1}=d_{(k-1,i_1)}$, \newline
\noindent $d_{(k-1,i_2)}, \dots, d_{(k-1,i_{2r-1})}$ are the $2r-1$ directions in the image of $Dg_{k-1,1}$ (none of which is $d^u_{k-1} = d^{pa}_{k-1}=d_{k-1,j}$).  $Dg_k$ is defined by $Dg_k(d^{pu}_{k-1})= d^a_k(=d_{k,j})$, mapping $d^{pu}_{k-1}$ to $d^a_k=d_{k,j}$ and $d_{k-1,i_t}$ to $d_{k,i_t}$  for $2 \leq t \leq 2r-1$ (replacing the index $i_1$ with the previously absent index $j$ and fixing all other indices).  Since $i_1 \neq i_t$ for $1 \neq t$, since $i_1$ is replaced by $j$, and since $d^u_k=d_{k,i_1}$, we can conclude that $d^u_k$ is not in the image of $Dg_{k,1}$.  Thus we have shown our two desired conclusions in this case also, i.e. that $g_{k,1}$ has $2r-1$ distinct gates and $d^u_k$ is not in the image of $Dg_{k,1}$, as desired.

We have thus completed the inductive step and consequently inductively proved the backward direction, completing the proof of the entire proposition.   \newline
\noindent QED.

\begin{cor} \textbf{(of Proposition \ref{P:PeriodicDirections})}  For each $k$, $t^R_k= \{\overline{d^a_k}, d^u_k\}$, must contain either $d^{pu}_k$ or $d^{pa}_k$. \end{cor}

\noindent \emph{Proof}:  We showed that, for each $1 \leq k \leq n$, the illegal turn $T_{k+1}=\{d^{pa}_k, d^{pu}_k \}$ always contains $d^u_k$.  At the same time, we know that $t^R_k= \{\overline{d^a_k}, d^u_k\}$, implying $t^R_k$ contains $d^u_k$ and thus either $d^{pa}_k$ or $d^{pu}_k$. \newline
\noindent QED.

\vskip5pt

\noindent We have shown that every ideally decomposed Type (*) representative satisfies
\begin{description}
\item[AM Property II:] At each graph $G_k$, the illegal turn $T_{k+1}$ for the generator $g_{k+1}$ exiting $G_k$  always contains the unachieved direction $d^u_k$ for the generator $g_k$ entering the graph $G_k$, i.e. either $d^u_k=d^{pa}_k$ or $d^u_k=d^{pu}_k$.
\end{description}

\subsection{The Nonperiodic Red Direction and AM Property III}{\label{S:AMPropertyIII}}

The conditions for this subsection are the same as described in the start of the section. \newline
\indent The following corollary of the proof of Proposition \ref{P:PeriodicDirections} gives \emph{AM Property III}.

\begin{cor}{\label{C:UnachievedDirection}} \textbf{(of Proof of Proposition \ref{P:PeriodicDirections})} For each $1 \leq k \leq n$, $d^u_k$ is the direction not fixed by $Df_k$, i.e. $d^u_k$ is not a periodic direction for $f_k$.  In particular, the vertex labeled by $d^u_k$ in $G_k$ is red and $[t^R_k]= [\overline{d^a_k} , d^u_k]$ is a red edge in $G_k$. \end{cor}

\noindent \emph{Proof}: Since $Dg_k$ is defined by $d^{pu}_{k-1} \mapsto d^a_k$ (and $Dg_k(d_{k-1,j})=d_{k-1,j}$ for all $j$ such that $e^{pu}_{k-1} \neq e_{k-1,j}$ and $e^{pu}_{k-1} \neq \overline{e_{k,j}}$), $d^u_k$ is not in the image of $Dg_k$.  Suppose for the sake of contradiction, that $d^u_k$ is a periodic direction for $f_k$.  Let $N$ be a sufficiently high power for $f_k$ so that $D(f_k^N)$ fixes all periodic directions of $f_k$.  Then $g_k$ would still be the final generator in the decomposition of $f_k$ and thus $d^u_k$ would also not be in the image of $D(f_k^N)$.  This contradicts $d^u_k$ being a periodic direction for $f_k$.  So $d^u_k$ is not a periodic direction for $f_k$ and hence labels a red vertex in $G_k$.

Before we can identify that $[t^R_k]$ is a red edge, we first need to show that $[t^R_k]$ is an edge of $LW(f_k)$.  In order to show that $[t^R_k]$  is an edge of $LW(f_k)$, it suffices to show that $[t^R_k]=[\overline{d^a_k}, d^u_k]$ is in $f_k(e^u_k)$.  Let $e^u_k=e_{k,l}$.  By Lemma \ref{L:PreLemma} we know that the edge path $g_{k-1,k+1} (e^u_k=e_{k,l})$ contains $e_{k-1,l}$.  Let $e_j$ be edges in $\Gamma_{l-1}$ such that $g_{k-1,k+1}(e^u_k)= e_1 \dots e_{q-1} e_{k-1,l} e_{q+1} \dots e_m$.  Then $f_k(e^u_k ) =g_{k,k+1} (e^u_k )= \gamma_1 \dots \gamma_{q-1} e^a_k e^u_k \gamma_{q+1}\dots \gamma_m$ where $\gamma_j =g_k(e_{i_j})$ for all $j$.  Thus $f_k(e^u_k)$ contains $\{ \bar d^a_k, d^u_k \}$, as desired and [$t^R_k$] is an edge of $LW(f_k)$.  Since $[\overline{d^a_k} , d^u_k]$ contains the red vertex $d^u_k$, $[\overline{d^a_k}, d^u_k]$ is a red edge in $G_k$. \newline
\noindent QED.

\begin{df}\label{D:CreatesEdge} As a consequence of the proof of Corollary \ref{C:UnachievedDirection}, we will say that $g_k$ \emph{creates the edge} $[\overline{d^a_k}, d^u_k]$ in $G_k$ (in the sense that $\{\bar d^a_k, d^u_k \}$ is a turn in the image of $g_{k,l}(e^u_k)$ for any $1 \leq l \leq n$ and $[\overline{d^a_k}, d^u_k]$ is in $G_k$ (and is, in fact, the red edge of $G_k$, as it is not in the image of $d^C(g_k)$, but in the image of the black edge [$e^u_k$] in $G_{k-1}$ under $dg_k^T$)).  Further details are discussed in Subsection \ref{S:AMPropertyV}.

As a consequence of Corollary \ref{C:UnachievedDirection}, we will also henceforth sometimes refer to $d^u_k$ as the \emph{(red) unachieved direction in $G_k$ (and $G_{k,l}$)}, $t^R_k=\{\overline{d^a_k}, d^u_k \}$ as the \emph{new red turn in $G_k$ (and $G_{k,l}$)}, and $e^R_k=$[$t^R_k$] as the \emph{red edge in $G_k$ (and $G_{k,l}$)}.  We will sometimes call $d^a_k$ the \emph{twice-achieved direction in $G_k$(and $G_{k,l}$)} for reasons ascertainable by analyzing the proof of Proposition \ref{P:PeriodicDirections}.
\end{df}

\begin{rk} Visually what we established in Corollary \ref{C:UnachievedDirection} for a Type (*) representative is that, in each augmented LTT structure $G_A(f_k)$, the intersection of the red edge and green edge is the red vertex.
\end{rk}

\noindent We have now shown that every ideally decomposed Type (*) representative satisfies
\begin{description}
\item[AM Property III:] The vertex labeled by $d^u_k$ is red in $G_k$ and $[t^R_k]=[d^u_k, \overline{d^a_k}]$ is a red edge in $G_k$.
\end{description}

\noindent Since we have what is necessary to do so at this point, and it will be used later, we will prove a final lemma about periodic directions here before shifting our focus.

\begin{lem}{\label{L:DirectImage}} Suppose that $g: \Gamma \to \Gamma$ is semi-ideally decomposed and has $2r-1$ periodic directions.  Then the image under $Dg$ of the $2r$ directions at the vertex $v \in \Gamma$ is precisely the set of the $2r-1$ periodic directions for $g$.
\end{lem}

\noindent \emph{Proof}:  Since $Dg_n$'s image is missing $d^u_n$, it is clear that the image of $Dg$ has at most $2r-1$ directions.  So we are left to show that $Dg$'s image cannot be missing a periodic direction for $g$.

For the sake of contradiction, let $d_k$ be a periodic direction not in the image of $Dg$. Then $d_k$ is also not in the image of any $Dg^n$ since $Dg^n=Dg \circ Dg^{n-1}$. Let $N$ be such that $Dg^N$ fixes every periodic direction. Then $d_k$ is still not in the image of $Dg^N$, so it cannot be one of the periodic directions.  This is a contradiction, meaning that the image of $Dg$ cannot be missing a periodic direction for $g$.  The lemma is proved. \newline
\noindent QED.

\vskip10pt

\subsection{$D^Cg_{k,l}$ Edge Images and AM Property IV}

\noindent The following lemma gives \emph{AM Property IV}.

\begin{lem}{\label{L:EdgeImage}} Let $g:\Gamma \to \Gamma$ be a semi-ideally decomposed representative of $\phi \in Out(F_r)$ with the standard notation.  If $[d_{(l,i)},d_{(l,j)}]$ is a purple or red edge in $G_l$, then $[D^tg_{k,l+1}(\{d_{(l,i)},d_{(l,j)}\})]$ is a purple edge in $G_k$.
\end{lem}

\noindent \emph{Proof}:  It suffices to show two things: \newline
\indent (1) $D^tg_{k,l+1}(\{d_{(l,i)},d_{(l,j)}\})$ is a turn in some edge path $f_l^p(e_{l,m})$ with $p \geq 1$ and \newline
\indent (2) $Dg_{k,l+1}(d_{l,i})$ and $Dg_{k,l+1}(d_{l,j})$ are periodic directions for $f_l$. \newline
\noindent We will proceed by induction and start with (1).  For the base case of (1) assume that the turn $\{d_{(k-1,i)}, d_{(k-1,j)}\}$ is represented by a purple or red edge in $G_{k-1}$.  Then $f_{k-1}^p(e_{k-1,t})=$ \newline
\noindent $s_1 \dots \overline{e_{(k-1,i)}} e_{(k-1,j)}... s_m$ for some edges $e_{(k-1,t)},s_1, \dots s_m \in \mathcal{E}_{k-1}$ and $p \geq 1$.  By Lemma \ref{L:PreLemma}, $e_{k-1,t}$ is contained in the edge path $g_{k-1} \circ \cdots \circ g_1 \circ g_n \circ \cdots \circ g_{k+1} (e_{k,t})$.  Thus, since $f_{k-1}^p(e_{k-1,t})=s_1 \dots \overline{e_{(k-1,i)}} e_{(k-1,j)} \dots s_m$ and no $g_{i,j}(e_{j-1,t})$ can have cancellation, $f_{k-1}^p \circ g_{k-1} \circ \cdots \circ g_1 \circ g_n \circ \cdots \circ g_{k+1} (e_{k,t})$ contains $s_1 \dots \overline{e_{(k-1,i)}} e_{(k-1,j)} \dots s_m$ as a subpath.  Applying $g_k$ to $f_{k-1}^p \circ g_{k-1} \circ \cdots \circ g_1 \circ g_n \circ \cdots \circ g_k (e_{k-1,t})$, we get $f_k^{p+1}(e_{k,t})$.

Suppose first that $Dg_k(e_{k-1,i})=e_{k,i}$ and $Dg_k(e_{k-1,j})=e_{k,j}$.  Then $g_k(\dots\overline{e_{(k-1,i)}} e_{(k-1,j)} \dots) =$ $\dots \overline{e_{(k,i)}} e_{(k,j)} \dots$, with possibly different edges before and after $\overline{e_{k,i}}$ and $e_{k,j}$ than before and after $\overline{e_{k-1,i}}$ and $e_{k-1,j}$.  Thus, in this case, $f^{p+1}_k( \dots \overline{e_{(k-1,i)}} e_{(k-1,j)} \dots)$ contains the turn $\{d_{(k,i)},d_{(k,j)}\}$, which in this case is $D^tg_k(\{d_{(k-1,i)}, d_{(k-1,j)}\})$.  So $D^tg_k(\{d_{(k-1,i)}, d_{(k-1,j)}\})$ is represented by an edge in $G_k$.

Now suppose that $g_k: e_{k-1,j} \mapsto e_{k,l} e_{k,j}$.  Then $g_k(\dots \overline{e_{(k-1,i)}} e_{(k-1,j)}\dots) = \dots \overline{e_{(k,i)}} e_{(k,l)}e_{(k,j)} \dots$,
(again with possibly different edges before and after $\overline{e_{k,i}}$ and $e_{k,j}$).  So $g_k(\dots \overline{e_{(k-1,i)}} e_{(k-1,j)} \dots)$ contains the turn $\{\overline{d_{(k,l)}},d_{(k,j)}\}$, which in this case is $D^tg_k(\{d_{(k-1,i)},d_{(k-1,j)}\})$, so $D^tg_k(\{d_{(k-1,i)},d_{(k-1,j)}\})$ is again represented by an edge in $G_k$.

Finally, suppose that $g_k$ is defined by $g_k: e_{k-1,j} \mapsto e_{k,j}e_{k,l}$.  Unless $\overline{e_{k-1,i}} = e_{(k-1,j)}$, \newline
\noindent $g_k(\dots \overline{e_{(k-1,i)}} e_{(k-1,j)} \dots)=\dots \overline{e_{(k,i)}} e_{(k,j)} e_{(k,l)} \dots$, which contains the turn $\{d_{(k,i)},d_{(k,j)}\} =$ \newline
\noindent $D^tg_k(\{d_{(k-1,i)},d_{(k-1,j)}\})$, implying that $D^tg_k(\{d_{(k-1,i)},d_{(k-1,j)}\})$ is represented by an edge in $G_k$ in this case also.

If $\overline{e_{k-1,i}} = e_{k-1,j}$, then we are actually in a reflection of the previous case.  The other cases ($g_k: \overline{e_{k-1,i}} \mapsto \overline{e_{k,i}} e_{k,l}$ and $g_k: \overline{e_{k-1,i}} \mapsto e_{k,l} \overline{e_{k,i}}$) follow similarly by symmetry.  We have thus completed the base case for our proof of (1).

We now must prove the base case for (2).  Since $[D^tg_k(\{d_{(k-1,i)},d_{(k-1,j)}\})] = [Dg_k(d_{(k-1,i)}),Dg_k(d_{(k-1,j)})]$, both vertices of $[D^tg_k(\{d_{(k-1,i)},d_{(k-1,j)}\})]$ are directions in the image of $Dg_k$.  By Lemma \ref{L:DirectImage}, combined with Lemma \ref{L:fk}, this means that both vertices represent periodic directions.  Thus, $[D^tg_k(\{d_{(k-1,i)},d_{(k-1,j)}\})]$ is actually a purple edge in $G_k$, concluding our proof of the base case.

Now suppose inductively that $[d_{(l,i)},d_{(l,j)}]$ is a purple or red edge in $G_l$ and \newline
\noindent $[D^tg_{k-1,l+1}(\{d_{(l,i)},d_{(l,j)}\})]$ is a purple edge in $G_{k-1}$.  Then the base case implies that \newline \noindent $[D^tg_k(D^tg_{k-1,l+1}(\{d_{(l,i)},d_{(l,j)}\})]$  is a purple edge in $G_k$.  But $D^tg_k(D^tg_{k-1,l+1}(\{d_{(l,i)},d_{(l,j)}\}))$ \newline
\noindent $= D^tg_{k,l+1}(\{d_{(l,i)},d_{(l,j)}\})$.  So the lemma is proved. \newline
\noindent QED.

\vskip8pt

\noindent We have now shown that every ideally decomposed Type (*) representative satisfies
\begin{description}
\item[AM Property IV:]  If $[d_{(l,i)}, d_{(l,j)}]$ is a purple or red edge in $G_l$, then [$D^Cg_{k,l+1}(\{ d_{(l,i)}, d_{(l,j)} \})$] is a purple edge in $G_k$.
\end{description}

\vskip10pt

\subsection{The Red Turn and AM Property V}{\label{S:AMPropertyV}}

The aim of this subsection is to better understand red edges, their properties, and how they are ``created.''  This subsection should also begin to shed light on how generic edges in an ideal Whitehead graph are ``created'' by generators in an ideal decomposition.\newline
\indent \emph{For this subsection, $g:\Gamma \to \Gamma$ is an ideally decomposed Type (*) representative of $\phi \in Out(F_r)$ with the standard ideal decomposition notation.} \newline
\indent As above, we say that $g_k$ creates the edge $e=[d_{(k,i)}, d_{(k,j)}]$ of $G_k$ if $g_k$ is defined by either $e_{k-1,i} \mapsto \overline{e_{k,j}} e_{k,i}$ or $e_{k-1,j} \mapsto \overline{e_{k,i}} e_{k,j}$.  The first and second of the following lemmas, together with Lemma \ref{L:NG}, tell us that $g_k$ ``creating'' $\{\overline{d^a_k}, d^u_k \}$ means what we intuitively want for it to mean.

\begin{lem}{\label{L:RedEdgeImage}} For each $1 \leq l,k \leq n$, $[D^tg_{l,k}(\{\overline{d^a_{k-1}}, d^u_{k-1} \})]$ is a purple edge in $G_l$.
\end{lem}

\noindent \emph{Proof}: By Property IV proved in Lemma \ref{L:EdgeImage}, it suffices to show that $[\overline{d^a_{k-1}}, d^u_{k-1}]$ is a colored edge of $G_{k-1}$.  This was shown in Corollary \ref{C:UnachievedDirection}. \newline
\noindent QED.

\begin{lem} $[\overline{d^a_k}, d^u_k]$ is not in $D^Cg_k(G_{k-1})$. \end{lem}

\noindent \emph{Proof}: By Lemma \ref{L:EdgeImage}, all purple and red edges of $G_{k-1}$ are mapped to purple edges in $G_k$.  On the other hand, $[\overline{d^a_k}, d^u_k]$ is a red edge in $G_k$.  Thus, $[\overline{d^a_k}, d^u_k]$ is not in $D^Cg_k(G_{k-1})$. \newline
QED.

\begin{rk} Notice that the above lemmas also show the uniqueness of $g_k$ once the red edge and red nonperiodic direction vertex of $G_k$ are known.  This is explained further in the next section.
\end{rk}

\smallskip

\noindent The following Lemma (together with Corollary \ref{C:UnachievedDirection}) gives \emph{AM Property V}.

\begin{lem} $LW(g)$ can have at most 1 edge segment connecting the nonperiodic red direction vertex to the set of purple periodic direction vertices. \end{lem}

\noindent \emph{Proof}: First notice that the nonperiodic direction vertex is the red vertex $d^u_k$ in $G_k$. If $g_k(e_{k-1,i})= e_{k,i}e_{k,j}$, then the red direction in $G_k$ is $\overline{d_{k,i}}$ (where $d_{k,i}=D_0(e_{k,i})$ and $d_{k,j}=D_0(e_{k,j})$). If $g_k$ is the final generator in the decomposition, then the vertex $\overline{d_{k,i}}$ will be adjoined to the vertex for $d_{k,j}$ and only $d_{k,j}$, as every occurrence of $e_{k-1,i}$ in the image under $g_{k-1,1}$ of any edge has been replaced by $e_{k,i}e_{k,j}$ and every occurrence of $\overline{e_{k,i}}$ has been replaced by $\overline{e_{k,i}} \overline{e_{k,j}}$, ie, there are no copies of $e_{k,j}$ without $e_{k,i}$ following them and no copies of $\overline{e_{k,i}}$ without $\overline{e_{k,j}}$ preceding them. \newline
\noindent QED.

\vskip10pt

\noindent We have now shown that every ideally decomposed Type (*) representative satisfies
\begin{description}
\item[AM Property V:] Each $C(G_k)$ can have at most one edge segment connecting the red (nonperiodic) vertex of $G_k$ to the set of purple (periodic) vertices of $G_k$.  This single edge is red and is in fact the edge $[t^R_k]=[d^u_k, \overline{d^a_k}]$.
\end{description}

\vskip10pt

\subsection{The Ingoing Nielsen Generator and AM Property VI}

\noindent Given an LTT structure $G_k$ in a Type (*) representative ideal decomposition (or even just given the red vertex or red edge), there is only one possibility for the generator $g_k$ entering $G_k$.  We will use this fact when constructing representatives yielding our desired ideal Whitehead graphs.

\emph{We continue to assume that $g:\Gamma \to \Gamma$ is an ideally decomposed Type (*) representative of $\phi \in Out(F_r)$ with the standard ideal decomposition notation.}

The following lemma gives AM Property VI.

\begin{lem}{\label{L:NG}} Let $g:\Gamma \to \Gamma$ be an ideally decomposed Type (*) representative of $\phi \in Out(F_r)$ with the standard notation.  Suppose that the unique red edge in $G_k$ is $[t^R_k]= [d_{(k,j)}, \overline{d_{(k,i)}}]$ and that the vertex representing $d_{k,j}$ is red.  Then $g_k(e_{k-1,j}) = e_{k,i} e_{k,j}$ and $g_k(e_{k-1,t})=e_{k,t}$ for $e_{k-1,t} \neq (e_{k-1,j})^{\pm 1}$, where $D_0(e_{s,t})=d_{s,t}$ and $D_0(\overline{e_{s,t}}) = \overline{d_{s,t}}$ for all $s$, $t$.  \end{lem}

\noindent \emph{Proof}:  By the definition of an ideal decomposition, $g_k$ must be of the form $g_k: e_{k-1,j} \mapsto e_{k,i} e_{k,j}$ ($g_k(e_{k-1,i})= e_{k,i}$ for $e_{k-1,i} \neq (e_{k-1,j})^{\pm 1}$ and $e_{k,i} \neq (e_{k,j})^{\pm 1}$).  Corollary \ref{C:UnachievedDirection} indicates that $D_0(e_{k,j})=d^u_k$, i.e. the direction associated to the red vertex of $G_k$.  In other words, the second index of $d^u_k$ uniquely determines the index $j$ and so $e_{k-1,j}= e^{pu}_{k-1}$ and $e_{k,i}=e^a_k$.  Additionally, the proof of Corollary \ref{C:UnachievedDirection} indicates that $[\overline{d_{(k,i)}}, d_{(k,j)}]$ is the red edge of $G_k$.  This means that we must have $e_{k,i}=e^a_k$.  $g_k$ has thus been determined to be $g_k : e^{pu}_{k-1} \mapsto e^a_k e^u_k$, i.e, $e_{k-1,j} \mapsto e_{k,i} e_{k,j}$, as desired. \newline
\noindent QED.

\begin{df} The $g_k$ in Lemma \ref{L:NG} will be called the \emph{ingoing Nielsen generator} for $G_k$. \end{df}

\vskip5pt

\noindent We have now shown that every ideally decomposed Type (*) representative satisfies
\begin{description}
\item[AM Property VI:] Given that $[t^R_k] = [d^u_k, \overline{d^a_k}]$ is the red edge of $G_k$ and $d^u_k$ labels the single red vertex of $G_k$, $g_k$ is defined by  $g_k(e^{pu}_{k-1})=e^a_k e^u_k$ and $g_k(e_{k-1,i})=e_{k,i}$ for $e_{k-1,i} \neq (e^{pu}_{k-1})^{\pm 1}$, where $D_0(e^u_k)=d^u_k$, $D_0(\overline{e^a_k})=\overline{d^a_k}$, $e^{pu}_{k-1}= e_{(k-1,j)}$, and $e^u_k=e_{k,j}$.
\end{description}

\vskip10pt

\subsection{Isomorphic Ideal Whitehead Graphs and AM Property VII}

\noindent The aim of this subsection is \emph{AM Property VII} (stated in Proposition \ref{P:SWsIsomorphic}), giving that representatives of the same outer automorphism have isomorphic ideal whitehead graphs.

\begin{prop}{\label{P:SWsIsomorphic}} Let $g:\Gamma \to \Gamma$ be an ideally decomposed Type (*) representative of $\phi \in Out(F_r)$ with the standard notation.  For each $0 \leq l,k \leq n$, $Dg_{l, k+1}$ induces an isomorphism from SW($f_k$) onto SW($f_l$). \end{prop}

\smallskip

\noindent The proof of the proposition will come after the following two lemmas used in the proof.  Notice that Lemma \ref{L:fk} implies that SW($f_k$) and SW($f_l$) are isomorphic and so the key point of the proposition is that this isomorphism is induced by $Dg_{l, k+1}$.

\begin{lem} Each $D^Cf_k$ maps the purple subgraph $PI(G_k)$ of $G_k$ isomorphically (as a graph) onto itself.  Further, the graph isomorphism preserves the vertex and edge labels.
\end{lem}

\noindent \emph{Proof}: Lemma \ref{L:EdgeImage} implies that $D^Cf_k$ maps the purple subgraph of $G_k$ into itself.  However, $Df_k$ fixes all directions corresponding to vertices of the purple graph.  Thus,  $D^Cf_k$ restricted to $PI(G_k)$ is a label-preserving graph isomorphism onto its image. \newline
\noindent QED.

\begin{lem}{\label{L:PurpleEdgeImages}} The set of purple edges of $G_{k-1}$ is mapped by $D^Cg_k$ injectively into the set of purple edges of $G_k$.  \end{lem}

\noindent \emph{Proof}:  Since $d^a_k$ is the only direction with more than one preimage of $Dg_k$ and these two preimages are $d^{pa}_{k-1}$ and $d^{pu}_{k-1}$, the only edges in $G_k$ with more than one preimage under $D^Cg_k$ are those of the form $[d_{(k,i)}, d^a_k]$ and the two preimages are the edges $[d_{(k-1,i)}, d^{pa}_{k-1}]$ and $[d_{(k-1,i)}, d^{pu}_{k-1}]$ in $G_{k-1}$.  However, by Proposition \ref{P:PeriodicDirections}, either $e^u_{k-1}= e^{pu}_{k-1}$ or $e^u_{k-1}=e^{pa}_{k-1}$, meaning that one of the preimages of $d^a_k$ is actually $d^u_{k-1}$, i.e. one of the preimage edges is actually $[d_{(k-1,i)}, d^u_{k-1}]$.  Since $[t^R_{k-1}]$ is the only purple or red edge of $G_{k-1}$ containing $d^u_{k-1}$, one of the preimages of $[d_{(k,i)}, d^a_k]$ must be $[e^R_{k-1}]$, leaving only one possible purple preimage. \newline
\noindent QED.

\vskip5pt

\noindent \emph{Proof of Proposition \ref{P:SWsIsomorphic}}: Since compositions of injective maps are injective, by Lemma \ref{L:PurpleEdgeImages}, the set of purple edges of $G_k$ is mapped injectively by $D^Cg_{l, k+1}$ into the set of purple edges of $G_l$.  Likewise, the set of purple edges of $G_l$ is mapped injectively by $D^Cg_{k, l+1}$ into $G_k$.  Additionally, by the first lemma proved above, $D^Cf_k=(D^Cg_{k, l+1}) \circ (D^Cg_{l, k+1})$ is a bijection.  Thus, since each of these sets of edges is a finite set, the map that $D^Cg_{l, k+1}$ induces on the set of purple edges of $G_k$ is a bijection.  It is only left to show that two purple edges share a vertex in $G_k$ if and only if their images under $D^Cg_{l, k+1}$ share a vertex in $G_l$.

Suppose that we have two purple edges $[x, d_1]$ and $[x, d_2]$ in $G_k$ sharing the vertex $x$.  Then $D^Cg_{l, k+1}([x, d_1])=[Dg_{l, k+1}(x), Dg_{l, k+1}(d_1)]$ and $D^tg_{l, k+1}([x, d_2])=$ \newline
\noindent $[Dg_{l, k+1}(x), Dg_{l, k+1}(d_2)]$ share the vertex $Dg_{l, k+1}(x)$.  This proves the forward direction.  To prove the other direction, observe that, if two purple edges $[w, d_3]$ and $[w, d_4]$ in $G_l$ share the vertex $w$, then $[D^tg_{k, l+1}(\{w, d_3 \})]=[Dg_{k, l+1}(w), Dg_{k, l+1}(d_3)]$ and $[D^tg_{k, l+1}(\{w, d_4 \})]=[Dg_{k, l+1}(w), Dg_{k, l+1}(d_4)]$ share the vertex $Dg_{k, l+1}(w)$ in $G$.  Since $D^Cf_l$ is an isomorphism on $PI(G_l)$, $D^Cg_{l, k+1}$ and $D^Cg_{k, l+1}$ act on inverses.  So the preimages of $[w, d_3]$ and $[w, d_4]$ under $D^Cg_{l, k+1}$ share a vertex in $G_l$. \newline
\noindent QED.

\begin{cor}  \textbf{(of Proposition \ref{P:SWsIsomorphic})} Purple edges of $G_k$ are images under $D^Cg_k$ of purple edges of $G_{k-1}$.
\end{cor}

\noindent \emph{Proof of Corollary}: From the proposition, we know that $D^Cg_k$ gives a bijection on the set of purple edges of $G_{k-1}$.  In particular, it is surjective, meaning that the purple edges of $G_k$ are all images under $Dg_k$ of purple edges of $G_{k-1}$, as desired. \newline
\noindent QED.

\vskip10pt

\noindent We have now shown that every ideally decomposed Type (*) representative satisfies:

\begin{description}
\item[AM Property VII:] $Dg_{l,k+1}$ induces an isomorphism from $SW(f_k)$ onto $SW(f_l)$ for $0 \leq l,k \leq n$.
\end{description}

\vskip10pt

\subsection{Irreducibility and AM Property VIII}

In order for a train track map to represent a fully irreducible outer automorphism, it certainly needs to be irreducible. We begin this subsection with several definitions.

\begin{df} The \emph{transition matrix} for an irreducible TT representative $g$ is the square matrix such that, for each $i$ and $j$, the $ij^{th}$ entry is the number of times $g(E_j)$ crosses $E_i$ in either direction.  A matrix $A=[a_{ij}]$ is an \emph{irreducible matrix} if each entry $a_{ij} \geq 0$ and if, for each $i$ and $j$, there exists a $k>0$ so that the $ij^{th}$ entry of $A^k$ is strictly positive.  If the same $k$ works for each index pair $\{i, j \}$, then the matrix is called \emph{aperiodic}.  If each sufficiently high $k$ works for all index pairs $\{i, j \}$, then the matrix is called \emph{Perron-Frobenius (PF)}. [BH92]
\end{df}

\begin{rk} PF matrices are part of the Full Irreducibility Criterion.  We collect here the following facts about transition matrices and PF matrices:
\begin{itemize}
\item [(1)] Any power of a Perron-Frobenius matrix is Perron-Frobenius and irreducible.
\item [(2)] A power of an irreducible matrix need not be irreducible.
\item [(3)] While aperiodic matrices are irreducible, the converse is not always true.
\item [(4)] A topological representative is irreducible if and only if its transition matrix is irreducible [BH92].
\end{itemize}
\end{rk}

The following three lemmas give properties stemming from irreducibility (though not proving irreducibility).  Together they comprise AM Property VIII.

\emph{We will assume that $g:\Gamma \to \Gamma$ is a semi-ideally decomposed train track representative of $\phi \in Out(F_r)$ with the standard notation}.

\begin{lem} For each $1 \leq j \leq r$, there exists a $k$ such that either $e^u_k=E_{k,j}$ or $e^u_k= \overline{E_{k,j}}$.
\end{lem}

\noindent \emph{Proof}: Suppose, for the sake of contradiction, that, if there is some $j$ so that $e^u_k \neq E_{k,j}$ and $e^u_k \neq \overline{E_{k,j}}$ for all $k$.  We will proceed by induction to show that $g(E_{0,j})=E_{0,j}$ and so $g$ is certainly reducible.  Induction will be on the $k$ in $g_{k-1,1}$.

For the base case, we need to show that $g_1(E_{0,j})=E_{1,j}$ if $e_1^u \neq E_{1,j}$ and $e_1^u \neq \bar{E_{1,j}}$.  $g_1$ is defined by $e^{pu}_0 \mapsto e^a_1 e^u_1$ and $g_1(e_{0,1})=e_{1,l}$ for all $e_{0,1} \neq (e^{pu}_0)^{\pm 1}$.  Since $e_1^u \neq E_{1,j}$ and $\overline{e_1^u} \neq \overline{E_{(1,j)}}$, $e^{pu}_0 \neq E_{(0,j)}$ and $e^{pu}_0 \neq \overline{E_{(0,j)}}$.  Thus, $g_1(E_{0,j})=E_{(1,j)}$, as desired.  Now suppose inductively that $g_{k-1,1}(E_{0,j})=E_{k-1,j}$ and that neither $e_k^u=E_{k,j}$ nor $\overline{e_k^u}= E_{k,j}$.
Then $e^{pu}_{k-1} \neq E_{k-1,j}$ and $e^{pu}_{k-1} \neq \overline{E_{k-1,j}}$.  Thus, since $g_k$ is defined by $e^{pu}_{k-1} \mapsto e^a_k e^u_k$ and $g_k(e_{k-1,l})= e_{k,l}$ for all $e_{k-1,l} \neq (e^{pu}_{k-1})^{\pm 1}$, $g_k(E_{k-1,j})=E_{k,j}$.  So $g_{k, 1}(E_{0,j})=g_k (g_{k-1,1}(E_{0,j}))=g_k (E_{k-1,j})= E_{(k,j)}$, as desired.  Inductively, this proves that $g(E_{0,j})=E_{0,j}$, we have our contradiction and the lemma is proved. \newline
\noindent QED

\begin{lem} For each $1 \leq j \leq r$, either $e^a_k=E_{k,j}$ or $e^a_k= \overline{E_{k,j}}$ for some $k$.
\end{lem}

\noindent \emph{Proof}: For the sake of contradiction, suppose that, for some $1 \leq j \leq r$, $e^a_k \neq E_{k,j}$ and $e^a_k \neq \overline{E_{k,j}}$ for each $k$.  The goal will be to inductively show that, for each $E_{0,i}$ with $E_{0,i} \neq E_{0,j}$ and $E_{0,i} \neq \overline{E_{(0,j)}}$, $g(E_{0,i})$ does not contain $E_{0,j}$ and does not contain $\overline{E_{0,j}}$ (contradicting irreducibility).

We start with the base case.  $g_1$ is defined by $e^{pu}_0 \mapsto e^a_1 e^u_1$ (and $g_1(e_{0,l})= e_{1,l}$ for all $e_{0,l} \neq (e^{pu}_{0})^{\pm 1}$).  First suppose that either $E_{0,j}=e^{pu}_0$ or $E_{0,j}= \overline{e^{pu}_0}$.  Then $e^{pu}_0 \neq E_{0,i}$ and $e^{pu}_0 \neq \overline{E_{0,i}}$ (since $E_{0,i} \neq E_{0,j}$ and $E_{0,i} \neq \overline{E_{(0,j)}}$) and so $g_1(E_{0,i})= E_{1,i}$, which does not contain $E_{1,j}$ or $\overline{E_{1,j}}$.  Now suppose that $E_{0,j} \neq e^{pu}_0$ and $E_{0,j} \neq \overline{e^{pu}_0}$  Then $e^a_1 e^u_1$ does not contain $E_{1,j}$ or $\overline{E_{1,j}}$ (since $e^a_k \neq (E_{k,j})^{\pm 1}$ by the assumption), which means that $E_{1,j}$ and $\overline{E_{1,j}}$ are not in the image of $E_{0,i}$ if $E_{0,i} = e^{pu}_0$ (since the image is of $E_{0,i}$ is then $e^a_1 e^u_1$) and are not in the image of $\overline {E_{0,i}}$ (since the image is $\overline{e^u_1} \overline{e^a_1}$) and are not in the image $E_{0,i}$ if $E_{0,i} \neq e^{pu}_0$ and $E_{0,i} \neq \overline{e^{pu}_0}$ (since the image is $E_{1,i}$, which does not equal $E_{1,j}$ or $\overline{E_{1,j}}$).  So the base case is proved.

Now inductively suppose that $g_{k-1,1}(E_{0,i})$ does not contain $E_{k-1,j}$ or $\overline{E_{k-1,j}}$.  A similar analysis to the above shows that $g_k(E_{k-1,i})$ does not contain $E_{k,j}$ or $\overline{E_{k,j}}$ for any $E_{k,i} \neq E_{k,j}$ and $E_{k,i} \neq \overline{E_{k,j}}$.  Since $g_{k-1,1}(E_{k-1,i})$ does not contain $E_{k-1,j}$ or $\overline{E_{k-1,j}}$, $g_{k-1, 1}(E_{0,i})= e_1 \dots e_m$ with each $e_i \neq E_{k-1,j}$ and $e_i \neq \overline{E_{k-1,j}}$.  Thus, no $g_k(e_i)$ contains $E_{k,j}$ or $\overline{E_{k,j}}$, which means that $g_{k, 1}(E_{0,i})=g_k(g_{k-1,1}(E_{0,i}))=g_k(e_1) \dots g_k(e_m)$ does not contain $E_{k,j}$ or $\overline{E_{k,j}}$.  This completes the inductive step and thus proves the lemma. \newline
\noindent QED.

\begin{rk} While the above lemmas are necessary for $g$ to be irreducible, they are not sufficient to prove the irreducibility of a semi-ideally decomposed representative. For example, the composition of $a \mapsto ab$, $b \mapsto ba$, $c \mapsto cd$, and $d \mapsto dc$ would satisfy these lemmas, but is clearly reducible.  On the other hand, Lemma \ref{L:PF} below gives a necessary and sufficient condition for irreducibility.
\end{rk}

\begin{df}  Let $g=g_n \circ \dots \circ g_1: \Gamma \to \Gamma$ be a semi-ideally decomposed Type (*) representative of $\phi \in Out(F_r)$ with the standard notation, except that we return to the convention of Lemma \ref{L:PreLemma} and index the generators in the decomposition of all powers $g^p$ of $g$ so that $g^p=g_{pn} \circ g_{pn-1} \circ \dots \circ g_{(p-1)n} \circ \dots \circ g_{(p-2)n} \circ \dots \circ g_{n+1} \circ g_n \circ \dots \circ g_1$ ($g_{mn+i}=g_i$, but we want to use the indices to keep track of a generator's place in the decomposition of $g^p$).  Again, with this notation, $g_{k,l}$ will mean $g_k \circ \dots \circ g_l$. We recursively define the \emph{edge containment sequence} for an edge $E_{j,0}$ of $\Gamma$ (or just for $j$).  For $1 \leq j \leq r$, the \emph{level-1 edge containment set for $j$}, denoted $\mathcal{C}^1_j$, contains each index $i$ such that, for some $k$, $e^{pu}_k=(E_{k,j})^{\pm 1}$ and $e^a_{k+1}=(E_{k+1,i})^{\pm 1}$.  Recursively define the \emph{level $k$ edge containment set for $j$}, denoted $\mathcal{C}^k_j$, as $\underset{i \in \mathcal{C}^{k-1}_j}{\cup} \mathcal{C}^1_i$ with duplicates of indices removed.  The \emph{edge containment sequence} for $E_{j,0}$ (or just $j$) is $\{\mathcal{C}^1_j, \mathcal{C}^2_j, \dots \}$.
\end {df}

\begin{lem}{\label{L:PF}} $g$ has a Perron-Frobenius transition matrix if and only if for each $1 \leq k,l \leq r$, we have $l \in \mathcal{C}^i_k$ for some $i$.
\end{lem}

\noindent \emph{Proof}: Suppose that for some $1 \leq k,l \leq r$, we have that $l$ is not in $\mathcal{C}^i_k$ for any $i$.  Let $H$ be the subgraph of $\Gamma$ that includes precisely the edges $E_t$ where $t \in \mathcal{C}^i_k$ for some $i$.  Then it is not too difficult to see that $H$ is a proper invariant subgraph (proper since it does not contain $E_l$).  This proves the that $g$ is not irreducible and, in particular, does not have a Perron-Frobenius transition matrix.

Now suppose that for each $1 \leq k,l \leq r$, we have that $l \in \mathcal{C}^i_k$ for some $i$.  This means that for each $1 \leq k,l \leq r$ some $g^{p(k,l)}(E_k)$ passes over $E_l$ (in some direction).  Let $p$ be the least common multiple of the $p(k,l)$. Then $M^{p(k,l)}$ is strictly positive where $M$ is the transition matrix for $g$.  And, in fact, $M^N$ is strictly positive for any $N \geq p(k,l)$, since $g$ maps each $E_l$ over itself.  This proves that $g$ has a Perron-Frobenius transition matrix and thus proves the reverse direction. \newline
\noindent QED.

\begin{rk} It will be relevant later that a semi-ideally decomposed train track representative satisfies that, for each $1 \leq k,l \leq r$, we have $l \in \mathcal{C}^i_k$ for some $i$, is not just irreducible, but actually has a Perron-Frobenius transition matrix.  Since this is a condition in the FIC, it is useful to have this way to check the condition.
\end{rk}

\vskip10pt

\noindent We have now shown that every ideally decomposed Type (*) representative satisfies:
\begin{description}
\item[AM Property VIII:] $g$ is irreducible (and, in fact, has a PF transition matrix), i.e.
\begin{itemize}
\item[(a)] for each $1 \leq j \leq r$, either $e^u_k=E_{k,j}$ or $e^u_k= \overline{E_{k,j}}$ for some $k$;
\item[(b)] for each $1 \leq j \leq r$, either $e^a_k=E_{k,j}$ or $e^a_k= \overline{E_{k,j}}$ for some $k$; and
\item[(c)] for each $1 \leq k,l \leq r$, we have that $l \in \mathcal{C}^i_k$ for any $i$.
\end{itemize}
\end{description}

\vskip10pt

\subsection{Admissible Map Properties Summarized}{\label{S:AMPropertiesSummarized}}

We proved in this section that a list of properties hold for any ideally decomposed Type (*) representative of a $\phi \in Out(F_r)$.  However, one can at least analyze whether they hold in any situation where $\Gamma = \Gamma_0 \xrightarrow{g_1} \Gamma_1 \xrightarrow{g_2} \cdots \xrightarrow{g_{n-1}} \Gamma_{n-1} \xrightarrow{g_n} \Gamma_n = \Gamma$ is an ideal decomposition of a TT representative $g$ such that $SW(g)$ is a Type (*) pIWG.

For the sake of clarity we list here the properties we proved hold for ideally decomposed Type (*) representatives and call them ``Admissible Map (AM) Properties''.  We use the standard ideal decomposition notation.

\begin{df}
Let $\mathcal{G}$ be a Type (*) pIWG.  Let $(g_{(i-k,i)}, G_{i-k-1}, G_i)$, with $k \geq 0$, be a triple such that $g_{i-k,i}$ can be decomposed as $\Gamma_{i-k-1} \xrightarrow{g_{i-k}} \Gamma_{i-k}
\xrightarrow{g_{i-k+1}} \cdots \xrightarrow{g_{i-1}}\Gamma_{i-1} \xrightarrow{g_i} \Gamma_i$,
with a sequence of LTT structures for $\mathcal{G}$ \newline
\indent $G_{i-k-1} \xrightarrow{D^T(g_{i-k})} G_{i-k} \xrightarrow{D^T(g_{i-k+1})} \cdots \xrightarrow{D^T(g_{i-1})} G_{i-1} \xrightarrow{D^T(g_i)} G_i$. \newline
We say $(g_{(i-k,i)}, G_{i-k-1}, G_i)$ satisfies the \emph{Admissible Map (AM) Properties} if it satisfies:
\begin{description}
\item[AM Property I:] Each LTT structure $G_j$, with $i-k-1 \leq j \leq i$, is birecurrent.
\item[AM Property II:] For each LTT structure $G_j$, with $i-k-1 \leq j \leq i$, the illegal turn $T_{j+1}$ for the generator $g_{j+1}$ exiting $G_j$ contains the unachieved direction $d^u_j$ for the generator $g_j$ entering the graph $G_j$, i.e. either $d^u_j=d^{pa}_j$ or $d^u_j=d^{pu}_j$.
\item[AM Property III:] In each LTT structure $G_j$, with $i-k-1 \leq j \leq i$, the vertex labeled $d^u_j$ and the edge $[t^R_j]=[d^u_j, \overline{d^a_j}]$ are both red.
\item[AM Property IV:]  For all $i-k-1 \leq j < m \leq i$, if $[d_{(j,i)}, d_{(j,l)}]$ is a purple or red edge in $G_j$, then $D^Cg_{m,j+1}$([$d_{(j,i)}, d_{(j,l)}$]) is a purple edge in $G_m$.
\item[AM Property V:] For each $i-k-1 \leq j \leq i$, $C(G_j)$ has precisely one edge segment containing the red (nonperiodic) vertex $d^u_j$ of $G_j$.  This single edge is red and is in fact $[t^R_j]=[d^u_j, \overline{d^a_j}]$.
\item[AM Property VI:] For each $i-k \leq j \leq i$, the generator $g_j$ is defined by $g_j:e^{pu}_{j-1} \mapsto e^a_j e^u_j$ (where $e^u_j=e_{j,m}$, $D_0(e^u_j)=d^u_j$, $D_0(\overline{e^a_j})=\overline{d^a_j}$, and $e^{pu}_{j-1}=e_{j-1,m}$).
\item[AM Property VII:] $Dg_{l,j+1}$ induces an isomorphism from $SW(f_j)$ onto $SW(f_l)$ for all $i-k-1 \leq j < l \leq i$.
\item[AM Property VIII:] $g$ is irreducible (and, in fact, has a Perron-Frobenius transition matrix), i.e.
\begin{itemize}
\item[(a)] for each $1 \leq j \leq r$, either $e^u_m=E_{m,j}$ or $e^u_m = \overline{E_{m,j}}$ for some $m$;
\item[(b)] for each $1 \leq j \leq r$, either $e^a_m=E_{m,j}$ or $e^a_m = \overline{E_{m,j}}$ for some $m$; and
\item[(c)] for each $1 \leq m,l \leq r$, we have that $l \in \mathcal{C}^i_m$ for any $i$.
\end{itemize}
\end{description}
\end {df}

\section{Lamination Train Track Structures are Lamination Train Track Structures}{\label{Ch:LTTMeansLTT}}

In this section we simply show that the LTT structures defined in Subsection \ref{SS:RealLTTs} are indeed abstract LTT structures.

\begin{lem}{\label{L:PF}} Let $g:\Gamma \to \Gamma$ be a Type (*) representative of $\phi \in Out(F_r)$ such that $IW(g) \cong \mathcal{G}$.  Then $G(g)$ is a Type (*) LTT structure with base graph $\Gamma$. Furthermore,
\begin{enumerate}
\item $PI(G(g)) \cong \mathcal{G}$ and
\item if $\Gamma = \Gamma_0 \xrightarrow{g_1} \Gamma_1 \xrightarrow{g_2} \cdots \xrightarrow{g_{n-1}} \Gamma_{n-1} \xrightarrow{g_n} \Gamma_n = \Gamma$ is an ideal decomposition of $g$ with the standard notation, then each $G_j=G(f_j)$ is a Type (*) LTT structure with base graph $\Gamma_j$ such that
\begin{itemize}
\item [a.] $PI(G_j) \cong \mathcal{G}$,
\item [b.] the vertex labeled $d^u_j$ is the red vertex of $G_j$, and
\item [c.] the red edge of $G_j$ is $[t^R_j]= [d^u_j, \overline{d^a_j}]$.
\end{itemize}
\end{enumerate}
\end{lem}

\noindent \emph{Proof}:  We first need that each $G_j$ is a Type (*) LTT structure with base graph $\Gamma_j$.  However, since each $f_j$ is also an ideally decomposed Type (*) representative with the same ideal Whitehead graph as $G(g)$ (and even the same ideal decomposition with simply a shifting of indices), it suffices to show that $G(g)$ is a Type (*) LTT structure with base graph $\Gamma$

For STTG1 to hold, we need that $G(g)$ has a colored edge containing each vertex, since each vertex labeled $d_i$ or $\overline{d_i}$ is contained in the black edge $[e_i]$.  Note that this would also prove STTG3 and LTT5.  Since $\mathcal{G}$ must have $2r-1$ vertices and $PI(G(g)) \cong \mathcal{G}$, there is at most one vertex without a colored edge containing it.  However, this vertex must be the red vertex contained in the red edge $[t^R_n]= [d^u_n, \overline{d^a_n}]$ by AM Property V and Corollary \ref{C:UnachievedDirection}, in particular.  We now prove STTG2.  Colored edges of $G(g)$ contain distinct vertices because they correspond to turns taken by images of edges. The black edges contain distinct vertices because they connect the directions corresponding to the initial and terminal directions in each edge of $\Gamma$, which are distinct.  This proves STTG2 and $G(g)$ is a smooth train track graph.

LTT1-LTT3 hold by construction (in the definition of $G(g)$).  That the edges of $G(g)$ are either black, purple, or red follows from the construction of the definition.  LTT4(Black Edges) holds by construction (in the definition of $G(g)$).  If an edge is red in $G(g)$, by how $G(g)$ is constructed, it means that the edge is in $LW(g)$, but not in $SW(g)$.  For this to be true, it must have a nonperiodic vertex, i.e. a red vertex.  This implies both LTT4(Red Edges) and LTT4(Purple Edges).  Since LTT5 was proved above, we are left to show that $\Gamma$ is a base for $G(g)$, but this can also easily be seen to be true by construction in the definition of $G(g)$.

LTT(*)1 holds because the fact that each $\Gamma_j$ is a rose means that each $G_j$ has $2r$ vertices and because AM Property VII implies that each $G_j$ has precisely $2r-1$ purple vertices.  LTT(*)2 holds by AM Property V

(1) is true by construction.  (2a) is true by (1) combined with AM Property VII.  (2b) and (2c) are true by AM Property III. \newline
\noindent QED

\section{Extensions, Switches, Construction Compositions, Switch Paths, and Peels}{\label{Ch:Peels}}

\noindent Let $\mathcal{G}$ be a Type (*) pIWG.  We saw in Section \ref{Ch:IdealDecompositions} that, if there is an ageometric, fully irreducible $\phi \in Out(F_r)$ with $IW(\phi) \cong \mathcal{G}$, then there is a Type (*) representative $g$ of a power of $\phi$.  By Section \ref{Ch:AMProperties}, such a representative would satisfy AM Properties I-VIII.  Thus, if we can show that a representative satisfying all these properties does not exist, then we have shown that the Type (*) representative cannot exist, and thus that there is no ageometric, fully irreducible $\phi \in Out(F_r)$ with $IW(\phi) \cong \mathcal{G}$.

We will show how to construct all ideally decomposed representatives satisfying the AM properties by determining, given knowledge of an LTT structure $G_k$ in the decomposition, all possibilities for $g_k$ and $G_{k-1}$ respecting the AM properties.  We prove in Proposition \ref{P:ExtensionsSwitches} that, if the structure $G_k$ and a purple edge $[d, d^a_k]$ in $G_k$ are set, then there is only one $g_k$ possibility and at most two $G_{k-1}$ possibilities (one possible triple $(g_k, G_{k-1}, G_k)$ will be called a ``switch'' and the other an ``extension'').

In Section \ref{Ch:AMDiagrams}, we construct the ``PreAdmissible Map (PreAM) Diagram" for $\mathcal{G}$ (denoted $PreAMD(\mathcal{G})$) from all ``permitted switches'' and ``permitted extensions.''  Then, from $PreAMD(\mathcal{G})$, we construct the ``Admissible Map (AM) Diagram" for $\mathcal{G}$ (denoted $AMD(\mathcal{G})$), in which any Type (*) representative $g$ with $IW(g) \cong \mathcal{G}$ will be ``realizable'' as a loop.  Thus, as a consequence of the above, if no loop in the $AMD(\mathcal{G})$ satisfies all the AM properties, then $\mathcal{G}$ is ``unachievable.'' The simplest ``unachievable'' examples arise when all loops in the $AMD(\mathcal{G})$ represent reducible maps.

One should note that, while we do not restrict the rank $r$, we only consider ageometric, fully irreducible $\phi \in Out(F_r)$ such that $IW(\phi)$ is a Type (*) pIWG.  The definitions below would need to be tailored for any other circumstance.

\subsection{Peels}{\label{S:Peels}}

\noindent As a warm-up for the following subsections, we describe here a geometric method for visualizing ``switches'' and ``extensions.''  It is sometimes useful to visualize switches and extensions as geometric moves (called ``peels'') transforming an LTT structure $G_i$ into an LTT structure $G_{i-1}$.  For such a peel, we call $G_i$ the \emph{source LTT structure} and $G_{i-1}$ the \emph{destination LTT structure}.

\begin{df}  Each of the two peels associated to a generator $g_i$ and source LTT structure $G_i$ involve three directed edges of $G_i$.  The \emph{first edge of the peel} begins with $d^u_i$ and ends with $\overline{d^a_i}$ (it is also known as the \emph{new red edge} in $G_i$).  The \emph{second edge of the peel} is the black edge from $\overline{d^a_i}$ to $d^a_i$ (the \emph{twice-achieved edge} in $G_i$).  The \emph{third edge of the peel} (the \emph{determining edge} $[d^a_i, d]$ for the peel), starts with $d^a_i$ and ends with some direction $d$, that will actually be labeled $\overline{d^a_{i-1}}$ in our peel's destination LTT structure $G_{i-1}$ (it will be the attaching vertex of the red edge in $G_i$).  For each choice of a determining edge, we arrive at one ``peel switch'' and one ``peel extension'' as follows. \newline
\indent Starting at the vertex $\overline{d^a_i}$, split open the black edge $[\overline{d^a_i}, d^a_i]$ and the third edge $[d^a_i, d]$, keeping $d$ fixed.  The first edge of the peel now concatenates with the black edge and third edge (all internal vertices of the edge formed by the concatenation are removed) to form $[d^u_i, d]$.  The split has created two edges: the third edge of the peel $[d^a_i, d]$ and an edge $[d^u_i, d]$.  In the peel's destination LTT structure, one of these edges will be purple and one will be red.  If $[d^u_i, d]$ is the red edge in $G_{i-1}$, then the peel is called a \emph{peel extension} (and the triple $(g_i, G_{i-1}, G_i)$ with $g_i$ as in AM Property VI, will be called the \emph{potential extension associated to $[d^a_i, d]$}).  A composition of peels will be called an \emph{extended peel}.  If $[d^a_i, d]$ is the red edge in $G_{i-1}$, then the peel will be called a \emph{peel switch} (and the triple $(g_i, G_{i-1}, G_i)$ with $g_i$ as in AM Property VI, will be called the \emph{potential switch associated to} $[d^a_i, d]$), but the structure must be further altered in this case.  One could image just shifting all edges containing $d^a_i$ to contain $d^{pu}_{i-1}$ instead.  Alternatively, one could imagine that, instead of the peel performed as above, one at a time, for each $[d^a_i, d']$ in $G_i$, the peeler peels a copy of $[d^a_i, \overline{d^a_i}, d^u_i]$ off to concatenate with $[d^a_i, d']$ and form the edge $[d^u_i, d']$.  Only in the case of $[d^a_i, d]$, is an edge (the new red edge) left behind by the process.
\end{df}

\begin{ex} The following illustrates a peel extension:

\vskip5pt

\begin{figure}[H]
\centering
\includegraphics[width=2.5in]{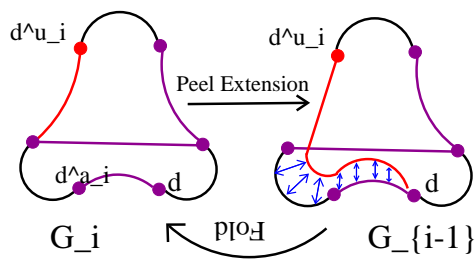}
\caption{{\small{\emph{Peel Extension: Note that the first, second, and third edges of the peel concatenate to $[d^u_i, d]$ after removal of the two interior vertices.  Here $[d^u_i, d]$ is red.}}}}
\label{fig:PeelExtension}
\end{figure}

\noindent And the following illustrates a peel switch, where $[d^a_i, d]$ is the only purple edge containing the vertex labeled $d^a_i$:

\vskip5pt

\begin{figure}[H]
\centering
\includegraphics[width=2.5in]{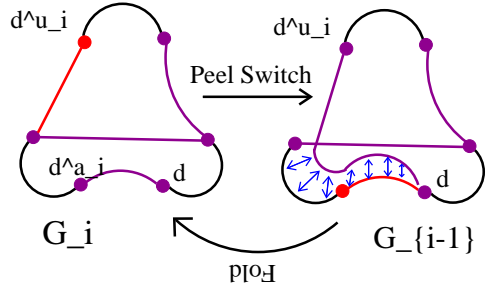}
\caption{{\small{\emph{Peel Switch: Again the first, second, and third edges of the peel concatenate to $[d^u_i, d]$ after removal of the\ interior vertices.  However, here, $[d^u_i, d]$ is purple and $[d^a_i, d]$ changes to red (as does the vertex $d^a_i$).}}}}
\label{fig:PeelSwitch}
\end{figure}

\noindent If there are other purple edges containing the vertex $d^a_i$, imagine first peeling each of them off as in Figure \ref{fig:PeelSwitch2} below (and then either splitting open $[d, d^a_i]$ and then $[d^a_i, \overline{d^a_i}]$ starting at $d$ or performing the peeling off of $[d, d^u_i]$, as demonstrated in Figure \ref{fig:PeelSwitch} above):

\vskip5pt

\begin{figure}[H]
\centering
\includegraphics[width=4.3in]{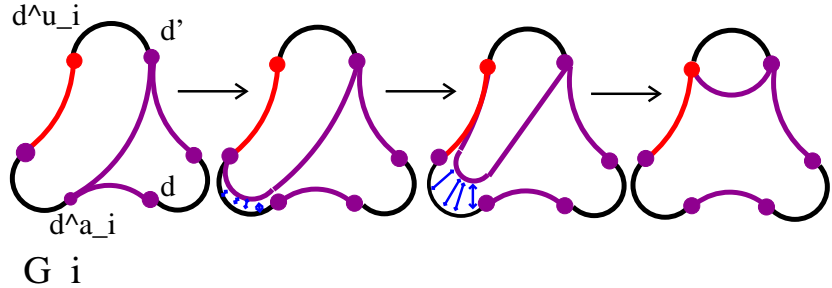}
\caption{{\small{\emph{Peel Switch Continued: For each purple edge $[d^a_i, d']$ in $G_i$, the peeler peels a copy of $[d^a_i, \overline{d^a_i}, d^u_i]$ off to concatenate with $[d^a_i, d']$ and form the purple edge $[d^u_i, d']$.}}}}
\label{fig:PeelSwitch2}
\end{figure}
\end{ex}

\begin{rk} We record here two properties of peels: \newline
\indent (1) A peel can be viewed as an extension to the LTT structure of a fold inverse.
\newline
\indent (2) In Subsection \ref{S:ConstructionCompositions} we see that a composition of extensions ``constructs'' a path (the ``construction path'') in an LTT structure.  This LTT structure is the source structure for the corresponding extended peel.  The extended peel splits open the structure along this path. \end{rk}

\subsection{Extensions and Switches}

Let $\mathcal{G}$ be a Type (*) pIWG.  In this subsection we define (permitted) extensions and switches ``entering'' a Type (*) admissible LTT structure $G_k$ for $\mathcal{G}$.  When composed, extensions ``construct'' a smooth path in an LTT structure. The colored edges in the image of the path are ``constructed'' in $\mathcal{G}$ (see Proposition \ref{P:PathConstruction}).  ``Switches'' change the LTT structure more dramatically, start the construction of a new path, and are necessary for reducibility.  Our goal in building TT representatives will be to maximize the number of extensions and minimize the number of switches, while still creating an irreducible representative yielding the entire graph.  This will allow us to more easily track our progress in ``constructing'' $\mathcal{G}$.

\vskip10pt

We collect and establish here the standard notation we will use for a Type (*) LTT structure $G_k$ for a Type (*) pIW graph $\mathcal{G}$ with rose base graph $\Gamma_k$.

\vskip10pt

\noindent \textbf{Standard Notation for a Type (*) LTT Structure:} Let $G_k$ be a Type (*) LTT structure for a Type (*) pIW graph $\mathcal{G}$ with rose base graph $\Gamma_k$.
\begin{itemize}
\item The Type (*) LTT structure definition gives a unique red vertex. Motivated by AM Property III, we denote this vertex by $d^u_k$ and call the corresponding direction the \emph{unachieved direction} in $G_k$.
\item The Type (*) LTT structure definition also gives a unique red edge (containing $d^u_k$). Motivated again by AM Property III, we denote this edge $e^R_k = [t^R_k]$.
\item The other endpoint of $e^R_k$ is a uniquely determined purple vertex. Also motivated by AM Property III, we denote this vertex by $\overline{d^a_k}$.
\item  In summary, the red edge can always be denoted $e^R_k = [t^R_k] = [d^u_k, \overline{d^a_k}]$.
\item  There is a uniquely determined black edge with endpoint $\overline{d^a_k}$, i.e. $[d^a_k, \overline{d^a_k}]$.  This edge is denoted $e^a_k$ and called the \emph{twice-achieved edge} in $G_k$.
\item  $d^a_k$ will be called the \emph{twice-achieved direction} in $G_k$ (see Corollary \ref{P:PeriodicDirections}).
\end{itemize}

\vskip18pt

\begin{lem}
Let $G_k$ be a Type (*) LTT structure for a Type (*) pIW graph $\mathcal{G}$ with rose base graph $\Gamma_k$ and the standard notation.  There exists a colored edge having an endpoint at $d^a_k$, so that it may be written $[d^a_k, d_{k,l}]$ and this edge must be purple.
\end{lem}

\noindent \emph{Proof}:  If $d^a_k$ were red, the red edge would be $[d^a_k, \overline{d^a_k}]$, violating [LTT(*)3].  If $d_{k,l}$ were red, i.e. $d_{k,l}=d^u_k$, then both $[d^u_k, \overline{d^a_k}]$ and $[d^u_k, d^a_k]$ would be red, violating [LTT(*)2].  So the edge must be purple.   At least one such purple edge must exist or $\mathcal{G}$ would not have $2r-1$ vertices. \newline
\noindent QED.

\begin{df}{\label{D:Extension}} Let $G_k$ be a Type (*) admissible LTT structure for a Type (*) pIW graph $\mathcal{G}$ with rose base graph $\Gamma_k$ and the standard notation.  An \emph{extension} associated to a purple edge $[d^a_k, d_{k,l}]$ is a triple $(g_k, G_{k-1}, G_k)$ satisfying each of the following:
\indent \begin{description}
\item[(EXTI):] $G_{k-1}$ is additionally
\begin{itemize}
\item[(a)] a Type (*) LTT structure for $\mathcal{G}$ and
\item[(b)] a Type (*) LTT structure with rose base graph (denoted $\Gamma_{k-1}$)
\end{itemize}
\item[(EXTII):]  $g_k: \Gamma_{k-1} \to \Gamma_{k-1}$ is the tight homotopy equivalence defined by $g_k(e_{k-1,j_k})=e_{k,i_k} e_{k,j_k}$ where:
\end{description}
\begin{itemize}
\item[(a)] The edge sets $\mathcal{E}(\Gamma_{k-1})= \mathcal{E}_{k-1}$ and $\mathcal{E}(\Gamma_{k})=\mathcal{E}_{k}$ are respectively denoted:
\end{itemize}

\noindent $\{E_{(k,1)}, \overline{E_{(k,1)}}, E_{(k,2)}, \overline{E_{(k,2)}}, \dots,  E_{(k,r)}, \overline{E_{(k,r)}}\}=\{e_{(k,1)}, e_{(k,2)}, \dots, e_{(k,2r-1)}, e_{(k,2r)} \}$ and \newline
\noindent $\{E_{(k-1,1)}, \overline{E_{(k-1,1)}}, \dots, E_{(k-1,2r)}, \overline{E_{(k-1,2r)}}\} =\{e_{(k-1,1)}, \dots, e_{(k-1,2r)} \}$ \newline
\indent (we let $d_{k,j}=D_0(e_{k,j}$), $\overline{d_{(k,j)}}=D_0(\overline{e_{k,j}}$), $d_{(k-1,j)}=D_0(e_{k-1,j}$), and  $\overline{d_{k-1,j}}=D_0(\overline{e_{k-1,j}}$));
\begin{itemize}
\item[(b)] $e_{k-1,j_k} \in \mathcal{E}_{k-1}$ and $e_{(k,i_k)}, e_{(k,j_k)} \in \mathcal{E}_k$ are such that $D_0(e_{k,i_k})=d_{k,i_k}=d^a_k$ and $D_0(e_{k,j_k})=d_{k,j_k}=d^u_k$;
\item[(c)] $g_k(e_{k-1,i})=e_{k,i}$ for all $e_{k-1,i} \neq e_{k-1,j_k}^{\pm 1}$.
\end{itemize}
\begin{description}
\item[(EXTIII):] $Dg_k$ induces an ornamentation-preserving graph isomorphism from $PI(G_{k-1})$ onto $PI(G_k)$ defined by sending each vertex labeled $d_{k-1,j}$ to the vertex labeled $d_{k,j}$ and extending linearly over edges.
\item[(EXTIV):] $d^u_{k-1}=d_{(k-1,j_k)}$, i.e. $d_{k-1,j_k}$ labels the single red vertex in $G_{k-1}$.
\item[(EXTV):] The single red edge of $G_{k-1}$ is $[d_{(k-1,l)},d_{(k-1,j_k)}]$.
\end{description}
\end{df}

\begin{rk} Using the definition of a peel extension in Subsection \ref{S:Peels}, one can always construct a ``potential extension'' $(g_k, G_{k-1}, G_k)$ for a purple edge $[d^a_k, d_{k,l}]$.  If it satisfies EXTI-EXTV, it will be the unique such extension.  This uniqueness will be proved in Lemma \ref{L:ExtensionUniqueness}.  Another method for obtaining the potential extension is described in Definition \ref{D:PotentialExtension} \end{rk}

\noindent \textbf{Standard Extension Terminology and Notation:} Let $G_k$ be a Type (*) admissible LTT structure for a Type (*) pIWG $\mathcal{G}$ with rose base graph $\Gamma_k$ and let $(g_k, G_{k-1}, G_k)$ be an extension, as in Definition \ref{D:Extension}.
\begin{enumerate}
\item An extension will be called \emph{admissible} if $G_k$ and $G_{k-1}$ are both birecurrent (and thus are actually Type (*) admissible LTT structures for $\mathcal{G}$).
\item We call $G_{k-1}$ the \emph{source LTT structure} and $G_k$ the \emph{destination LTT structure}.
\item The single red vertex $d_{k-1,j_k}$ in $G_{k-1}$ will be denoted by both $d^{pu}_{k-1}$ and $d^{u}_{k-1}$ (as before, ``p'' is for ``pre'').
\item $\overline{d^{a}_{k-1}}$ denotes $d_{k-1,l}$ (and $d^{a}_{k-1}$ denotes $\overline{d_{k-1,l}}$).  Consequently, the red edge $e^{R}_{k-1}$ in $G_{k-1}$ can be written, among other ways, as $[d^u_{k-1}, \overline{d^{a}_{k-1}}]$ or $[d^{pu}_{k-1}, d_{(k-1,l)}]$.
\item $e^{pa}_{k-1}$ denotes $e_{k-1,i_k}$ (again ``p'' is for ``pre'').
\item{\label{N:IngoingGeneratorTerminology}}  $g_k$ will be called the \emph{ingoing generator} and, motivated by AM Property VI, is sometimes written $g_k: e^{pu}_{k-1} \mapsto e^a_k e^u_k$.
\item We will call $[d^a_k, d_{(k,l)}]$ the \emph{(purple) edge determining $(g_k, G_{k-1}, G_k)$}.
\end{enumerate}

\vskip20pt

As explained in Subsection \ref{S:Peels}, but with this subsection's notation, an extension transforms LTT structures as follows:

\vskip5pt

\begin{figure}[H]
\centering
\includegraphics[width=5in]{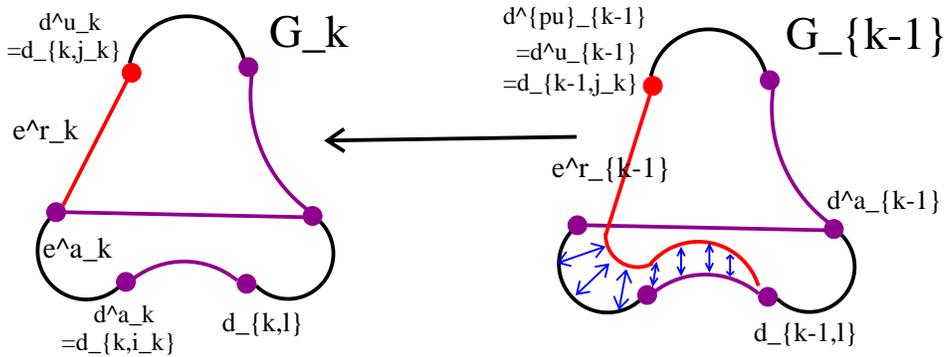}
\caption{{\small{\emph{Extension}}}}
\label{fig:ExtensionDiagram}
\end{figure}

\begin{df}{\label{D:PotentialExtension}} For a Type (*) admissible LTT structure $G_k$ with Type (*) pIWG $\mathcal{G}$ and rose base graph, we obtain the \emph{potential extension} $(g_k, G_{k-1}, G_k)$ associated to a purple edge $[d^a_k, d_{k,l}]$ by the following procedure (justification in Lemma \ref{L:ExtensionUniqueness}):
\begin{enumerate}
\item Replace each vertex label $d_{k,i}$ with $d_{k-1,i}$ and each vertex label $\overline{d_{k,i}}$ with $\overline{d_{k-1,i}}$.
\item Remove the interior of the red edge from $G_k$.
\item Add a red edge connecting the red vertex $d_{k-1,j_k}$ to $d_{k-1,l}$.
\end{enumerate}
\end{df}

\vskip10pt

\begin{rk}{\label{R:ExtensionRestriction}} It can be noted that we cannot have $d_{k,l}= \overline{d^u_k}$ (similarly $e_{k,i_k} \neq (e_{k,j_k})^{\pm 1}$). If $d_{k,l}= \overline{d^u_k}$, then the red edge of $G_{k-1}$ would be $[d^{pu}_{k-1}, \overline{d^{pu}_{k-1}}]= [d^{u}_{k-1}, \overline{d^{u}_{k-1}}]$, which would contradict that $G_{k-1}$ is a Type (*) LTT structure because of [LTT(*)4].  The following two reasons motivate this requirement.  First, $e^{R}_{k-1}= [d^{u}_{k-1}, \overline{d^{u}_{k-1}}]$ would force the ingoing generator for $G_{k-1}$ to be $e^{pu}_{k-2} \mapsto e^u_{k-1} e^u_{k-1}$, which is not a generator.  Second, $[d^{u}_{k-1}, \overline{d^{u}_{k-1}}]$ would be the only edge in $G_{k-1}$ containing the vertex labeled $d^{u}_{k-1}$, forcing $G_{k-1}$ to be nonbirecurrent (this actually motivates [LTT(*)4]).
\end{rk}

\vskip10pt

\begin{df}{\label{D:Switch}}  Let $G_k$ be a Type (*) admissible LTT structure for a Type (*) pIWG $\mathcal{G}$ with rose base graph $\Gamma_k$ with the standard notation.  The \emph{switch} associated to a purple edge $[d^a_k, d_{k,l}]$ in $G_k$ is a triple $(g_k, G_{k-1}, G_k)$ satisfying the properties (EXTI) and (EXTII) of Definition \ref{D:Extension} (with the notation of (EXTI) and (EXTII)), as well as:
\indent \begin{description}
\item[(SWITCHIII):] $Dg_k$ induces an isomorphism from $PI(G_{k-1})$ to $PI(G_k)$ defined by
    $$PI(G_{k-1}) \xrightarrow{d_{k-1, j_k} \mapsto d^a_k=d_{k, i_k}} PI(G_k)$$
    ($d_{k-1,t} \mapsto d_{k,t}$ for $d_{k-1,t} \neq d_{k-1, j_k}$) and extended linearly over edges.
\item[(SWITCHIV):] $d_{k-1, i_k}$ labels the red nonperiodic vertex in $G_{k-1}$.
\item[(SWITCHV):] The single red edge of $G_{k-1}$ is $[d_{(k-1,i_k)}, d_{(k-1,l)}]$.
\end{description}
\end{df}

\vskip10pt

\noindent \textbf{Standard Switch Terminology and Notation:} Let $G_k$ be a Type (*) admissible LTT structure for a Type (*) pIW graph $\mathcal{G}$ with rose base graph $\Gamma_k$ and let $(g_k, G_{k-1}, G_k)$ be a switch, as in Definition \ref{D:Switch}.
\begin{itemize}
\item As with extensions, a switch is \emph{admissible} if $G_k$ and $G_{k-1}$ are both birecurrent (and thus actually Type (*) admissible LTT structures for $\mathcal{G}$).
\item Again $G_{k-1}$ is the \emph{source LTT structure} and $G_k$ the \emph{destination LTT structure}.
\item Again we call $[d^a_k,d_{k,j}]$ is called the \emph{(purple) edge determining the switch $(g_k, G_{k-1}, G_k)$}.
\item $d^{pu}_{k-1}$ again denotes $d_{k-1, j_k}$, though it does not label the red vertex of $G_{k-1}$ in the case of a switch (as it had in the case of an extension).
\item $d_{k-1, i_k}$ is denoted by both $d^{pa}_{k-1}$ and $d^u_{k-1}$ as, in this case, $d_{k-1, i_k}$ is also the label on the red nonperiodic (unachieved) direction vertex of $G_{k-1}$.
\item The red edge $e^R_{k-1}=[t^R_{k-1}]$ in $G_{k-1}$ can thus be written as $[d^u_{k-1}, d^a_{k-1}]$ or $[d^{pa}_{k-1}, d^a_{k-1}]$ or $[d^{pa}_{k-1}, d_{k-1, l}]$, among other ways.
\item $g_k$ is still called the \emph{ingoing generator} and can again be written $g_k: e^{pu}_k \mapsto e^a_k e^u_k$.
\end{itemize}

\vskip15pt

As explained in Subsection \ref{S:Peels}, but with this subsection's notation, a switch transforms LTT structures as follows:

\vskip5pt

\begin{figure}[H]
\centering
\includegraphics[width=5in]{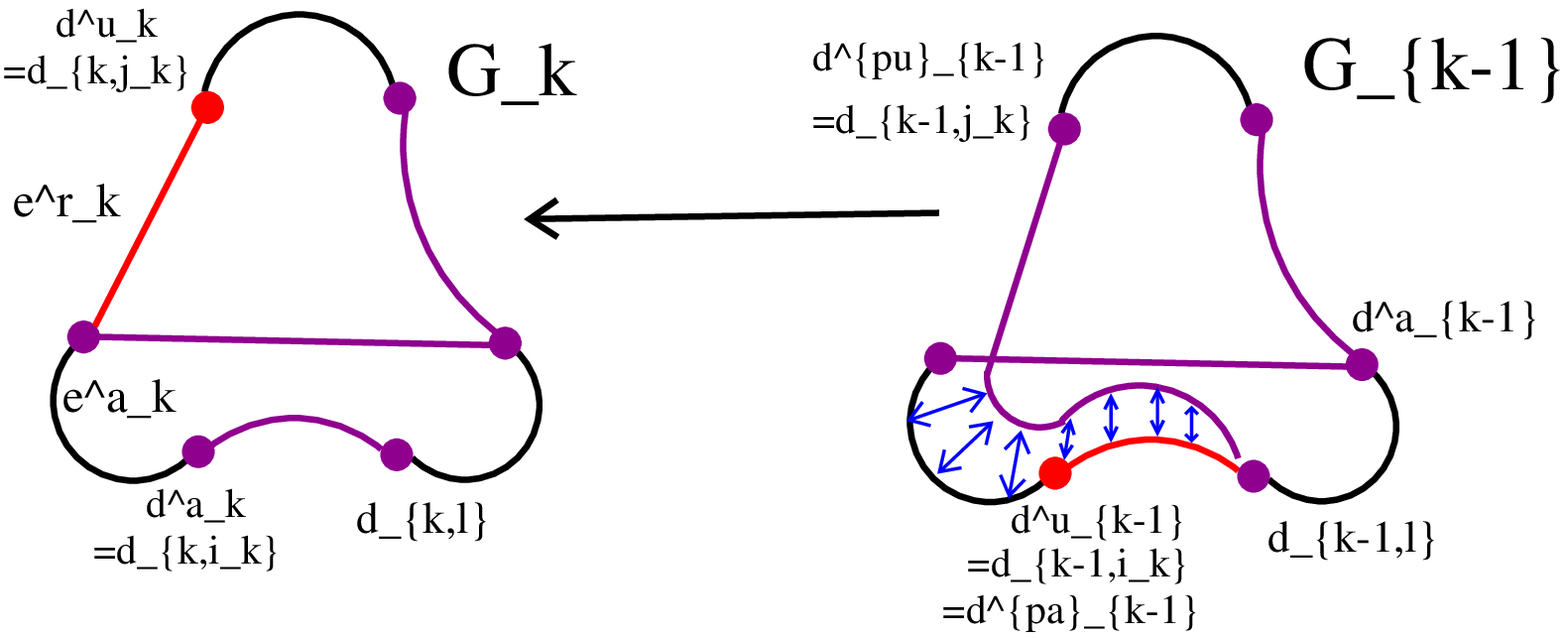}
\caption{{\small{Switch}}}
\label{fig:SwitchDiagram}
\end{figure}

\begin{rk}{\label{R:SwitchRestriction}} For similar reasons as to why we could not have $d_{k,l}= \overline{d^u_k}$ for an extension, we cannot have $d_{k,l}= d^a_k$ for a switch.
\end{rk}

\begin{rk}{\label{R:SwitchUniqueness}} One can realize the uniqueness of the potential switch associated to a purple edge $[d^a_k, d_{k,l}]$ as follows.  $PI(G_{k-1})$ is determined by the isomorphism in (SWITCHIII).  Since $G_{k-1}$ must be a Type (*) LTT structure, [LTT(*)1] implies it has only a single red vertex and [LTT(*)2] implies it has a unique red edge. The label on the red vertex is dictated by (SWITCHIV) and the vertices of the red edge are dictated by (SWITCHV).  Since the black edges of an LTT structure connect precisely vertices with inverse labels ((EXTII) dictates inverse by vertex label indices), and we have already determined the colored edges and vertex labels, the LTT structure $G_{k-1}$ is uniquely determined by Definition \ref{D:Switch}. $g_k$ is uniquely determined by (EXTII).
\end{rk}

\begin{df}{\label{D:PotentialSwitch}} For a Type (*) admissible LTT structure $G_k$ with Type (*) pIWG $\mathcal{G}$ and rose base graph, we obtain the \emph{potential switch} $(g_k, G_{k-1}, G_k)$ associated to a purple edge $[d^a_k, d_{k,l}]$ by the following procedure (justification in Lemma \ref{R:SwitchUniqueness}):
\begin{enumerate}
\item Start with $PI(G_k)$.
\item Replace each vertex label $d_{k,i}$ with $d_{k-1,i}$ and vertex label $\overline{d_{k,i}}$ with $\overline{d_{k-1,i}}$.
\item Switch the attaching (purple) vertex of the red edge to be $d_{k-1,l}$.
\item Switch the labels $d_{(k-1,j_k)}$ and $d_{(k-1,i_k)}$, so that the red vertex of $G_{k-1}$ will be $d_{k-1,i_k}$ and the red edge of $G_{k-1}$ will be $[d_{(k-1,i_k)}, d_{(k-1,l)}]$.
\item Include black edges connecting inverse pair labeled vertices (there is a black edge $[d_{(k-1,i)}, d_{(k-1,j)}]$ in $G_{k-1}$ if and only if there is a black edge $[d_{(k,i)}, d_{(k,j)}]$ in $G_k$).
\end{enumerate}
\end{df}

\begin{df}{\label{D:PotentialSwitch}} If $G_k$ is a Type (*) admissible LTT structure for a Type (*) pIW graph $\mathcal{G}$ with rose base graph, then we call the triple $(g_k, G_{k-1}, G_k)$ obtained in this manner the \emph{potential switch} associated to the purple edge $[d^a_k, d_{k,l}]$.
\end{df}

\noindent We will call a composition of switches and extensions an ``(admissible) composition'':

\begin{df} An \emph{(admissible) composition} for a Type (*) pIW graph $\mathcal{G}$ is a pair $(g_{i-k}, \dots, g_i, G_{i-k-1}, \dots, G_i)$, with $k \geq 0$, such that $g_{i-k,i}$ can be written
$$\Gamma_{i-k-1} \xrightarrow{g_{i-k}} \Gamma_{i-k}
\xrightarrow{g_{i-k+1}} \cdots \xrightarrow{g_{i-1}}\Gamma_{i-1} \xrightarrow{g_i} \Gamma_i,$$
with an \emph{associated sequence of LTT structures} for $\mathcal{G}$
$$G_{i-k-1} \xrightarrow{D^T(g_{i-k})} G_{i-k} \xrightarrow{D^T(g_{i-k+1})} \cdots \xrightarrow{D^T(g_{i-1})} G_{i-1} \xrightarrow{D^T(g_i)} G_i$$
where, for each $i-k-1 \leq j < i$,
\begin{itemize}
\item[(1)] $(g_{j+1}, G_j, G_{j+1})$ is either an (admissible) extension or switch for $\mathcal{G}$ and
\item[(2)] either $d^u_j=d^{pa}_j$ or $d^u_j=d^{pu}_j$.
\end{itemize}
\end{df}

\begin{prop}{\label{P:ExtensionsSwitches}} (Admissible) switches and extensions for a Type (*) pIWG $\mathcal{G}$ satisfy AM Properties I-VII. Conversely, the (admissible) switches and extensions for a Type (*) pIWG $\mathcal{G}$ are the only possible triples $(g_k, G_{k-1}, G_k)$ for $\mathcal{G}$ satisfying:
\begin{enumerate}
\item $G_{k-1}$ and $G_k$ are admissible Type (*) LTT structures for $\mathcal{G}$ with respective base graphs $\Gamma_{k-1}$ and $\Gamma_k$;
\item AM Properties I-VII; and
\item the vertex labels are mapped according to $Dg_k$.
\end{enumerate}
\noindent In particular, in the circumstance where $d^u_{k-1}= d^{pa}_{k-1}$, the triple is a switch and, in the circumstance where $d^u_{k-1}= d^{pu}_{k-1}$, the triple is an extension.
\end{prop}

\noindent \emph{Proof}: We start with the forward direction.

Since we have required that the extensions and switches be admissible, $G_{k-1}$ and $G_k$ are both birecurrent.  This gives us AM Property I.

Notice that the first and second parts of AM Property II are equivalent and that the second part holds by (EXTIIb) and (EXTIV) in the case of an extension and (EXTIIb), together with (SWITCHIV), in the case of a switch.  That there is only a single red unachieved direction vertex labeled $d^u_k$ in $G_k$ and that there is only a single red unachieved direction vertex labeled $d^u_{k-1}$ in $G_{k-1}$ follows from the requirement in (EXTI) that $G_k$ and $G_{k-1}$ are Type (*) LTT structures (the notation of (EXTII) makes this notation consistent with that in the AM properties).

What is left of AM Property III is that the edge $[t^R_k]= [d^u_k, \overline{d^a_k}]$ in $G_k$ and the edge $[t^R_{k-1}]= [d^u_{k-1}, \overline{d^a_{k-1}}]$ in $G_{k-1}$ are both red.  This follows from the notation of (EXTII) with (EXTV) in the extension case and (SWITCHV) in the switch case.

We now prove AM Property IV.  First notice that $d_{k,j_k} \neq d_{k,i_k}$ (as this would make the red edge of $G_k$, [$d_{(k,j_k)}, \overline{d_{(k,i_k)}}$], equal to [$d_{(k,i_k)}, \overline{d_{(k,i_k)}}$] contradicting the birecurrency of $G_k$) and $d_{k,l} \neq d_{k,i_k}$ (as then the purple edge determining the extension would be [$d_{(k,i_k)}, d_{(k,i_k)}$], contradicting SGTT2) and $d_{k,l} \neq d_{k,j_k}$ (as then the edge determining the extension would be [$d_{(k,i_k)}, d_{(k,j_k)}$], which would be red and not purple).

\indent Consider the extension case.  (EXTV) implies that the red edge of $G_{k-1}$ is $[d_{(k-1,l)},d_{(k-1,j_k)}]$.  By (EXTIIc) and the fact that $d_{k,l} \neq d_{k,i_k}$, since $Dg_k(d_{k-1,m})=d_{k,m}$ for all $m \neq j_k$, we know that $D^Cg_{k}$([$d_{(k-1,l)}, d_{(k-1,j_k)}$]) $=$ [$d_{(k,l)},d_{(k,i_k)}$].  This edge must be purple in $G_k$ since $d_{k,j_k}$ labels the only red vertex of $G_k$ and we showed above that $d_{k,j_k} \neq d_{k,i_k}$ and $d_{k,l} \neq d_{k,j_k}$.

\indent Now consider the switch case.  (SWITCHV) implies that the red edge of $G_{k-1}$ is $[d_{(k-1,l)},d_{(k-1,i_k)}]$.  Since (EXTII) still implies $Dg_k(d_{k-1,m})=d_{k,m}$ for all $m \neq j_k$ and we still have that $d_{k,l} \neq d_{k,i_k}$, we know $D^Cg_k($[$d_{(k-1,l)}, d_{(k-1,i_k)}$]$) = $[$d_{(k,l)},d_{(k,i_k)}$].  This edge must be purple, exactly as in the extension case.

For AM Property V, notice that AM Property III implies that $e^R_k$ is a red edge containing the red vertex $d^u_k$.  LTT(*)1 and LTT(*)2 imply the uniqueness of the red edge and red direction.

Since AM Property VI follows from (EXTII) and AM Property VII follows from (EXTIII) in the extension case and (SWITCHIII) in the switch case, we have proved the forward direction.

We now prove the converse.  Consider a triple $(g_k, G_{k-1}, G_k)$ satisfying AM Properties I-VII, as well as the other conditions in the proposition statement.  We will first show that the triple is either a switch or an extension, as $G_{k-1}$ and $G_k$ are birecurrent by AM Property I.

Assumption (1) in the proposition statement implies (EXTI).

By AM Property VI, $g_k$ is defined by  $g_k:e^{pu}_{k-1} \mapsto e^a_k e^u_k$ (where $g_k(e_{k-1,i})=e_{k,i}$ for $e_{k-1,i} \neq (e^{pu}_{k-1})^{\pm 1}$, $D_0(e^u_k)= d^u_k$, $D_0(\overline{e^a_k})= \overline{d^a_k}$, and $e^{pu}_{k-1}= e_{(k-1,j)}$, where $e^u_k=e_{k,j}$).  This gives us (EXTII).

By AM Property VII, $Dg_k$ induces on isomorphism from $SW(G_{k-1})$ to $SW(G_k)$.  Since the only direction whose second index is not fixed by $Dg_k$ is $d^{pu}_{k-1}$, the only vertex label of $SW(G_{k-1})$ that is not determined by this isomorphism is the preimage of $d^a_k$ (which AM Property IV dictates to be either $d^{pu}_{k-1}$ or $d^{pa}_{k-1}$).  When the preimage is $d^{pa}_{k-1}$, this gives us (EXTIII).  When the preimage is $d^{pu}_{k-1}$, this gives (SWITCHIII).  For the isomorphism to extend linearly over edges, we need that the image of an edge in $G_{k-1}$ is an edge in $G_k$, i.e. $[Dg_k(d_{k-1,i}), Dg_k(d_{k-1,j})]$ is an edge in $G_k$ for each edge $[d_{(k-1,i)}, d_{(k-1,j)}]$ in $G_{k-1}$.  This follows from AM Property IV.

AM Property II tells us that either $d^u_{k-1}= d^{pa}_{k-1}$ or $d^u_{k-1}=d^{pu}_{k-1}$.  When we are in the switch case, the above arguments tell us that $d^{pu}_{k-1}$ labels a purple periodic vertex, so we must have that $d^u_{k-1}=d^{pa}_{k-1}$ (since AM Property III tells us $d^u_{k-1}$ is red).  This gives us (SWITCHIV) once one appropriately coordinates the notation.  In the extension case, the above arguments tell us that instead $d^{pa}_{k-1}$ labels a purple periodic vertex, meaning that $d^u_{k-1}=d^{pu}_{k-1}$ (again since AM Property III tells us $d^u_{k-1}$ is red).  This gives us (EXTIV).  We are now only left with (EXTV) and (SWITCHV).

By AM Property V, $G_{k-1}$ has a single red edge $[t^R_{k-1}] = [\overline{d^a_{k-1}}, d^u_{k-1}]$.  By AM Property IV, the image of $[t^R_{k-1}]$ is a purple edge in $G_k$.  First consider what we established is the switch case, i.e. assume $d^u_{k-1}=d^{pa}_{k-1}$.  The goal is to determine that $[t^R_{k-1}]$ is $[d_{(k-1,i_k)}, d_{(k-1,l)}]$, where $d^a_k=d_{k,i_k}$ ($d_{k-1,i_k}=d^{pa}_{k-1}$) and $[d^a_k, d_{k,l}]$ is a purple edge in $G_k$ (making $(g_k, G_{k-1}, G_k)$ the switch determined by $[d^a_k, d_{k,l}]$).  Since $d^u_{k-1}=d^{pa}_{k-1}$, we know $[t^R_{k-1}]= [\overline{d^a_{k-1}}, d^u_{k-1} ]=[\overline{d^a_{k-1}}, d^{pa}_{k-1}]$.  We know $\overline{d^a_{k-1}} \neq d^{pa}_{k-1}$ (since (STTG2) implies $\overline{d^a_{k-1}} \neq d^u_{k-1}$, which equals $d^{pa}_{k-1}$).  Thus, AM Property VI says $D^C([t^R_{k-1}]) = D^C([\overline{d^a_{k-1}}, d^{pa}_{k-1}]) = [d_{(k,l)}, d^{a}_{k}]$ where $\overline{d^a_{k-1}}=e_{k-1,l}$.  So $[d_{(k,l)}, d^{a}_{k}]$ is a purple edge in $G_k$.  We thus have (SWITCHV).  Now consider what we established is the extension case, i.e. assume $d^u_{k-1}= d^{pu}_{k-1}$.  We need that the red edge $[t^R_{k-1}]$ is $[d_{(k-1,j_k)}, d_{(k-1,l)}]$, where $d^u_{k-1}=d_{k-1,j_k}$ and $[d^a_k, d_{k,l}]$ is a purple edge in $G_k$ (making $(g_k, G_{k-1}, G_k)$ the extension determined by $[d^a_k, d_{k,l}]$).  Since $d^u_{k-1}= d^{pu}_{k-1}$, we know $[t^R_{k-1}]= [\overline{d^a_{k-1}}, d^u_{k-1} ]= [\overline{d^a_{k-1}}, d^{pu}_{k-1}]$.  We know $\overline{d^a_{k-1}} \neq d^{pu}_{k-1}$ (since (STTG2) implies $\overline{d^a_{k-1}} \neq d^u_{k-1}$, which equals $d^{pu}_{k-1}$).  Thus, by AM Property VI, $D^C([t^R_{k-1}]) = D^C([\overline{d^a_{k-1}}, d^{pu}_{k-1}]) = [d_{(k,l)}, d^{a}_{k}]$, where $\overline{d^a_{k-1}}= e_{k-1,l}$.  We thus have (EXTV). \newline
\noindent QED.

\vskip8pt

\begin{df} An (admissible) switch or extension ($g_i$, $G_{i-1}$, $G_i$) such that $G_{i-1}$ and $G_i$ are Type (*) (admissible) LTT structures for $\mathcal{G}$ will be called an \emph{(admissible) generator triple} for $\mathcal{G}$.  (If both $G_{i-1}$ and $G_i$ are birecurrent, then they will actually be admissible Type (*) LTT structures for $\mathcal{G}$ and ($g_i$, $G_{i-1}$, $G_i$) will be an admissible generator triple).

Two generator triples ($g_i$, $G_{i-1}$, $G_i$) and ($g_i'$, $G_{i-1}'$, $G_i'$) will be considered \emph{equivalent} if they are equivalent as generating triples in the sense of Definition \ref{D:GeneratorExtendsToLTTStructures}.

\indent If ($g_i$, $G_{i-1}$, $G_i$) is an (admissible) generator triple, then we call both $g_i$ and its corresponding automorphism an \emph{(admissible) generator}.
\end{df}

\begin{lem}{\label{L:ExtensionUniqueness}} The extension $(g_k, G_{k-1}, G_k)$ associated to a purple edge $[d^a_k, d_{k,l}]$ is unique up to generator triple equivalence.  $G_{k-1}$ can be obtained from $G_k$ by:
\begin{enumerate}
\item Replacing each vertex label $d_{k,i}$ with $d_{k-1,i}$ and each vertex label $\overline{d_{k,i}}$ with $\overline{d_{k-1,i}}$.
\item Removing the interior of the red edge from $G_k$.
\item Adding a red edge $e^R_k$ connecting the red vertex to $d_{k-1,l}$.
\end{enumerate}
\end{lem}

\noindent Proof: $PI(G_{k-1})$ is uniquely determined by the isomorphism in (EXTIII) to differ from $PI(G_k)$ by the relabeling of vertices described in (1).  Since $G_{k-1}$ must be a Type (*) LTT structure, [LTT(*)1] implies it has precisely one red vertex and [LTT(*)2] implies it has a unique red edge. The label on the red vertex is dictated by (EXTIV) to be $d_{(k-1,j)}$, where $d^u_k=d_{(k,j)}$, and the red edge is dictated to be $[d_{(k-1,j)},d_{(k-1,l)}]$ by (EXTV).  Since the black edges of an LTT structure connect precisely vertices with inverse labels, and we have already determined the colored edges and vertex labels, the LTT structure $G_{k-1}$ is uniquely determined by Definition \ref{D:Extension}. $g_k$ is uniquely determined by (EXTII). It is clear that the procedure gives us the structure $G_{k-1}$ described.  \newline
\noindent QED.

\vskip10pt

\subsection{Construction Compositions}{\label{S:ConstructionCompositions}}

\noindent It is not enough for a representative $g$ to be composed of permitted switches and extensions, for it to start and end with the same LTT structure, and for its transition matrix to be Perron-Frobenius.  We also need that $IW(g) \cong \mathcal{G}$.  To ensure that $IW(g) \cong \mathcal{G}$, we use ``building block'' compositions of extensions called ``construction compositions:''

\begin{df}  An \emph{(admissible) construction composition} for a Type (*) pIW graph $\mathcal{G}$ is a pair $(g_{i-k}, \dots, g_i; G_{i-k-1}, \dots, G_i)$, together with a decomposition
$$\Gamma_{i-k-1} \xrightarrow{g_{i-k}} \Gamma_{i-k} \xrightarrow{g_{i-k+1}} \cdots \xrightarrow{g_{i-1}}\Gamma_{i-1} \xrightarrow{g_i} \Gamma_i,$$
and sequence of LTT structures for $\mathcal{G}$ \newline
\indent $G_{i-k-1} \xrightarrow{D^T(g_{i-k})} G_{i-k} \xrightarrow{D^T(g_{i-k+1})} \cdots \xrightarrow{D^T(g_{i-1})} G_{i-1} \xrightarrow{D^T(g_i)} G_i,$ where:
\begin{enumerate}
\item each $(g_j, G_j, G_{j+1})$ with $i-k \leq j \leq i$ is an (admissible) extension,
\item $(g_{i-k}, G_{i-k-1}, G_{i-k})$ is an (admissible) switch, and
\item $PI(G_j)=\mathcal{G}$ for each $i-k-1 \leq j \leq i$.
\end{enumerate}
\noindent We call the composition of generators $g_{i,i-k}=g_i \circ \dots \circ g_{i-k}$ a \emph{construction automorphism}, $G_{i-k-1}$ the \emph{source LTT structure} and $G_k$ the \emph{destination LTT structure}.  \newline
\indent If $k=1$, the composition is simply an (admissible) switch.  If we leave out the switch, we get a \emph{purified construction automorphism} $g_p=g_i \circ \dots \circ g_{i-k+1}$ and pair $(g_{i-k+1}, \dots, g_i; G_{i-k}, \dots, G_i)$ called a \emph{purified construction composition}.  Now $G_{i-k}$ is the \emph{source LTT structure} and $G_{i-k}$ is the \emph{destination LTT structure}.\newline
\indent An (admissible) construction composition $(g_{i-k}, \dots, g_i; G_{i-k-1}, \dots, G_i)$ is said to be \emph{realized} if there exists an ideally decomposed Type (*) representative $g:\Gamma' \to \Gamma'$ of $\phi \in Out(F_r)$ decomposed as $\Gamma' = \Gamma_0' \xrightarrow{g_1'} \Gamma_1' \xrightarrow{g_2'} \cdots \xrightarrow{g_{n-1}'} \Gamma_{n-1}' \xrightarrow{g_n'} \Gamma_n' = \Gamma'$ with the sequence of LTT structures $G_0' \xrightarrow{D^T(g_1')} G_1' \xrightarrow{D^T(g_2')} \cdots \xrightarrow{D^T(g_{n-1}')} G_{n-1}' \xrightarrow{D^T(g_n')} G_n'$, such that $g_j \cong g_j'$ for all $i-k \leq j \leq i$ and $G_j \cong G_j'$ for all $i-k-1 \leq j \leq i$.
\end{df}

\vskip8pt

A construction automorphism will always have the form of a Dehn twist automorphism $e^{pu}_{i-k-1} \mapsto w e^u_{i-k}$, where $w=e^a_{i-k} \dots e^a_i$, and one can view the composition as twisting the edge corresponding to $e^{pu}_{i-k-1}$ around the path corresponding to $w$ in the destination LTT structure.  Since construction compositions are in ways analogous to Dehn twist mapping class group representatives and since many construction methods for pseudo-Anosov  mapping classes (including those of Penner in [P88]) used Dehn twists, it was natural to look into properties of construction compositions.  Their properties and connections to Dehn Twists certainly warrant further investigation and Clay and Pettet use them in their fully irreducible construction methods.  However, our methods here utilize a single special construction compositions property, i.e. that they, in some sense, ``construct'' smooth paths in the destination LTT structures (see Proposition \ref{P:PathConstruction}).  Since we use the paths in our procedure for constructing (ideal Whitehead graph)-yielding representatives, we include their definition.

However, before giving the definition, we first note that we abuse notation throughout this section by dropping indices.  While not necessary, this abuse may make the visual aspects of the properties and procedures much clearer, as well as reduce the potential for confusion over indices.

\begin{df} A \emph{construction path} associated to a construction composition \newline
\noindent $(g_1, \dots, g_k; G_1, \dots, G_k)$ is a path in $G_k$ starting with the red vertex $d^u_k$, transversing the red edge $[d^u_k, \overline{d^a_k}]$ from the red vertex $d^u_k$ to the vertex $\bar d^a_k$, transversing the black edge $[\overline{d^a_k}, d^a_k]$ from the vertex $\overline{d^a_k}$ to the vertex $d^a_k$, transversing the purple edge $[d^a_k, d_k]=[d^a_k, \overline{d^a_{k-1}}]$ from $d^a_k$ to $d_k= \overline{d^a_{k-1}}$, transversing the black edge $[\overline{d^a_{k-1}}, d^a_{k-1}]$ from the vertex $\overline{d^a_{k-1}}$ to the vertex $d^a_{k-1}$, transversing the purple edge $[d^a_{k-1}, d_{k-1} ]= [d^a_{k-1}, \overline{d^a_{k-2}}]$ from the vertex $d^a_{k-1}$ to the vertex $d_{k-1}= \overline{d^a_{k-2}}$, transversing the black edge $[\overline{d^a_{k-2}}, d^a_{k-2}]$ from the vertex $\overline{d^a_{k-2}}$ to the vertex $d^a_{k-2}$, continues as such through the purple edges determining each $g_i$, and finally ending by transversing $[d^a_2, d_2]= [d^a_2, \overline{d^a_1}]$ and the black edge $[\overline{d^a_1}, d^a_1]$ from the vertex  $\overline{d^a_1}$ to the vertex $d^a_1$.
\end{df}

\begin{lem}{\label{L:ConstructionPathSmooth}} The construction path associated to a realized construction composition \newline
 $(g_{i-k}, \dots, g_i; G_{i-k-1}, \dots, G_i)$ for a Type (*) pIW graph $\mathcal{G}$, with
\begin{itemize}
\item decomposition $\Gamma_{i-k-1} \xrightarrow{g_{i-k}} \Gamma_{i-k} \xrightarrow{g_{i-k+1}} \cdots \xrightarrow{g_{i-1}}\Gamma_{i-1} \xrightarrow{g_i} \Gamma_i$ and
\item LTT structures $G_{i-k-1} \xrightarrow{D^T(g_{i-k})} G_{i-k} \xrightarrow{D^T(g_{i-k+1})} \cdots \xrightarrow{D^T(g_{i-1})} G_{i-1} \xrightarrow{D^T(g_i)} G_i$
\end{itemize}
is the smooth path $[d^u_i, \overline{d^a_i}, d^a_i, \overline{d^a_{i-1}}, d^a_{i-1}, \dots, d^a_{s+1}, \overline{d^a_{i-k}}, d^a_{i-k}]$ in the LTT structure $G_i$. \end{lem}

\noindent \emph{Proof}:  We will proceed by induction for decreasing $s$ values.  Since nothing in the proof will rely on $G_{i-k-1}$ (which is the only thing that distinguishes that $(g_{i-k}, G_{i-k-1}, G_{i-k})$ is a switch instead of an extension), proof by induction is valid here.

\indent For the base case realize that $[d^u_i, \overline{d^a_i}]$ is the red edge in $G_i$.  So $[d^u_i, \overline{d^a_i}, d^a_i]$ is a path in $G_i$ and is smooth because it alternates between colored and black edges ($[d^u_i, \overline{d^a_i}]$ is colored and $[\overline{d^a_i}, d^a_i]$ is black).  For the sake of induction assume that, for $i > s > i-k$, \newline
\noindent $[d^u_i, \overline{d^a_i}, d^a_i, \overline{d^a_{i-1}}, d^a_{i-1}, \dots, d^a_{s+1}, \overline{d^a_s}, d^a_s]$ is a smooth path in $G_i$ (ending with the black edge $[\overline{d^a_s}, d^a_s]$).

\indent  The red edge that $g_{s-1}$ creates in $G_{s-1}$ is $[d^u_{s-1}, \overline{d^a_{s-1}}]$ (see Definition \ref{D:CreatesEdge} the proof of Corollary \ref{C:UnachievedDirection}).  By Lemma \ref{L:PurpleEdgeImages}, $D^Cg_s([d^u_{s-1}, \overline{d^a_{s-1}}])= [d^a_s, \overline{d^a_{s-1}}]$ is a purple edge in $G_s$. Since purple edges are always mapped to themselves by extensions (in the sense that $D^C$ preserves the second index of their vertex labels) and $D^Cg_{s}([d^u_{s}, \overline{d^a_{s-1}}])=[d^a_{s}, \overline{d^a_{s-1}}]$ is a purple edge in $G_s$, $D^Cg_{n,s}(\{d^u_{s-1}, \overline{d^a_{s-1}} \}) = D^Cg_{n,s+1}(D^Cg_{s}([d^u_{s-1}, \overline{d^a_{s-1}}]))$ $= D^Cg_{n,s+1}([d^a_s, \overline{d^a_{s-1}}])=[d^a_s, \overline{d^a_{s-1}}]$ is a purple edge in $G_i$.  Thus, including the purple edge $[d^a_s, \overline{d^a_{s-1}}]$ in the smooth path $[d^u_i, \overline{d^a_i}, d^a_i, \overline{d^a_{i-1}}, d^a_{i-1}, \dots, d^a_{s+1}, \overline{d^a_s}, d^a_s]$ gives the smooth path $[d^u_i, \overline{d^a_i}, d^a_i, \overline{d^a_{i-1}}, d^a_{i-1}, \dots, d^a_{s+1}, \overline{d^a_s}, d^a_s, \overline{d^a_{s-1}}].$  (This path is smooth because we added a colored edge to a path with edges alternating between colored and black that ended with a black edge).  By including the black edge $[\overline{d^a_{s-1}}, d^a_{s-1}]$ we get the construction path \newline
\noindent $[d^u_i, \overline{d^a_i}, d^a_i, \overline{d^a_{i-1}}, d^a_{i-1}, \dots, d^a_s, \overline{d^a_{s-1}}, d^a_{s-1}].$
(Again this path is smooth because we added a black edge to a path with edges alternating between colored and black that ended with a colored edge).

This concludes the inductive step and hence the proof. \newline
\noindent QED.

\begin{df} Let $G$ be an admissible Type (*) LTT structure with red vertex $d^u$. The \emph{construction subgraph} $G_C$ is constructed from $G$ via the following procedure:
\begin{itemize}
\item Start by removing the interior of the black edge $[e^u]$, the purple vertex $\overline{d^u}$, and the interior of any purple edges containing the vertex $\overline{d^u}$.  Call the graph with these edges and vertices removed $G^1$.
\item Given $G^{j-1}$, recursively define $G^j$ as follows: Let $\{ \alpha_{j-1,i} \}$ be the set of vertices in $G^{j-1}$ not contained in any colored edge of $G^{j-1}$.  $G^j$ will be the subgraph of $G^{j-1}$ obtained by removing all black edges containing a vertex $\alpha_{j-1,i} \in \{ \alpha_{j-1,i} \}$, as well as the interior of each purple edge containing a vertex of the form $\overline{\alpha_{j-1,i}}$.
\item $G_C=\underset{j}{\cap} G^j$.
\end{itemize}
\end{df}

\vskip10pt

\begin{lem}{\label{L:ConstructionAutomorsphismFromPath}} Let $G$ be an admissible Type (*) LTT structure.  Consider a smooth path $\gamma = [d^u, \overline{x_1}, x_1, \overline{x_2}, x_2, \dots, x_{k+1}, \overline{x_{k+1}}]$ in the construction subgraph $G_C$ starting with $e^R$ (oriented from $d^u$ to $\overline{d^a}$) and ending with the black edge $[x_{k+1}, \overline{x_{k+1}}]$.

Mark $r$-petaled roses $\Gamma_{i-k-1}, \dots, \Gamma_i$ with edge sets $\mathcal{E}_k = \mathcal{E}(\Gamma_k)$ denoted \newline
\noindent $\{E_{(k,1)},
\overline{E_{(k,1)}}, E_{(k,2)}, \overline{E_{(k,2)}}, \dots , E_{(k,r)}, \overline{E_{(k,r)}} \} =
\{e_{(k,1)}, e_{(k,2)}, \dots, e_{(k,2r-1)}, e_{(k,2r)}\}$ \newline
\noindent so that, for each $i$, $G$ can be viewed as having base graph $\Gamma_i$ where $e^u=e_{(i,s)}$ and $d^u=D_0(e^u)$.

Define the homotopy equivalences $\Gamma_{i-k-1} \xrightarrow{g_{i-k}} \Gamma_{i-k} \xrightarrow{g_{i-k+1}} \cdots \xrightarrow{g_{i-1}}\Gamma_{i-1} \xrightarrow{g_i} \Gamma_i$
by $g_l: e_{l-1,s} \mapsto e_{l,t_l} e_{l,s}$, where $D_0(e_{l,t_l})= \overline{x_{i-l+1}}$, and $g_l(e_{l-1,j})=e_{l,j}$ for $e_{l-1,j} \neq e_{l-1,s}^{\pm}$.

Define the LTT structures with respective base graphs $\Gamma_j$ (and maps between) \newline
\noindent $G_{i-k-1} \xrightarrow{D^T(g_{i-k})} G_{i-k} \xrightarrow{D^T(g_{i-k+1})} \cdots \xrightarrow{D^T(g_{i-1})} G_{i-1} \xrightarrow{D^T(g_i)} G_i$ by having
\indent \begin{description}
\item 1. each $PI(G_l)$ isomorphic to $PI(G_i)$ via an isomorphism preserving the second indices of the vertex labels,
\item 2. the second index of the vertex label on the single red vertex in each $G_l$ be ``s'' (the same as in $G_i$), and
\item 3. the single red edge in $G_l$ be $[d_{l,s}, \overline{d_{l,t_l}}]$.
\end{description}

If each $G_j$ is a Type (*) LTT structure for $\mathcal{G}$ with base graph $\Gamma_j$, then \newline
\noindent $(g_{i-k}, \dots, g_i; G_{i-k-1}, \dots, G_i)$ is a purified construction composition.  In fact, it is the unique realized purified construction composition with $\gamma$ as its construction path.
\end{lem}

\noindent \emph{Proof}: We need to show that $(g_{i-k}, \dots, g_i; G_{i-k-1}, \dots, G_i)$ is indeed a construction composition, that its construction path is $[d^u_i, \overline{d^a_i}, d^a_i, \overline{d^a_{i-1}}, d^a_{i-1}, \dots, d^a_{i-k+1}, \overline{d^a_{i-k}}, d^a_{i-k}]$, and that it is the unique construction composition with that path.

We first show that each $(g_l, G_{l-1}, G_l)$ is an extension.  (EXT I) is ensured by our requirement that each $G_j$ is a Type (*) LTT structure with rose base graph.  The $G_l$ are all Type (*) LTT structures for $PI(G)$ by (1)-(3) in the lemma statement. The second index of the single red vertex label is the same in each $G_l$ as in $G_i$, giving (EXT IV).  (EXT II)(a) holds by how we determined our notation.  (EXT II)(b) holds by how we determined our notation and the construction in the lemma statement.  (EXT III) is true by construction (by (1), in particular).  (EXT V) is true by construction (by (3), in particular).

The construction path is $[d^u_i, \overline{d^a_i}, d^a_i, \overline{d^a_{i-1}}, d^a_{i-1}, \dots, d^a_{i-k+1}, \overline{d^a_{i-k}}, d^a_{i-k}]$ by Lemma \ref{L:ConstructionPathSmooth}.

\indent That the $G_l$ must be as stated follows from
\noindent \begin{description}
\item 1. each $PI(G_l)$ being isomorphic to $PI(G_i)$ via an isomorphism preserving the second indices of the vertex labels in order for the $(g_l, G_{l-1}, G_l)$ to be extensions
\item 2. the $G_l$ being Type (*) LTT structures,
\item 3. the second index of the red vertices being the same, making each $(g_l, G_{l-1}, G_l)$ an extension, and
\item 4. knowing, by Lemma \ref{L:ConstructionPathSmooth}, that the attaching vertex for $e^R_l$ in $G_l$ must be $x_{i-l+1}$.
\end{description}
Once each $G_l$ is determined, that $g_l$ must be as stated follows from AM Property VI. \newline
\noindent QED.

\begin{df} We call $\Gamma_{i-k-1} \xrightarrow{g_{i-k}} \Gamma_{i-k} \xrightarrow{g_{i-k+1}} \cdots \xrightarrow{g_{i-1}}\Gamma_{i-1} \xrightarrow{g_i} \Gamma_i$, together with its sequence of LTT structures $G_{i-k-1} \xrightarrow{D^T(g_{i-k})} G_{i-k} \xrightarrow{D^T(g_{i-k+1})} \cdots \xrightarrow{D^T(g_{i-1})} G_{i-1} \xrightarrow{D^T(g_i)} G_i$, as in the lemma, the \emph{construction composition corresponding} to the path $\gamma = [d^u, \overline{x_1}, x_1, \overline{x_2}, x_2, \dots, x_{k+1}, \overline{x_{k+1}}]$.
\end{df}

\begin{ex}{\label{E:ConstructionAutomorphism}} In the following LTT structure, $G$, for Graph XX, the numbered edges give a construction path associated to the construction automorphism $a \mapsto ab\bar{c}\bar{c}bbcb$ (all other edges are fixed by the automorphism) in the LTT structure.  In the following figure, $A$ denotes $\bar{a}$, etc. We will continue with this notation throughout the document.

\begin{figure}[H]
\centering
\includegraphics[width=1.5in]{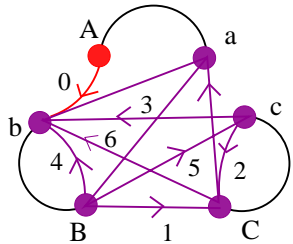}
\caption{{\small{\emph{Construction path associated to construction automorphism $a \mapsto ab\bar{c}\bar{c}bbcb$}}}}
\label{fig:ConstructionPathExample}
\end{figure}

In Figure \ref{fig:ConstructionCompositionSequence} we show the construction composition corresponding to the construction path of Figure \ref{fig:ConstructionPathExample}. The source LTT structure of the switch is left out to highlight the fact that it does not affect the construction path.

\noindent \begin{figure}[H]
\centering
\noindent \includegraphics[width=6.3in]{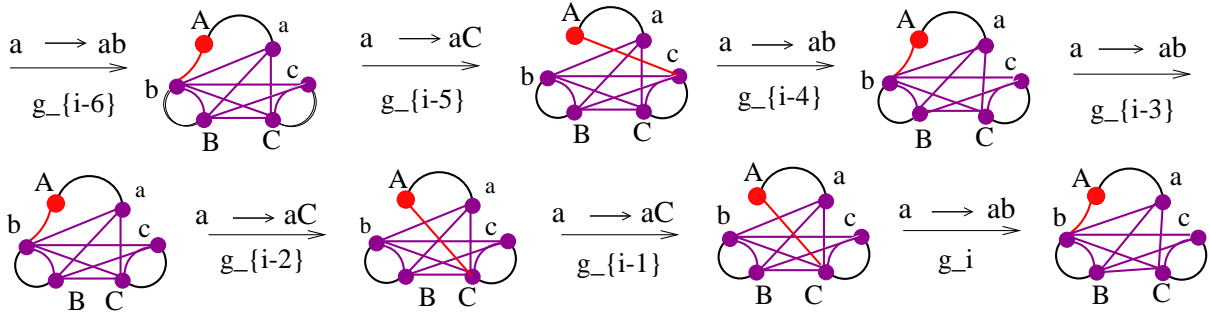}
\caption{{\small{\emph{Construction path associated to Figure \ref{fig:ConstructionPathExample} construction composition}}}}
\label{fig:ConstructionCompositionSequence}
\end{figure}

We determine the sequence of LTT structures in the construction composition by attaching the red edge in $G_{i-k}$ to the terminal vertex of edge $k$ in the construction path.  The generator can be determined by the red edge in its destination LTT structure: If the red vertex of $G_j$ is $d_s$ and the red edge is $[d_s,d_t]$, then $g_j$ is defined by $e_s \mapsto \bar{e_t}e_s$.
\end{ex}

\vskip10pt

\noindent The following proposition tells us that construction paths are ``built'' by construction compositions.  By saying a turn is \emph{taken by $g_{(k,l)}$}, we will mean that the turn is taken by some $g_{k,l}(e_{l-1,i})$.

\begin{prop}{\label{P:PathConstruction}}  Let $g:\Gamma \to \Gamma$ be an ideally decomposed Type (*) representative of $\phi \in Out(F_r)$ with $IW(\phi)=\mathcal{G}$.  Suppose that $g$ is decomposable as $\Gamma = \Gamma_0 \xrightarrow{g_1} \Gamma_1 \xrightarrow{g_2} \cdots \xrightarrow{g_{n-1}}\Gamma_{n-1} \xrightarrow{g_n} \Gamma_n = \Gamma$,
with the sequence of LTT structures for $\mathcal{G}$: \newline
\indent $G_{i-k-1} \xrightarrow{D^T(g_{i-k})} G_{i-k} \xrightarrow{D^T(g_{i-k+1})} \cdots \xrightarrow{D^T(g_{i-1})} G_{i-1} \xrightarrow{D^T(g_i)} G_i$. \newline
\noindent Assume the Standard Notation for a Type (*) LTT Structure.  If $g'=g_n \circ \dots \circ g_{k+1}$ is a construction composition, then $\mathcal{G}$ contains as a subgraph the purple edges in the construction path for $g'$. \end{prop}

\noindent \emph{Proof}: We will proceed by induction for decreasing $k$ values.  Since nothing in the proof will rely on $G_k$ (which is the only thing that distinguishes that $(g_k, G_k, G_{k+1})$ is a switch instead of an extension), proof by induction is valid here.

\indent For the base case consider $g_n \circ g_{n-1}$. By the  Corollary \ref{C:UnachievedDirection} proof, $g_{n-1}$ creates the red edge $[d^u_{n-1}, \overline{d^a_{n-1}}]$ in $G_{n-1}$.  We know that $g_n$ is defined by $g_n$:$e^{pu}_{n-1} \mapsto e^a_n e^u_n$ and $g_n(e_{n-1,l})= e_{n,l}$ for all $e_{n-1,l} \neq (e^{pu}_{n-1})^{\pm 1}$.  Thus, since $d^{pu}_{n-1}= d^u_{n-1} \neq \overline{d^a_{n-1}}$, we know that $Dg_n(\overline{d^a_{n-1}}) = \overline{d^a_{n-1}}$. So $D^Cg_n([d^u_{n-1}, \overline{d^a_{n-1}}]) = D^Cg_n([d^{pu}_{n-1}, \overline{d^a_{n-1}}]) = [d^a_n, \overline{d^a_{n-1}}]$ and, since $D^Cg_n$ images of purple and red edges of $G_{n-1}$ are purple edges of $G_n$, $[d^a_n, \overline{d^a_{n-1}}]$ is a purple edge in $G_n$.  The base case is proved.

\indent For the inductive step assume that, for $n > s > k+1$, $G_n$ contains as a subgraph the purple edges in the construction path associated to $g_{n,s}$.  The red edge that $g_{s-1}$ creates in $G_{s-1}$ is $[d^u_{s-1}, \overline{d^a_{s-1}}]$ (see the proof of Corollary \ref{C:UnachievedDirection}).  As above, $D^Cg_s([d^u_{s-1}, \overline{d^a_{s-1}}])= [d^a_s, \overline{d^a_{s-1}} ]$ is represented by a purple edge in $G_s$. Since purple edges are always  mapped to themselves by extensions and $D^Cg_{s}([d^u_{s}, \overline{d^a_{s-1}}])=[d^a_{s}, \overline{d^a_{s-1}}]$, $D^Cg_{n,s}([d^u_{s-1}, \overline{d^a_{s-1}}]) = D^Cg_{n,s+1}(D^tg_{s}([d^u_{s-1}, \overline{d^a_{s-1}}])) = D^Cg_{n,s+1}([d^a_s, \overline{d^a_{s-1}}])=[d^a_s, \overline{d^a_{s-1}}]$, proving the inductive step. The proposition is proved.  \newline
\noindent QED.

\subsection{Switch Paths}

While they do not give insight into the progress of building $\mathcal{G}$ and have more restrictions, switch sequences also have associated paths.  The usefulness of switch paths lies in the aid they give in constructing the switch sequences required in our methods below.  This subsection focuses on switch sequences and their associated switch paths.

\emph{We will continue in this section with the abuse of notation from the previous subsection (this mainly consists of ignoring second indices).}

\vskip8pt

\begin{df} A \emph{switch sequence} for a Type (*) pIW graph $\mathcal{G}$ is a pair \newline
\noindent $(g_{i-k}, \dots, g_i; G_{i-k-1}, \dots, G_i)$ with a decomposition of $g_{(i,i-k)}$ as \newline
\indent $\Gamma_{i-k-1} \xrightarrow{g_{i-k}} \Gamma_{i-k} \xrightarrow{g_{i-k+1}} \cdots \xrightarrow{g_{i-1}}\Gamma_{i-1} \xrightarrow{g_i} \Gamma_i$, \newline
\noindent and associated sequence of LTT structures \newline
\indent $G_{i-k-1} \xrightarrow{D^T(g_{i-k})} G_{i-k} \xrightarrow{D^T(g_{i-k+1})} \cdots \xrightarrow{D^T(g_{i-1})} G_{i-1} \xrightarrow{D^T(g_i)} G_i$, where:
\indent{\begin{description}
\item  [(SS1)] each  $(g_j, G_{j-1}, G_j)$ with $i-k \leq j \leq i$ is a permitted switch,
\item  [(SS2)] $PI(G_j)=\mathcal{G}$ for each $i-k-1 \leq j \leq i$, and
\item  [(SS3)] $d^a_{n+1}=d^u_n \neq d^u_l=d^a_{l+1}$ and $\overline{d^a_l} \neq
d^u_n=d^a_{n+1}$ for all $i \geq n > l \geq i-k$.
\end{description}}
\noindent Sometimes we will just call $\Gamma_{i-k-1} \xrightarrow{g_{i-k}} \Gamma_{i-k} \xrightarrow{g_{i-k+1}} \cdots \xrightarrow{g_{i-1}}\Gamma_{i-1} \xrightarrow{g_i} \Gamma_i$ a switch sequence when either the LTT structures should be clear from the decomposition or are unnecessary for discussion. \newline
\indent We call the composition of generators $g_{i,i-k}=g_i \circ \dots \circ g_{i-k}$ a \emph{switch sequence automorphism} with \emph{source LTT structure} $G_{i-k-1}$ and \emph{destination LTT structure} $G_i$.
\end{df}

\begin{rk} (SS3) is not implied by (SS1) and (SS2) and is necessary for a switch path to indeed be a path.  Certain statements in the proof of Lemma \ref{L:SwitchPathSmooth} below (where we show that the switch path corresponding to a switch sequence is realized as a smooth path in the destination LTT structure) would be incorrect without (3).
\end{rk}

\begin{df} A \emph{realized switch sequence} is a pair $(g_{i-k}, \dots, g_i; G_{i-k-1}, \dots, G_i)$ such that there exists an ideally decomposed Type (*) representative $\Gamma = \Gamma_0 \xrightarrow{g_1} \Gamma_1 \xrightarrow{g_2} \cdots \xrightarrow{g_{n-1}}\Gamma_{n-1} \xrightarrow{g_n} \Gamma_n = \Gamma$ of $\phi \in Out(F_r)$ having the sequence of LTT structures $G_{0} \xrightarrow{D^T(g_{1})} G_{1} \xrightarrow{D^T(g_{2})} \cdots \xrightarrow{D^T(g_{n-1})} G_{n-1} \xrightarrow{D^T(g_n)} G_n$ and satisfying that, for each $i-k < j \leq i$, $(g_j, G_{j-1}, G_j)$ is a switch.

Again we call the composition of generators $g_{i,i-k} = g_i \circ \dots \circ g_{i-k}$ a \emph{switch sequence automorphism} with \emph{source LTT structure} $G_{i-k-1}$ and \emph{destination LTT structure} $G_i$.
\end{df}

\begin{df} A \emph{switch path} associated to a (realized) switch sequence $(g_j, \dots, g_k; G_{j-1}, \dots, G_k)$ is a path in the destination LTT structure $G_k$ for $g_k$ that starts with the red vertex $d^u_k$, transverses the red edge [$d^u_k, \overline{d^a_k}$] for $g_k$ from the red vertex $d^u_k$ to the vertex $\overline{d^a_k}$, transverses the black edge [$\overline{d^a_k}, d^a_k$] from the vertex $\overline{d^a_k}$ to the vertex $d^a_k$, transverses what is the red edge [$d^u_{k-1}, \overline{d^a_{k-1}}$] = [$d^a_k, \overline{d^a_{k-1}}$] in $G_{k-1}$ (and a purple edge in $G_k$) from $d^a_k=d^u_{k-1}$ to $\overline{d^a_{k-1}}$, transverses the black edge [$\overline{d^a_{k-1}}, d^a_{k-1}$] from the vertex $\overline{d^a_{k-1}}$ to the vertex $d^a_{k-1}$, continues as such through all of the new red edges for the $g_i$ with $j \leq i \leq j$, and ends by transversing the black edge [$\overline{d^a_{j+1}}, d^a_{j+1}$] from the vertex $\overline{d^a_{j+1}}$ to the vertex $d^a_{j+1}$ and what is the red edge [$d^u_j, \overline{d^a_j}$] = [$d^a_{j+1}, \overline{d^a_j}$] in $G_j$ (purple edge in $G_k$), and then the black edge [$\overline{d^a_j}, d^a_j$] from the vertex  $\overline{d^a_j}$ to the vertex $d^a_j$.

In other words, a switch path alternates between the red edges (oriented from the unachieved direction $d^u_j$ to $\overline{d^a_j}$) for the $G_j$ (for descending $j$) and the black edges between.
\end{df}

\begin{rk} We clarify here some ways switch paths and construction paths differ.
\begin{description}
\item[(1)] Switch paths look like construction paths but, while the purple edges in the construction path for a construction composition $(g_{i-k}, \dots, g_i; G_{i-k-1}, \dots, G_i)$ are purple in each $G_l$ with $i-l \leq l < i$, for a switch path, they will be red edges in the structure $G_l$ they are created in and then will not exist at all in the structures $G_m$ with $m<l$.  The change of color of red edges and then disappearance of edges is the reason for (3) in the switch sequence definition.
\item[(2)] Unlike constructions paths, switch paths do not give subpaths of lamination leaves.
\end{description}
\end{rk}

\begin{lem}{\label{L:SwitchPathSmooth}} The switch path associated to a realized switch sequence $(g_j, \dots, g_k; G_{j-1}, \dots, G_k)$ forms a smooth path in the LTT structure $G_k$. \end{lem}

\noindent \emph{Proof}: The red edge in $G_k$ is [$d^u_k, \overline{d^a_k}$].  We are left to show (by induction) that:

(1) For each $1 \leq l < k$, [$d^u_{l}, \overline{d^a_{l}}]=[d^a_{l+1}, \overline{d^a_{l}}$] is a purple edge of $G_k$ and

(2) the purple edges [$d^a_{l+1}, \overline{d^a_{l}}$] (together with the black edges in the switch sequence) form a smooth path in $G_k$.

Start with the base case.  By the switch properties, the red edge in $G_{k-1}$ is [$d^u_{k-1}, \overline{d^a_{k-1}}]=[d^a_k, \overline{d^a_{k-1}}$].  Since $d^a_k \neq d^u_k$ and $\overline{d^a_{k-1}} \neq d^u_k$ (by the switch sequence definition), $D^tg_k(\{d^a_k, \overline{d^a_{k-1}}\})=\{d^a_k, \overline{d^a_{k-1}}\}$.  Thus, by Lemma \ref{L:RedEdgeImage}, $[d^a_k, \overline{d^a_{k-1}}$], is a purple edge in $G_k$.  The red edge in $G_k$ is $[d^u_k, \overline{d^a_k}]$.  So, by including the black edge [$\overline{d^a_k}, d^a_k$], we have a path $[d^u_k, \overline{d^a_k}, d^a_k, \overline{d^a_{k-1}} ]$ in $G_k$.  This path is smooth since it alternates between colored and black edges.  So our proof of the base case is complete.

We now prove the inductive step.  By the inductive hypothesis we assume that the sequence of switches associated to $g_k,\dots ,g_{k-i}$ gives us a smooth path $[d^u_k, \dots, \overline{d^a_{k-i}}]$ in $G_k$ ending with a purple edge with ``free'' vertex $\overline{d^a_{k-i-1}}$.  We know that the red edge in $G_{k-i-1}$ is $[d^u_{k-i-1}, \overline{d ^a_{k-i-1}}]= [d^a_{k-i}, \overline{d^a_{k-i-1}}$].  As long as we do not have $d^u_l=d^a_{k-i}$ or $d^u_l= \overline{d^a_{k-i-1}}$ for $k-i \leq l \leq k$ (which holds again by the definition of a switch sequence), $D^tg_{k,k-i}(\{d^u_{(k-i-1)}, \overline{d^a_{(k-i-1)}}\})= D^tg_{k,k-i}(\{d^a_{(k-i)}, \overline{d^a_{(k-i-1)}}\}) = \{d^a_{(k-i)}, \overline{d^a_{(k-i-1)}}\}$.  This, as above, makes $[d^a_{(k-i)}, \overline{d^a_{(k-i-1)}}]$ a purple edge in $G_k$ by Lemma \ref{L:RedEdgeImage}.

Since $[d^u_k, \dots, \overline{d^a_{k-i}}]$ is a smooth path in $G_k$ ending with a black edge, \newline
\noindent $[d^u_k, \dots, \overline{d^a_{k-i}}, d^a_{k-i}, \overline{d^a_{k-i-1}}]$ is also a smooth path in $G_k$, [$\overline{d^a_{k-i}},d^a_{k-i}$] is a black edge in $G_k$ and as [$d^a_{k-i}, \overline{d^a_{k-i-1}}$] is a purple edge in $G_k$.\newline
\noindent QED.

\begin{ex}{\label{E:SwitchPath}}  We return to the LTT structure $G$ (for Graph XX) of Example \ref{E:ConstructionAutomorphism} and number below the colored edges of a switch path.

\begin{figure}[H]
\centering
\includegraphics[width=1.3in]{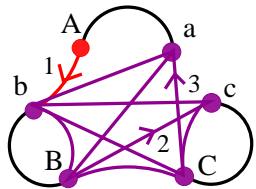}
\caption{{\small{\emph{Switch path in $G$ with colored edges numbered}}}}
\label{fig:SwitchPath}
\end{figure}

\noindent A switch sequence for $G$ constructed from the switch path is:

\begin{figure}[H]
\centering
\noindent \includegraphics[width=6in]{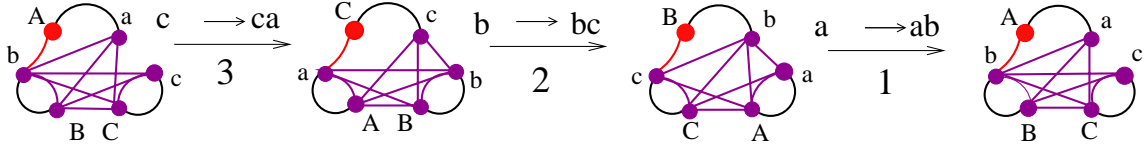}
\caption{{\small{\emph{Switch sequence constructed from the switch path in Figure \ref{fig:SwitchPath}}}}}
\label{fig:SwitchSequence}
\end{figure}
\end{ex}

Notice how the red edge in the destination LTT structure for the generator labeled by (1) is the first colored edge in the switch path, the red edge in the destination LTT structure for (2) is the second colored edge in the switch path, etc.

\section{AM Diagrams}{\label{Ch:AMDiagrams}}

In this section we describe how to construct the ``AM Diagram'' for a Type (*) pIWG, as well as prove that Type (*) representatives are realized as loops in these diagrams.

\begin{df} Let $\mathcal{G}$ be a given Type (*) pIW graph.  A \emph{PreAdmissible Map Diagram (PreAM Diagram) for $\mathcal{G}$ or \emph{$PreAMD(\mathcal{G})$})} is the directed graph where
\begin{itemize}
\item [1.] the nodes correspond to equivalence classes of admissible LTT structures for $\mathcal{G}$
\item [2.] for each admissible generator triple ($g_i$, $G_{i-1}$, $G_i$) for $\mathcal{G}$, there exists a directed edge $E(g_i, G_{i-1}, G_i)$ from the node [$G_{i-1}$] to the node [$G_i$].
\end{itemize}
The disjoint union of the maximal strongly connected subgraphs of $PreAMD(\mathcal{G})$ will be called the \emph{Admissible Map Diagram for $\mathcal{G}$} (or \emph{$AMD(\mathcal{G})$}).
\end{df}

\indent The following proposition shows that we can restrict our search for representatives to loops in AM Diagrams.  For the proposition and henceforth after, $E(g_1, G_{0}, G_1)$ will denote the oriented edge from [$G_{0}$] to [$G_1$] labeled by $g_1$ and $E(g_1, G_{0}, G_1) * \dots * E(g_k, G_{k-1}, G_k)$ will denote the oriented path in $AMD(\mathcal{G})$ from [$G_0$] to [$G_k$] that transverses the edges corresponding to the generators in the order they are composed (starting with the edge $E(g_1, G_{0}, G_1)$ and concluding with the $E(g_k, G_{k-1}, G_k)$).  We say that $E(g_1, G_{0}, G_1) * \dots * E(g_k, G_{k-1}, G_k)$ \emph{realizes} $g$ or that $g$ is \emph{realized} by $E(g_1, G_{0}, G_1) * \dots * E(g_k, G_{k-1}, G_k)$.  In particular, we say that the decomposition $\Gamma_{0} \xrightarrow{g_{1}} \Gamma_{1} \xrightarrow{g_{2}} \cdots \xrightarrow{g_{k-1}}\Gamma_{k-1} \xrightarrow{g_k} \Gamma_k$ of $g$ is \emph{realized} by $E(g_1, G_{0}, G_1) * \dots * E(g_k, G_{k-1}, G_k)$ in $AMD(\mathcal{G})$.

\begin{prop}{\label{P:ReferenceLoop}} If $g=g_{k} \circ \cdots \circ g_1: \Gamma \to \Gamma$ is an ideally decomposed Type (*) representative of $\phi \in Out(F_r)$ such that $IW(\phi)=\mathcal{G}$, with corresponding sequence of LTT structures $G_{0} \xrightarrow{D^T(g_{1})} G_{1} \xrightarrow{D^T(g_{2})} \cdots \xrightarrow{D^T(g_{k-1})} G_{k-1} \xrightarrow{D^T(g_k)} G_k$, then $E(g_1, G_{0}, G_1) * \dots * E(g_k, G_{k-1}, G_k)$ forms an oriented loop in $AMD(\mathcal{G})$. \end{prop}

\noindent \emph{Proof}: Suppose that $g$ is such a representative.  We showed in Proposition \ref{P:ExtensionsSwitches} that $g$ is a composition of permitted switches $g_i$ and permitted extensions $g_i$, and thus of admissible maps.  This tells us that $E(g_1, G_{0}, G_1) * \dots * E(g_k, G_{k-1}, G_k)$ forms an oriented loop in $PreAMD(\mathcal{G})$.  Since all loops of a graph are contained in the union of the maximal strongly connected subgraphs of the graph, we know that  $E(g_1, G_{0}, G_1) * \dots * E(g_k, G_{k-1}, G_k)$ is actually in $AMD(\mathcal{G})$, proving the proposition. \newline
\noindent QED.

\begin{df} We denote the loop $E(g_1, G_{0}, G_1) * \dots * E(g_k, G_{k-1}, G_k)$ by $L(g_1, \dots, g_k)$.
\end{df}

\vskip5pt

\begin{cor} \textbf{(of Proposition \ref{P:ReferenceLoop})} If no loop in $AMD(\mathcal{G})$ gives a Type (*) representative of an ageometric, fully irreducible outer automorphism $\phi \in Out(F_r)$ such that $IW(\phi) = \mathcal{G}$, then such a $\phi$ does not exist.  In particular, any of the following properties of an AM Diagram would prove that such a representative does not exist:
\begin{itemize}
\item [(1)] There is at least one edge direction pair $\{d_i, \overline{d_i}\}$, where $e_i \in \mathcal{E}(\Gamma)$, such that no red vertex in $AMD(\mathcal{G})$ is labeled by either $d_i$ or $\overline{d_i}$.
\item [(2)] The representative corresponding to each loop in $AMD(\mathcal{G})$ has a PNP.
\end{itemize}
\end{cor}

\noindent \emph{Proof}:  Proposition \ref{P:IdealDecomposition} says that such a $\phi$ would have an ideally decomposed Type (*) representative $g$ and Proposition \ref{P:ReferenceLoop} shows that any ideally decomposed Type (*) representative would be realized by a loop in $AMD(\mathcal{G})$. Thus, $g$ has a realization $L(g_1, \dots, g_m)$ in $AMD(\mathcal{G})$.

If, for some $1 \leq i \leq r$, $L(g_1, \dots, g_m)$ did not contain an LTT structure $G_k$ where either $d^u_k=d^i$ or $d^u_k=\overline{d^i}$, then the corresponding automorphism would fix the generator of $F_r$ corresponding to $E_i$, which would make $g$ reducible, contradicting that $\phi \in Out(F_r)$ is fully irreducible.  So (1) is proved.

Since Type (*) representatives must be PNP-free, if no loop in $AMD(\mathcal{G})$ realizes a PNP-free automorphism, then no Type (*) representative exists.  This proves (2) and thus the entire corollary since the first sentence is a direct consequence of the Proposition \ref{P:ReferenceLoop}. \newline
\noindent QED. \newline

\section{Full Irreducibility Criterion}{\label{Ch:FIC}}

The main goal of this section is the proof of a ``Folk Lemma'' giving a criterion, the ``Full Irreducibility Criterion (FIC),'' for an irreducible train track map to represent a fully irreducible outer automorphism.  Our original approach to proving the criterion involved the ``Weak Attraction Theorem,'' several notions of train tracks, laminations, and the basin of attraction for a lamination.  However, Michael Handel graciously provided a way to finish the proof making much of our initial work unnecessary.  The proof of the criterion we give here uses Michael Handel's recommendation.

\subsection{Free Factor Systems, Filtrations, and RTTs}{\label{S:RTTs}}

The following definitions are necessary to understand the definition of a relative train track representative for an outer automorphism.  While [BH92] gives that we always have train track representatives for irreducible outer automorphisms, this is not true for reducible outer automorphisms.  Relative train tracks were invented by Bestvina and Handel to approximate train tracks as best possible in this circumstance.  We use relative train tracks in our proof of the Full Irreducibility Criterion.

We begin by defining a ``free factor system'' for a free group, $F_r$, of rank $r$.

\begin{df} [BFH00]
$\mathcal{F}=\{[[F^1]], \dots, [[F^k]]\}$  is a \emph{free factor system} for $F_r$ if $F^1*F^2* \dots *F^k$ is a free factor of $F_r$ and each $F^i$ is nontrivial.  For free factor systems $\mathcal{F}_1$ and $\mathcal{F}_2$, we say that $\mathcal{F}_1 \sqsubset \mathcal{F}_2$ when, for each $[[F^i]] \in \mathcal{F}_1$, there exists some $[[F^j]] \in \mathcal{F}_2$ such that $[[F^i]] \sqsubset [[F^j]]$, i.e. $F^i$ is conjugate to a free factor of $F^j$.
\end{df}

A distinguishing characteristic of reducible outer automorphism representatives is the existence of proper nontrivial invariant subgraphs.  A relative train track representative of such an outer automorphism will have a ``filtration'' of invariant subgraphs ``realizing'' a nested sequence of free factors.  Over the course of the next few definitions we describe what this means.

\begin{df}  [BH92]
For a topological representative $g: \Gamma \to \Gamma$ of $\phi \in Out(F_r)$, an increasing sequence of $g$-invariant subgraphs $\emptyset =\Gamma_0 \subset \Gamma_1 \subset \dots \subset \Gamma_k = \Gamma$ such that each component of each subgraph contains at least one edge is called a \emph{filtration}.  For such a filtration, the closure $H_t$ of $\Gamma_t - \Gamma_{t-1}$ is called the \emph{$t^{th}$ stratum}.  Let $\mathcal{E}^+_t= \{E^t_1, \dots, E^t_{n_t}\}$ denote the set of edges of $H_t$ with some prescribed orientation and let $\mathcal{E}_t= \{E^t_1, \overline{E^t_1}, \dots, E^t_{n_t}, \overline{E^t_{n_t}}\}$.  The \emph{transition submatrix} for the stratum $H_t$ is the square matrix $M_t$ such that, for each $i$ and $j$, the $ij^{th}$ entry is the number of times $g(E^t_j)$ crosses over $E^t_i$ in either direction.  Strata with zero matrices as their transition submatrices  are called \emph{zero strata}.  Let $\lambda_t$ denote the \emph{Perron-Frobenius eigenvalue} for $M_t$, i.e. the real eigenvalue of largest norm.  Then $H_t$ is called \emph{exponentially growing (EG)} if $\lambda_t > 1$ and a \emph{nonexponentially growing (NEG)} if $\lambda_t=1$.
\end{df}

\noindent We are now ready for the relative train track representative definition, as defined in [BH92].

\begin{df} [BH92] A \emph{Relative Train Track (RTT) Representative} of an outer automorphism $\phi \in Out(F_r)$ is a topological representative $g: \Gamma \to \Gamma$ and filtration $\emptyset =\Gamma_0 \subset \Gamma_1 \subset \dots \subset \Gamma_k = \Gamma$ such that:
\begin{description}
\item[(1)] Each $M_t$ is either the zero matrix or is irreducible;
\item[(2)] all vertices have valence greater than one; and
\item[(3)] each EG-stratum satisfies:
\begin{itemize}
\item[(a)] for each edge $E \in \mathcal{E}_t$, the first edge of $g(E)$ is in $\mathcal{E}_t$,
\item[(b)] $g_{\#}(\beta)$ is nontrivial for each nontrivial path $\beta \subset \Gamma_{t-1}$ having endpoints in $\Gamma_{t-1} \cap H_t$, and
\item[(c)] $g(\gamma) \subset \Gamma_r$ is a \emph{t-legal path} (i.e. $\Gamma_{t-1}$ contains each of its illegal turns) for each legal path $\gamma \subset H_t$.
\end{itemize}
\end{description}
\end{df}

\noindent The following are needed for understanding the revised versions of RTTs used in proving the Full Irreducibility Criterion.

\begin{df} Suppose $\Gamma$ is a marked graph and $\Gamma_i$ is a subgraph with non-contractible components $C_1, \dots, C_k$. Then $\mathcal{F}(\Gamma_i)=\{\pi(C_1), \dots, \pi(C_k) \}$ is called the free factor system $\mathcal{F}(\Gamma_i)$ \emph{realized} by $\Gamma_i$.  A nested sequence of free factor systems $\mathcal{F}^1 \sqsubset \mathcal{F}^2 \sqsubset \dots \sqsubset \mathcal{F}^m$ is said to be \emph{realized} by an RTT $g: \Gamma \to \Gamma$ and filtration $\emptyset =\Gamma_0 \subset \Gamma_1 \subset \dots \subset \Gamma_k = \Gamma$ if each $\mathcal{F}^j$ is realized by some $F_{i_j}$. [BH92]

A topological representative $g: \Gamma \to \Gamma$ and filtration $\emptyset =\Gamma_0 \subset \Gamma_1 \subset \dots \subset \Gamma_k = \Gamma$ are called \emph{reduced} if each $H_t$ satisfies: For any $l>0$ and each $\phi^l$-invariant free factor system $\mathcal{F}$ such that $\mathcal{F}(\Gamma_{i-1}) \sqsubset \mathcal{F} \sqsubset \mathcal{F}(\Gamma_i)$, either $\mathcal{F}=\mathcal{F}(\Gamma_{i-1})$ or $\mathcal{F}= \mathcal{F}(\Gamma_i)$. [BFH00]
\end{df}

We will use a correspondence proved in [BFH00] between attracting laminations for an outer automorphism $\phi \in Out(F_r)$ and the EG-strata of an RTT representative $g: \Gamma \to \Gamma$ of  $\phi$: For each EG stratum $H_t$ of $g$, there exists a unique attracting lamination (denoted by $\Lambda_t$) having $H_t$ as the highest stratum crossed by the realization $\lambda \subset \Gamma$ of a $\Lambda_t$-generic line.  $H_t$ is called the \emph{EG-stratum determined by} $\Lambda_t \in \mathcal{L}(\phi)$.

\indent We will remind the reader of the definition of a revised version of a relative train track called a ``complete split relative train track (CT).''  These train tracks are defined by M. Feighn and M. Handel in [FH09].  However, we first give the definition of a ``complete splitting.''  Both these definitions are specialized definitions for the case of ageometric outer automorphisms (where there are, in particular, no closed PNPs).

\begin{df} [FH09]
\indent \emph A nontrivial path or circuit $\gamma$ is \emph{completely split} if it has a \emph{complete splitting}, i.e. can be written as $\gamma = \dots \gamma_{l-1}\gamma_l \dots$ where
\noindent \begin{itemize}
\item [(1)] each $\gamma_i$ is either a single edge in an irreducible stratum, an iNP, or a taken connecting path in a zero stratum and
\item [(2)] $g_{\#}^k(\gamma) = \dots g_{\#}^k(\gamma_{l-1}) g_{\#}^k(\gamma_l) \dots$.
\end{itemize}
\noindent In the case of a complete splitting we write $\gamma = \dots \gamma_{l-1} \bullet \gamma_l \dots$.
\end{df}

\begin{df}  [FH09]
An RTT representative $g: \Gamma \to \Gamma$ of $\phi \in Out(F_r)$, together with its filtration $\mathcal{F}$ given by $\emptyset = \Gamma_0 \subset \Gamma_1 \subset \dots \subset \Gamma_k=\Gamma$, is a \emph{Completely Split Relative Train Track (CT)} if it satisfies all of the following:
\begin{description}
\item [(CT1)] $g$ is rotationless.
\item [(CT2)] $g$ is completely split, that is:
\begin{itemize}
\item [(a)] $g(E)$ is completely split for each edge $E$ in each irreducible stratum and
\item [(b)] $g_{\#}(\sigma)$ is completely split for each taken connecting path $\sigma$ in a zero stratum.
\end{itemize}
\item [(CT3)] $\mathcal{F}$ is reduced and the cores of the $\Gamma_i$ are also filtration elements. (Recall that the \emph{core} of a finite graph $K$ is the subgraph consisting of all edges of $K$ crossed by some circuit in $K$).
\item [(CT4)] The endpoints of all iNPs are vertices (necessarily principal).  For each NEG-stratum $H_i$ and nonfixed edge of $H_i$, the terminal endpoint is principal (and hence fixed).
\item [(CT5)] Periodic edges are fixed and the endpoints of fixed edges are principal.  For a fixed stratum $H_t$ with unique edge $E_t$, either $E_t$ is a loop or each end of $E_t$ is in $\Gamma_{t-1}$ and $\Gamma_{t-1}$ is a core graph.
\item [(CT6)] For each zero stratum $H_i$, there is an EG-stratum $H_t$ (with $t>i$) such that:
\begin{itemize}
\item [(a)] $H_i$ is enveloped by $H_t$, i.e.
{\begin{itemize}
\item [(i)] for some $u<i<t$, $H_u$ is irreducible,
\item [(ii)] no component of $G_t$ is contractible,
\item [(iii)] $H_i$ is a component of $G_{t-1}$, and
\item [(iv)] each $H_i$-vertex is of valence greater than one in $G_t$.
\end{itemize}}
\item [(b)] Each edge of $H_i$ is \emph{$t$-taken} (i.e. is a maximal subpath of $g^k_{\#}(E)$ in $H_i$ for some $k>0$ and edge $E$ in $H_t$)
\item [(c)] $H_t$ contains every vertex of $H_i$ and the link of each vertex of $H_i$ is contained in $H_i \cup H_t$.
\end{itemize}
\item [(CT7)] There are no \emph{linear edges} (i.e. edges $E_i$ of NEG-strata $H_i$ such that $g(E_i)=E_iu_i$ for some nontrivial NP $u_i \in \Gamma_{i-1}$).
\item [(CT8)] The highest edges of iNPs do not belong to NEG-strata.
\item [(CT9)] Suppose that $H_t$ is an EG-stratum and that $\rho$ is a height-$t$ iNP.  Then the restriction of $g$ to $\Gamma_t$ is $\theta \circ g_{t-1} \circ g_t$ where:
\begin{itemize}
\item [(a)] $g_t: \Gamma_t \to \Gamma^1$ can be decomposed into proper extended folds defined by iteratively folding $\rho$,
\item [(b)] $g_{t-1}: \Gamma^1 \to \Gamma^2$ can be decomposed into folds involving edges in $\Gamma_{t-1}$, and
\item [(c)] $\theta:\Gamma^2 \to \Gamma_t$ is a homomorphism.
\end{itemize}
\end{description}
\end{df}

\begin{rk} If $\phi \in Out(F_r)$ is forward rotationless and $\mathcal{C}$ is a nested sequence of $\phi$-invariant free factor systems, then $\phi$ is represented by a CT $g: \Gamma \to \Gamma$ and filtration $\mathcal{F}$ that realizes $\mathcal{C}$ [FH09, Theorem 4.29].
\end{rk}

Before we can finally give our Fully Irreducibility Criterion proof, we need to remind the reader of the following.

\begin{df} [BFH00]
The \emph{complexity} of the free factor system $\mathcal{F}=\{[[F^1]], \dots, [[F^k]]\}$  is defined to be zero if $\mathcal{F}$ is trivial and is otherwise defined to be the non-increasing sequence of positive integers obtained by rearranging the set $\{rank(F^1), \dots, rank(F^k)\}$.  The set of all complexities is given a lexicographic ordering.

The \emph{free factor support} for a subset $B \subset \mathcal{B}$ is defined in [BFH00, Corollary 2.6.5] to be the unique free factor system of minimal complexity carrying every element of $\mathcal{B}$.
\end{df}

The only relevant information about the free factor support for our proof of the FIC is that, if a lamination is carried by a proper free factor, then its support is a proper free factor.  If this were not the case, then the free factor support would have to have rank $r$ (and thus complexity $\{r\}$), while the free factor carrying the lamination (because it is proper) must be of rank less than $r$, giving it complexity less than the free factor support and hence contradicting that a free factor support is of minimal complexity.

\smallskip

We would like to credit Michael Handel for his contributions to the proof of the following ``Folk lemma,'' which we will call the \emph{Full Irreducibility Criterion (FIC)}.

\begin{lem}{\label{L:FIC}} (The Full Irreducibility Criterion)
Let $g: \Gamma \to \Gamma$ be an irreducible train track representative of $\phi \in Out(F_r)$.  Suppose that
\begin{itemize}
\item [(I)] $g$ has no PNPs,
\item [(II)] the transition matrix for $g$ is Perron-Frobenius, and
\item [(III)] all $LW(x;g)$ for $g$ are connected.
\end{itemize}
\indent Then $\phi$ is a fully irreducible outer automorphism.
\end{lem}

\noindent \emph{Proof:} Suppose that $g: \Gamma \to \Gamma$ is an irreducible TT representative of $\phi \in Out(F_r)$ with Perron-Frobenius transition matrix, with connected local Whitehead graphs, and with no PNPs.  Since $g$ has a Perron-Frobenius transition matrix, as an RTT, it has precisely one stratum and that stratum is EG. Hence, it has precisely one attracting lamination [BFH00].  Since the number of attracting laminations belonging to a TT representative of $\phi$ is independent of the choice of representative, any representative of $\phi$ would also have precisely one attracting lamination.

Suppose, for the sake of contradiction, that $\phi$ were not fully irreducible.  Then some power $\phi^k$ of $\phi$ would be reducible.   If necessary, take an even higher power so that $\phi$ will also be rotationless (this does not change the reducibility).  Notice that, since $\mathcal{L}(\phi)$ is $\phi$-invariant, any representative of $\phi^k$ would also have precisely one attracting lamination.

Since $\phi^k$ is reducible (and rotationless), there exists a completely split train track representative $h: \Gamma' \to \Gamma'$ of $\phi^k$ with more than one stratum [FH09, Theorem 4.29].  Since $\phi^k$ has precisely one attracting lamination, $h$ will have precisely one EG-stratum $H_t$.  Each stratum $H_i$, other than $H_t$ and any zero strata (if they exist), would be an NEG-stratum consisting of a single edge $E_i$ [FH09, Lemma 4.22].  We will consider separately the cases where $t=1$ and where $t>1$.

Notice that, since any zero stratum has zero transition matrix (and thus must have every edge mapped to a lower filtration element by $h$), a zero stratum could not be $H_1$.  Thus, if $t>1$, then $H_1$ is NEG and must consist of a single edge $E_1$.  Since $H_1$ is bottom-most, it would have to be fixed, as there are no lower strata for its edge to be mapped into.  But then, according to (CT5), $E_1$ would have to be an invariant loop, which would mean that $\phi^k$ would have a rank-1 invariant free factor.  However, $g$ was PNP-free, which means that $g^k$ was PNP-free and thus that $\phi^k$ should not be able to have a rank-1 invariant free factor. We have thus reached a contradiction for the case where $t>1$.

Now assume that $t=1$.  This would imply that $\Lambda(\phi^k)(=\Lambda(\phi))$ is carried by a proper free factor.  Proposition 2.4 of [BFH97] states that, if a finitely generated subgroup $A \subset F_r$ carries $\Lambda_{\phi}$, then $A$ has finite index in $F_r$. The necessary conditions for this proposition are actually only: (1) the transition matrix of $g$ is irreducible and (2) at each vertex of $\Gamma$, the local Whitehead graph is connected. (Up to the contradiction in the proof of Proposition 2.4 of [BFH97], the only properties used in the proof are that the support is finitely generated, proper, and carries the lamination.  The contradiction uses Lemma 2.1 of [BFH97], which simply proves that properties (1) and (2) carry over to lifts of $g$ to finite-sheeted covering spaces, using no properties other than properties (1) and (2).)  The assumptions (1) and (2) are assumptions in the hypotheses of our criteria and $\Lambda$ is still the attracting lamination for $g$ and so we can apply the proposition to create a contradiction with the fact that $\Lambda$ has proper free factor support.  Applying the proposition, since proper free factors have infinite index, the support must be the whole group. This contradicts that the EG-stratum is $H_1$ and that there must be more than one stratum.

We have thus shown that we cannot have more than one stratum with $t=1$ or $t >1$.  So all powers of $\phi$ must be irreducible and thus $\phi$ is fully irreducible, as desired. \newline
\noindent QED.

\begin{rk} To apply this lemma, we need a procedure for proving the nonexistence of PNPs, as stated in Proposition \ref{P:NPIdentification} of the following section.
\end{rk}

\section{Nielsen Path Identification}{\label{Ch:NPIdentification}}

In this section we give a method for finding all iPNPs, thus PNPs, in our circumstance.  While we have not yet proved the method's finiteness, its application ended quickly in all examples thus far.

\begin{ex} As a warm-up to the procedure, we give the following example of the procedure applied to the construction composition $g$:

\noindent \begin{figure}[H]
\centering
\includegraphics[width=6.2in]{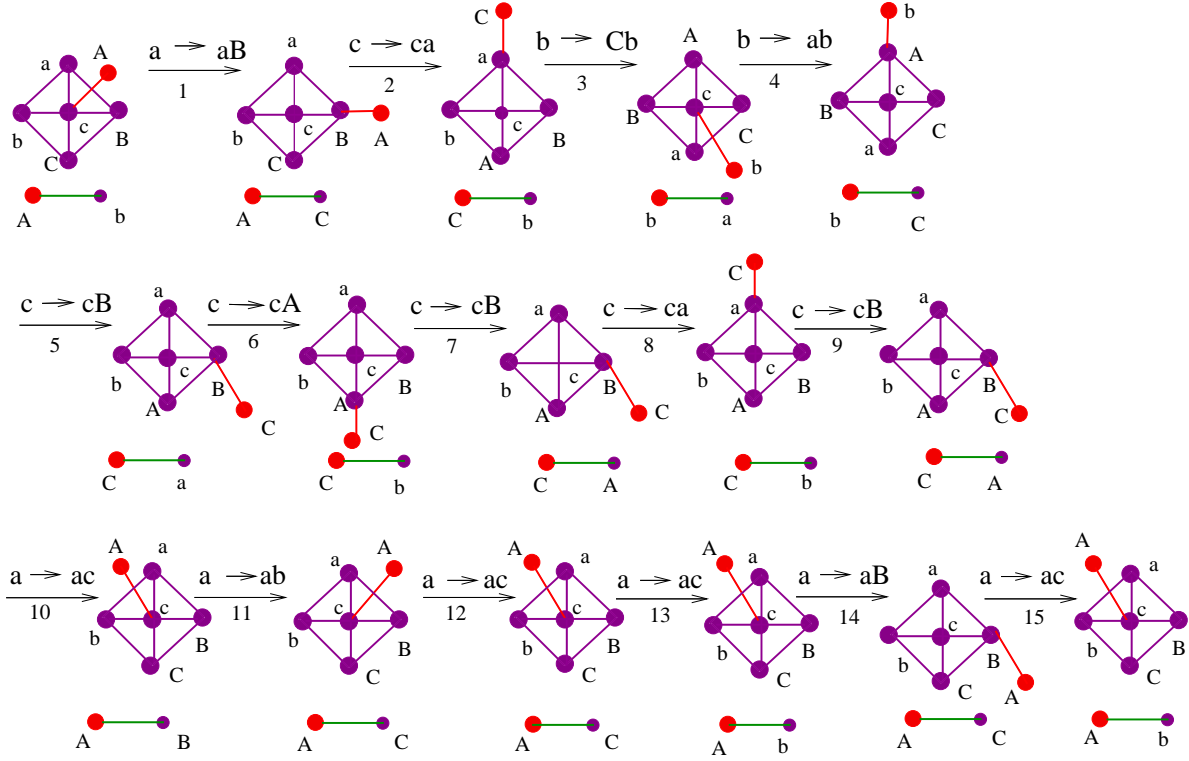}
\caption{{\small{\emph{(underneath each graph we included the green illegal turn from the augmented LTT structure)}}}}
\label{fig:NPAlgorithmExample}
\end{figure}

The green illegal turn for $g$ is $\{\overline{a}, b\}$.  We will try to build an iPNP $\overline{\rho_1}\rho_2$ where $\rho_1= \overline{a}\dots$ and $\rho_2=b\dots$ are legal paths.

Since $g_1(b)=b$ is a subpath of $g_1(\overline{a})= b\overline{a}$, we must add another edge after $b$ in $\rho_2$.  And, since $D_0(\overline{c})$ is red in $G_1$, the edge added after $b$ will be such that its initial direction is a preimage of $D_0(\overline{c})$ (the other direction in $T_2$) under $Dg_1$.  The only such preimage is $D_0(\overline{c})$.  Thus, the next edge of $\rho_2$ would have to be $\overline{c}$.

$g_{2,1}(\overline{a})=g_2(b)\overline{a}$ is a subpath of $g_{2,1}(b\overline{c})= g_2(b)\overline{a}\overline{c}$, so we must add another edge after $\overline{a}$ in $\rho_1$.  Since $D_0(\overline{a})$ is red in $G_2$, the edge added after $\overline{a}$ will be such that  its initial direction is a preimage of $D_0(b)$ (the other direction in $T_3$) under $Dg_{2,1}$.  There are two such preimages: Case 1 will be where the initial direction is $D_0(\overline{a})$ and Case 2 will be where the initial direction is $D_0(b)$. \newline
\noindent We try each of these as the next edge of $\rho_1$.

For Case 1, suppose that the next edge of $\rho_1$ is $\overline{a}$.  Since $g_{3,1}(b\overline{c})= g_{3,2}(b)g_3(\overline{a})\overline{c}$ is a subpath of $g_{3,1}(\overline{a}\overline{a})= g_{3,2}(b)g_3(\overline{a})\overline{c}b\overline{a}$, it follows that we must add another edge after $\overline{c}$ in $\rho_2$.  Since $D_0(b)$ is red in $G_3$, the new edge will have to be such that its initial direction is a preimage of $D_0(a)$ (the other direction in $T_4$) under $Dg_3$.  The only such preimage is $D_0(a)$.  So the next edge of $\rho_2$ would have to be $a$.

Now, $g_{4,1}(b\overline{c}a)= g_{4,2}(b)g_{4,3}(\overline{a})g_4(\overline{c})a\overline{b}\overline{a}c$ and $g_{4,1}(\overline{a}\overline{a})= g_{4,2}(b)g_{4,3}(\overline{a})g_4(\overline{c})ab\overline{a}$.  Since $\{b, \overline{b}\} \neq T_5$ (and $b \neq \overline{b}$), $\overline{a}\overline{a}$ and $b\overline{c}a$ could not start $\rho_1$ and $\rho_2$, respectively.

For Case 2, suppose that the next edge of $\rho_1$ is $b$.  Since $g_{3,1}(b\overline{c})= g_{3,2}(b)g_3(\overline{a})\overline{c}$ is a subpath of $g_{3,1}(\overline{a}b)= g_{3,2}(b)g_3(\overline{a})\overline{c}b$, it follows that we must add another edge of $\rho_2$ after $\overline{c}$.  As above, $D_0(b)$ is red in $G_3$, so we can follow the logic above and see that the next edge of $\rho_2$ would have to be $a$.

$g_{4,1}(b\overline{c}a)= g_{4,2}(b)g_{4,3}(\overline{a})g_4(\overline{c})a\overline{b}\overline{a}c$ and $g_{4,1}(\overline{a}b)= g_{4,2}(b) g_{4,3}(\overline{a}) g_4(\overline{c})ab$, which again leaves us with $\{b, \overline{b}\}$.  So, as above, $\overline{a}b$ and $b\overline{c}a$ could not start $\rho_1$ and $\rho_2$, respectively.

This rules out all possibilities for $\overline{\rho_1}\rho_2$ and so $g$ has no iPNPs, and thus no PNPs, as desired.
\end{ex}

\smallskip

\noindent \emph{Throughout this section, $g:\Gamma \to \Gamma$ will be a semi-ideally decomposed TT representative of $\phi\in Out(F_r)$ with the same notation as that given at the end of Section \ref{Ch:IdealDecompositions}.  We will additionally require that $g=g_n \circ \cdots \circ g_1$ satisfies AM Properties (I)-(VIII) and, in particular, is a permitted composition.  Let $f_k=g_k \circ g_1 \circ g_n \circ \cdots \circ g_{k+1}$ and $G_k=G(f_k)$.}

\vskip10pt

\begin{prop}{\label{P:NPIdentification}} There exists a procedure for determining all iPNPs $\rho= \overline{\rho_1} \rho_2$ for $g$, where $\rho_1 = e_1\dots e_m$; $\rho_2 = e_1'\dots e_{m'}'$; $e_1,\dots, e_m, e_1',\dots,e_{m'}' \in \mathcal{E}(\Gamma)$; and $\{D_0(e_1), D_0(e_1')\}=\{d_1, d_1'\}$ is the unique illegal turn of $\rho$. \newline
\indent Let $\rho_{1,k}=e_1\dots e_k$ and $\rho_{2,l}= e_1'\dots e_l'$ throughout the following procedure.
\begin{description}
\item[(I)] Apply generators $g_1$, $g_2$, etc, to $e_1$ and to $e_1'$ until $Dg_{j,1}(e_1')=Dg_{j,1}(e_1)$.  Either $g_{j,1}(e_1)$ is a subpath of $g_{j,1}(e_1')$ or vice versa.  Without loss of generality assume that $g_{j,1}(e_1')$ is a subpath of $g_{j,1}(e_1)$ so that $g_{j,1}(e_1)=g_{j,1}(e_1')t_2 \dots$, for some edge $t_2$.  (Otherwise just switch all $e_i$ and $e_i'$, $\rho_1$ and $\rho_2$, etc, in the following arguments).  Then, $\rho_2$ must contain another edge $e_2'$. \newline
\indent We explain in (III) how to find all possibilities for this edge.
\item[(II)] Inductively, assume that $g_{j,1}(\rho_{1,k})=g_{j,1}(\rho_{2,s})t_{s+1}\dots$ (or again switch $e_i$ for $e_i'$, $\rho_1$ for $\rho_2$, and so on).
\item[(III)] We must add another edge $e_{s+1}'$ to $\rho_2$.  There are two cases to consider:
\begin{itemize}
\item[(a)] If $D_0(t_{s+1})=d^{u}_{j}$, then the different possibilities for $e_{s+1}'$ are determined by the directions $d_{s+1}'$ such that $T_{j+1}=\{Dg_{j,1}(d_{s+1}'), D_0(t_{s+1})\}$ where $D_0(e_{s+1}')=d_{s+1}'$.  (In this case there is a green segment between $Dg_{j,1}(d_{s+1}')$ and $t_{s+1}$ in $G_{TA}(f_k)$, and $Dg_{j,1}(d_{s+1}')$ is the non-red direction of $T_{j+1}$).
\item[(b)] If $D_0(t_{s+1}) \neq d^{u}_{j}$, then the different possibilities for $e_{s+1}'$ are all edges $e_{s+1}'$ such that $Dg_{j,1}(d_{s+1}')=D_0(t_{s+1})$ where $D_0(e_{s+1}')=d_{s+1}'$. \newline
\indent After throwing out any choices for $d_{s+1}'$ such that $T_0=\{\overline{d_s'}, d_{s+1}'\}$ is the green illegal turn for $g$, each remaining $d_{s+1}'$ in (a) or (b) gives another prospective iPNP that we must continue applying the algorithm to.
\end{itemize}
\item[(IV)] Suppose that at some step we do not have that $g_{j,l}(\rho_{2,k})$ is a subpath of $g_{j,l}(\rho_{1,s})$ (or vice versa).  Then we compose with generators $g_i$ until either:
\begin{itemize}
\item[(a)] We composed with enough $g_i$ to obtain some $g^{p'}$ such that  $g^{p'}(\rho_{2,k})=\tau'e_1'\dots$ and $g^{p'}(\rho_{1,s})=\tau'e_1\dots$ for some legal path $\tau'$  (in this case proceed to (V) ),
\item[(b)] $g_{j,l}(\rho_{2,k})$ is a subpath of $g_{j,l}(\rho_{1,s})$ or vice versa (in this case, return to \indent (II) and continue with the algorithm as before), or
\item[(c)] some $g_{l,1}\circ g^{p'}(\rho_{2,k})=\tau'\gamma_{2,k}$ and $g_{l,1}\circ g^{p'}(\rho_{1,s})=\tau'\gamma_{1,s}$ where $\{D_0(\gamma_{2,k}), D_0(\gamma_{1,s})\}$ is a legal turn in $G_l$.  In this case there cannot be an iPNP with $\overline{\rho_{2,k}}\rho_{1,s}$ as a subpath (proceed to (VIII)).
\end{itemize}
\item[(V)] For each $1 \leq p'$ such that $g^{p'}(\rho_{2,m})=\tau'e_1\dots$ and $g^{p'}(\rho_{1,n})=\tau'e_1'\dots$ for some legal path $\tau'$ (for the appropriate m and n), check if $g^{p'}_{\#}(\overline{\rho_{1,n}}\rho_{2,m}) \subset \overline{\rho_{1,n}}\rho_{2,m}'$ or vice versa and follow the appropriate step (among (a)-(d) below).
\begin{itemize}
\item[(a)] If, for some $1 \leq p'$, $g^{p'}_{\#}(\overline{\rho_{1,n}}\rho_{2,m}) = \overline{\rho_{1,n}}\rho_{2,m}$, then $\overline{\rho_{1,n}}\rho_{2,m}$ is the only possible iPNP for $g$.
\item[(b)] For each $1 \leq p'$ such that $g^{p'}_{\#}(\overline{\rho_{1,n}}\rho_{2,m}) \subset \overline{\rho_{1,n}}\rho_{2,m}$ (where containment is proper), proceed to (VII).
\item[(c)] If $\overline{\rho_{1,n}}\rho_{2,m} \subset g^{p'}_{\#}(\overline{\rho_{1,n}}\rho_{2,m})$ (where containment is proper), proceed to (VI).
\item[(d)] If we do not have $g^{p'}_{\#}(\overline{\rho_{1,n}}\rho_{2,m}) \subset \overline{\rho_{1,n}}\rho_{2,m}'$ or vice versa for any $1 \leq p' \leq b$, then there is only one circumstance where we can possibly have an iPNP with $\overline{\rho_{2,m}}\rho_{1,n}$ as a subpath.  This is the case where $g^{p'}_{\#}(\overline{\rho_{1,n}}\rho_{2,m}) =\overline{\gamma_{1,n}}\gamma_{2,m}$ where either $\gamma_{1,n} \subset \rho_{1,n}$ and $\rho_{2,m} \subset \gamma_{2,m}$ or $\gamma_{2,m} \subset \rho_{2,m}$ and $\rho_{1,n} \subset \gamma_{1,n}$.  In this case, apply (VII) to the side that is too short.  Otherwise, there cannot be an iPNP with $\overline{\rho_{2,m}}\rho_{1,n}$ as a subpath, so proceed to (VIII).
\end{itemize}
\item[(VI)] We assume here that $ \bar \rho_{1,n}\rho_{2,m} \subset g^{p'}(\overline{\rho_{1,n}}\rho_{2,m})$.  Consider the final occurrence of $e_n$ in the copy of $\rho_{1,n}$ in $g^{p'}(\overline{\rho_{1,n}}\rho_{2,m})$ and the final occurrence of $e_m'$ in the copy of $\rho_{2,m}$ in $g^{p'}(\overline{\rho_{1,n}}\rho_{2,m})$.  This final occurrence of $e_n$ must have come from $g^{p'}(e_n)$ and this final occurrence of $e_m'$ must have come from $g^{p'}(e_m')$.  This means that we have fixed points in $e_n$ and $e_m'$.  Replace $\overline{\rho_{1,n}}\rho_{2,m}$ with $\overline{\rho_{1,n}'} \rho_{2,m}'$ where $\overline{ \rho_{1,n}'}\rho_{2,m}'$ is the same as $\overline{\rho_{1,n}}\rho_{2,m}$, except that $e_n$ and $e_m'$  are replaced with some partial edges ending at the fixed points.  Repeat this process until some $\overline{ \rho_{1,n}'}\rho_{2,m}'$ is an iPNP.
\item[(VII)] We assume  here that $g^{p'}(\bar \rho_{1,n}\rho_{2,m}) \subset \overline{\rho_{1,n}}\rho_{2,m}$ (where containment is proper).  Without loss of generality, assume that there exists a $t_{m+1}$ such that  $g^{p'}(\overline{\rho_{1,n}}\rho_{2,m})t_{m+1} \subset \overline{\rho_{1,n}}\rho_{2,m}$.  For each direction $d_i$ such that $Dg^{p'}(d_i)=D_0(t_{m+1})$ and such that $\{D_0(\overline{e_{i-1}}), D_0(e_i) \}$ is not the green illegal turn $\{D_0(e_1), D_0(e_1') \}$ for $g$, return to (V) with $\rho_{2,m+1}$ where $D_0(e_{m+1})=d_i$.
\item[(VIII)] Continue to rule out the other possible subpaths that have arisen via this procedure (by different choices of $d_i$, as in (III) or (VII)).  If there are no other possible subpaths, then we have shown there are no iPNPs for $g_j$, and thus no PNPs for $g$.
\end{description}
\end{prop}

\bigskip

\noindent We will need the following Lemma(s) for the proof of this proposition.

\begin{lem} Subpaths of legal paths are legal.
\end{lem}

\noindent \emph{Proof of Lemma:}  The set of turns of the subpath is a subset of the set of turns of the path.  Since all turns of the path are legal, all turns of the subpath must also be legal.  So the subpath must also be a legal path.  \newline
\noindent QED.

\vskip5pt

\noindent  \emph{Proof of Proposition:}  We begin with an argument that will be used throughout the proof.  Since $\rho=\overline{\rho_1} \rho_2$ is an iPNP, $\rho_1$ and $\rho_2$ are both legal paths.  Since subpaths of legal paths are legal and since the images under $g_{k,1}$ of legal paths are legal, the paths $g_{k,1}(e_1 \dots e_l)$ and $g_{k,1}(e_1' \dots e'_{l'})$ are legal for each $1 \leq k \leq n$, $1 \leq l \leq m$, and $1 \leq l' \leq n'$.

Since $\{D_0(e_1), D_0(e_1')\}$ is an illegal turn, for some $j$, $Dg_{j,1}(e_1')=Dg_{j,1}(e_1)$.  We need to prove that either $g_{j,1}(e_1)$ is a subpath of $g_{l,1}(e_1')$ or vice versa.  Let $d_1=D_0(e_1)$ and $d_1'=D_0(e_1')$.  Since $\{D_0(e_1), D_0(e_1')\}$ is an illegal turn for $g$ and the only illegal turn for $g$ is $T_1=\{d^{pu}_0, d^{pa}_0\}$, $\{d_1, d_1'\}=\{d^{pu}_0, d^{pa}_0\}$.  Without loss of generality suppose that $d_1=d^{pu}_0$ (and $d_1'=d^{pa}_0$).  Since $g_1$ is defined by $e^{pu}_0 \mapsto e^a_1 e^u_1$, it follows that $g_1(e_1)=g_1(e^{pu}_0)=e^a_1e^u_1$ and $g_1(e_1')=g_1(e^{pa}_0)=e^a_1$.  So $g_1(e_1')$ is a subpath of $g_1(e_1)$.  Since $g_{j,2}$ is an automorphism and also takes legal paths to legal paths, $g_{j,2}(g_1(e_1'))=g_{j,2}(e^a_1)$ is a subpath of $g_{j,2}(g_1(e_1)) = g_{j,2}(e^a_1 e^u_1) = g_{j,2}(e^a_1) g_{j,2}(e^u_1)$, as desired.

The final thing to show for (I) is that $\rho_2$ must contain a second edge $e_2'$.  Suppose $\rho_2$ did not contain a second edge.  Then tightening would cancel out all of $g_{j,1}(\rho_2)$ with an initial subpath of $g_{j,1}(\rho_1)$ and so certainly would cancel out all of $g_{j,1}(\rho_2)$ with an initial subpath of $g_{j,1}(\rho_1)$.  Thus, $(g_{j,1})_{\#}(\rho) = (g_{j,1})_{\#}(\overline{\rho_1}\rho_2)$ would be a subpath of $g_{j,1}(\rho_2)$ and hence would be legal.  So $g^p_{\#}(\rho)=(g^{p-1} \circ g_{n,j+1})_{\#}((g_{j,1})_{\#}(\rho))$ is legal for all $p$, contradicting that some $g^p_{\#}(\rho)=\rho$, which has an illegal turn.

To start off proving (III) we need a similar argument, as in the previous paragraph, to show that $\rho_2$ would have to have another edge $e_{s+1}'$.  For the sake of contradiction, suppose that $\rho_2$ ended with $e_s'$.  Then $(g_{j,1})_{\#}(\rho)= (g_{j,1})_{\#}(\overline{\rho_1}\rho_2)$ would be a subpath of $g_{j,1}(\rho_1)$ (for similar reasons as above), which leads to a contradiction as in the argument above.

We now prove the claims of (IIIa) and (IIIb). Since $g_{j,1}(\rho_{1,k})= g_{j,1}(\rho_{2,s})t_{s+1} \dots$, \newline \noindent $(g_{j,1})_{\#}(\overline{\rho_{1,k}} \rho_{2,s}) = \dots \overline{t_{s+1}} = \overline{\gamma}$ for some legal path $\gamma$ ($\gamma$ is legal because it is a subpath of the image of the legal path $\rho_{1,k}$).  Additionally, $g_{j,1}(e_{s+1}')$ will be legal since $e_{s+1}'$ is a legal path.  Thus, $(g_{j,1})_{\#}(\overline{\rho_{1,k}} \rho_{2,s}) = (\overline{\gamma} g_{j,1}(e_{s+1}'))_{\#}$, which is just $\overline{\gamma} g_{j,1}(e_{s+1}')$ unless $Dg_{j,1}(d_{s+1}')= D_0(t_{s+1})$,  and then is legal unless  $\{D_0(\gamma), D_0(g_{j,1}(e_{s+1}')) \}$ is an illegal turn, i.e. $T_{j+1}=\{Dg_{j,1}(d_{s+1}'),D_0( t_{s+1}) \}$ with $D_0(e_{s+1}')=d_{s+1}'$.

Suppose first that $D_0(t_{s+1})=d^u_j$ (as in (III)(a)).  Notice that, in this case, $D_0(t_{s+1})$ is not in the image of $Dg_j$ and thus is not in the image of $Dg_{j,1}=D(g_j \circ g_{j-1,1})$ and so $Dg_{j,1}(d_{s+1}') \neq D_0( t_{s+1})$.  This tells us that $(\overline{\gamma} g_{j,1}(e_{s+1}'))_{\#} = \overline{\gamma} g_{j,1}(e_{s+1}')$, which will be a legal path unless $T_{j+1}=\{Dg_{j,1}(d_{s+1}'),D_0( t_{s+1})\}$.  However, if $(g_{j,1})_{\#}(\overline{\rho_{1,k}} \rho_{2,s}$) $= \overline{\gamma} g_{j,1}(e_{s+1}')$, then $(g_{j,1})_{\#}(\rho) = (g_{j,1})_{\#}(\overline{\rho_1} \rho_2)=$ \newline
\noindent $(g_{j,1})_{\#} (\overline{e_m} \dots \overline{e_{k+1}})\overline{\gamma}g_{j,1}(e_{s+1}') (g_{j,1})_{\#}(e_{s+2}'\dots e_m')=$ $g_{j,1}(\overline{e_m} \dots \overline{e_{k+1}}) \overline{\gamma} g_{j,1}(e_{s+1}') g_{j,1}(e_{s+2}' \dots e_m')$, since $\rho_1$ and $\rho_2$ are legal paths and the images of edges are legal.  But $g_{j,1}(\overline{e_m}\dots\overline{e_{k+1}}) \overline{\gamma}$ is a subpath of $g_{j,1}(\overline{\rho_{1}})$, so is legal, and $g_{j,1}(e_{s+1}') g_{j,1}(e_{s+2}' \dots e_m')$ is a subpath of $g_{j,1}(\rho_2)$, so is legal, and we still have that $\overline{\gamma} g_{j,1}(e_{s+1}')$ is legal, which together would make $(g_{j,1})_{\#}(\rho)$ legal.  This contradicts that some $g^p_{\#}(\rho)$ $= (g^{p-1}\circ g_{j,n+1})_{\#} (g_{j,1})_{\#}(\rho))$ must be $\rho$, which has an illegal turn.  So, $T_{j+1} = \{Dg_{j,1}(d_{s+1}'), D_0(t_{s+1})\}$, as desired.

Suppose now (as in (III)(b)) that $D_0(t_{s+1}) \neq d^u_j$.  For the sake of contradiction suppose that $Dg_{j,1}(d_{s+1}') \neq D_0(t_{s+1})$, where $D_0(e_{s+1}') = d_{s+1}'$.  First off, notice that this means that again $(\overline{\gamma} g_{j,1}(e_{s+1}'))_{\#} = \overline{\gamma} g_{j,1}(e_{s+1}')$.  Also, since $Dg_{j,1}(d_{s+1}')$ cannot be $d^u_j$ (see reasoning above) and $D_0(t_{s+1}) \neq d^u_j$, we cannot have $T_{j+1} = \{Dg_{j,1}(d_{s+1}'),D_0( t_{s+1}) \}$.  This would make $\overline{\gamma} g_{j,1}(e_{s+1}')$ legal, which leads to a contradiction as above.  So $Dg_{j,1}(d_{s+1}') = D_0(t_{s+1})$, as desired.

The final observation about (III) is that choices for $e_{s+1}'$ such that $T_0=\{D_0(e_{s}'), D_0(e_{s+1}')\}$ must be thrown out since $\rho_2$ must be a legal path.

We need to show for (IVc) that, if $g_{l,1}\circ g^{p'}(\rho_{2,k}) = \tau' \gamma_{2,k}$ and $g_{l,1}\circ g^{p'}(\rho_{1,s})=\tau'\gamma_{1,s}$ where $\{D_0(\gamma_{1,s}), D_0(\gamma_{2,k})\}$ is a legal turn in $G_l$, then there cannot be an iPNP with $\overline{\rho_{2,k}}\rho_{1,s}$ as a subpath. We prove this now.  Under the stated conditions, $g_{l,1}\circ g^{p'}(\rho_2)=\tau' \gamma_2$ and $g_{l,1}\circ g^{p'}(\rho_{1})=\tau'\gamma_{1}$ where $\gamma_{2,k}$ is an initial subpath of $\gamma_2$ (both of which are legal) and $\gamma_{1,s}$ is an initial subpath of $\gamma_{1}$ (both of which are legal).  Since $\{D_0(\gamma_{1}), D_0(\gamma_{2})\} = \{D_0(\gamma_{1,s}), D_0(\gamma_{2,k})\}$ is a legal turn, $(g_{l,1}\circ g^{p'})_{\#}(\rho) = (g_{l,1}\circ g^{p'})_{\#}(\overline{\rho_{1}} \rho_2) = \bar\gamma_1\gamma_2$, which is a legal path.  Let $p$ be such that $g^p_{\#}(\rho))=\rho$.  (Without loss of generality we can assume that $p > p'$ by replacing $p$ by a multiple of $p$ if necessary).  Then, $g^p_{\#}(\rho) = ((g^{p-p'-1} \circ g_{n, l+1}) \circ (g_{l,1} \circ g^{p'}))_{\#}(\rho) =  (g^{p-p'-1} \circ g_{n, l+1})_{\#} ((g_{l,1} \circ g^{p'})_{\#}(\rho)) = (g^{p-p'-1} \circ g_{n,l})_{\#}(\overline{\gamma_1}\gamma_2) = (g^{p-p'-1} \circ g_{n,l})(\overline{\gamma_1}\gamma_2)$, since $\overline{\gamma_1}\gamma_2$ is a legal path.  This makes $g^p_{\#}(\rho)$ legal since images under permitted compositions of legal paths are legal.  And this contradicts that $g^p_{\#}(\rho)=\rho$, which is not a legal path.  We have now verified everything needing verification in (IV).

As in (V), suppose that $g_{l,1} \circ  g^{p'}(\rho_{2,m})=\tau'e_1'\dots$ and $g_{l,1} \circ  g^{p'}(\rho_{1,n}) = \tau' e_1\dots$ for some legal path $\tau'$ (for the appropriate $m$ and $n$).  (Va) is true by definition.  Since (Vb) and (Vc) just refer us to later steps, we can just focus on (Vd) for now.  The first thing that we need to prove for (Vd) is that there is only one circumstance where we can possibly have an iPNP with $\overline{\rho_{1,n}}\rho_{2,m}$ as a subpath.  Suppose that, for no power $p'$ do we ever have $g^{p'}_{\#}(\overline{\rho_{1,n}}\rho_{2,m}) = \overline{\rho_{1,n}}\rho_{2,m}$, $g^{p'}_{\#}(\overline{\rho_{1,n}}\rho_{2,m}) \subset \overline{\rho_{1,n}}\rho_{2,m}$, $\overline{\rho_{1,n}}\rho_{2,m} \subset g^{p'}_{\#}(\overline{\rho_{1,n}}\rho_{2,m})$, or $g^{p'}_{\#}(\overline{\rho_{1,n}}\rho_{2,m}) =\overline{\gamma_{1,n}}\gamma_{2,m}$ where either $\gamma_{1,n} \subset \rho_{1,n}$ and $\rho_{2,m} \subset \gamma_{2,m}$ or $\gamma_{2,m} \subset \rho_{2,m}$ and $\rho_{1,n} \subset \gamma_{1,n}$.  Now, for the sake of contradiction, suppose that some $\overline{\rho_{1,n+k}}\rho_{2,m+l}$ containing $\overline{\rho_{1,n}}\rho_{2,m}$ is an iPNP of period $p$.  Since $\overline{\rho_{1,n+k}}\rho_{2,m+l}$ is an iPNP, $\rho_{1,n+k}$ and $\rho_{2,m+l}$ are both legal paths (as are the subpath $\rho_{1,n}$ and $\rho_{2,m}$).  This tells us that $g^{p}_{\#}(\overline{\rho_{1,n+k}}\rho_{2,m+l})= g^{p}(\overline{e_{n+k}'} \dots \overline{e_{n+1}})g^{p'}_{\#}(\overline{\rho_{1,n}}\rho_{2,m})g^{p}(e_{m+1}\dots e_{m+l})$.  So $g^{p'}_{\#}(\overline{\rho_{1,n}}\rho_{2,m}) \subset g^{p}_{\#}(\overline{\rho_{1,n+k}}\rho_{2,m+l}) = \overline{\rho_{1,n+k}}\rho_{2,m+l}$.  But this lands us in one of the situations we said could not occur, which is a contradiction.

The verification of (VI) is left to the reader and there is nothing really to prove in (VII) since the conditions for (V) still hold. \newline
\noindent QED.

\section{Representative Loops}{\label{Ch:RepresentativeLoops}}

The goal of this section is to prove, for a Type (*) pIW graph $\mathcal{G}$, that representatives coming from the loops in the $AMD(\mathcal{G})$ and satisfying certain prescribed properties, as mentioned before, are indeed TT representatives of ageometric, fully irreducible $\phi \in Out(F_r)$ such that $IW(\phi)=\mathcal{G}$.  This result is given in Proposition \ref{P:RepresentativeLoops}.

The following three lemmas are used in the proof of Proposition \ref{P:RepresentativeLoops}.

\begin{lem}{\label{L:GreenIllegalTurn}} For each $k$, $T_k$ is not represented by any edge in $G_{k-1}$.
\end{lem}

\noindent \emph{Proof}: For the sake of contradiction suppose that $T_k= \{d^{pa}_{k-1}, d^{pu}_{k-1} \}$ were represented by an edge in $G_{k-1}$.  Then $\{d^{pa}_{k-1}, d^{pu}_{k-1} \}$ would be a turn in some $f_{k-1}^p(e_{k-1,i})$.  By Lemma \ref{L:PreLemma}, the edge path $g_{k-1,1}(e_{0,i})$ contains $e_{k-1,i}$.  Thus, the edge path $f_{k-1}^p\circ g_{k-1,1}(e_{k-1,i})$ would also contain the turn $\{d^{pa}_{k-1}, d^{pu}_{k-1} \}$.  Since $g_{k}(e^{pa}_{k-1})=g_{k}(e^{pu}_{k-1})$, this would contradict that $g_{k} \circ f_{k-1}^p\circ g_{k-1,1}(e_{0,i})=g_{k-1,1}^p(e_{0,i})$ cannot have cancelation.  So $T_k$ cannot be represented by any edge in $G_{k-1}$, as desired. \newline
\noindent QED.

\begin{lem} Since one vertex of the green illegal turn $T_{k+1}$ of $G_k$ will be $d^u_k$, $T_{k+1}$ cannot also be a purple edge in $G_k$.  Also, $d^u_k$ must be a vertex of the red edge $[t^R_k]$ of $G_k$.  \end{lem}

\noindent \emph{Proof}: Since $d^u_k$ has no preimage under $Dg_k$, it cannot be a vertex of a purple edge of $G_k$, as the purple edges of $G_k$ are images of purple and red turns of $G_{k-1}$ under $Dg_k^C$.  The red edge in $G_k$ is $[t^R_k]=[\overline{d^a_k}, d^u_k]$, which contains $d^u_k$.  \newline
\noindent QED.

\begin{lem} For permitted compositions, images of legal paths and turns are legal.
\end{lem}

\noindent \emph{Proof of Lemma:}  Suppose that $\gamma$ is a legal path and suppose that $g$ is a permitted composition.  Since permitted compositions are train track maps, the image under $g$ of any edge of $\gamma$ is legal.  Thus, we only need to be concerned about what happens with the turns of $\gamma$.  Since $\gamma$ is legal, all turns of $\gamma$ are legal.  Since images of legal turns are legal, the images of all turns of $\gamma$ are legal.  Thus, the image of $\gamma$ is legal, as desired. \newline
\noindent QED.

\begin{prop}{\label{P:RepresentativeLoops}} Suppose that $\mathcal{G}$ is a Type (*) pIW graph and that $L(g_1, \dots, g_k)=E(g_1, G_{0}, G_1)  * \dots *  E(g_k, G_{k-1}, G_k)$ is a loop in $AMD(\mathcal{G})$ satisfying:
\begin{description}
\item 1. Each purple edge of $G(g)$ correspond to a turn taken by some $g^k(E_j)$ where $E_j \in \mathcal{E}(\Gamma)$;
\item 2. for each $1 \leq i,j \leq q$, there exists some $k \geq 1$ such that $g^k(E_j)$ contains either $E_i$ or $\bar E_i$; and
\item 3. $g$ has no periodic Nielsen paths.
\end{description}
\noindent Then $g: \Gamma \to \Gamma$ is a train track representative of an ageometric fully irreducible $\phi \in Out(F_r)$ such that $IW(\phi)=\mathcal{G}$.
\end{prop}

\noindent  \emph{Proof of Proposition:}  By the FIC, we only need to show that $g$ is a train track map, the transition matrix of $g$ is Perron-Frobenius, and $IW(\phi)=\mathcal{G}$.  Property (2) of this proposition is the same as AM Property (VIc) and so that the transition matrix is Perron-Frobenius follows from Lemma \ref{L:PF}.  That $g$ is a train track map follows from Lemma \ref{L:PF}.  Since $g$ has no PNPs, $IW(g)=SW(g)$, where $SW(g)= \underset{\text{vertices v} \in \Gamma}{\bigcup}LSW(v;g)$.  By the definition of $LSW(v;g)$, $\underset{\text{vertices v} \in \Gamma}{\bigcup}LSW(v;g)$ edges correspond precisely with turns crossed over by some $g^k(E_j)$ were $E_j \in \mathcal{E}(\Gamma)$.  By the definition of $G(g)$ the edges of $\underset{\text{vertices v} \in \Gamma}{\bigcup}LSW(v;g)$ correspond precisely with the purple edges of $G(g)$.

This completes the proof. \newline
\noindent QED.

\section{Procedure for Building Ideal Whitehead Graphs}{\label{Ch:Procedure}}

In this section we give three different methods for constructing train track representatives that have the potential to be of Type (*) with a Type (*) pIW graph ideal Whitehead graph $\mathcal{G}$.  Different methods work better in different circumstances.  For example, if most of the LTT structures $G$ with $PI(G)=\mathcal{G}$ are birecurrent, then Methods II and III are better suited (as $AMD(\mathcal{G})$ may be very large and impractical to construct).  On the other hand, if only a few of the LTT structures $G$ with $PI(G)=\mathcal{G}$ are birecurrent, then constructing $AMD(\mathcal{G})$ is much simpler than using ``guess and check'' methods and so Method I generally proves more practical.

\vskip7pt

\noindent \emph{In all figures of this section we continue with the convention that Y denotes $\bar{y}$, etc.}

\vskip10pt

\subsection{Method I (Building the Entire AM Diagram)}

\indent Let $\mathcal{G}$ be a Type (*) pIW graph.  We explain in Steps 1-6 below a procedure for constructing $AMD(\mathcal{G})$. Once $AMD(\mathcal{G})$ has been built, one still needs to find an appropriate loop in $AMD(\mathcal{G})$.  We explain how to do this in Step 7.  The final step will be to test the representative constructed from the loop to ensure that it is PNP-free.  The procedure for identifying PNPs was explained and proved in Section \ref{Ch:NPIdentification}.

\vskip10pt

\noindent \textbf{\underline{STEP 1}: \emph{CONSTRUCT AN LTT CHART FOR $\mathcal{G}$}}

\begin{description}
\item[A.] An LTT Chart for $\mathcal{G}$ will contain precisely one column for each way of labeling the vertices of $\mathcal{G}$ with $\{x_1, \overline{x_1}, x_2, \overline{x_2}, \dots, x_{r-1}, \overline{x_{r-1}}, x_r \}=\{X_1, X_2, \dots, X_{2r-2}, X_{2r-1} \}$, except that:
\begin{itemize}
\item[1.] We consider labelings of the vertices of $\mathcal{G}$ equivalent when there is a permutation of the indices $1 \leq i \leq r-1$ and a permutation of the two elements of each pair $\{x_i, \overline{x_i}\}$ making the labelings identical (when we have such a permutation, we say that the labeled graphs are \emph{Edge Pair Permutation (EPP) Isomorphic}, even in a graph where we also have a vertex labeled by $x_r$ or include black edges, making the graph an LTT structure) \newline

\begin{ex} The following are two EPP-isomorphic graphs (in this example and the following examples, $z$ denotes $x_1$, $Z$ denotes $\overline{x_1}$, $y$ denotes $x_2$, $Y$ denotes $\overline{x_2}$, $x$ denotes $x_3$, and $X$ denotes $\overline{x_3}$):

\begin{figure}[H]
\centering
\includegraphics[width=2.2in]{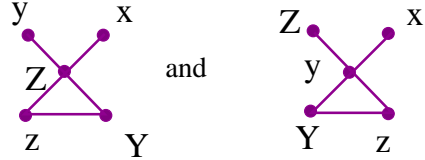}
\caption[width=4.5in]{{\small{\emph{Graphs EPP Isomorphic by the permutation taking $y$ to $\overline{z}$ and $\overline{y}$ to $z$}}}}
\label{fig:EPPIsomorphicGraphs}
\end{figure}
\end{ex}

\item[2.] We leave out any labeling causing $\mathcal{G}$ to have two edges of the form $[x_i, \overline{x_i}]$, each with a valence-one vertex (we will call an edge containing a valence-one vertex a \emph{valence-one edge} and a pair of vertices of the form $\{x_i, \overline{x_i}\}$ an \emph{edge pair}, making the statement say that we leave out any labeling causing $\mathcal{G}$ to have two distinct pairs of edge-pair labeled vertices, each connected by a valence-one edge). \newline

\begin{ex} Here we show a graph not included as a result of having distinct pairs of edge-pair labeled vertices, each connected by a valence-1 edge.
\begin{figure}[H]
\centering
\includegraphics[width=.9in]{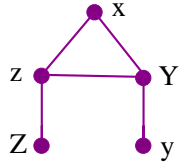}
\caption[width=4.5in]{{\small{\emph{Not in LTT Chart since $[y, \overline{y}]$ and $[z, \overline{z}]$ valence-one edges}}}}
\label{fig:NotInLTTChart}
\end{figure}
\end{ex}
\end{itemize}

\item[B.] Graphs determining the columns of an LTT Chart are colored purple and labeled from left to right by $SW_I$, $SW_{II}$, etc.  We call these graphs the \emph{Determining SW-Graphs} for the columns.

\begin{rk}
There is no canonical way of ordering the graph that will be $SW_I$, $SW_{II}$, etc.  There is also no canonical way of choosing a labeled graph representing an EPP-isomorphic equivalence class.  Thus, there will be multiple possible LTT Charts for $\mathcal{G}$, all of which will lead to the same diagram $AMD(\mathcal{G})$.
\end{rk}

\begin{ex} In the case of Graph VII, one set of representatives we could use to determine the columns would be:

\begin{figure}[H]
\centering
\includegraphics[width=4.5in]{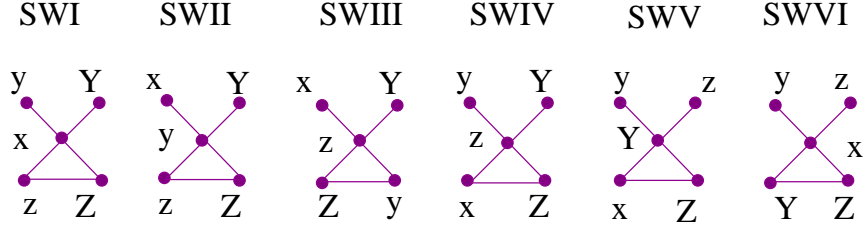}
\caption[width=4.8in]{{\small{\emph{Determining SW-Graphs for Graph VII}}}}
\label{fig:DeterminingSWGraphs}
\end{figure}
\end{ex}

\item[C.] Each column of the LTT Chart will contain a graph for each way of attaching (at a single vertex, called the \emph{attaching vertex}) a red edge to the column's determining SW-graph so that:
\begin{itemize}
\item [1.] the valence-one vertex of the red edge (which is colored red and called the \emph{free vertex}) is labeled by $\overline{x_r}$ and
\item [2.] no valence-one edge connects edge pair vertices (there are no edges of the form $[x_i, \overline{x_i}]$ where either $x_i$ has valence one or $\overline{x_i}$ has valence one).

\begin{ex} We show here two graphs that would not be included in an LTT Chart for Graph VII because they each have a valence-one edge connecting edge-pair vertices.

\begin{figure}[H]
\centering
\includegraphics[width=2.1in]{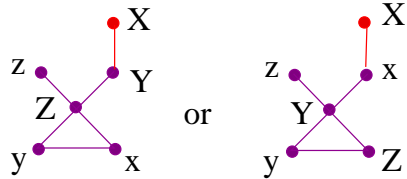}
\caption[width=4.5in]{{\small{\emph{Graphs not in LTT Chart because $[z, \overline{z}]$ and $[\overline{x}, x]$ are valence-one edges}}}}
\label{fig:ContainingValenceOneEdgePairEdge}
\end{figure}
\end{ex}

\end{itemize}
\item[D.] We label graphs in the first ($SW_I$) column from top to bottom $I_a$, $I_b$, $I_c$, $\dots$; graphs in the second ($SW_II$) column from top to bottom $II_a$, $II_b$, $II_c$, $\dots$; etc.
\end{description}

\smallskip

\begin{ex}{\label{E:MethodIIILTTChart}} \emph{An LTT Chart for Graph VII}

\smallskip

\noindent In the following chart, the graphs we leave out are either EPP-isomorphic to one among those we included (see Figure \ref{fig:EPPIsomorphicGraphs}) or violate one of the conditions for inclusion (as in Figure \ref{fig:ContainingValenceOneEdgePairEdge}).  We continue to denote $z$ by $x_1$, we denote $Z$ by $\overline{x_1}$, we denote $y$ by $x_2$, we denote $Y$ by $\overline{x_2}$, we denote $x$ by $x_3$, and we denote $X$ by $\overline{x_3}$.

\begin{figure}[H]
\centering
\includegraphics[width=4.3in]{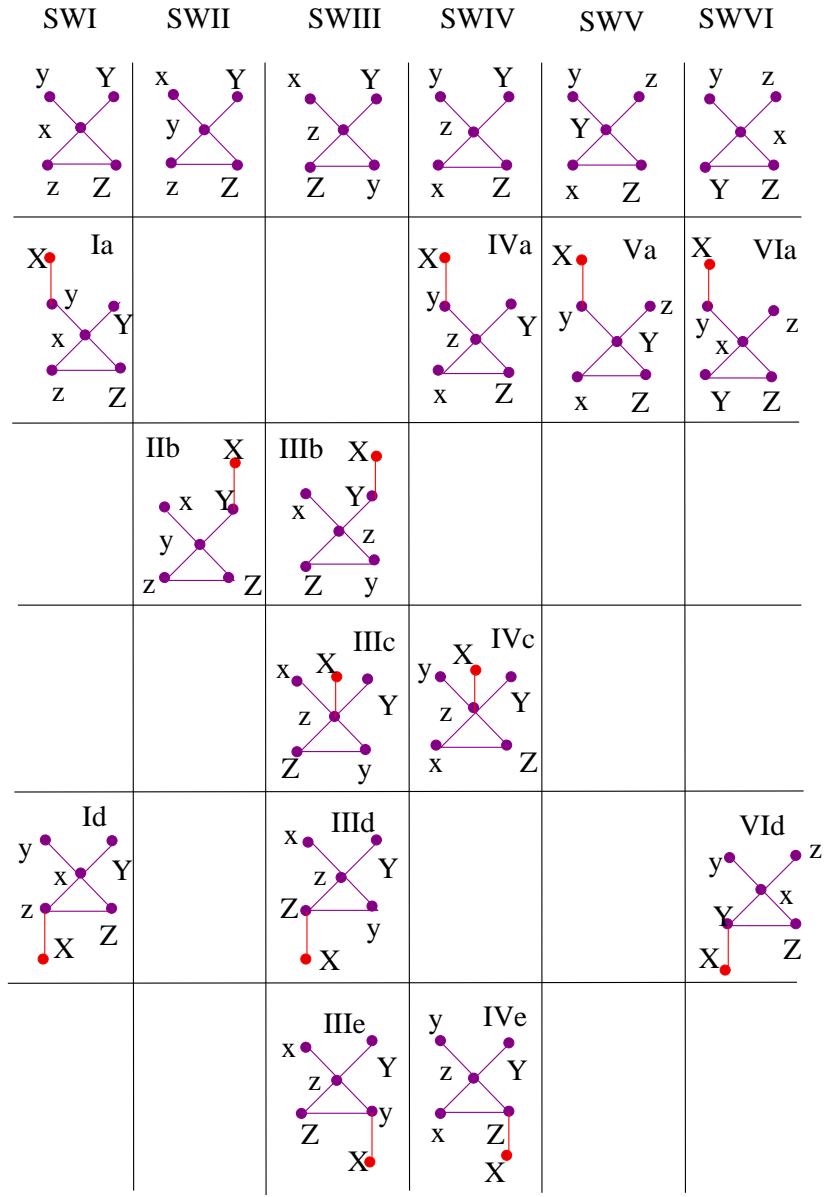}
\caption{{\small{\emph{The LTT Chart for Graph VII}}}}
\label{fig:LTTChart}
\end{figure}
\end{ex}

\newpage

\noindent \textbf{\underline{STEP 2}: \emph{BIRECURRENCY}}

\vskip5pt

For each graph in the LTT Chart for $\mathcal{G}$, and each $1 \leq i \leq r$, add a black edge $[x_i, \overline{x_i}]$ connecting $x_i$ and $\overline{x_i}$ (notice that the graphs are now smooth train track graphs).  Check each of these graphs for birecurrency, put a box around each birecurrent graph in the LTT chart, and cross out each nonbirecurrent graph in the LTT chart.

\begin{ex}{\label{E:MethodIIIBirecurrency}} \emph{Birecurrency for Graph VII}

\smallskip

We leave out here labels on the vertices of each LTT structure to highlight that the only significant information is vertex color and whether two vertices are part of an edge pair (forming an edge pair with the free vertex will have added significance).  The fact that black edges connect edge-pair vertices encodes the information needed.

\begin{figure}[H]
\centering
\includegraphics[width=4in]{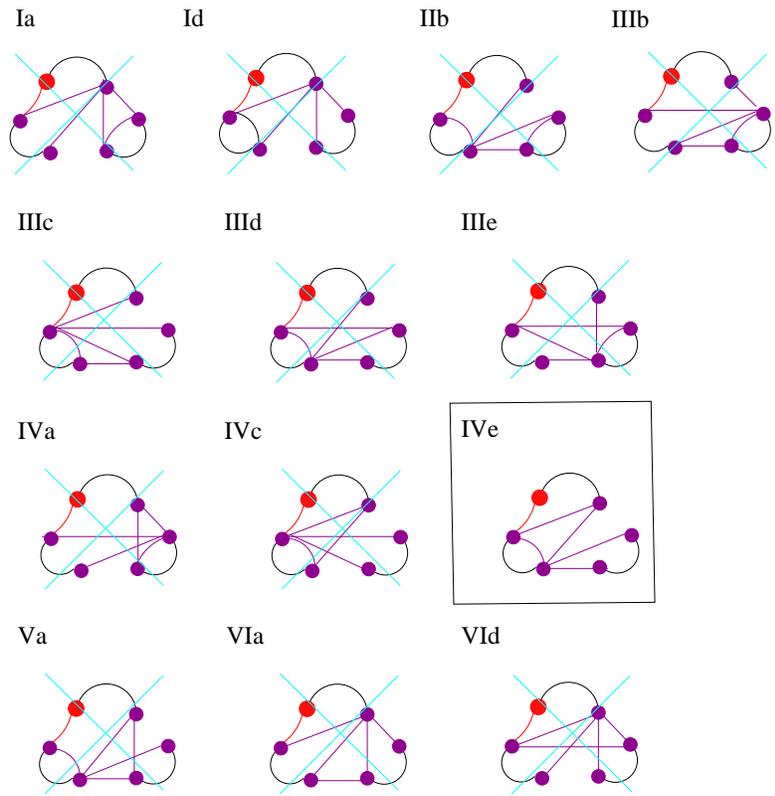}
\caption{{\small{\emph{Birecurrency Analysis for Graph VII: We drew the LTT structure for each labeled graph in the LTT Chart and checked the structures for birecurrency.  The only birecurrent LTT structure was that for IVe, which is why it is the single ``boxed'' graph.}}}}
\label{fig:Birecurrency}
\end{figure}

\noindent From the above figure one can see that the only LTT chart graph with a birecurrent LTT structure is IVe.  We cross out and ``box'' the corresponding graphs in the LTT Chart to get:

\begin{figure}[H]
\centering
\includegraphics[width=4.3in]{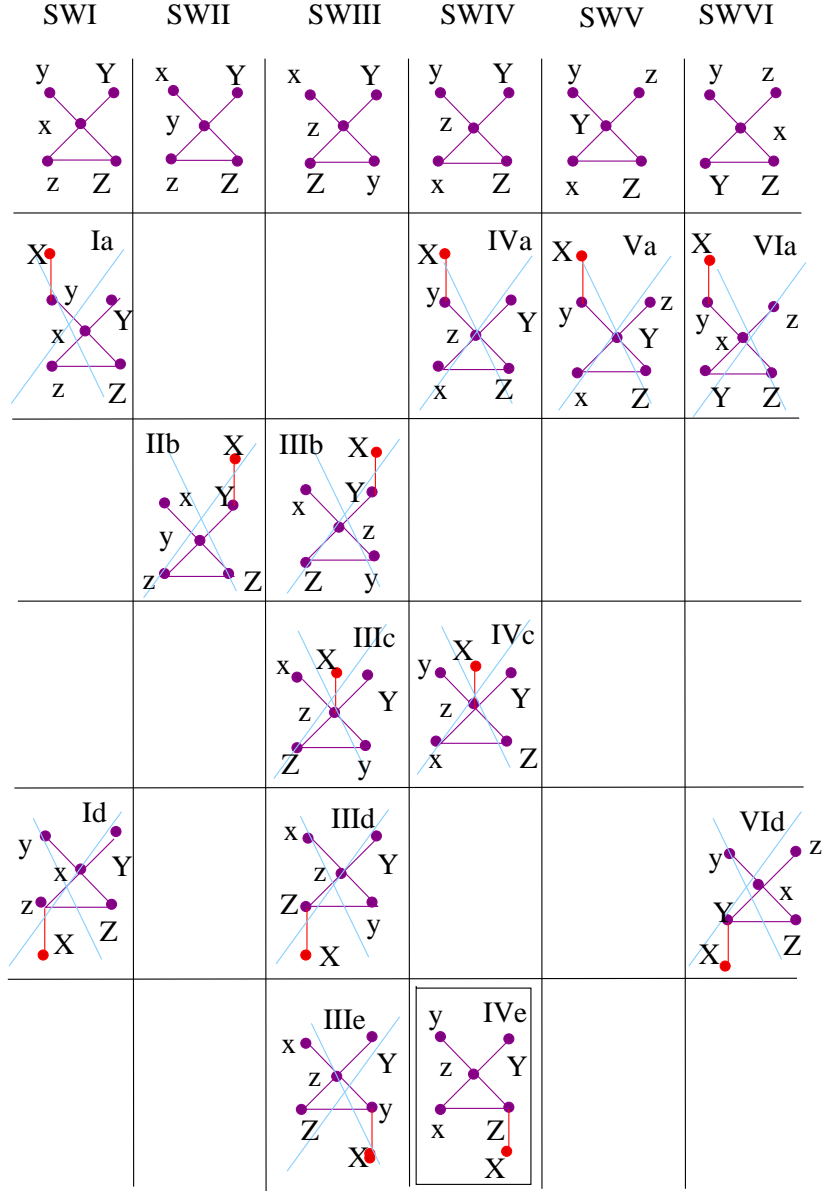}
\caption{{\small{\emph{We crossed out and ``boxed'' the graphs in the LTT Chart for Graph VII that correspond to the LTT structures crossed out and ``boxed'' in Figure \ref{fig:Birecurrency}}}}}
\label{fig:LTTChart2}
\end{figure}

\end{ex}

\bigskip

\noindent \textbf{\underline{STEP 3}: \emph{A PERMITTED EXTENSION/SWITCH WEB $\mathcal{W}(\mathcal{G})$}}

\vskip5pt

\noindent Consider the LTT structure $G$ for any boxed graph $C(G)$ from the LTT chart.  In what follows, for each $1 \leq i \leq r$, $v_i$ will be used to denote both $d_{k,i}$ and $d_{k-1,i}$.  The attaching vertex of $G$ (the purple vertex of the red edge of $G$) will be denoted by $v_a$ (the ``a'' here is for ``attaching'').  Additionally, in the figures, we will use $V_i$ to denote $\overline{v_i}$.

\begin{description}
\item[A.] We determine all permitted extensions and permitted switches $(g_{k}, G_{k-1}, G_k)$ with $G$ as their destination LTT structure (i.e. $G=G_k$).
\begin{description}
\item[1.] Remove the interior of all black edges from $G$ to retrieve the colored subgraph $C(G)$, called the \emph{CLW-graph} for $G$, from the LTT chart.

\item[2.] For each vertex  $v_s$ distance-one in $C(G)$ from $\overline{v_a}$ we determine two potential ``ingoing LTT Structures'' (a switch and an extension).  In LTT structure notation (where $G=G_k$), $\overline{v_a}$ is denoted $d^a_k$ and the label on the free vertex is denoted $d^u_k$.
\begin{itemize}
\item[(a)] Associated to the distance-one vertex $v_s$, the \emph{potential ingoing extension graph}, $(C(G))^e_{v_s}$, will have:
{\begin{itemize}
\item the same purple subgraph as that of $C(G)$
\item the second index of the label on the free vertex will remain the same, i.e. $d^{pu}_{k-1}=d^{u}_{k-1}$ and, if we use the same letters to label the vertices in $G_k$ and $G_{k-1}$, the same letter remains on the free vertex
\item the red edge will now be attached at $v_s$, instead of at $v_a$
\end{itemize}}

\begin{figure}[H]
\centering
\includegraphics[width=4in]{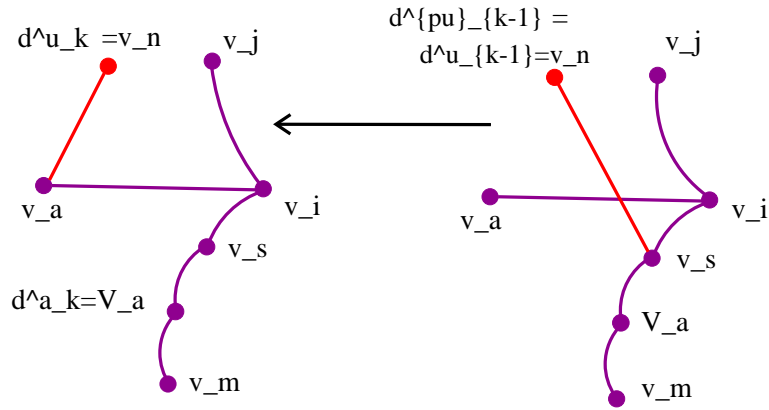}
\caption{{\small{\emph{Potential Ingoing Extension Graph for $v_s$: Notice that the label on the free vertex remains the same ($v_n$) but that the attaching vertex of the red edge changes from $v_a$($=\overline{d^a_k}$) to the distance-one vertex $v_s$.}}}}
\label{fig:PotentialIngoingExtensionGraph}
\end{figure}

\item[(b)] The \emph{potential ingoing switch graph}, $(C(G))^s_{v_s}$, will have an isomorphic purple subgraph to that of $C(G)$.  The labeling is almost the same labeling except that the second index of $d^u_k$ in $G_k$ is the second index of the vertex in $G_{k-1}$ that is mapped by the isomorphism to the vertex labeled with $d^a_k$.  (If we use the same letters to label the vertices in $G_k$ and $G_{k-1}$, it will look as if the label on the red vertex in $C(G)$ has moved to the position in $(C(G))^s_{v_s}$ of what had been $d^a_k$ in $C(G)$ (the inverse of the attaching vertex)).
{\begin{itemize}
\item The label on the red vertex in $(C(G))^s_{v_s}$ is $d^a_k$.
\item The red edge in $(C(G))^s_{v_s}$ is attached at $v_s$.
\item The label on the free vertex is now $d^a_k$.
\end{itemize}}
\end{itemize}

\begin{figure}[H]
\centering
\includegraphics[width=3.8in]{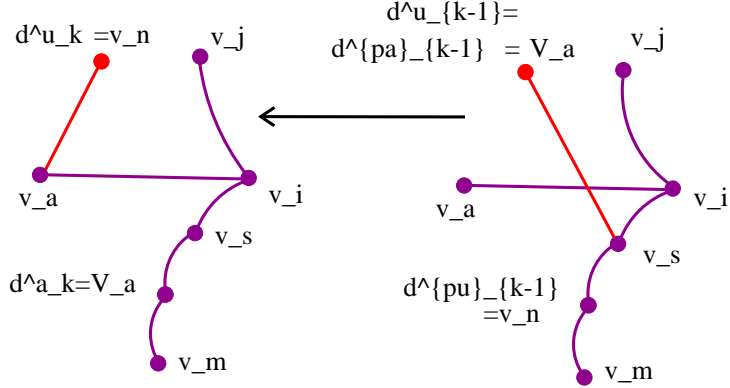}
\caption{{\small{\emph{Potential Ingoing Switch Graph for $v_s$:  Notice how the purple graph with its labels looks the same except for movement of $V_a$ and $v_n$.  Also notice how the red edge shifts.}}}}
\label{fig:PotentialIngoingSwitchGraph}
\end{figure}

\item[3.] Construct in this way the CLW-graph for each potential ingoing LTT structure corresponding to each vertex $v_s$ that is distance-one in $C(G)$ from $\overline{v_a}$.
\item[4.] Label each CLW-graph $(C(G))^s_{v_s}$ and $(C(G))^e_{v_s}$ with the label of the LTT Chart EPP-isomorphic graph(labels will be of the form $(Roman Numeral)_{Letter}$).
\item[5.] If a CLW-graph is not EPP-isomorphic to a graph in the LTT Chart (i.e. if there is an edge of the from $[x_i, \overline{x_i}]$ where either $x_i$ or $\overline{x_i}$ has valence-one), cross it out.  Also cross out any graph EPP-isomorphic to a graph crossed out in the LTT chart.
\item[6.] Each remaining potential ingoing LTT structure CLW-graph, $H_j$, will be EPP-isomorphic to a graph with a box around it in the LTT Chart.  Box the label of the graph $H_j$ in the web and draw an arrow in the web from each such $H_j$ to $C(G)$.
\item[7.] We now have a graph $C(G)$ with a number of arrows entering it.  Each arrow originates at a potential ingoing LTT structure CLW-graph $H_j$ that is labeled with the label of the EPP isomorphic boxed graph in the LTT chart. The label of each of these $H_j$ is boxed in the web.  If some $H_j$ has the same LTT Chart label as $C(G)$, we box the entire graph $H$ in the web.  Such a boxed graph $H_j$ will be called a \emph{branch end}.
\item[8.] For each potential ingoing LTT structure CLW-graph for $C(G)$ with a box around its label, carry out the same procedure to find all of its potential ingoing LTT structure CLW-graphs, their labels, whether they have their labels boxed, and their arrows.
\begin{itemize}
\item In this stage (and all subsequent stages), box graphs (not crossed out) EPP-isomorphic to graphs that have already occurred elsewhere in the web.  Such boxed graphs are again called \emph{branch ends} and do not need their potential ingoing LTT structures determined (as they are determined elsewhere in the web).
\end{itemize}
\item[9.] Continue this process recursively.
\item[10.] When the recursion ends (see Remark \ref{R:RecursionEnds} for why it must end), one portion of the web is complete.  This process will be repeated to create further portions that will need to be glued together in the end:
\indent \begin{itemize}
\item[(a)] If there is some boxed graph $G_1$ in the LTT chart not EPP-isomorphic to any graph in the portion of the web obtained by starting with $C(G)$, then use any such $G_1$ to start the construction of another portion of the web in the same way we constructed the portion starting with $C(G)$.
\item[(b)] Keep constructing portions of the web as such until all boxed graphs in the LTT chart have arisen in at least one portion of the web constructed.
\end{itemize}
\end{description}
\end{description}

\begin{rk}{\label{R:RecursionEnds}}
The recursion process ends, as the LTT chart is finite (so there are a finite number of ``portions'' of the web) and a branch can only contain each LTT Chart element once.
\end{rk}

\begin{rk}{\label{R:RefiningWeb}}
A refined version of the web could be obtained by gluing together distinct pieces of the web at graphs that are both the starting graph of one piece and (up to EPP-isomorphism) a branch end in a separate piece.  If the graphs glued are only EPP-isomorphic (and not identical), this would also involve an EPP isomorphism of one of the entire pieces glued (this is possible since we are never gluing a piece to itself).  None of this refinement is necessary to properly perform the remaining steps of Method I and so one can just leave their web in a nonminimal number of pieces.
\end{rk}

\begin{ex}{\label{E:MethodIIIWeb}} \emph{A Permitted Extension/Switch Web for Graph VII}

\vskip4pt

\indent Above we saw that the only LTT chart graph with a birecurrent LTT structure is IVe.  $\overline{v_a}=z$ and so the distance-1 vertices are those labeled $y, \overline{y}, \overline{x}$, and $\overline{z}$.  We get four potential ingoing extension graphs by attaching the red edge to each of these four vertices.  However, the only one having a birecurrent LTT structure is that where the red edge was not moved (and thus is still attached to the vertex labeled $z$).  Therefore, the only permitted ingoing extension is the self-map of IVe.  To get the four potential ingoing switch graphs for IVe, we again attach the red edge to each of the four vertices $y, \overline{y}, \overline{x}$, and $\overline{z}$, but this time also switch the labels so that the free vertex is labeled $z$ and what was labeled $z$ before is now labeled $x$.  The only of these potential ingoing switch graphs having a birecurrent LTT structure is the one EPP-isomorphic to IVe, i.e. the one where the attaching vertex is $\overline{x}$.  Thus, the Permitted Extension/Switch Web looks like:

\begin{figure}[H]
\centering
\includegraphics[width=3.4in]{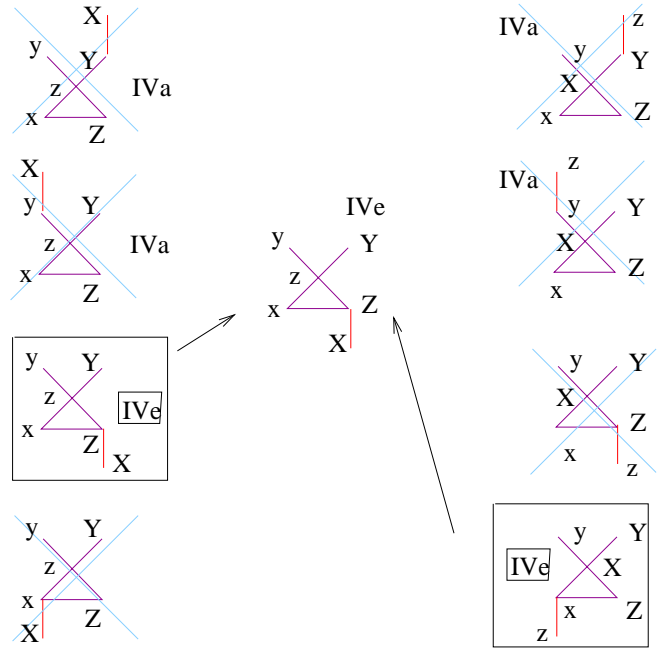}
\caption{{\small{\emph{Permitted Extension/Switch Web for Graph VII (the permitted ingoing extensions are on the left and the permitted ingoing switches are on the right)}}}}
\label{fig:Web}
\end{figure}

\end{ex}

\bigskip

\noindent \textbf{\underline{STEP 4}: THE \emph{(REFINED) SCHEMATIC EXTENSION/SWITCH WEB}}

\vskip5pt

Construct a directed graph $W(\mathcal{G})$ from the information in the Permitted Extension/Switch Web by drawing:

\begin{itemize}
\item a vertex for each boxed graph in the LTT chart labeled as the boxed graph in the LTT chart is labeled

\item a directed edge from a vertex $V_1$ of $W(\mathcal{G})$ to a vertex $V_2$ of $W(\mathcal{G})$ for each arrow in the web $\mathcal{W}(\mathcal{G})$ of Step 3 pointing from a graph with the labeling of $V_1$ to a graph with the labeling of $V_2$.
\end{itemize}

\begin{rk}
Notice that this schematic web is missing information about an edge pair permutation that may occur in its loops, but we currently do not need that information and can/will retrieve it from the web of Step 3 when necessary.
\end{rk}

Before proceeding we trim the schematic web, leaving only its maximal strongly connected components.  This can be done by inductively removing any vertex with no colored edges entering it as well as removing the interiors of edges that exit a removed vertex. We call this ``trimmed'' web a \emph{Refined Schematic Extension/Switch Web} and denote it here by $\mathcal{W}'(\mathcal{G})$. \newline

\begin{ex}{\label{E:MethodIIIRefinedSchematicWeb}} \emph{(Refined) Schematic Extension/Switch Web for Graph VII}

\vskip3pt

\indent Since we have precisely a single extension and switch mapping IVe to itself, the Schematic Permitted Extension/Switch Web is just:

\vskip3pt

\begin{figure}[H]
\centering
\includegraphics[width=1.8in]{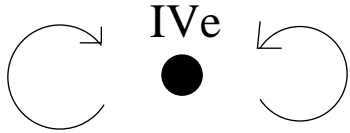}
\caption{{\small{\emph{Schematic Permitted Extension/Switch Web for Graph VII}}}}
\label{fig:SchematicWeb}
\end{figure}

\noindent There was no trimming necessary to find the maximal strongly connected subgraph.  Thus, the Schematic Refined Permitted Extension/Switch Web is the same. \newline
\end{ex}

\bigskip

\noindent \textbf{\underline{STEP 5}: \emph{ILLUSTRATIVE AM DIAGRAM $\mathcal{I}AMD(\mathcal{G})$}}

\vskip5pt

\indent The illustrative AM Diagram, $\mathcal{I}AMD(\mathcal{G})$, for a Type (*) pIW graph $\mathcal{G}$ is constructed as follows:

\begin{description}
\item[(A)] Choose some vertex $V_i$ in the Refined Schematic Web $\mathcal{W}'(\mathcal{G})$.  Let $G(V_i)$ be the graph in the LTT chart labeled the same as $V$.
\item[(B)] Build the corresponding LTT Structure $LTT(V_i)$ from $G(V_i)$ by adding a black edge to $G(V_i)$ for each set of edge-pair vertices $\{x_j, \overline{x_j}\}$ in $G(V_i)$.
\item[(C)] For each directed edge entering $V_i$ in $\mathcal{W}'(\mathcal{G})$ there will be an arrow in $\mathcal{I}AMD(\mathcal{G})$ entering $G(V_i)$.
\item[(D)] To determine the LTT structures at the initial ends of the arrows entering $V_i$:
{\begin{description}
\item 1. Find the corresponding arrow in the web $\mathcal{W}(\mathcal{G})$ of Step 3.
\item 2. If the terminal LTT structure $LTT(V_i)'$ in the web $\mathcal{W}(\mathcal{G})$ differs from $LTT(V_i)$ by an EPP-isomorphism, apply the permutation of vertex labels realizing $LTT(V_i)$ from $LTT(V_i)'$ to the arrow's initial LTT structure in the web $\mathcal{W}(\mathcal{G})$  will give the appropriate initial structure for the arrow in the illustrative AM Diagram.
\item 3. If $LTT(V_i)=LTT(V_i)'$, include only one copy of $LTT(V_i)$, drawing the arrow as a loop.
\item 4. Differentiate between different EPP-isomorphic graphs and include them separately when they occur.
\end{description}}
\item[(E)] Label all arrows terminating at $LTT(V_i)$ by $x \mapsto xy$ where the red vertex of $LTT(V_i)$ is labeled by $x$ and the red edge of $LTT(V_i)$ is $[x, \overline{y}]$.
\item[(F)] For each initial LTT structure for each arrow constructed thus far, carry out the same process as was carried out for $LTT(V_i)$ in (D) and (E). However, in no step should the same LTT structure be drawn twice.  If the initial LTT structure for an arrow is already in the diagram just start the arrow at the copy of the LTT structure already in the diagram.
\item[(G)] Recursively follow this process until every structure in the diagram has the same number of arrows entering it as the corresponding vertex in the schematic diagram had entering it.  If the schematic web had more than one component, the whole process must be carried out for each component of the schematic web. (Several components in the schematic web may end up as part of the same component in the $\mathcal{I}AMD(\mathcal{G})$ if components of the schematic web were not glued in every possible circumstance.)
\end{description}

\begin{ex}{\label{E:MethodIIIillustrativeAMD}} \emph{Illustrative AM Diagram for Graph VII}

\smallskip

\indent By referencing information in the full permitted extension/switch web we can obtain the illustrative AM Diagram.  (There are more components, but each is the same up to EPP-isomorphism, so we show here only a single component.)  Recall that the red edges determine the maps labeling the arrows.

\begin{figure}[H]
\centering
\includegraphics[width=5in]{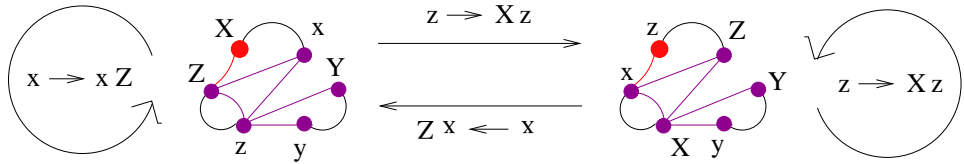}
\caption{{\small{\emph{Illustrative AM Diagram for Graph VII}}}}
\label{fig:IllustrativeAMDiagram}
\end{figure}

\end{ex}

\bigskip

\noindent \textbf{\underline{STEP 6}: \emph{THE AM DIAGRAM $AMD(\mathcal{G})$}}

\vskip5pt

\noindent $AMD(\mathcal{G})$ is obtained from $\mathcal{I}AMD(\mathcal{G})$ by
\begin{description}
\item[(1)] replacing LTT Structures with nodes (labeled by those LTT structures they replaced) and
\item[(2)] replacing the arrows with directed edges (the maps in $\mathcal{I}AMD(\mathcal{G})$ are the labels on the directed edges in $AMD(\mathcal{G})$).
\end{description}

\vskip10pt

\begin{ex}{\label{E:MethodIIIAMD}} \emph{AM Diagram for Graph VII}

\vskip5pt

\indent This is only a single component of the diagram, but including one component suffices to understand the entire diagram, as the others are all the same up to EPP-isomorphism.

\begin{figure}[H]
\centering
\includegraphics[width=4.8in]{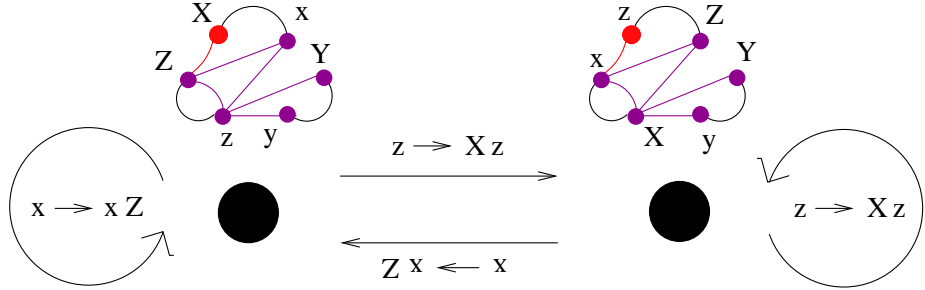}
\caption{{\small{\emph{AM Diagram for Graph VII: Obtained from the illustrative AM Diagram by replacing the LTT structures with nodes (labeled by those LTT structures they replaced) and replacing the arrows with directed edges (the maps in the illustrative AM diagram are the labels on the directed edges in the AM diagram).}}}}
\label{fig:AMD}
\end{figure}

\end{ex}

\vskip7pt

\noindent \textbf{\underline{STEP 7}: \emph{FINDING $L(g_1, \dots, g_k)$}}

\vskip5pt

\begin{description}
\item[(A)] Check Irreducibility Potential: Check whether, in each connected component of $AMD(\mathcal{G})$, for each edge vertex pair $\{d_i, \overline{d_i}\}$, there exists a node in the component such that either $d_i$ or $\overline{d_i}$ labels the red vertex in the LTT structure labeling the node.  There needs to be at least one component where this holds and we only needs to check for $L(g_1, \dots, g_k)$ in such components.  If there is no component for which it holds, then we have shown that the graph $\mathcal{G}$ is unachievable.
\end{description}

\begin{ex}{\label{E:GraphVIIReducibility}} \emph{AM Diagram for Graph VII}

\vskip3pt

\indent We return to the AM Diagram for Graph VII given in Figure \ref{fig:AMD}.  Since $AMD(\mathcal{G})$ contains only red vertices labeled $z$ and $\bar{x}$ (and thus leaves out the edge vertex pair $\{y, \overline{y}\}$) unless some other component contains all three edge vertex pairs ($\{x, \overline{x}\}$, $\{y, \overline{y}\}$, and $\{z, \overline{z}\}$), Graph VII would be deemed unachievable at this stage.  Since no other component does contain all three edge vertex pairs (all components are EPP-isomorphic), the analysis has led us to the conclusion that Graph VII is indeed unachievable. (We will later use these same arguments to deem Graph V unachievable.)
\end{ex}

\begin{description}
\item[(B)] To clarify the discussion here, we refer to the AM Diagram for Graph II.
\end{description}

\begin{figure}[H]
\centering
\includegraphics[width=5in]{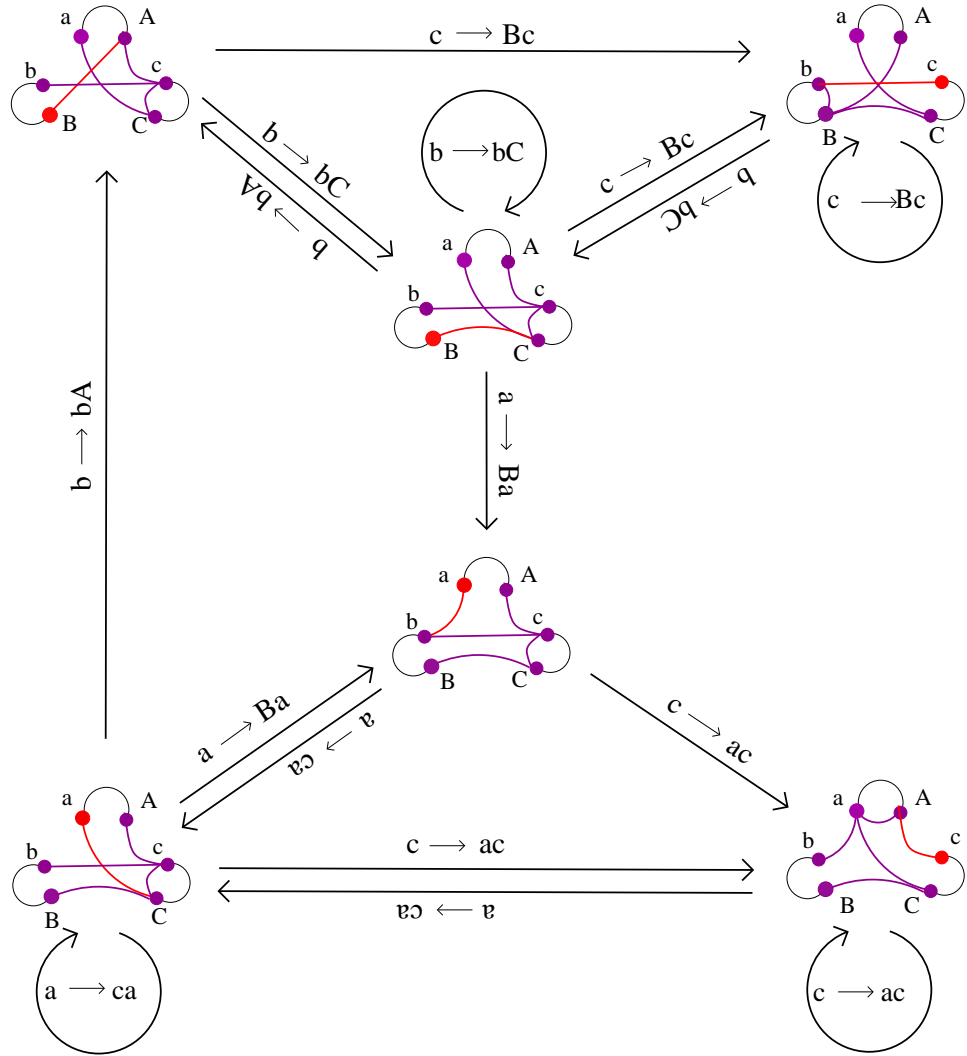}
\caption{{\small{\emph{ $AMD(\mathcal{G})$ where $\mathcal{G}$ is Graph II}}}}
\label{fig:GraphIIAMD}
\end{figure}

Each directed edge in $AMD(\mathcal{G})$ corresponds to either a switch or an extension.  Consider the subdiagram $(AMD(\mathcal{G}))_e$ of $AMD(\mathcal{G})$ consisting precisely of the directed edges corresponding to extensions (and their source and destination LTT structure nodes).

\begin{figure}[H]
\centering
\includegraphics[width=4.5in]{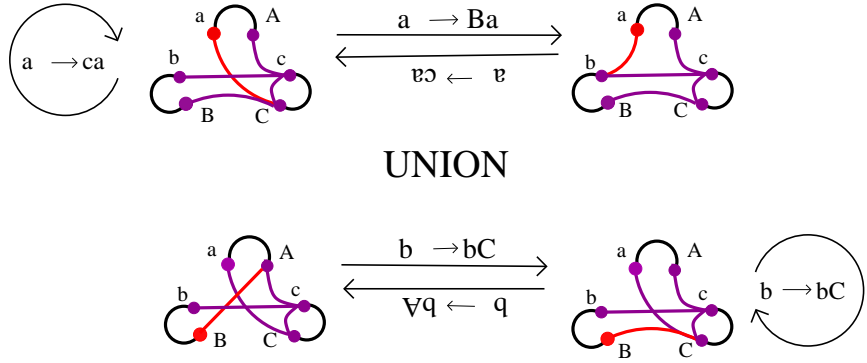}
\caption{{\small{\emph{The entire subdiagram $(AMD(\mathcal{G}))_e$ of $AMD(\mathcal{G})$ where $\mathcal{G}$ is Graph II ($(AMD(\mathcal{G}))_e$ contains only one component here up to EPP-isomorphism)}}}}
\label{fig:GraphIIAMDExtSubdiagram}
\end{figure}

Notice that all of the LTT structures (source and destination LTT structures for the extensions) labeling nodes in a connected component of $(AMD(\mathcal{G}))_e$ share the same purple subgraph (including the same labels on vertices).  We call this purple subgraph the \emph{potential composition PI subgraph} for the component.

\begin{figure}[H]
\centering
\includegraphics[width=3.8in]{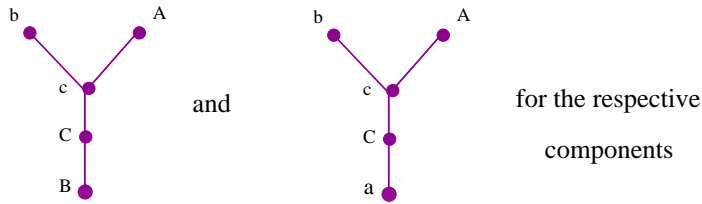}
\caption{{\small{\emph{Graph II Potential Composition PI Subgraph}}}}
\label{fig:GraphIIPotentialCompositionPISubgraph}
\end{figure}

If we add black edges connecting edge pair vertices in the potential composition PI subgraph and then recursively remove any valence-one edges (leaving the vertex with valence greater than one each time we remove a valence-one edge), we get the \emph{potential composition subgraph} for the connected component of $(AMD(\mathcal{G}))_e$.

\begin{figure}[H]
\centering
\includegraphics[width=4in]{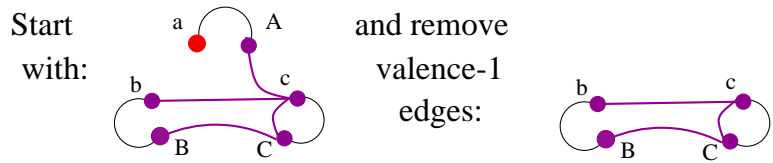}
\caption{{\small{\emph{Graph II Potential Composition Subgraph}}}}
\label{fig:GraphIIPotentialCompositionSubgraph}
\end{figure}

A directed smooth path $[d^a_i, d_{i, j_1}, \overline{d_{i, j_1}}, \dots, d_{i, j_n}, \overline{d_{i, j_n}}]$ (where the potential composition subgraph is viewed as a subgraph of an LTT structure $G_i$) is called a \emph{potential composition path}.

\begin{ex} We give here an example of a potential composition path.

\begin{figure}[H]
\centering
\includegraphics[width=1.4in]{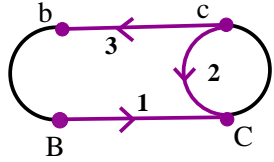}
\caption{{\small{\emph{The numbered colored edges, combined with the black edges between give a Graph II potential composition path. (Note: This path is not used to compute the representative below.)}}}}
\label{fig:GraphIIPotentialCompositionPath}
\end{figure}
\end{ex}

\begin{description}
\item[(C)] There are multiple techniques for finding $L(g_1, \dots, g_k)$.  Here are several:
\begin{itemize}
\item[(1)] One can use potential composition paths to build portions of the graph $\mathcal{G}$ (following progress using preimage subgraphs), see Method III. \newline
\item[(2)] One can test the TT map $g$ corresponding to a loop in $AMD(\mathcal{G})$ for being PNP-free, having Perron-Frobenius transition matrix, and having $IW(g) \cong \mathcal{G}$.  If such is not the case, one can ``attach'' small loops to the initial loop in $AMD(\mathcal{G})$ until the map satisfies those three necessary properties. If $IW(g) \cong \mathcal{G}$, the small loops attached can be determined by potential composition paths to ensure inclusion of necessary remaining edges (keeping in mind that the direction map will map purple edges of the construction path into the destination LTT structure).
    \end{itemize}
\end{description}

\begin{ex} In this example we look at how to find a representative for Graph I.  We start with the AM Diagram:

\noindent \begin{figure}[H]
\centering
\includegraphics[width=5.8in]{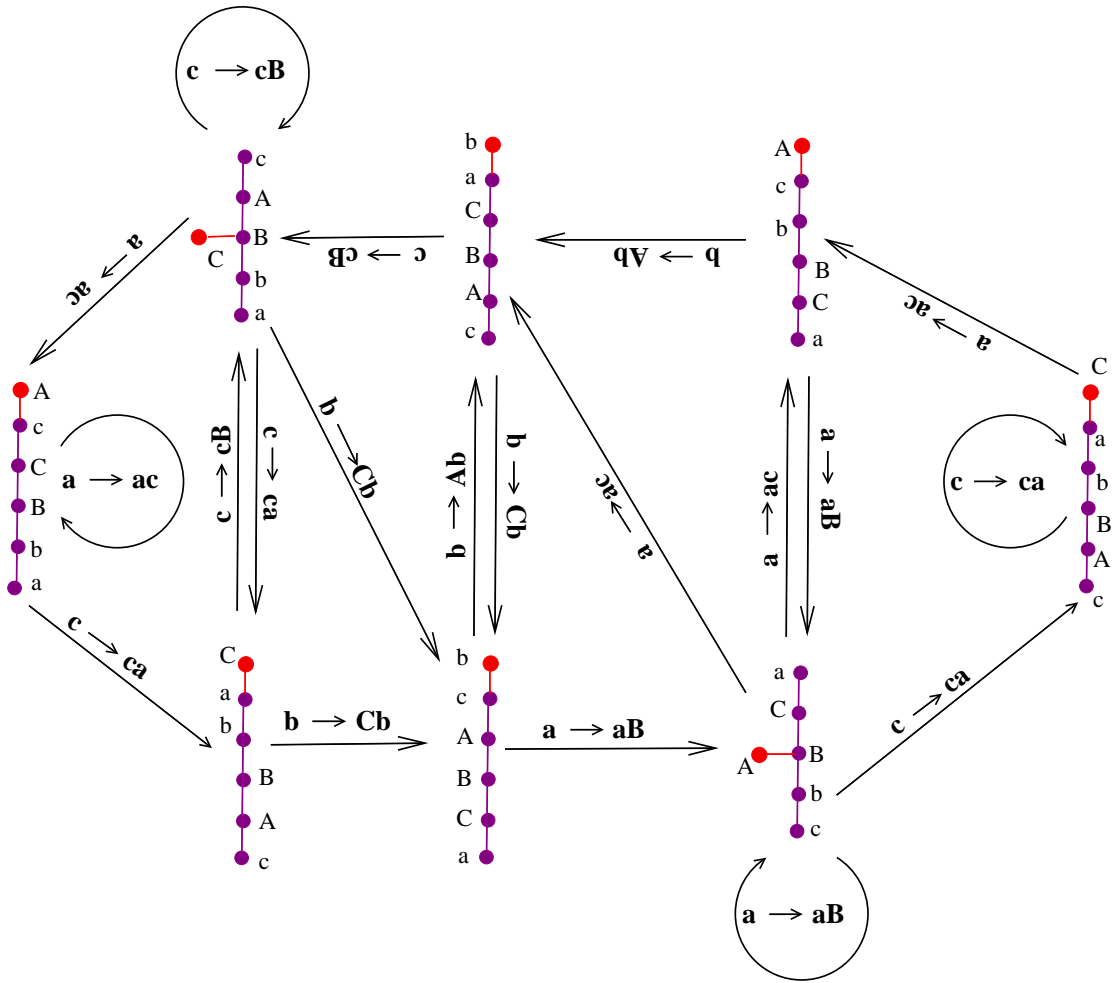}
\caption{{\small{\emph{$AMD(\mathcal{G})$ where $\mathcal{G}$ is Graph I (except that we leave out the black edges in the LTT structures for the sake of simplicity).}}}}
\label{fig:LineAMD}
\end{figure}

\noindent We find a loop to test in $AMD(\mathcal{G})$:

\noindent \begin{figure}[H]
\centering
\includegraphics[width=5.8in]{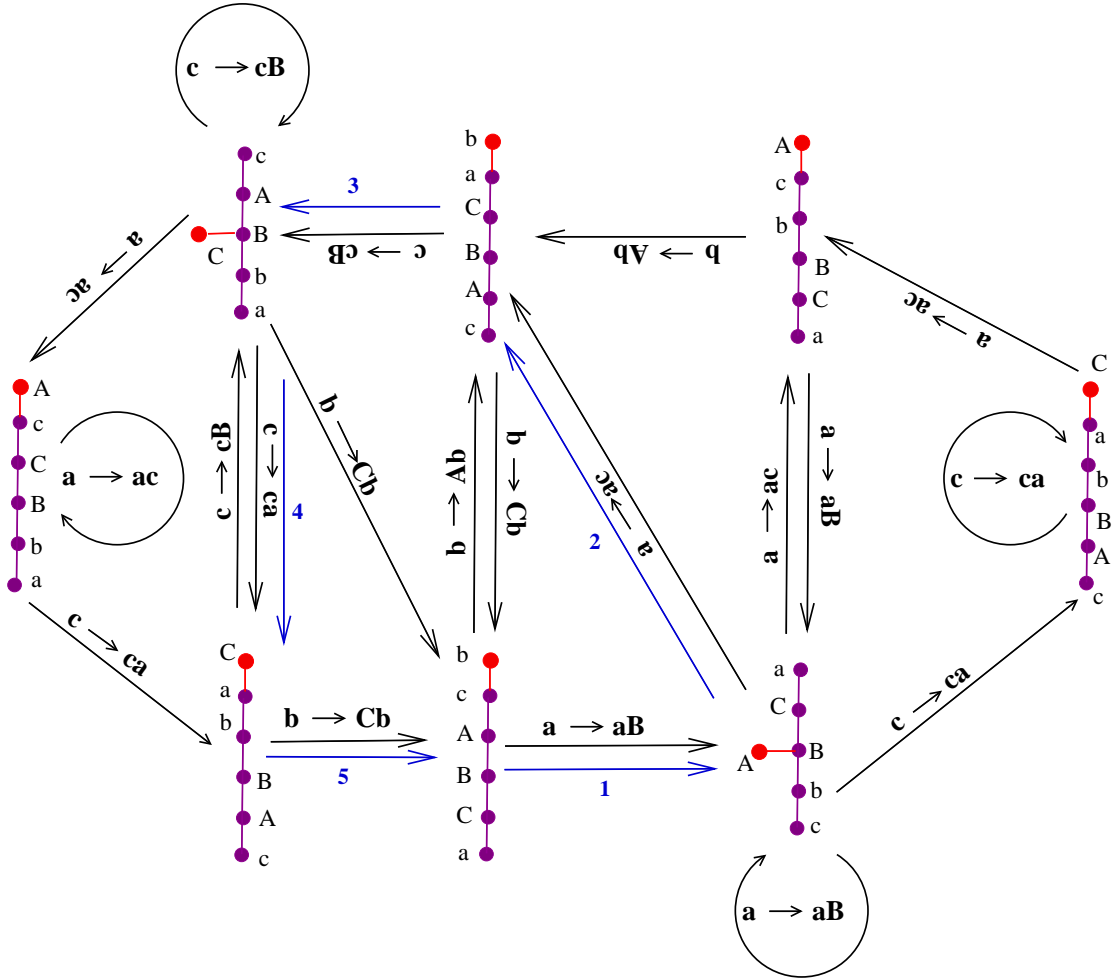}
\caption{{\small{\emph{The blue directed edges together give the loop we will test in $AMD(\mathcal{G})$ where $\mathcal{G}$.}}}}
\label{fig:LineAMDwithLoop}
\end{figure}

\noindent The loop gives:

\noindent \begin{figure}[H]
\centering
\includegraphics[width=4.5in]{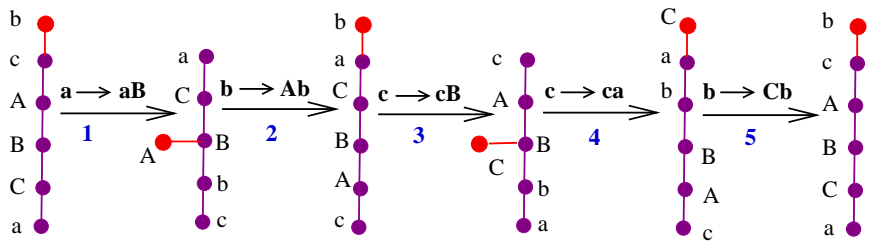}
\caption{{\small{\emph{Corresponding map for loop in $AMD(\mathcal{G})$}}}}
\label{fig:LineLoop1}
\end{figure}

\noindent Since, the loop is not enough (we do not get the edge $[\bar{b}, \bar{c}]$), we add a second loop to it:

\noindent \begin{figure}[H]
\centering
\includegraphics[width=5.8in]{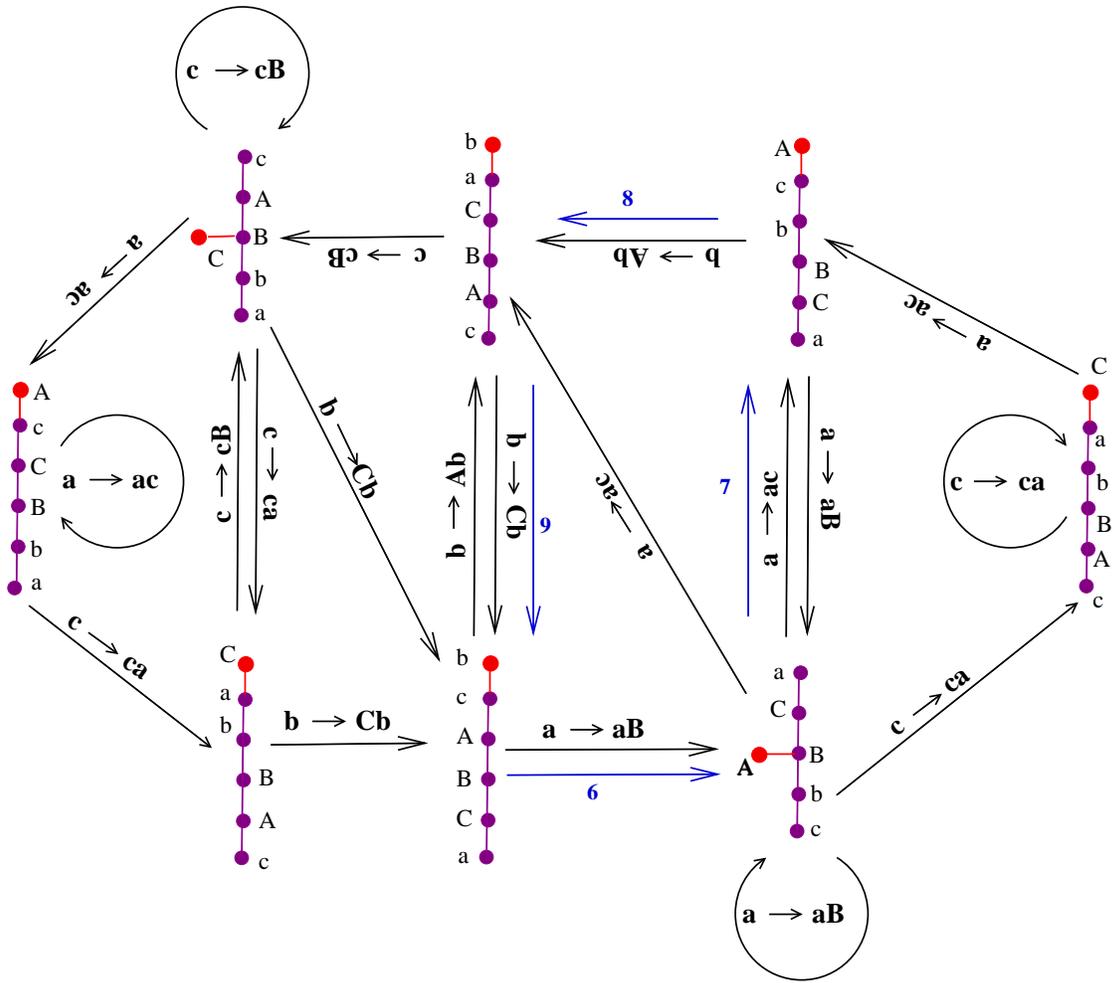}
\caption{{\small{\emph{We find another loop in $AMD(\mathcal{G})$ to add to the first loop.}}}}
\label{fig:LineAMDwithLoop2}
\end{figure}

\noindent The loop gives:

\noindent \begin{figure}[H]
\centering
\includegraphics[width=4in]{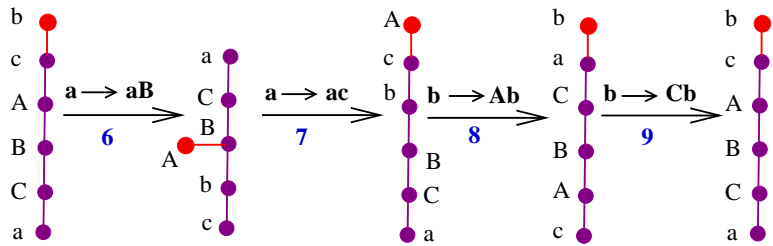}
\caption{{\small{\emph{What the second loop gives}}}}
\label{fig:LineLoop2}
\end{figure}

\noindent Combining the two loops we get the representative yielding Graph I (the line): \newline

\vskip5pt

\noindent \begin{figure}[H]
\centering
\includegraphics[width=5.5in]{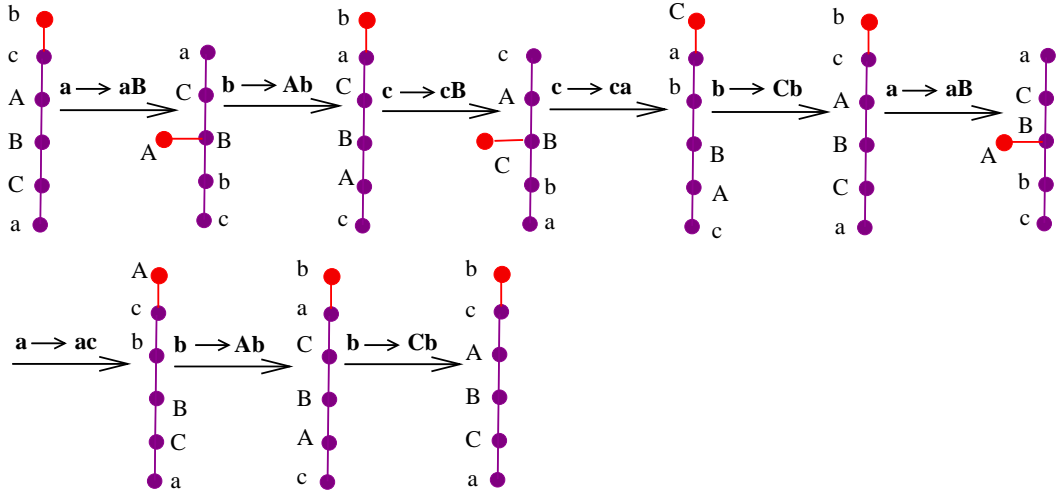}
\caption{{\small{\emph{Ideal Decomposition for the representative yielding Graph I (the line)}}}}
\label{fig:Example2New}
\end{figure}
\end{ex}

\vskip10pt

\noindent \textbf{\underline{STEP 8}: \emph{FINAL CHECKS}}

\vskip5pt

\noindent We need to check that:
\begin{itemize}
\item[(1)] There are $2r-1$ fixed directions.
\item[(2)] The map constructed is PNP-free.
\item[(3)] It is also important to check that the entire Type (*) pIW graph is ``built.''  We can do this by looking at the graph that is the union of the $Dg_{k+1,n}(t^R_k)$.

\begin{ex}{\label{E:GraphBuilding}} We show here an example of how to check that the entire Type (*) pIW graph is ``built'' (we iteratively take the image under each $Dg_k$ of the edges ``created'' thus far):

\begin{figure}[H]
\centering
\noindent \includegraphics[width=5.3in]{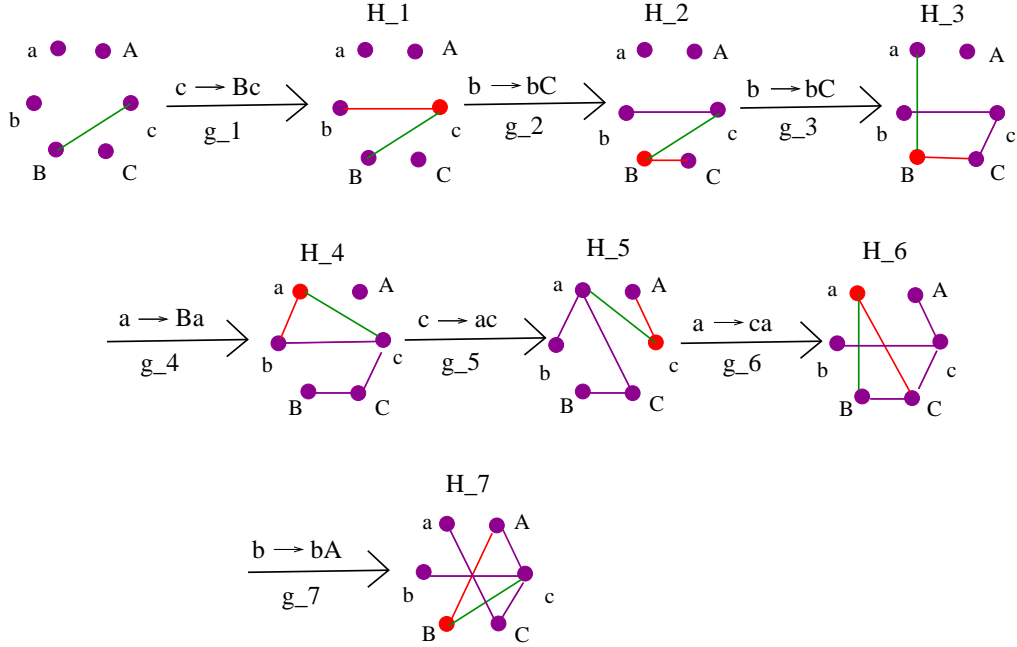}
\caption{{\small{\emph{``Graph Building''}}}}
\label{fig:Union}
\end{figure}

We include subgraphs $H_i$ of the LTT structures $G_i$ to track how edges are ``built.''  The first generator $g_1$ is defined by $c \mapsto \bar{b}c$.  Thus, the red edge $t^R_1$ in the destination LTT structure $G_1$ for $g_1$ will be $[c,b]$, where the red periodic direction vertex is labeled $c$.  The second generator $g_2$ is defined by $b \mapsto b\bar{c}$.  Thus, the red edge $t^R_2$ in the destination LTT structure $G_2$ for $g_2$ will be $[\bar{b},\bar{c}]$, where the red periodic direction vertex is labeled $\bar{c}$. This LTT structure will also contain the image $[c,b]$ of the red edge $[c,b]$ under the direction map $Dg_2:\bar{b} \mapsto c$.  The third generator $g_3$ is defined again by $b \mapsto b\bar{c}$.  Thus, the red edge $t^R_3$ in the destination LTT structure $G_3$ for $g_3$ will be again $[\bar{b},\bar{c}]$, where the red periodic direction vertex is labeled $\bar{c}$. This LTT structure will also contain the image $[c,\bar{c}]$ of the red edge $[\bar{b},\bar{c}]$ and the image $[c,b]$ of the purple edge $[c,b]$ under the direction map $Dg_3:\bar{b} \mapsto c$.  The fourth generator $g_4$ is defined by $a \mapsto \bar{b}a$.  Thus, the red edge $t^R_4$ in the destination LTT structure $G_4$ for $g_4$ will be $[a,b]$, where the red periodic direction vertex is labeled $a$. This LTT structure will also contain the image $[\bar{b},\bar{c}]$ of the red edge $[\bar{b},\bar{c}]$ and the images $[c,b]$ and $[c,\bar{c}]$ of the purple edges $[c,b]$ and $[c,\bar{c}]$ under the direction map $Dg_4:a \mapsto \bar{b}$.  The remaining $H_i$ are constructed in a similar fashion.
\end{ex}
\end{itemize}

\subsection{Method II}{\label{S:MethodI}}

\indent Again let $\mathcal{G}$ be a Type (*) pIWG.  And let $G$ be a Type (*) admissible LTT structure with $PI(G)=\mathcal{G}$ and the standard Type (*) admissible LTT structure notation.

\begin{df}\emph{$(G)_{ep}$} will denote the subgraph of $G$ containing all construction paths corresponding to construction compositions with destination LTT structure $G$.
\end{df}

\vskip10pt

\noindent \textbf{\underline{STEP 1}: FIRST BUILDING SUBGRAPH}

The first step we take in ``building'' a representative $g_{\mathcal{G}}$ with $IW(g_{\mathcal{G}})= \mathcal{G}$ will be to create a subgraph of $G$, which we denote by $G_{ep}'$, that contains and approximates $G_{ep}$.

\begin{df} The \emph{first building subgraph} \emph{$G_{ep}'$} for a Type (*) admissible LTT structure $G$ is obtained from $G$ as follows: \begin{description}
\item 1. Remove the interior of the black edge $[e^u]$, the vertex labeled $\overline{d^u}$, and any purple edges containing the vertex labeled $\overline{d^u}$.  Call the remaining graph $G^1$.
\item 2. Given $G^{j-1}$, recursively define $G^j$:  Let $\{ \alpha_{j-1,i} \}$ be the set of vertices in $G^{j-1}$ not contained in any purple edges in $G^{j-1}$.  $G^j$ will be the subgraph of $G^{j-1}$ obtained by removing all black and purple edges containing some vertex of the form $\overline{\alpha_{j-1,i}}$.
\item 3. $G_{ep}'=\underset{j}{\cap} G^j$.
\end{description}
\end{df}

\begin{ex} We find the first building subgraph $G_{ep}'$ for a Graph XIII LTT structure $G$. \newline

\indent Start with the following LTT Structure for Graph XIII:

\begin{figure}[H]
\centering
\noindent \includegraphics[width=1.2in]{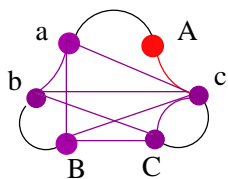}
\caption{{\small{\emph{LTT Structure for Graph XIII}}}}
\label{fig:firstbuildingsubgraph1}
\end{figure}

\indent Remove the interior of the black edge $[\bar a, a]$:

\begin{figure}[H]
\centering
\noindent \includegraphics[width=1.2in]{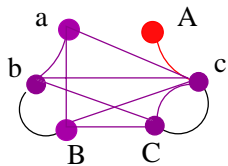}
\caption{{\small{\emph{Having removed the black edge $[\bar a, a]$ interior from the Graph XIII LTT Structure of Figure \ref{fig:firstbuildingsubgraph1}}}}}
\label{fig:firstbuildingsubgraph2}
\end{figure}

\indent Remove $a$ and the interior of all purple edges containing $a$:

\begin{figure}[H]
\centering
\noindent \includegraphics[width=1.6in]{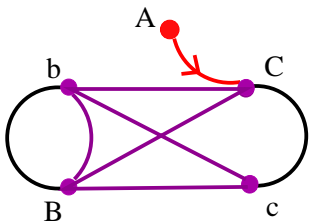}
\caption{{\small{\emph{First Building Subgraph for the LTT Structure in Figure \ref{fig:firstbuildingsubgraph1}}}}}
\label{fig:firstbuildingsubgraph3}
\end{figure}
\end{ex}

\vskip10pt

\noindent Recall that construction compositions behave like Dehn twists and are used to ``construct'' subgraphs of ideal Whitehead graphs. Thus, given the first building subgraph, we look for paths aspiring to be construction paths for construction compositions.

\begin{df}  By a \emph{potential construction path} in a first building subgraph for a Type (*) LTT structure $G$ we will mean a smooth oriented path \newline
$[d^u_i, \overline{d^a_i}, d^a_i, \overline{x_2}, x_2, \dots, x_{n-1}, \overline{x_n}, x_n]$ in $G$ that:
\begin{description}
\item 1. starts with the red edge of $G$ (oriented from $d^u$ to $\overline{d^a}$);
\item 2. is entirely contained in the first building subgraph after the initial red edge and subsequent black edge; and
\item 3. satisfies the following:  Each $G_t$ is an LTT structure (and, in particular, is birecurrent), where $G_t$ is obtained from $G$ by moving the red edge of $G$ to be attached in $G_t$ at $\overline{x_t}$.
\end{description}
\end{df}

\begin{ex}
\noindent A potential composition construction path in the first building subgraph of Figure \ref{fig:firstbuildingsubgraph3} is given by the numbered colored edges and black edges between in:

\begin{figure}[H]
\centering
\noindent \includegraphics[width=1.6in]{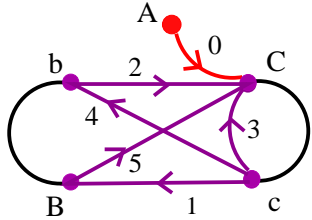}
\caption{{\small{\emph{A potential composition construction path in the first building subgraph of Figure \ref{fig:firstbuildingsubgraph3}}}}}
\label{fig:MethodIIPotentialCompositionConstructionPath}
\end{figure}

\indent The purple edges left after the construction path in Figure \ref{fig:MethodIIPotentialCompositionConstructionPath}:

\begin{figure}[H]
\centering
\noindent \includegraphics[width=.8in]{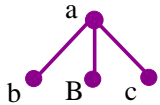}
\caption{{\small{\emph{Purple edges left after construction path in Figure \ref{fig:MethodIIPotentialCompositionConstructionPath}}}}}
\label{fig:PurpleEdgesLeftAfterConstructionPath}
\end{figure}

\end{ex}

\vskip10pt

\noindent \textbf{\underline{STEP 2}: PURIFIED CONSTRUCTION COMPOSITION $h^p_1$}

\vskip3pt

Choose a potential construction path $\gamma = [d^u, \overline{d^a_i}, d^a_i, \overline{x_2}, x_2, \dots, x_{k+1}, \overline{x_{k+1}}]$ in $G_{ep}'$ (a good choice would be one of minimal length among all potential construction paths transversing the maximum number of edges of $G_{ep}'$).  Let $h^p_1$ be the corresponding purified construction composition, if it exists, as in Lemma \ref{L:ConstructionAutomorsphismFromPath}.  If the corresponding construction composition does not exist (for example, if one of the $G_t$ is not birecurrent), then try other potential construction paths until one exists with a corresponding construction composition.  If none can be found, then one can use Method I to find $AMD(\mathcal{G})$ and determine whether $g_{\mathcal{G}}$ exists at all.

We denote the decomposition of $h^p_1$ by
$$\Gamma_{i-k} \xrightarrow{g_{i-k+1}} \cdots \xrightarrow{g_{i-1}}\Gamma_{i-1} \xrightarrow{g_i} \Gamma_i$$
and the corresponding sequence of LTT structures by
$$G_{i-k} \xrightarrow{D^T(g_{i-k+1})} \cdots \xrightarrow{D^T(g_{i-1})} G_{i-1} \xrightarrow{D^T(g_i)} G_i.$$

\begin{ex}{\label{ex:PurifiedConstructionComposition}} Purified construction composition corresponding to the potential composition construction path in Figure \ref{fig:MethodIIPotentialCompositionConstructionPath}:

\begin{figure}[H]
\centering
\noindent \includegraphics[width=5.8in]{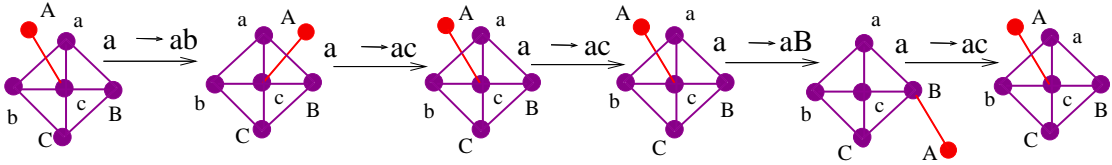}
\caption{{\small{\emph{Purified construction composition corresponding to the potential composition construction path in Figure \ref{fig:MethodIIPotentialCompositionConstructionPath}: We left out the black edges in the LTT structures, as they are not necessary to understand what is going on.}}}}
\label{fig:PurifiedConstructionComposition}
\end{figure}

\noindent It should be noted that all relevant LTT structures are Type (*) admissible LTT structures for $\mathcal{G}$ (and are birecurrent, in particular).

\end{ex}

\vskip5pt

\noindent \textbf{\underline{STEP 3}: SWITCH $s_1$}

\vskip3pt

In this step we determine the switch $(g_{i-k}, G_{i-k-1}, G_{i-k})$ that will precede $h^p_1$ in the decomposition of $g_{\mathcal{G}}$.  To determine choices that may give the switch, one has to look at the source LTT structure $G_{j_1}=G_{i-k}$ for the first generator in the purified construction composition.  There is one potential switch for each purple edge $[d^a_{j_1}, d]$  of $G_{j_1}$ such that $d \neq \overline{d^u_{k_1}}$ in $G_{j_1}$.  Disregard potential switches with source LTT structures that are not birecurrent (or are for other reasons not admissible Type (*) LTT structures).  Choose one of the remaining switches and call it $s_1$.  Denote the source LTT structure $G_{i-k-1}$ by $G_{j_1'}$. \newline

\begin{ex} The two options for the switch proceeding the pure construction composition of Example \ref{ex:PurifiedConstructionComposition} can be summed up giving their source LTT structures (as the generator is determined to be $a \mapsto ac$ by the red edge $[\bar a, c]$ in $G_{i-k}$.  The two source LTT structures are (the black edges in the LTT structures are left out, as they are easily ascertained):

\begin{figure}[H]
\centering
\noindent \includegraphics[width=3.1in]{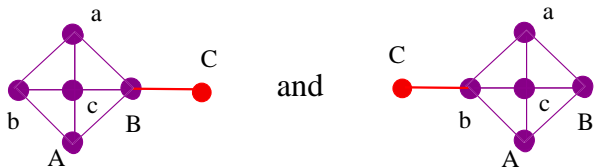}
\caption{{\small{\emph{Options for source LTT structures for switch proceeding pure construction composition of Example \ref{ex:PurifiedConstructionComposition} (in short hand)}}}}
\label{fig:SwitchChoices}
\end{figure}

\noindent Both of the LTT structures are admissible Type (*) LTT structures, so are options.
\end{ex}

\vskip15pt

\noindent \textbf{\underline{STEP 4}: RECURSIVE CONSTRUCTION COMPOSITION BUILDING}

\vskip3pt

\indent In order to track our progress in ensuring that all edges of $\mathcal{G}$ are actually in the ideal Whitehead graph for $g_{\mathcal{G}}$ we establish the notion of a ``preimage subgraph.''

\begin{df} Let ($g_i$, $G_i$, $G_{i-1}$) be a switch.  Recall that there exists an isomorphism from $PI(G_{i-1})$ to $PI(G_i)$ that sends the vertex labeled $d^{pu}_{i-1}$ to the vertex labeled $d^a_i$ and fixes the second index of the label of all other vertices of $PI(G_i)$ (it sends the vertex labeled $d_{i-1,j}$ in $PI(G_{i-1})$ to the vertex labeled $d_{i,j}$ in $PI(G_i)$ for all $d_{i-1,j} \neq d^{pu}_{i-1}$). The isomorphism extends naturally from the vertices to the edges. \newline
\indent The \emph{preimage subgraph} under ($g_i$, $G_i$, $G_{i-1}$) for a subgraph $H \subset PI(G_i)$ will be denoted $H^{-g_i}$ and is obtained from $H$ by replacing each edge of $H$ with its preimage under the isomorphism.  In other words, for each edge $[d_{(i,j)},d^a_i]$ in $H$ there is an edge $[d_{(i-1,j)},d^{pu}_i]$ in $H^{-g_i}$ and for each edge $[d_{(i,j)}, d_{(i,j')}]$ in $H$, where $d^a_i \neq d_{i,j}$ and $d^a_i \neq d_{i,j'}$, there is an edge $[d_{(i-1,j)}, d_{(i-1,j')}]$ in $H^{-g_i}$.
\end{df}

\begin{rk} One could obtain the preimage subgraph $H^{-g_i}$ from $H$ simply by changing the label in $H$ of $d^a_i$ to $d^{pu}_{i-1}$, as well as the labels $d_{i,j}$ to $d_{i-1,j}$ for $d_{i,j} \neq d^a_i$.
\end{rk}

\indent We define here further notation that will also be used to track our progress in ensuring that all edges of $\mathcal{G}$ are actually in the ideal Whitehead graph for $g_{\mathcal{G}}$.

\begin{df} For the purified construction composition $h^p_1$ with destination LTT structure $G$ we define \emph{$G^a_1$} as the subgraph of $G$ consisting of precisely the purple edges in the construction path for $h_1=h^p_1 \circ s_1$.  Let $P(\gamma_{h_n})$ denoted the set of purple edges in the construction path $\gamma_{h_n}$ for $h_n=h^p_n \circ s_n$.  Then $G^a_1=P(\gamma_{h_1})$ and we inductively define $G^a_n$ as $P(\gamma_{h_n}) \cup (G^a_{n-1})^{-s_{n-1}}$.
\end{df}

\vskip10pt

Now that we have established the notation to do so, we describe the recursive process of construction composition building.

\vskip5pt

\noindent \textbf{Recursive Process:} \newline
\noindent The following steps are repeated recursively until $G^a_N =PI(G_{j_N})$ for some $N$.
\indent \begin{description}
\item[I.] Determine the first building subgraph $(G_{j_n'})'_{ep}$ for $G_{j_n'}$.
\item[II.] Find a potential construction path in $(G'_{ep})_n$ (an ``optimal strategy,'' similar to that in Step 2, may involve choosing the path to be of minimal length among all potential construction paths transversing the maximum number of colored edges of $(G'_{ep})_n-G^a_{n+1}$).  Let $h^p_n$ be the corresponding purified construction composition, if it exists, as in Lemma \ref{L:ConstructionAutomorsphismFromPath}.  If the corresponding construction composition does not exist (for example, if one of the $G_t$ in the decomposition is not birecurrent), then try other potential construction paths until one has a valid corresponding construction composition.  If no valid corresponding construction composition can be found, then try using different construction compositions in the previous steps.  If this does not work, one can use Method I to find $AMD(\mathcal{G})$ and determine whether $g_{\mathcal{G}}$ exists at all.

\indent We denote the decomposition of $h^p_1$ by
$$\Gamma_{(i_n-k_n)} \xrightarrow{g_{(i_n-k_n+1)}} \cdots \xrightarrow{g_{(i_n-1)}}\Gamma_{(i_n-1)} \xrightarrow{g_{i_n}} \Gamma_{i_n}$$
and the corresponding sequence of LTT structures by
$$G_{(i_n-k_n)} \xrightarrow{D^T(g_{(i_n-k_n+1)})} \cdots \xrightarrow{D^T(g_{(i_n-1)})} G_{(i_n-1)} \xrightarrow{D^T(g_{i_n})} G_{i_n}.$$

\item[III.] Determine $s_n$: There is one potential switch for each purple edge $[d^a_{j_{n+1}},d]$ of $G_{j_{n+1}}=G_{(i_n-k_n)}$ such that $d \neq \overline{d^u_i}$ in $G_{j_{n+1}}$.  Disregard potential switches with source LTT structures that are not birecurrent (or are for other reasons not admissible Type (*) LTT structures).  Choose one of the remaining switches and call it $s_{n+1}$.   Denote the source LTT structure for $s_{n+1}$ by $G_{j_{n+1}'}$.
\item[IV.] Repeat (I)-(III)  recursively until $G^a_N =PI(G_{j_N})$ for some $N$.
\end{description}

\vskip10pt

\begin{ex}{\label{E:GraphXIII}}  We continue with the example for Graph XIII \newline

\noindent \includegraphics[width=6.2in]{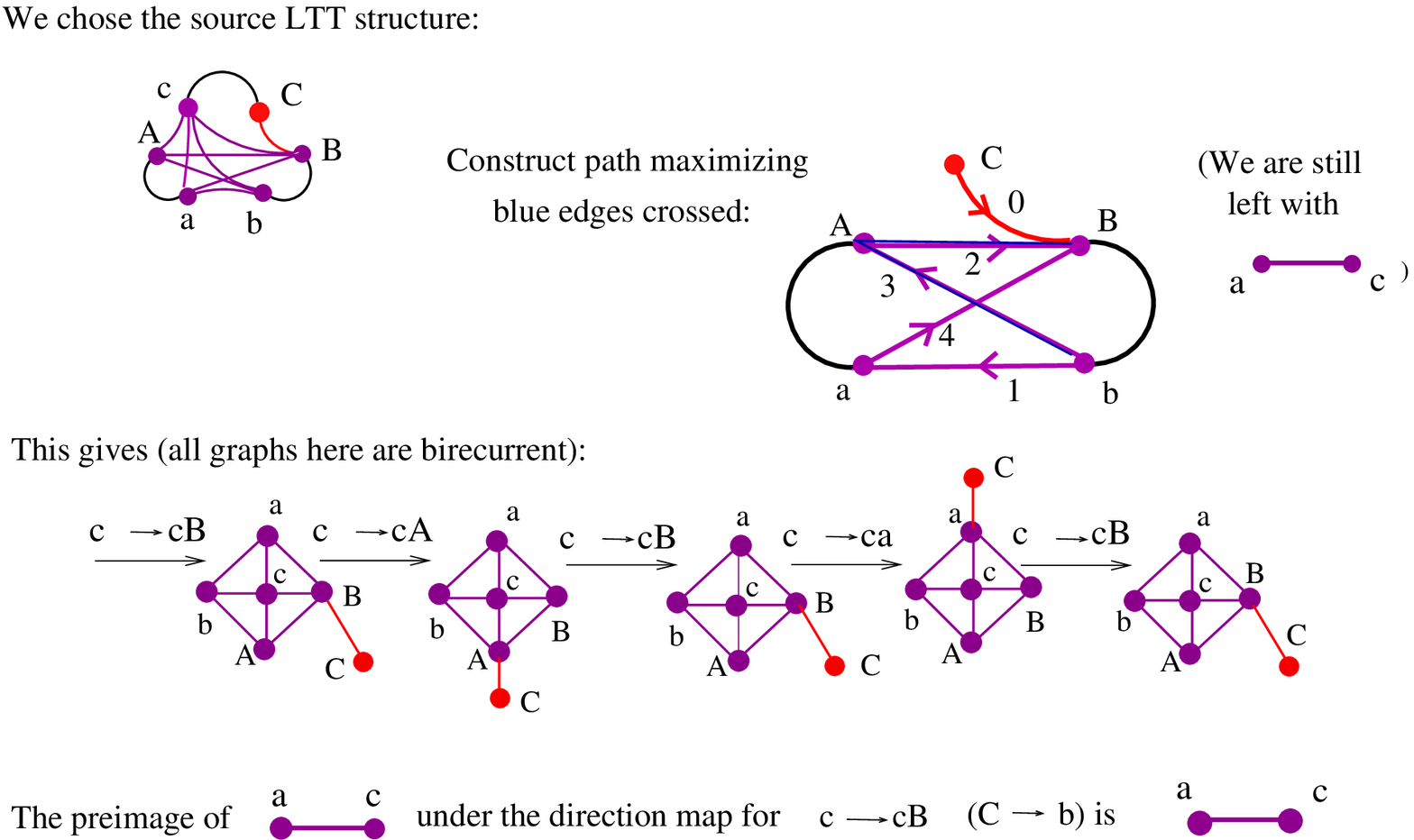} \newline

\vskip5pt

\noindent \includegraphics[width=6.9in]{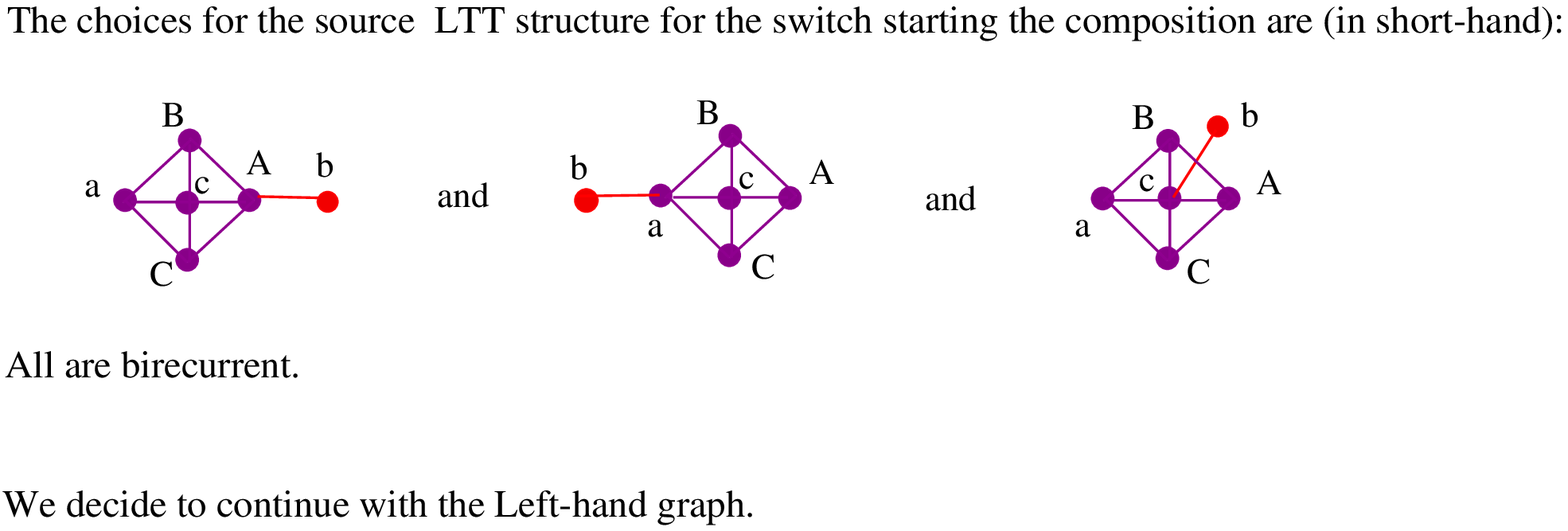} \newline

\vskip5pt

\noindent \includegraphics[width=4.7in]{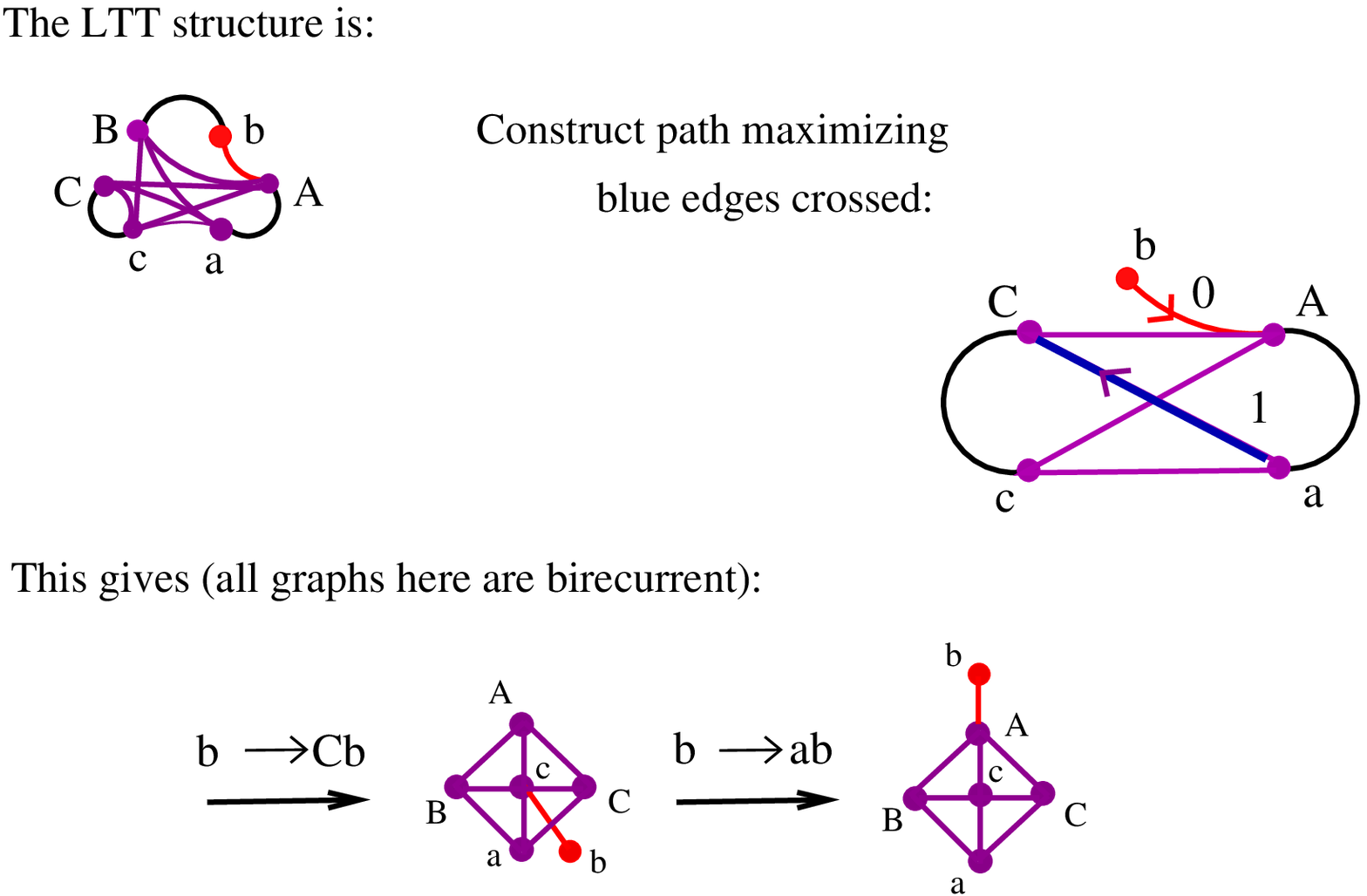}

\end{ex}

\vskip13pt

\noindent \textbf{\underline{STEP 5}: CONCLUDING SWITCH SEQUENCE}

\vskip3pt

Once we have that $G^a_N =PI(G_{j_N})$, we find the shortest possible switch sequence
$$G=G_{(i_N-k_N)} \xrightarrow{g_{(i_N-k_N+1)}} \cdots \xrightarrow{g_{(i_N-1)}} G_{(i_N-1)} \xrightarrow{g_{i_N}} G_{i_N}=G_{j_N}$$
with $G$ as the source LTT structure and $G_{j_N}$ as the destination LTT structure.  A switch path in $G_{j_N}$ may be used for this purpose, though it will be necessary to check that each $G_j$ with $i_N-k_N \leq j \leq i_N = j_N$ is an admissible Type (*) LTT structure for $\mathcal{G}$ (and, in particular, is birecurrent).

If it is not possible to get a pure sequences of switches, then one can try any permitted composition with $G$ as its source LTT structure (and $G_{j_N}$ its destination LTT structure) or, if necessary, find a path in $AMD(\mathcal{G})$ from $G$ to $G_{j_N}$ (see Method I).  It may be possible to find the path in $AMD(\mathcal{G})$ without actually building the entire diagram by instead just looking at the portion of the permitted extension/switch web constructed starting with $G_{j_N}$ (see Method I). \newline

\begin{ex} We continue with the example for Graph XIII.
\begin{figure}[H]
\centering
\noindent \includegraphics[width=3.9in]{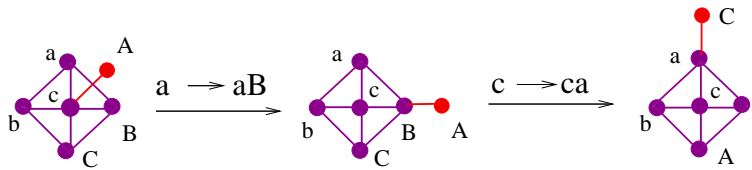}
\caption{{\small{\emph{Concluding sequence of generators for Graph XIII example}}}}
\label{fig:FinalGenerators}
\end{figure}

\noindent At this point we have the final map and get:

\begin{figure}[H]
\centering
\noindent \includegraphics[width=5.8in]{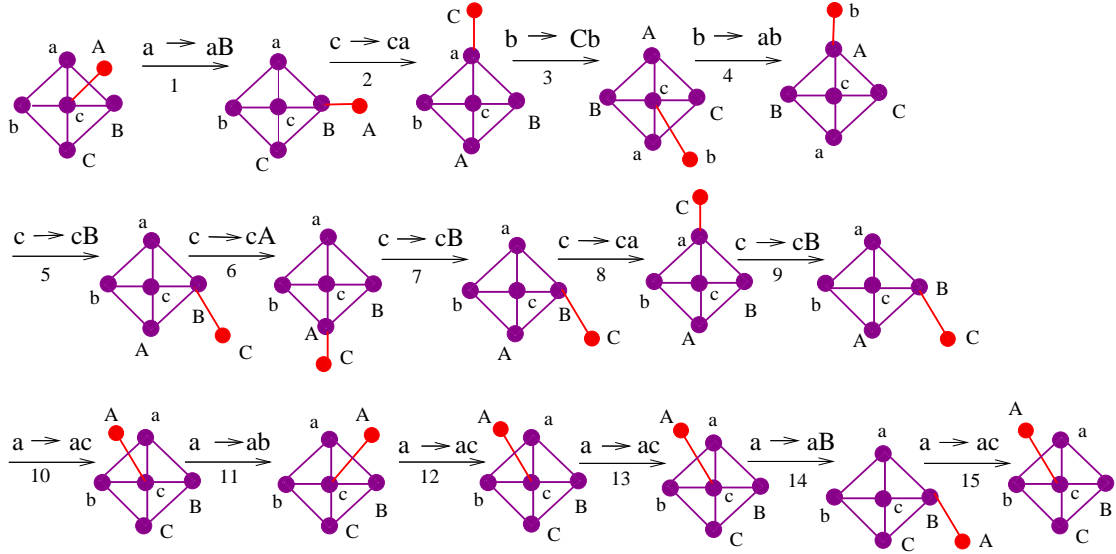} \newline
\caption{{\small{\emph{The entire representative for Graph XIII}}}}
\label{fig:FinalGenerators}
\end{figure}

\vskip5pt

\noindent We showed that this map does not have any PNPs in Section \ref{Ch:NPIdentification}.

\end{ex}

\vskip15pt

\begin{rk} There are choices for the switches at the beginning of each construction composition $h^p_n$ that make for shorter and simpler representatives.  The following guidelines are helpful.
\begin{itemize}
\item[(1)] We want that the switch changes the unachieved direction to be in a pair $\{ d^u_k, \overline{d^u_k} \}$ that does not contain the unachieved direction for $h_m$ for any $m \geq n$.
\item[(2)] We want to maximize the number of edges of $\mathcal{G}$ that $h_n$ constructs that have not yet been constructed by the $h_m$ with $m \geq n$.  To do this it helps that the switch does not result in an unachieved direction pair $\{d^u_k, \overline{d^u_k} \}$ where a lot of the edges in $(G'_{ep})_n-G^a_{n-1}$ that are not transversed by the construction path $\gamma_{h_n}$ contain the vertex labeled $\overline{d^u_k}$.
    \end{itemize}
\end{rk}

\vskip10pt

\noindent \textbf{\underline{STEP 6}: FINAL CHECKS}

\vskip3pt

\noindent The map is not acceptable if any of the following holds:
\begin{itemize}
\item[(1)] for some vertex edge pair $\{d_i, \overline{d_i}\}$, neither $d_i$ nor $\overline{d_i}$ is the red vertex in any LTT structure in the decomposition;
\item[(2)] there are not $2r-1$ fixed directions; or
\item[(3)] the map constructed has a PNP.
\end{itemize}

Check (1) visually, check (2) by composing directions maps of the generators, and then check (3) via the procedure in Section \ref{Ch:NPIdentification}.  If (1) fails, then one can try finding an alternative concluding switch sequence resolving the problem.  If (2) fails, one can simply take a power of the map so that all periodic directions are fixed.

\vskip15pt

\subsection{Method III}

\vskip10pt

\begin{description}
\item[(A)] Find a switch sequence $(g_{(i,i-k)}, G_{i-k-1}, G_i)$ with $G_{i-k-1}=G_i$ such that, for each vertex pair $\{d_i, \overline{d_i}\}$, either $d_i$ or $\overline{d_i}$ is the red vertex in some LTT structure in the sequence.  (Such a composition would be represented by a loop in $AMD(\mathcal{G})$ and can be found as a loop in $AMD(\mathcal{G})$ if not by trial and error.  It would also work to use a loop in $AMD(\mathcal{G})$ that does not represent a switch sequence, but the condition on vertex pairs still holds.)

\item[(B)] As in Method II, find a construction path in $(G_i)'_{ep}$ transversing as many edges of $(G_i)'_{ep}$ as possible, except that we now have the added condition that the corresponding purified construction composition must start and end with the same LTT structure.

\item[(C)] Proceed as in Method II with the added condition of (B) and with the condition that the switches between the purified construction compositions are determined by the switch sequence $(g_{(i,i-k)}, G_{i-k-1}, G_i)$.

\item[(D)] The map constructed is not acceptable if it has any PNPs and so the procedure of Section \ref{Ch:NPIdentification} must be used to check that it does not have any PNPs.  It is additionally still important that all periodic directions are fixed and that the map constructed has all of $\mathcal{G}$ as its ideal Whitehead graph.  Reference Method II for how to ensure this is the case.
    \end{description}

\begin{ex}  We return to Graph XX:

\smallskip

A switch sequence for this graph is given in Example \ref{E:SwitchPath}.  Our first construction composition was given in Example \ref{E:ConstructionAutomorphism}.  What was still needed after that composition was:

\begin{figure}[H]
\centering
\includegraphics[width=.9in]{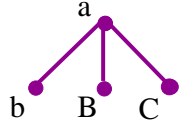}
\caption{{\small{\emph{Edges still needed after the composition given in Example \ref{E:ConstructionAutomorphism} (in Figure \ref{fig:ConstructionCompositionSequence})}}}}
\label{fig:WhatLeft}
\end{figure}

\noindent We take the preimage under the direction map for the final switch and get:

\begin{figure}[H]
\centering
\includegraphics[width=.9in]{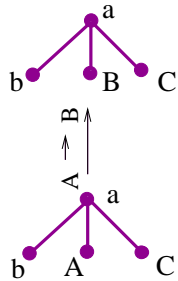}
\caption{{\small{\emph{Preimage of edges left (under the direction map for the switch in Example \ref{E:SwitchPath})}}}}
\label{fig:preimage}
\end{figure}

\noindent Since we could not obtain all of these edges from a single construction composition, we take another preimage:

\begin{figure}[H]
\centering
\includegraphics[width=.9in]{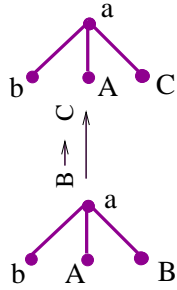}
\caption{{\small{\emph{Preimage of edges left (under the direction map of a second switch in the switch sequences of Example \ref{E:SwitchPath})}}}}
\label{fig:anotherpreimage}
\end{figure}

\noindent We use the construction composition for the following construction path to obtain these edges:

\begin{figure}[H]
\centering
\includegraphics[width=1.6in]{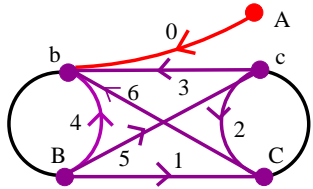}
\caption{{\small{\emph{Construction path in the Graph XX LTT structure used to obtain the remaining edges (given in Figure \ref{fig:anotherpreimage})}}}}
\label{fig:constructionpath}
\end{figure}

\noindent When composed we get:

\begin{figure}[H]
\centering
\includegraphics[width=5.3in]{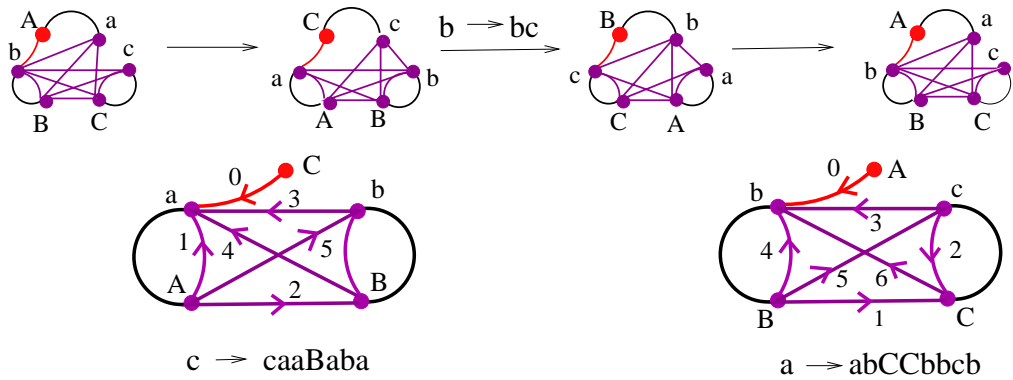}
\caption{{\small{\emph{The combination of everything we have so far in realizing Graph XX}}}}
\label{fig:ConstructionCompositions}
\end{figure}

\noindent The automorphism we have obtained is:

$$
h =
\begin{cases} a \mapsto ab \bar{c} \bar{c} bbcb \\
b \mapsto bc \\
c \mapsto cab \bar{c} \bar{c} bbcbab \bar{c} \bar{c} bbcb \bar{c} \bar{c} \bar{b} ab \bar{c} \bar{c} bbcbbccab \bar{c} \bar{c} bbcb
\end{cases},
$$

\noindent The periodic directions for this map are not fixed.  However, they are fixed when we compose $h$ with itself, so we take $g=h^2$.
\end{ex}

\section{Unachievable Ideal Whitehead Graphs}{\label{Ch:UnachievableIWGs}}

The unachievability of Graph VII was shown in Section \ref{Ch:Procedure}.  In this section we show that Graph II and Graph V are also unachievable.  \emph{In all figures of this section, we continue with the convention that X represents $\bar{x}$, etc.}

\subsection{Four Edges Sharing a Vertex (Graph II)}

The first unachievable 5-vertex Type (*) pIW graph is the graph $\mathcal{G}$ consisting of four edges adjoined at a single vertex.  For this graph we use Proposition \ref{P:Birecurrent}.  What we need is that every Type (*) Admissible LTT Structure for $\mathcal{G}$ is not birecurrent.  Up to EPP-isomorphism, there are two such LTT structures to consider neither of which is birecurrent):

\begin{figure}[H]
\centering
 \includegraphics[width=2.7in]{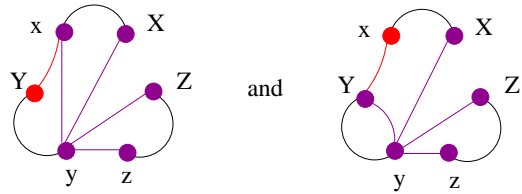}
\caption{{\small{\emph{Potential LTT structures for four edges adjoined at a single vertex}}}}
\label{fig:NotBirecurrent}
\end{figure}

\noindent These are the only structures worth considering as follows: \newline
\indent Either three of the valence-one vertices of $\mathcal{G}$ belong to different edge pairs or the valence-one vertices are labeled with two sets of edge pairs.  First consider the case where the valence-one vertices are labeled with two sets of edge pairs.  The red edge cannot be attached in such a way that it is labeled with an edge pair and all the other resulting LTT structures are the same as the first structure up to EPP-isomorphism.  Now consider the case where three of the valence-1 vertices of $\mathcal{G}$ belong to different edge pairs.  One of these three have the label of the inverse of the valence-four vertex. The red edge can only be attached at one vertex choice and without causing an edge pair labeled vertex set connected by a valence-1 edge in the colored subgraph of the LTT structure.  Up to EPP-isomorphism, this just leaves us with the second structure.

\subsection{Graph V}

\noindent We draw the illustrative AM diagram here without labels on the edges because it is clear even from this much that no map represented by a loop in this diagram would be irreducible (the only edge pairs labeling red vertices are $\{x, \overline{x} \}$ and $\{z, \overline{z} \}$):

\begin{figure}[H]
\centering
\includegraphics[width=4in]{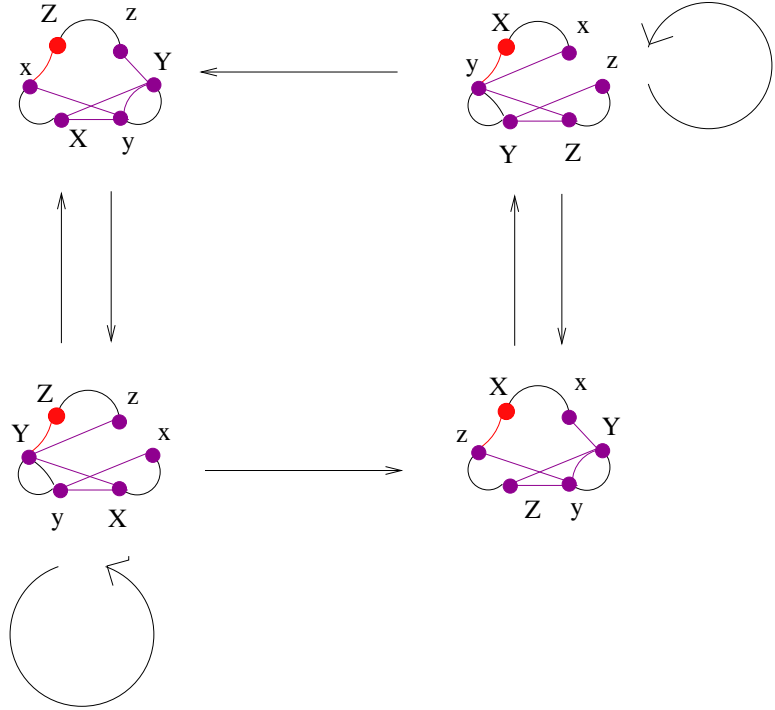}
\caption{{\small{\emph{A component of the Illustrative AM Diagram (all others are the same up to EPP-isomorphism)}}}}
\label{fig:IllustrativeAMDExample}
\end{figure}

\section{Achievable 5-Vertex Type (*) pIW Ideal Whitehead Graphs in Rank 3}{\label{Ch:Achievable}}

This section includes the main theorem of this document.  For the theorem we use our methods to determine which $5$-vertex Type (*) pIW graphs arise as $IW(\phi)$ for ageometric, fully irreducible $\phi \in Out(F_3)$. Since there are precisely twenty-one $5$-vertex Type (*) pIW graphs (see Figure \ref{fig:21Graphs} for a complete list), we can handle them on a case-by-case basis.  \emph{For all figures of this section we continue with the convention that A notates $\bar{a}$.}

\begin{thm}{\label{T:MainTheorem}} Precisely eighteen of the twenty-one 5-vertex Type (*) pIW graphs are ideal Whitehead graphs for ageometric, fully irreducible $\phi\in Out(F_3)$. \end{thm}

\noindent \emph{Proof}: The unachievable graphs (Graph II, Graph V, and Graph VII) were already handled in sections \ref{Ch:Procedure} and \ref{Ch:UnachievableIWGs}.  We now give representatives for the remaining graphs, leaving it to the reader to prove that these representatives are PNP-free (using the procedure of Section \ref{Ch:NPIdentification}), have Perron-Frobenius transition matrices, and have the appropriate ideal Whitehead graph.  Proposition \ref{P:RepresentativeLoops} then gives that they are representatives of ageometric, fully irreducible $\phi\in Out(F_r)$ with the desired ideal Whitehead graphs.

\newpage

\noindent GRAPH I (The Line):

\vskip7pt

\noindent  We give here the representative $g$ whose ideal whitehead graph, $\mathcal{G}$, is GRAPH I:
$$
g =
\begin{cases} a \mapsto ac \bar{b} ca \bar{b} cacac \bar{b} ca  \\
b \mapsto \bar{a} \bar{c} b \bar{c} \bar{a} \bar{c} \bar{a} \bar{c} b  \\
c \mapsto cac \bar{b} ca \bar{b} cac
\end{cases}
$$
\smallskip

\noindent Our ideal decomposition for $g$ is described by the following figure: \newline

\smallskip

\noindent \includegraphics[width=5.5in]{Example2New.eps} \newline

For this we constructed a component of $AM(\mathcal{G})$ and used Method I.  This method made the most sense here as there were only a few birecurrent LTT structures to be included in $AM(\mathcal{G})$.

\newpage

\noindent GRAPH II:

\vskip7pt

\noindent The representative $g$ whose ideal whitehead graph, $\mathcal{G}$, is GRAPH II was again constructed using Method I.  We started with a path in $AM(\mathcal{G})$.

\begin{figure}[H]
\centering
\includegraphics[width=4.8in]{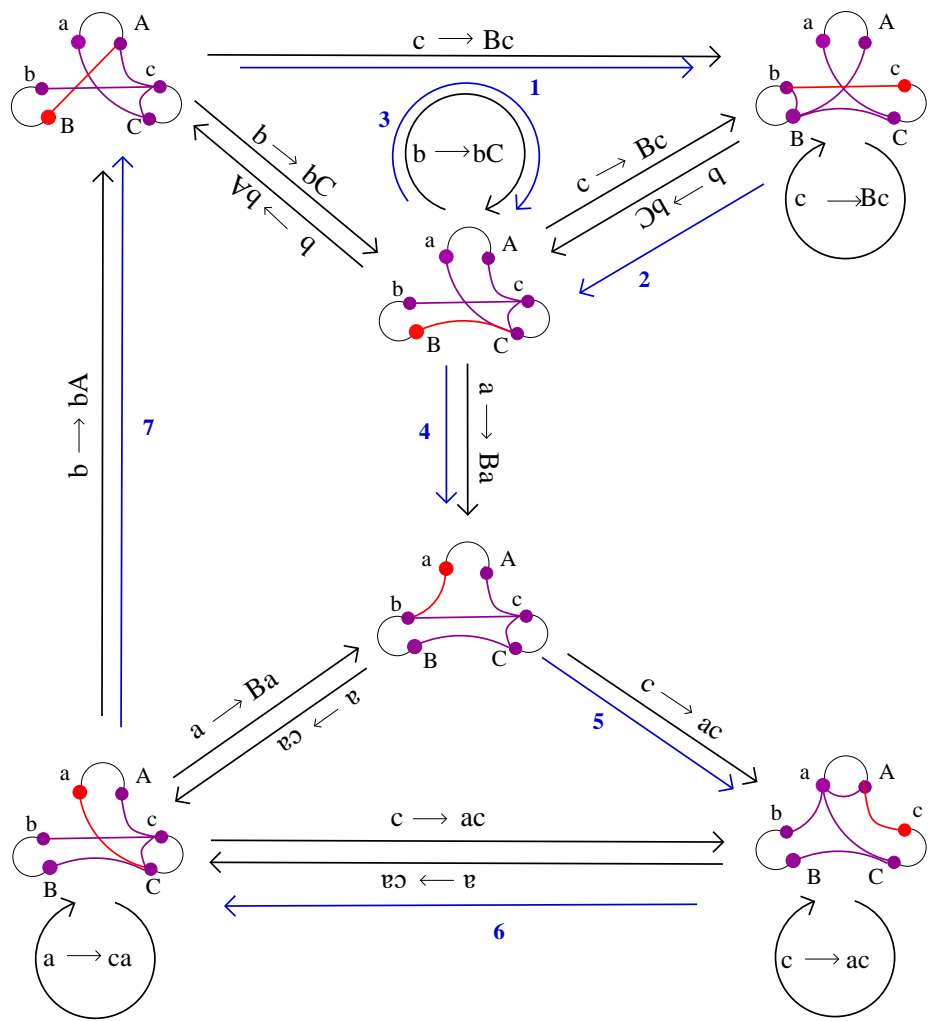}
\caption{{\small{\emph{The blue path in $AM(\mathcal{G})$ gives an ideal decomposition $g$.}}}}
\label{fig:GraphIIAMDPath}
\end{figure}

\vskip10pt

\noindent The path in $AM(\mathcal{G})$ corresponds to the ideal decomposition: \newline

\smallskip

\noindent \includegraphics[width=6.2in]{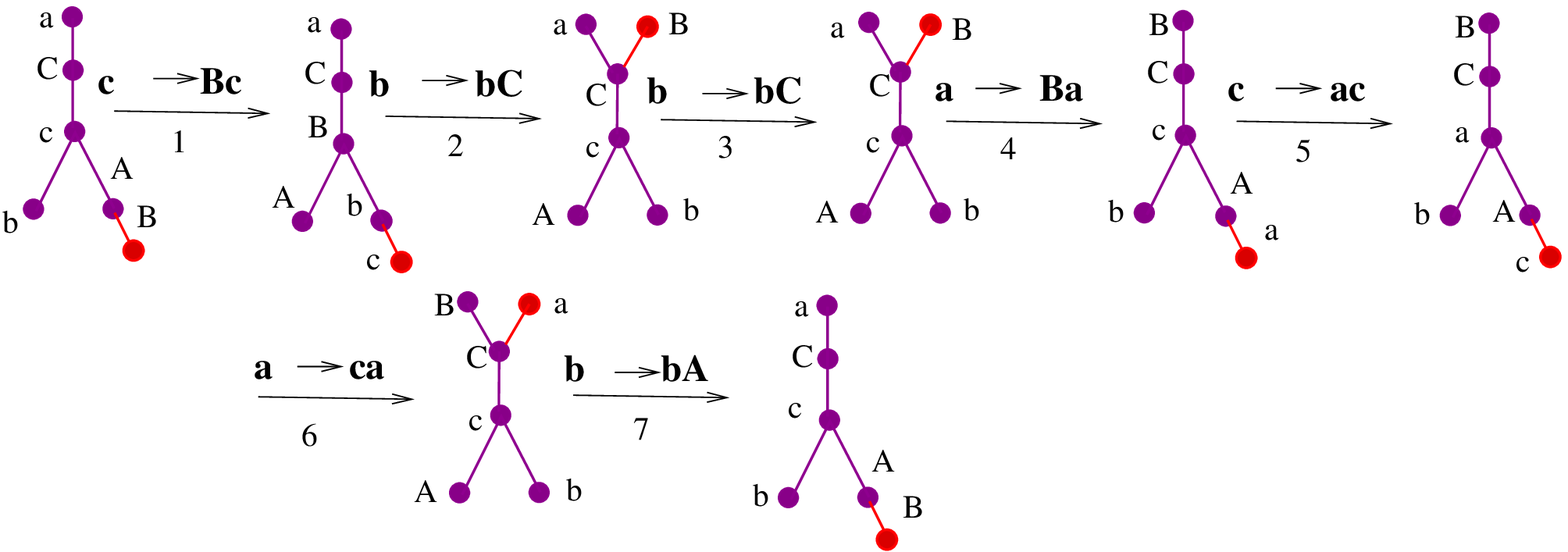}

\vskip5pt

\noindent And our representative is:
$$
g =
\begin{cases} a \mapsto a \bar{b} ca  \\
b \mapsto b \bar{a} \bar{c} \bar{a} \bar{c} \bar{c} \bar{a} \bar{c}  \\
c \mapsto caccaca \bar{b} cac
\end{cases}
$$

\vskip15pt

\noindent GRAPH IV:

\vskip7pt

\noindent The representative $g$ whose ideal whitehead graph, $\mathcal{G}$, is GRAPH IV is:
$$
g =
\begin{cases} a \mapsto c \bar{b} a  \\
b \mapsto bc \bar{a} b  \\
c \mapsto c \bar{b} abc  \bar{a} bc
\end{cases}
$$

\smallskip

\noindent We used Method I here. Constructing $AM(\mathcal{G})$ was exceptionally easy because of the symmetry in the graph.  Our ideal decomposition for $g$ is: \newline

\smallskip

\noindent \includegraphics[width=5in]{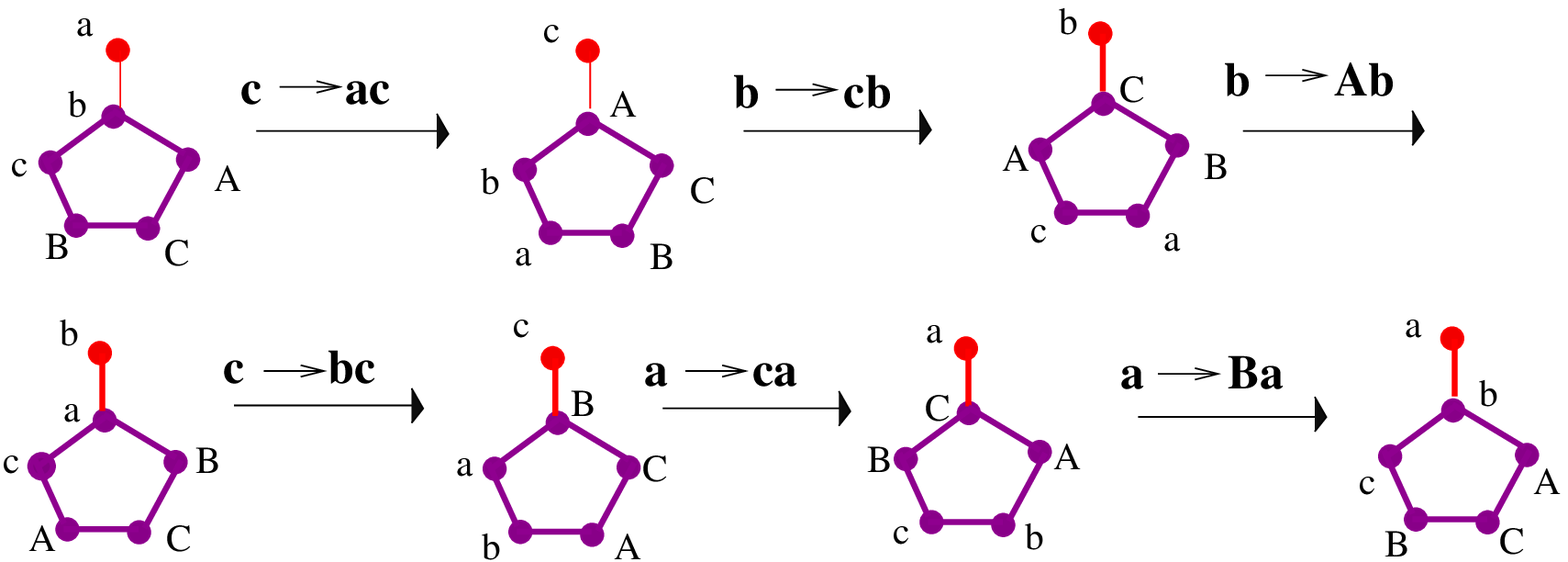} \newline

\newpage

\noindent GRAPH VI:

\vskip5pt

\noindent The representative $g$ whose ideal whitehead graph is GRAPH VI is:
$$
g =
\begin{cases} a \mapsto abacbaba \bar{c} abacbaba  \\
b \mapsto ba \bar{c}  \\
c \mapsto c \bar{a} \bar{b} \bar{a} \bar{b} \bar{a} \bar{b} \bar{c} \bar{a} \bar{b} \bar{a} c  \end{cases}
$$

\smallskip

\noindent Our ideal decomposition for $g$ is described by the following figure: \newline

\smallskip

\noindent \includegraphics[width=6.5in]{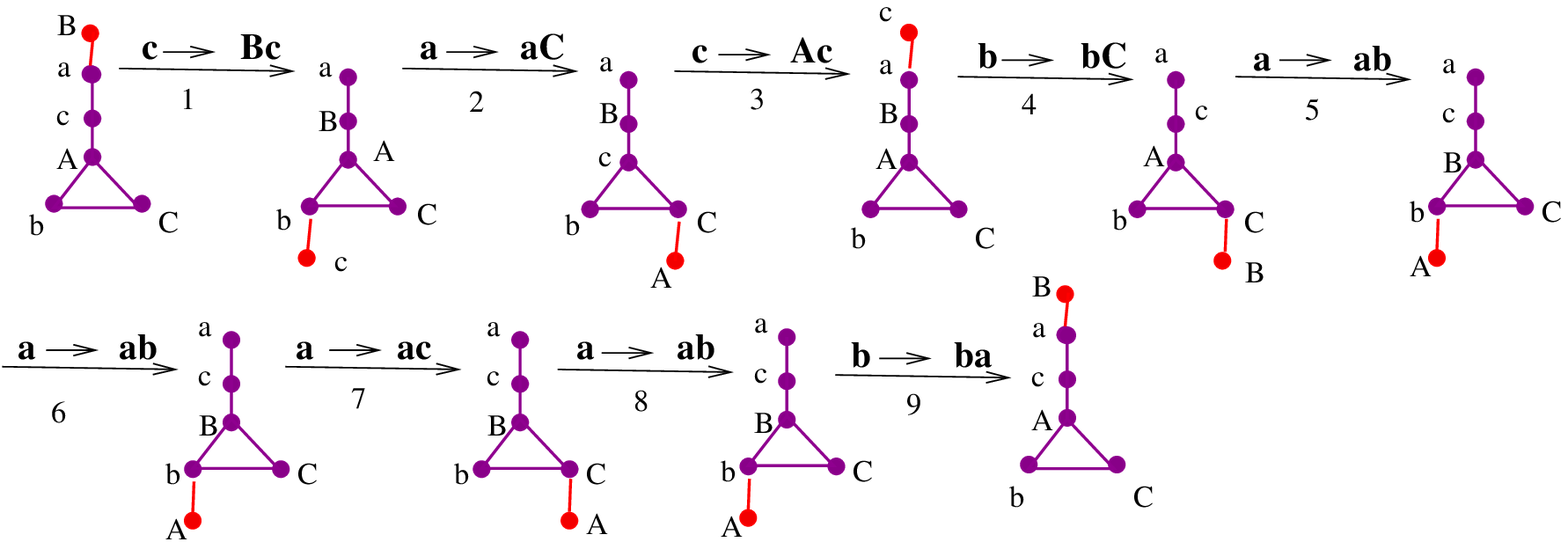}

\vskip15pt

\noindent GRAPH VIII:

\vskip7pt

\noindent The representative $g$ whose ideal whitehead graph, $\mathcal{G}$, is GRAPH VIII is:
$$
g =
\begin{cases} a \mapsto a \bar{c} aa \bar{b} a \bar{c} b \bar{a} \bar{a} ca \bar{c} aa \bar{b} a \bar{c} a  \\
b \mapsto b \bar{a} \bar{a} c  \\
c \mapsto c \bar{a} b \bar{a} \bar{a} c \bar{a} b \bar{a} \bar{a} c
\end{cases}
$$

\noindent Our ideal decomposition for $g$ is described by the following figure: \newline

\smallskip

\noindent \includegraphics[width=6.5in]{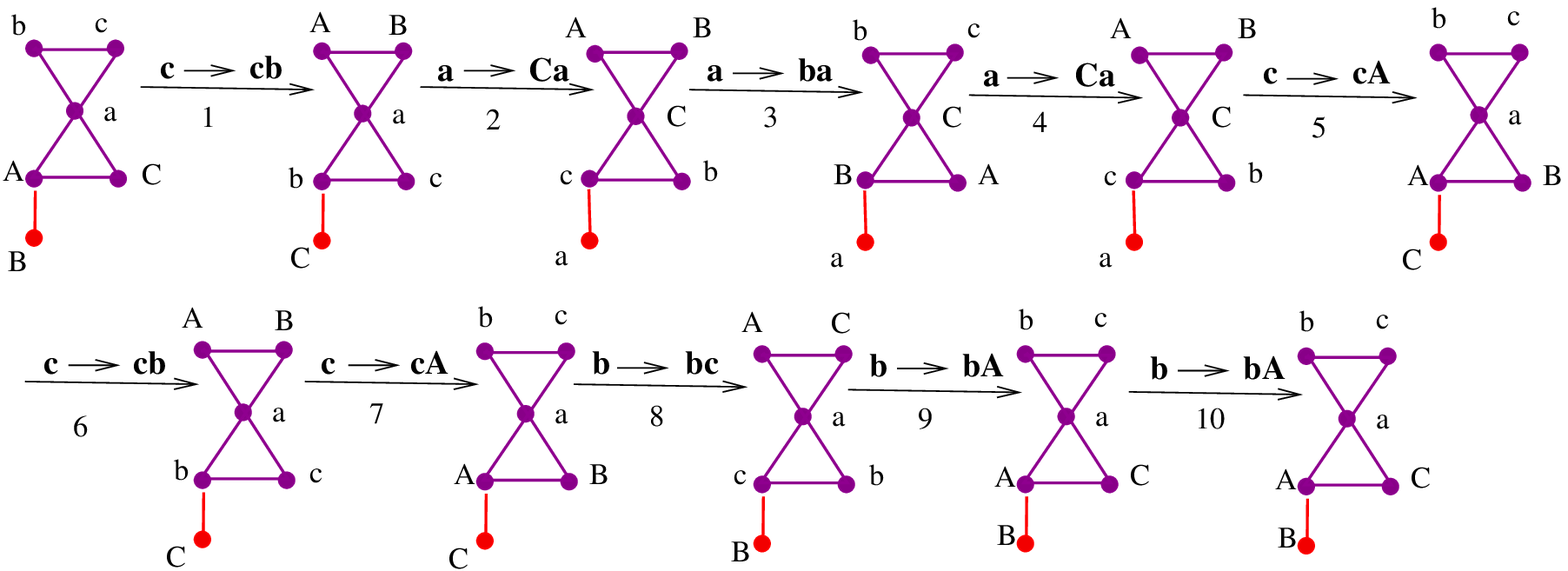} \newline

\smallskip

For this we constructed a component of $AM(\mathcal{G})$ and used Method I.  This method made the most sense here as there were only a few birecurrent LTT structures to be included in $AM(\mathcal{G})$.

\vskip15pt

\noindent GRAPH IX:

\vskip7pt

\noindent The representative $g$ whose ideal whitehead graph, $\mathcal{G}$, is GRAPH IX is:
$$
g =
\begin{cases} a \mapsto ab \bar{c} b \bar{c} ab \bar{c}b \bar{c} \bar{b} c \bar{b} \bar{a}\bar{c} \bar{b} ab \bar{c} b \bar{c} \\
b \mapsto bcab \bar{c} bc \bar{b} c \bar{b} \bar{a} cab \bar{c} b  \\
c \mapsto c \bar{b} c \bar{b} \bar{a}
\end{cases}
$$

\noindent Our ideal decomposition for $g$ is described by the following figure: \newline

\smallskip

\noindent \includegraphics[width=6.5in]{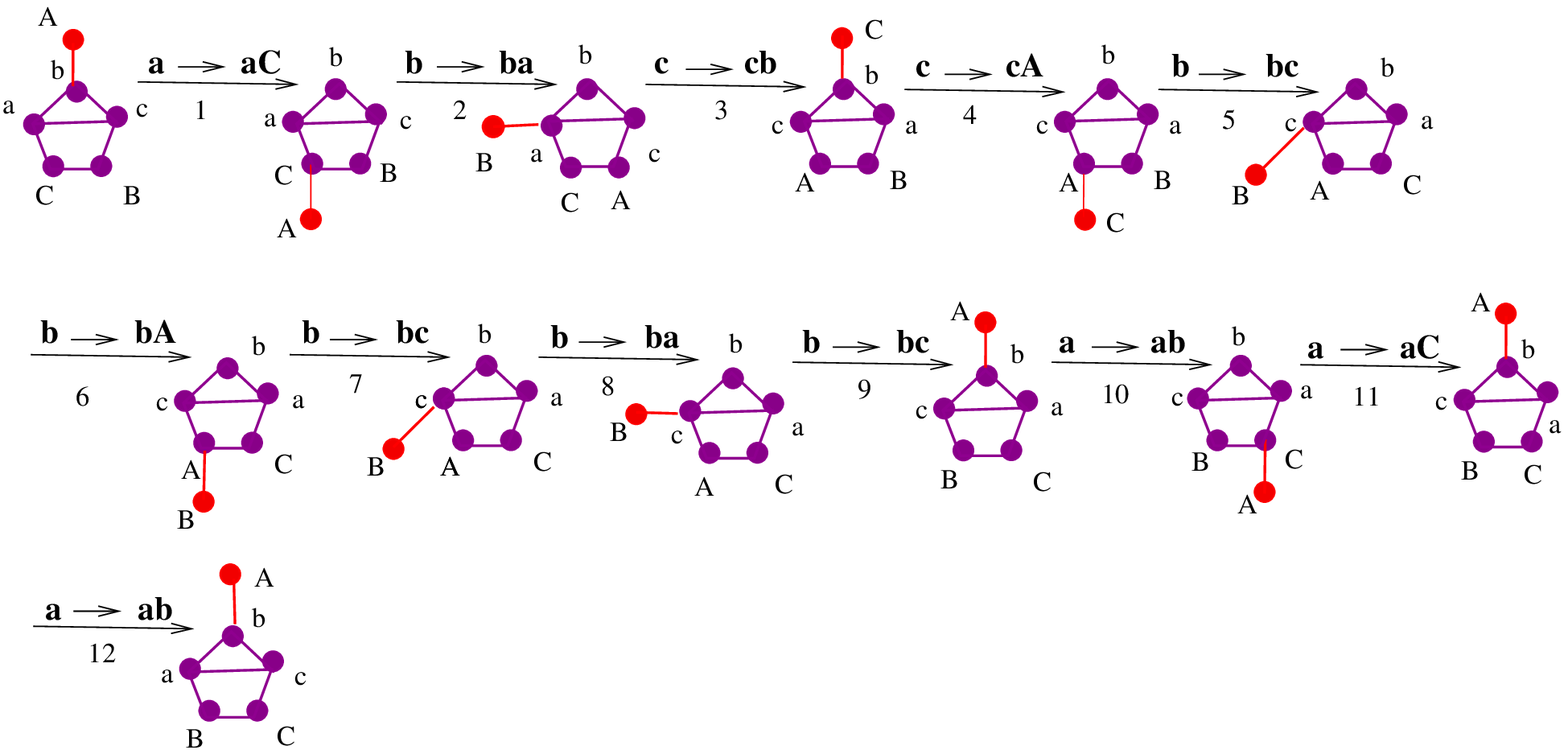}

\vskip15pt

\noindent GRAPH X:

\vskip7pt

\noindent The representative $g$ whose ideal whitehead graph, $\mathcal{G}$, is GRAPH X is:
$$
g =
\begin{cases} a \mapsto abacbabac \bar{a} b \bar{c} \bar{a} \bar{b} \bar{a} \bar{c} \bar{a} \bar{b} \bar{a} \bar{b} \bar{c} \bar{a} \bar{b} \bar{a} babac \bar{b} abacbabacabac \bar{b} a  \\
b \mapsto babac \bar{a} b \bar{c} \bar{a} \bar{b} \bar{a} \bar{c} \bar{a} \bar{b} \bar{a} \bar{b} \bar{c} \bar{a} \bar{b} \bar{a} babac \bar{a} b \bar{c} \bar{a} \bar{b} \bar{a} \bar{c} \bar{a} \bar{b} \bar{a} \bar{b} \bar{c} \bar{a} \bar{b} \bar{a} b  \\
c \mapsto babac \bar{a} b \bar{c} \bar{a} \bar{b} \bar{a} \bar{c} \bar{a} \bar{b} \bar{a} \bar{b} \bar{c} \bar{a} \bar{b} \bar{a} babac \bar{a} b \bar{c} \bar{a} \bar{b} \bar{a} \bar{c} \bar{a} \bar{b} \bar{a} \bar{b} \bar{c} \bar{a} \bar{b} \bar{a} babacbabac \bar{a} b \bar{c} \bar{a} \bar{b} \bar{a} \bar{c} \bar{a} \bar{b} \bar{a} \bar{b} \bar{c} \bar{a} \bar{b} \bar{a} babac
\end{cases}
$$

\smallskip

\noindent Instead of giving our entire ideal decomposition here, we give a condensed decomposition where construction compositions starting and ending at a graph are shown as paths below it. This is an example of a case where Method III is used to find our desired representative. If you leave out the initial generator (the upper left-most) and the pure construction compositions corresponding to the paths indicated, we have a switch sequence. \newline

\smallskip

\noindent \includegraphics[width=5.5in]{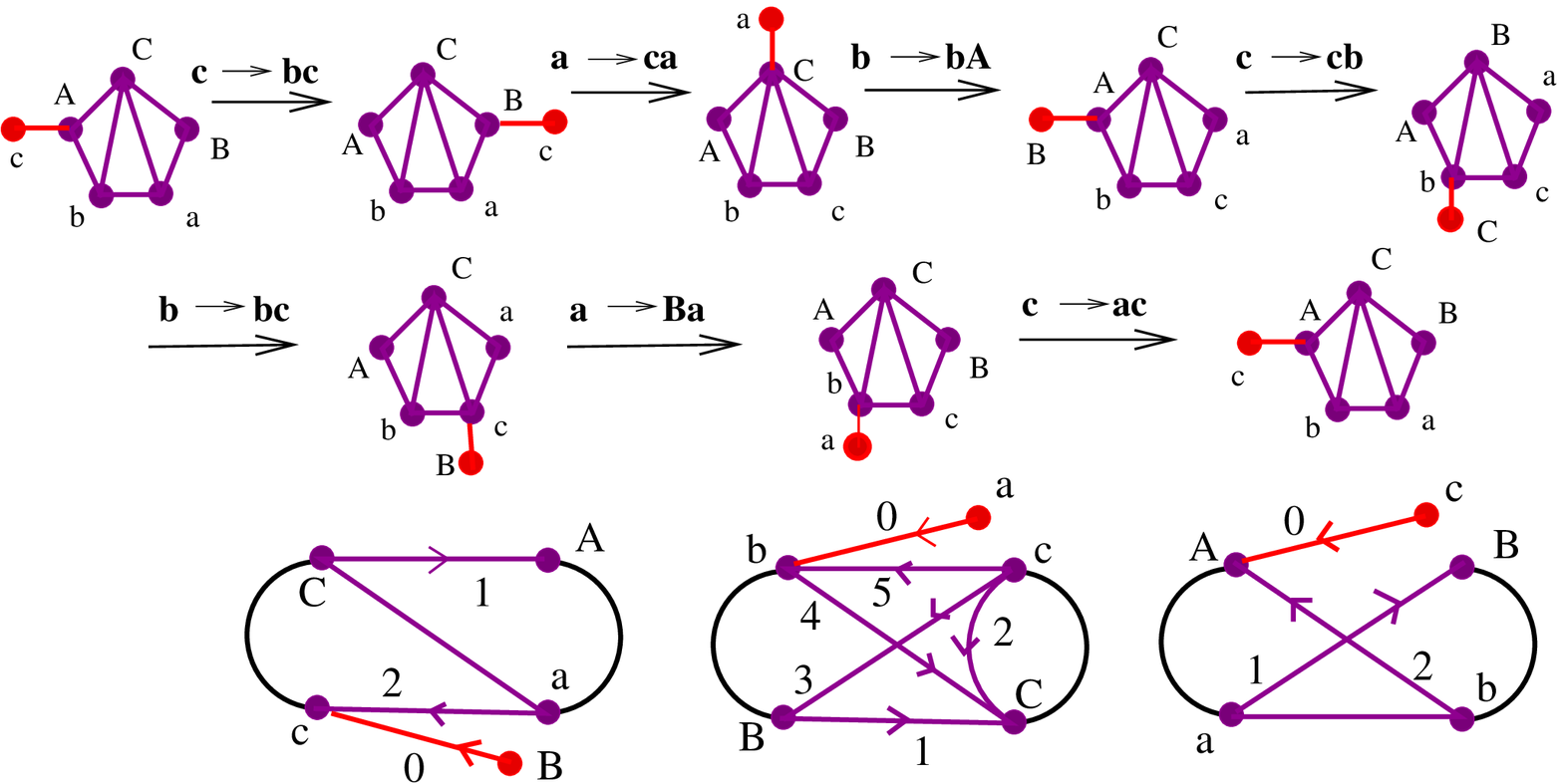}

\vskip15pt

\noindent GRAPH XI:

\vskip7pt

\noindent The representative $g$ whose ideal whitehead graph, $\mathcal{G}$, is GRAPH XI is:
$$
g =
\begin{cases} a \mapsto a \bar{b} \bar{b} \bar{c} b \bar{c} \bar{b} c \bar{b} cbc \bar{a} bc \bar{b} cbc \bar{a} \bar{b} c \bar{b}  \\
b \mapsto b \bar{c} ba \bar{c} \bar{b} \bar{c} b \bar{c} \bar{b} a \bar{c} \bar{b} \bar{c} b \bar{c} b  \\
c \mapsto c \bar{b} cbc \bar{a} \bar{b} c \bar{b} cbc \bar{a} bc \bar{b} cbc \bar{a} \bar{b} c  \end{cases}
$$

\smallskip

\noindent Again, instead of giving our entire ideal decomposition here, we give a condensed decomposition where construction compositions starting and ending at a graph are shown as paths below it.  However, the graph in the lower right actually gives a construction composition with source LTT structure above it to the left and destination LTT structure above it to the right. This representative thus actually uses a variant of Method III where we allow this.  (We in fact use a combination of Method II and Method III). \newline

\smallskip

\includegraphics[width=4.5in]{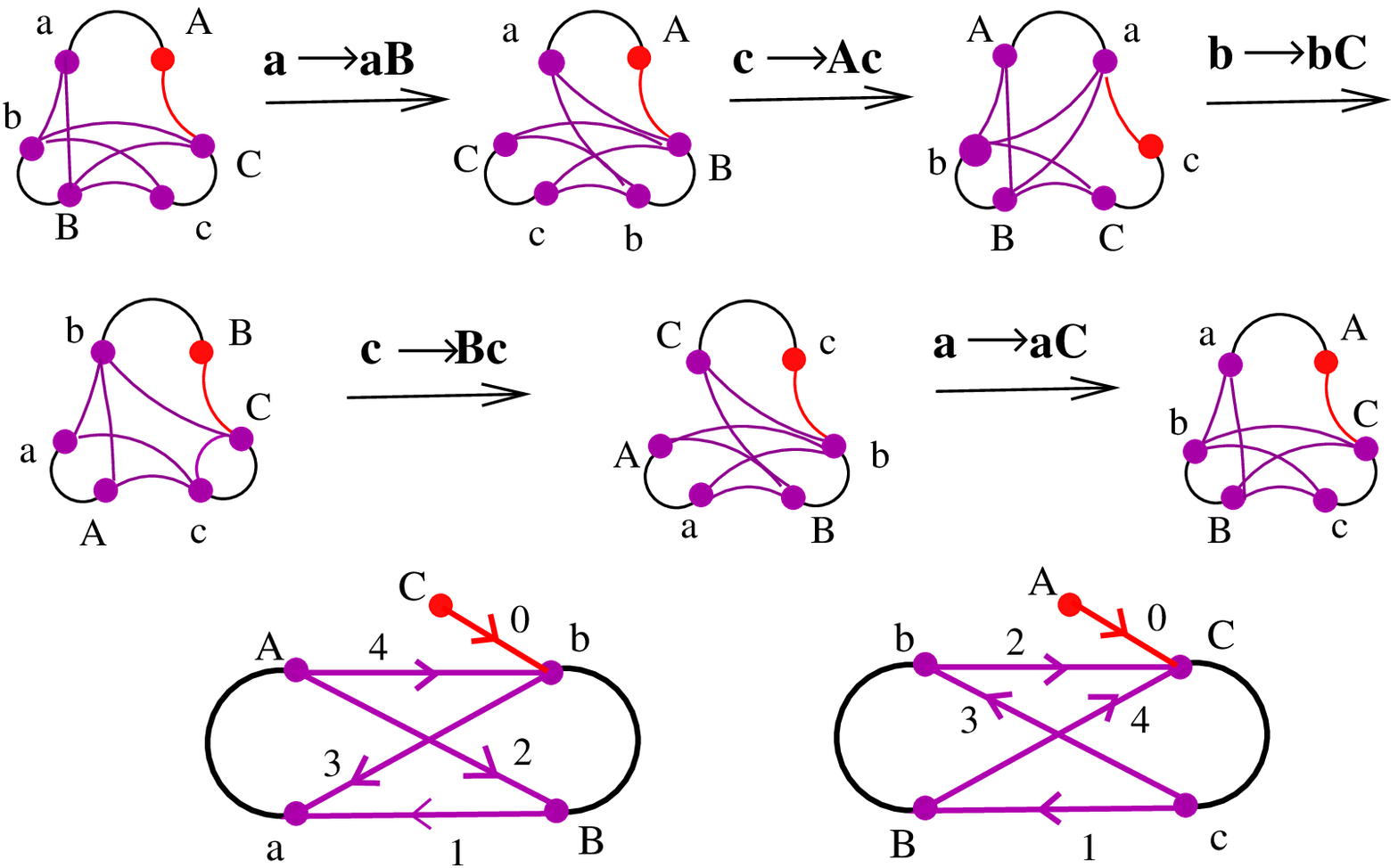}

\vskip15pt

\noindent GRAPH XII:

\vskip7pt

\noindent The representative $g$ whose ideal whitehead graph is GRAPH XII is:
$$
g =
\begin{cases} a \mapsto a \bar{c} \bar{b} \bar{b} \bar{c} b \bar{c} \bar{b} c \bar{b} cbcc \bar{a} bc \bar{b} cbcc \bar{a} \bar{b} c \bar{b}  \\
b \mapsto b \bar{c} ba \bar{c}\bar{c} \bar{b} \bar{c} b \bar{c} \bar{b} a \bar{c} \bar{c} \bar{b} \bar{c} b \bar{c} b  \\
c \mapsto c \bar{b} cbcc \bar{a} \bar{b} c \bar{b} cbcc \bar{a} bc \bar{b} cbcc \bar{a} \bar{b} c  \end{cases}
$$

\smallskip

\noindent Again, instead of giving our entire ideal decomposition here, we give a condensed decomposition where construction compositions starting and ending at a graph are shown as paths below it. \newline

\smallskip

\includegraphics[width=4.5in]{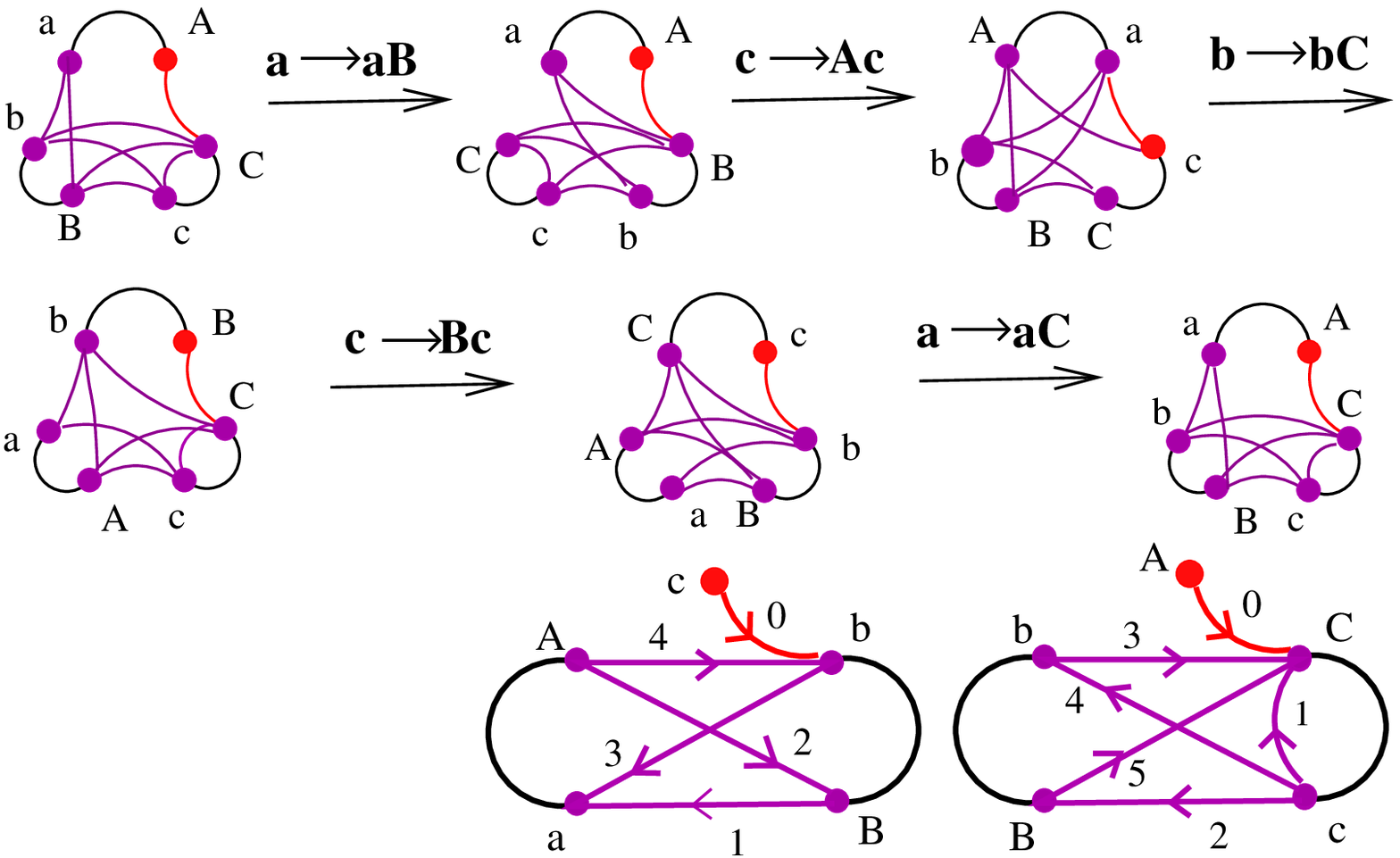}

\vskip15pt

\noindent GRAPH XIII:

\vskip7pt

\noindent The representative $g$ whose ideal whitehead graph is GRAPH XIII is:
$$
g =
\begin{cases} a \mapsto ac \bar{b} ccbc \bar{b} \bar{c} \bar{b} \bar{c} \bar{c} b \bar{c} \bar{a} c \bar{b} ac \bar{b} ccbc \bar{b} \bar{c} \bar{b} \bar{c} \bar{c} b \bar{c} \bar{a} \bar{b}  \\
b \mapsto bac \bar{b} ccbcb \bar{c} \bar{b} \bar{c} \bar{c} b \bar{c} \bar{a} b \bar{c} ac \bar{b} ccbcb  \\
c \mapsto c \bar{b} ac \bar{b} ccbc \bar{b} \bar{c} \bar{b} \bar{c} \bar{c} b \bar{c} \bar{a} \bar{b} ac \bar{b} ccbc
\end{cases}
$$

\smallskip

\noindent Our ideal decomposition for this representative and further explanation were given in Example \ref{E:GraphXIII}.

\vskip15pt

\noindent GRAPH XIV:

\vskip7pt

\noindent The representative $g$ whose ideal whitehead graph is GRAPH XIV is:
$$
g =
\begin{cases} a \mapsto acabaabcabaa  \\
b \mapsto \bar{a} \bar{a} \bar{b} \bar{a} \bar{c} \bar{b} \bar{a} \bar{a} \bar{b} \bar{a} \bar{c} \bar{a} b  \bar{a} \bar{a} \bar{b} \bar{a} \bar{c} acabaabcabaab \bar{a} \bar{a} \bar{b} \bar{a} \bar{c} \bar{a} \bar{a} \bar{b} \bar{a} \bar{c} \bar{b} \bar{a} \bar{a} \bar{b} \bar{a} \bar{c} \bar{a} b  \\
c \mapsto cabaa \bar{b}
\end{cases}
$$

\smallskip

\noindent Our ideal decomposition for $g$ is described by the following figure: \newline

\smallskip

\noindent \includegraphics[width=6.5in]{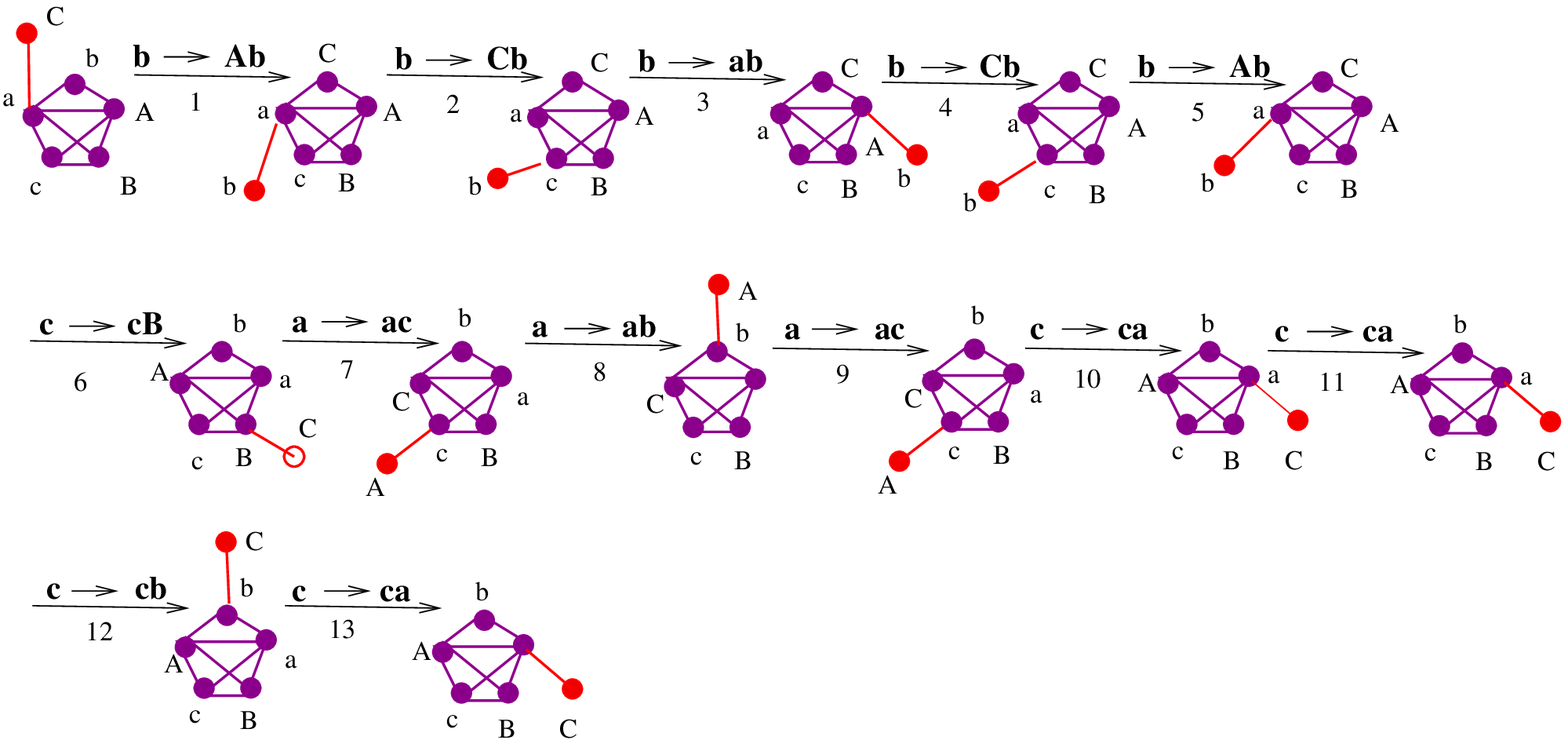} \newline

\vskip15pt

\noindent GRAPH XV:

\vskip7pt

\noindent The representative $g$ whose ideal whitehead graph, $\mathcal{G}$, is GRAPH XV is:
$$
g =
\begin{cases} a \mapsto a \bar{c} \bar{b} \bar{b} \bar{c} b \bar{c} \bar{b} c \bar{b} cbbc \bar{a} bc \bar{b} cbbc \bar{a} \bar{b} c \bar{b}  \\
b \mapsto b \bar{c} ba \bar{c} \bar{b} \bar{b} \bar{c} b \bar{c} \bar{b} a \bar{c} \bar{b} \bar{b}  \bar{c} b \bar{c} b  \\
c \mapsto c \bar{b} cbbc \bar{a} \bar{b} c \bar{b} cbbc \bar{a} bc \bar{b} cbbc \bar{a} bc  \end{cases}
$$

\smallskip

\noindent Again, instead of giving our entire ideal decomposition here, we give a condensed decomposition where construction compositions starting and ending at a graph are shown as paths below it.  The similarities between this construction and that of Graph XI are not a coincidence.  Since XI is a subgraph of XV missing only a single edge, once the representative for XI was constructed, we could alter the representative by adding the edge [$b, \bar{b}$] to the final (right-most) construction path to add that edge to $\mathcal{G}$. It must also be checked, however, that, if we add the preimages of this edge into the previous LTT structures that they are still birecurrent, that we still have a composition of switches and extensions, that we still have no PNPs, and that our initial and terminal LTT structures for the entire train track are the same. \newline

\smallskip

\includegraphics[width=4.4in]{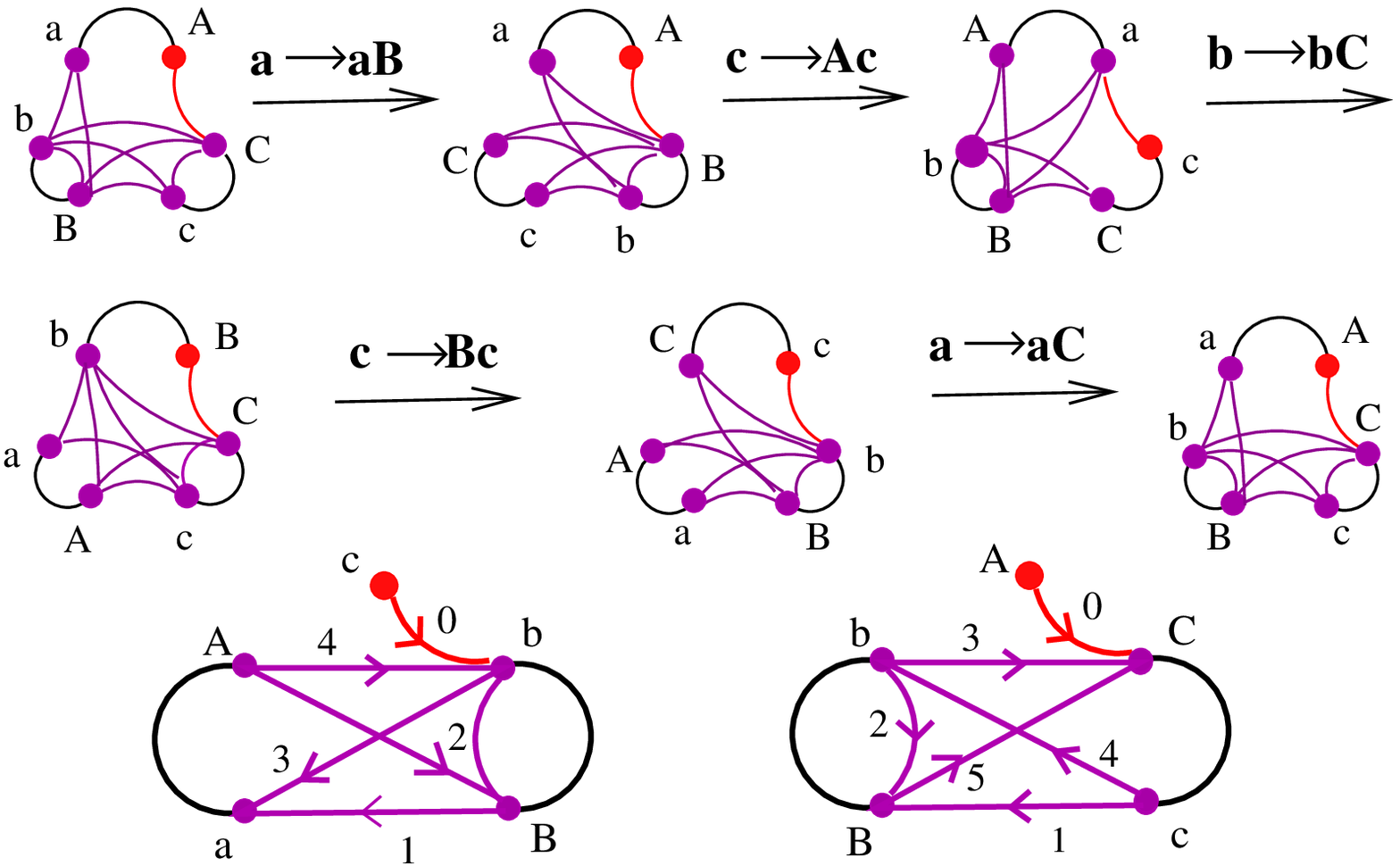}

\vskip15pt

\noindent GRAPH XVI:

\vskip7pt

\noindent The representative $g$ whose ideal whitehead graph is GRAPH XVI is:
$$
g =
\begin{cases} a \mapsto a \bar{b} ccbc \bar{b} c  \\
b \mapsto b \bar{c} \bar{b} \bar{c} \bar{c} b \bar{a} \bar{c} b  \\
c \mapsto ca \bar{b} ccbc \bar{b} c
\end{cases}
$$

\smallskip

\noindent Our ideal decomposition for $g$ is described by the following figure: \newline

\smallskip

\includegraphics[width=5.3in]{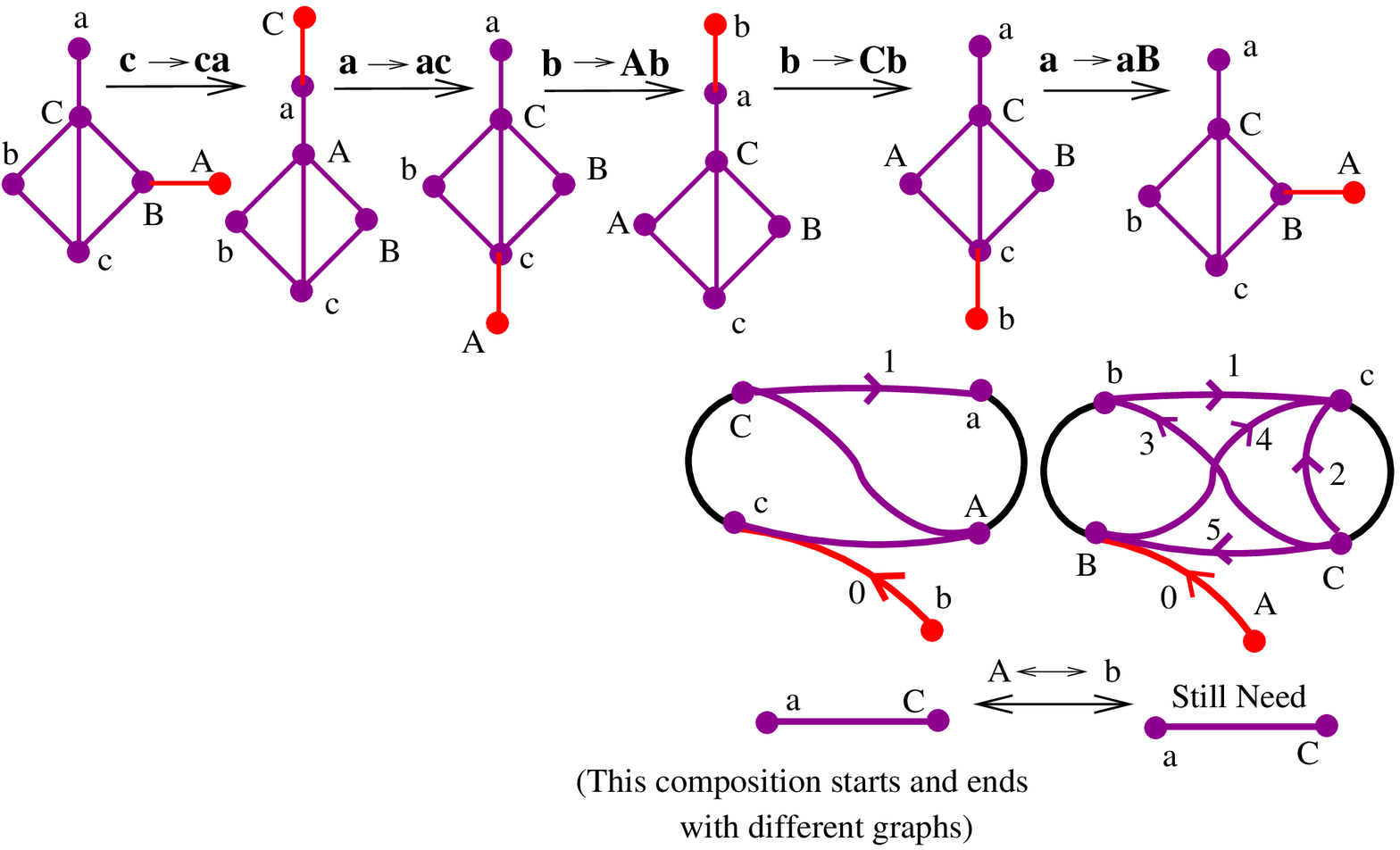}

\newpage

\noindent GRAPH XVII:

\vskip7pt

\noindent The representative $g$ whose ideal whitehead graph is GRAPH XVII is:
$$
g =
\begin{cases} a \mapsto acbcc \bar{b} c \bar{b} acbcc \bar{b} acbcc  \\
b \mapsto b \bar{c} \bar{c}\bar{b}\bar{c}\bar{a} b \bar{c} \bar{c} \bar{c}\bar{b}\bar{c}\bar{a} b \bar{c} \bar{c}\bar{b}\bar{c}\bar{a} b \bar{c} b  \\
c \mapsto c \bar{b} acbcc \bar{b} \bar{b} c \bar{b} acbcc \bar{b} acbccc \bar{b} acbcc \bar{b} acbcc \bar{b} c \bar{b} acbcc \bar{b} acbcc
\end{cases}
$$

\smallskip

\noindent Our ideal decomposition for $g$ is described by the following figure: \newline

\smallskip

\noindent \includegraphics[width=6.5in]{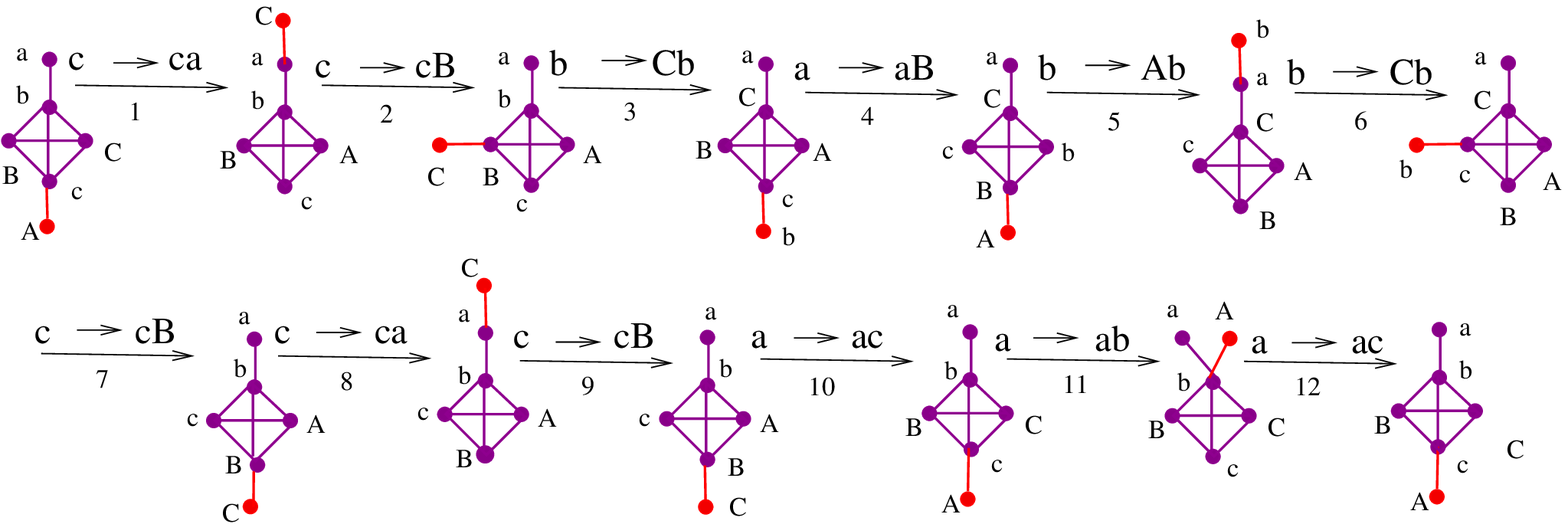}

\vskip15pt

\noindent GRAPH XVIII:

\vskip7pt

\noindent The representative $g$ whose ideal whitehead graph is GRAPH XVIII is:
$$
g =
\begin{cases} a \mapsto a \bar{b} c \bar{b} a \bar{b}  \bar{c}  \\
b \mapsto b \bar{a} b  \bar{c} b \bar{a} b  \bar{a} b  \\
c \mapsto cb \bar{a} b \bar{c} b \bar{a} b \bar{a} b \bar{c} b \bar{a} b \bar{a} bc
\end{cases}
$$

\smallskip

\noindent Our ideal decomposition for $g$ is described by the following figure: \newline

\smallskip

\includegraphics[width=5.8in]{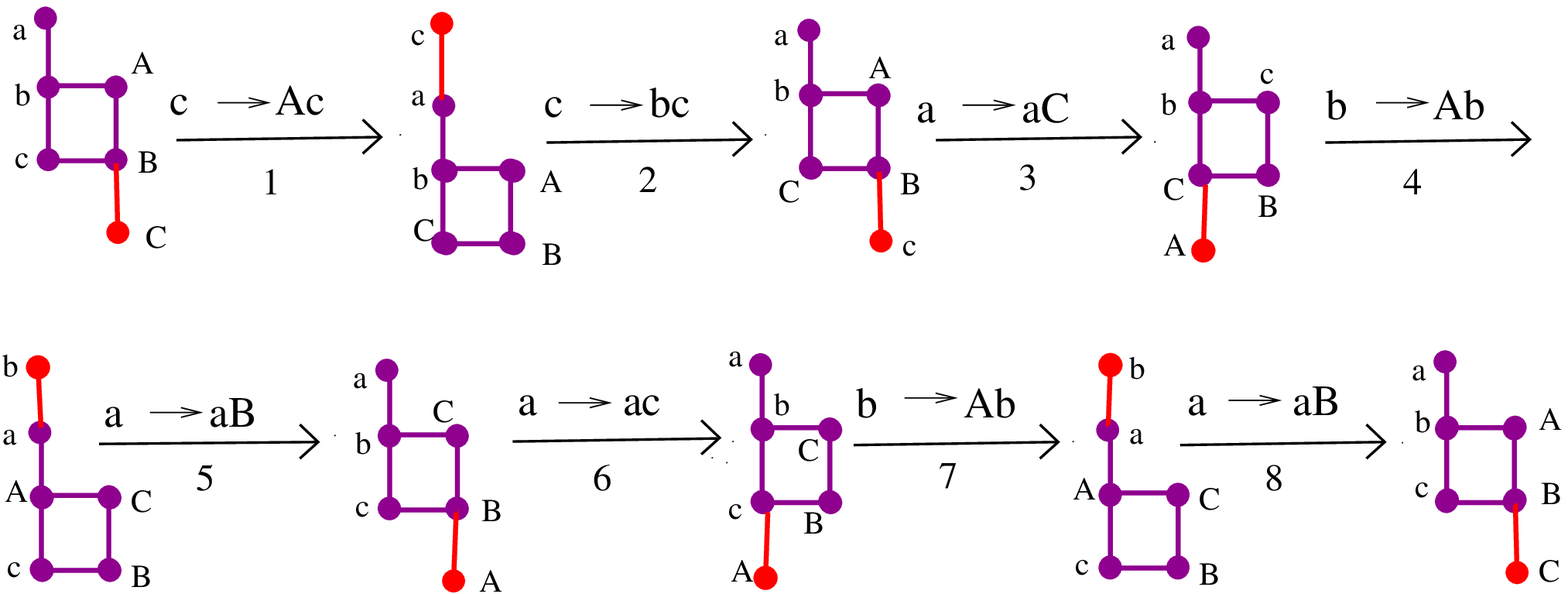}

\newpage

\noindent GRAPH XIX:

\vskip7pt

\noindent The representative $g$ whose ideal whitehead graph is GRAPH XIX is:
$$
g =
\begin{cases} a \mapsto acc \bar{b} cbc  \\
b \mapsto b \bar{c} \bar{b} \bar{c} b \bar{c} \bar{c} \bar{a} b \bar{c} b  \\
c \mapsto c \bar{b} acc \bar{b} cbc \bar{b} acc \bar{b} cbc
\end{cases}
$$

\smallskip

\noindent Again, instead of giving our entire ideal decomposition here, we give a condensed decomposition where construction compositions starting and ending at a graph are shown as paths below. \newline

\smallskip

\includegraphics[width=5.5in]{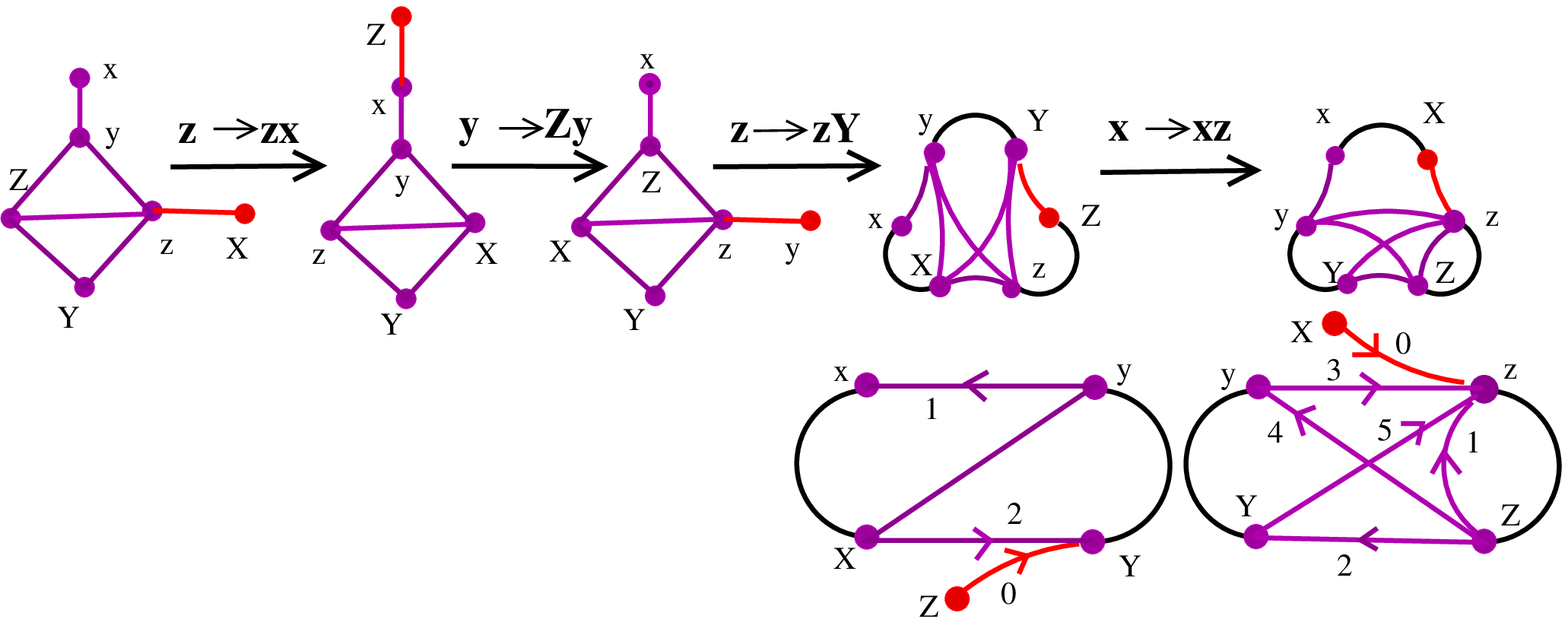}

\vskip15pt

\noindent GRAPH XX:

\vskip7pt

\noindent The representative $g=h^2$ having ideal Whitehead graph GRAPH XX, where
$$
h =
\begin{cases} a \mapsto ab \bar{c} \bar{c} bbcb  \\
b \mapsto bc  \\
c \mapsto cab \bar{c} \bar{c} bbcbab \bar{c} \bar{c} bbcb \bar{c} \bar{c} \bar{b} ab \bar{c} \bar{c} bbcbbccab \bar{c} \bar{c} bbcb
\end{cases},
$$
was constructed in the examples above.

\vskip15pt

\noindent GRAPH XXI (Complete Graph):

\vskip7pt

\noindent The representative $g$ whose ideal whitehead graph is GRAPH XXI is:
$$
g =
\begin{cases} a \mapsto aba \bar{b} aac \bar{b} aba \bar{b} aacbaba \bar{b} aacaba \bar{b} aac \bar{b} a  \\
b \mapsto baba \bar{b} aac \bar{a} b \bar{c} \bar{a} \bar{a} b \bar{a} \bar{b} \bar{a} \bar{c} \bar{a} \bar{a} b \bar{a} \bar{b} \bar{a} \bar{b} \bar{c} \bar{a} \bar{a} b \bar{a} \bar{b} \bar{a} b  \\
c \mapsto aba \bar{b} aac
\end{cases}
$$

Below we include an ideal decomposition of the representative.  We only label the vertices of the red edge in each graph of this example since the remainder of the graph is completely symmetric and so any permutation of the remaining labels gives exactly the same graph.  This is an example of a case where Method I would have been extremely impractical and Method II was particularly easy to apply.  Similar methods as used to construct this representative could also be used to construct the representative whose ideal Whitehead graph is the complete graph in any odd rank greater than five.  \newline

\smallskip

\noindent \includegraphics[width=6.3in]{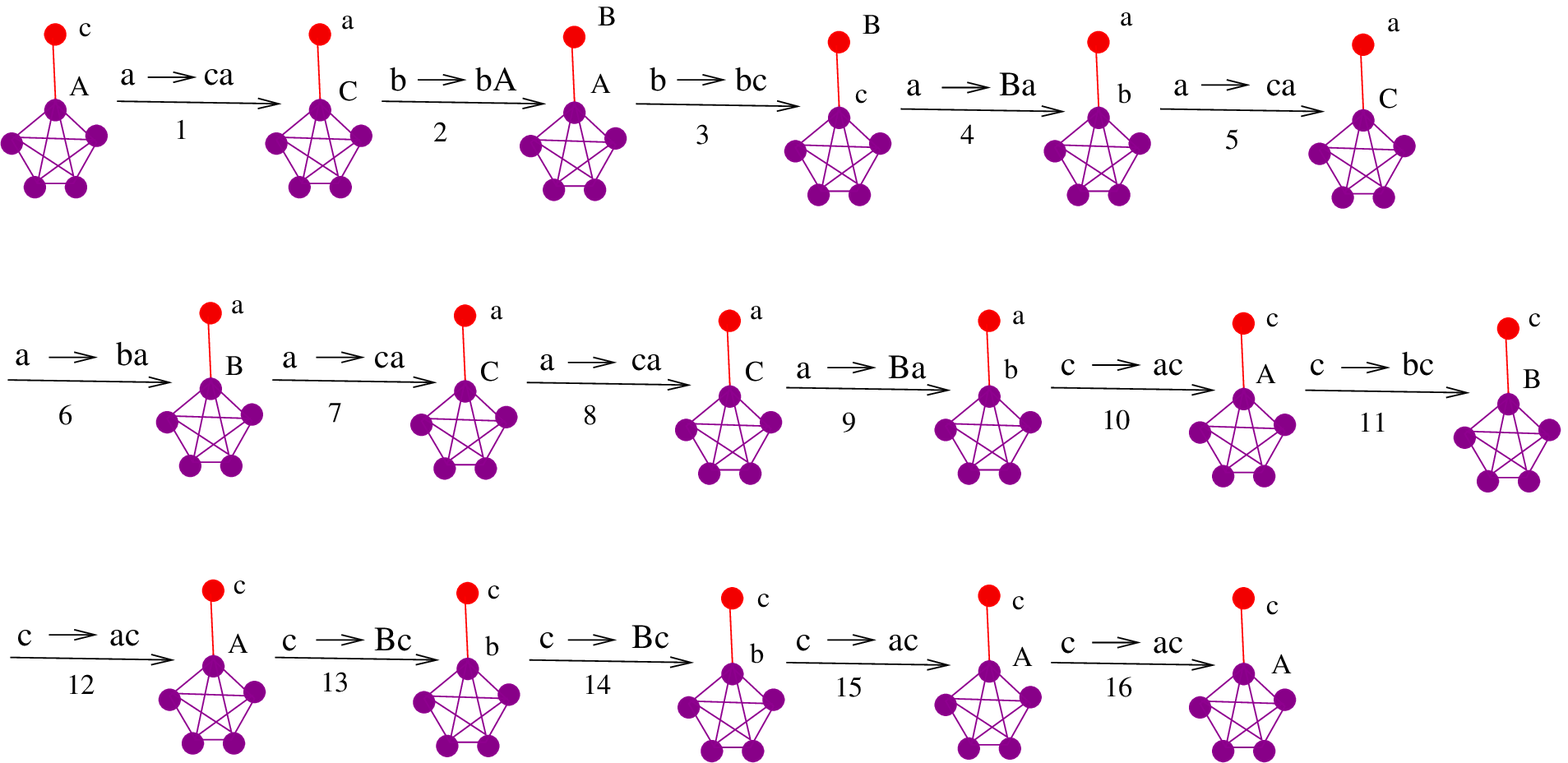}

\vskip15pt

\noindent Since we have either given representatives yielding or shown that they cannot exist for all twenty-one Type (*) pIW graphs with five vertices, we have completed the proof. \newline
\noindent QED.

\newpage

\section*{References}
\noindent [BF94] 	M. Bestvina, M. Feighn, M. Handel, \emph{Outer Limits}, preprint, 1993; available at  \newline
\indent http://andromeda.rutgers.edu/~feighn/papers/outer.pdf. \newline
\noindent [BFH00] 	M. Bestvina, M. Feighn, M. Handel, \emph{The Tits alternative for $Out(F_n)$ I: Dynamics of exponentially-growing automorphisms}, Ann. Math. 151 (2000), 517-623. \newline
\noindent [BFH97] 	M. Bestvina, M. Feighn, M. Handel, \emph{Laminations, trees, and irreducible automorphisms of free groups}, Geom. Funct. Anal. 7 (1997), 215-244. \newline
\noindent [BFH97] 	M. Bestvina, M. Feighn, M. Handel, \emph{Erratum to: Laminations, trees, and irreducible automorphisms of free groups}, Geom. Funct. Anal. 7 (1997),1143. \newline
\noindent [BH92]     M. Bestvina, M. Handel, \emph{Train tracks and automorphisms of free groups}, Ann. Math. 135 (1992), 1-51. \newline
\noindent [CP10]    M. Clay, A. Pettet, \emph{Twisting Out Fully Irreducible Automorphisms}, Geom. Funct. Anal. 20 (2010), no. 3, 657-689.  \newline
\noindent [CV86]    M. Culler and K. Vogtmann, \emph{Moduli of Graphs and automorphisms of free groups}, Invent. Math. 84 (1986), 91-119. \newline
\noindent [CP84]   D. Cvetkovic , M. Petric , \emph{A table of connected graphs on six vertices}, Discrete Math. 50 (1984) 37 49. \newline
\noindent [FB86]     Farb, Benson and Margalit, Dan. \emph{	
A Primer on Mapping Class Groups}, Princeton: Princeton University Press, 2012. \newline
\noindent [GJLL98] 	D. Gaboriau, A. Jaeger, G. Levitt, and M. Lustig, \emph{An Index for counting fixed points of automorphisms of free groups}, Duke Math. J. 93 (1998), no. 3, 425-452. \newline
\noindent [HM11] 	M. Handel and L. Mosher, \emph{Axes in outer space}, American Mathematical Society, 2011. \newline
\noindent  [LL03] Levitt, Gilbert; Lustig, Martin, \emph{Irreducible automorphisms of Fn have north-south dynamics on compactified outer space}. J. Inst. Math. Jussieu 2 (2003), no. 1, 59–72. \newline
\noindent [MS93]	H. Masur and J. Smillie, \emph{Quadratic differentials with prescribed singularities and pseudo-Anosov diffeomorphisms}, Comment. Math. Helv. 68 (1993), no.2, 289-307. \newline
\noindent [P88]	R.C. Penner. \emph{A construction of pseudo-Anosov homeomorphisms}, Trans. Amer. Math. Soc., 310 (1988) No 1, 179-197.
\noindent [N86]	\emph{Jakob Nielsen: collected mathematical papers. Vol. 1.} Edited and with a preface by Vagn Lundsgaard Hansen. Contemporary Mathematicians. Birkhäuser Boston, Inc., Boston, MA, 1986. xii + 459 pp. \newline
\noindent [St83] 	J. R. Stallings, \emph{Topology of Finite Graphs}, Inv. Math. 71 (1983), 552-565. \newline
\noindent [W36]  J. H. C. Whitehead, \emph{On certain sets of elements in a free group,} Proc. London Math. Soc., 41 (1936), pp. 48–56.

\end{document}